\documentclass[a4paper, 12pt]{amsart}

\usepackage{amsmath, amsthm, amssymb, amscd, enumitem, nccmath, framed, stmaryrd, mathtools}

\usepackage{etoolbox}

\makeatletter
\let\old@tocline\@tocline
\let\section@tocline\@tocline
\newcommand{\subsection@dotsep}{4.5}
\newcommand{\subsubsection@dotsep}{4.5}
\patchcmd{\@tocline}
  {\hfil}
  {\nobreak
     \leaders\hbox{$\m@th
        \mkern \subsection@dotsep mu\hbox{.}\mkern \subsection@dotsep mu$}\hfill
     \nobreak}{}{}
\let\subsection@tocline\@tocline
\let\@tocline\old@tocline

\patchcmd{\@tocline}
  {\hfil}
  {\nobreak
     \leaders\hbox{$\m@th
        \mkern \subsubsection@dotsep mu\hbox{.}\mkern \subsubsection@dotsep mu$}\hfill
     \nobreak}{}{}
\let\subsubsection@tocline\@tocline
\let\@tocline\old@tocline

\let\old@l@subsection\l@subsection
\let\old@l@subsubsection\l@subsubsection

\def\@tocwriteb#1#2#3{%
  \begingroup
    \@xp\def\csname #2@tocline\endcsname##1##2##3##4##5##6{%
      \ifnum##1>\c@tocdepth
      \else \sbox\z@{##5\let\indentlabel\@tochangmeasure##6}\fi}%
    \csname l@#2\endcsname{#1{\csname#2name\endcsname}{\@secnumber}{}}%
  \endgroup
  \addcontentsline{toc}{#2}%
    {\protect#1{\csname#2name\endcsname}{\@secnumber}{#3}}}%

\newlength{\@tocsectionindent}
\newlength{\@tocsubsectionindent}
\newlength{\@tocsubsubsectionindent}
\newlength{\@tocsectionnumwidth}
\newlength{\@tocsubsectionnumwidth}
\newlength{\@tocsubsubsectionnumwidth}
\newcommand{\settocsectionnumwidth}[1]{\setlength{\@tocsectionnumwidth}{#1}}
\newcommand{\settocsubsectionnumwidth}[1]{\setlength{\@tocsubsectionnumwidth}{#1}}
\newcommand{\settocsubsubsectionnumwidth}[1]{\setlength{\@tocsubsubsectionnumwidth}{#1}}
\newcommand{\settocsectionindent}[1]{\setlength{\@tocsectionindent}{#1}}
\newcommand{\settocsubsectionindent}[1]{\setlength{\@tocsubsectionindent}{#1}}
\newcommand{\settocsubsubsectionindent}[1]{\setlength{\@tocsubsubsectionindent}{#1}}

\renewcommand{\l@section}{\section@tocline{1}{\@tocsectionvskip}{\@tocsectionindent}{}{\@tocsectionformat}}%
\renewcommand{\l@subsection}{\subsection@tocline{2}{\@tocsubsectionvskip}{\@tocsubsectionindent}{}{\@tocsubsectionformat}}%
\renewcommand{\l@subsubsection}{\subsubsection@tocline{3}{\@tocsubsubsectionvskip}{\@tocsubsubsectionindent}{}{\@tocsubsubsectionformat}}%
\newcommand{\@tocsectionformat}{}
\newcommand{\@tocsubsectionformat}{}
\newcommand{\@tocsubsubsectionformat}{}
\expandafter\def\csname toc@1format\endcsname{\@tocsectionformat}
\expandafter\def\csname toc@2format\endcsname{\@tocsubsectionformat}
\expandafter\def\csname toc@3format\endcsname{\@tocsubsubsectionformat}
\newcommand{\settocsectionformat}[1]{\renewcommand{\@tocsectionformat}{#1}}
\newcommand{\settocsubsectionformat}[1]{\renewcommand{\@tocsubsectionformat}{#1}}
\newcommand{\settocsubsubsectionformat}[1]{\renewcommand{\@tocsubsubsectionformat}{#1}}
\newlength{\@tocsectionvskip}
\newcommand{\settocsectionvskip}[1]{\setlength{\@tocsectionvskip}{#1}}
\newlength{\@tocsubsectionvskip}
\newcommand{\settocsubsectionvskip}[1]{\setlength{\@tocsubsectionvskip}{#1}}
\newlength{\@tocsubsubsectionvskip}
\newcommand{\settocsubsubsectionvskip}[1]{\setlength{\@tocsubsubsectionvskip}{#1}}

\patchcmd{\tocsection}{\indentlabel}{\makebox[\@tocsectionnumwidth][l]}{}{}
\patchcmd{\tocsubsection}{\indentlabel}{\makebox[\@tocsubsectionnumwidth][l]}{}{}
\patchcmd{\tocsubsubsection}{\indentlabel}{\makebox[\@tocsubsubsectionnumwidth][l]}{}{}

\newcommand{\@sectypepnumformat}{}
\renewcommand{\contentsline}[1]{%
  \expandafter\let\expandafter\@sectypepnumformat\csname @toc#1pnumformat\endcsname%
  \csname l@#1\endcsname}
\newcommand{\@tocsectionpnumformat}{}
\newcommand{\@tocsubsectionpnumformat}{}
\newcommand{\@tocsubsubsectionpnumformat}{}
\newcommand{\setsectionpnumformat}[1]{\renewcommand{\@tocsectionpnumformat}{#1}}
\newcommand{\setsubsectionpnumformat}[1]{\renewcommand{\@tocsubsectionpnumformat}{#1}}
\newcommand{\setsubsubsectionpnumformat}[1]{\renewcommand{\@tocsubsubsectionpnumformat}{#1}}
\renewcommand{\@tocpagenum}[1]{%
  \hfill {\mdseries\@sectypepnumformat #1}}

\let\oldappendix\appendix
\renewcommand{\appendix}{%
  \leavevmode\oldappendix%
  \addtocontents{toc}{%
    \protect\settowidth{\protect\@tocsectionnumwidth}{\protect\@tocsectionformat\sectionname\space}%
    \protect\addtolength{\protect\@tocsectionnumwidth}{2em}}%
}
\makeatother

\makeatletter
\settocsectionnumwidth{2em}
\settocsubsectionnumwidth{2.5em}
\settocsubsubsectionnumwidth{3em}
\settocsectionindent{1pc}%
\settocsubsectionindent{\dimexpr\@tocsectionindent+\@tocsectionnumwidth}%
\settocsubsubsectionindent{\dimexpr\@tocsubsectionindent+\@tocsubsectionnumwidth}%
\makeatother

\settocsectionvskip{10pt}
\settocsubsectionvskip{0pt}
\settocsubsubsectionvskip{0pt}
    
\settocsectionformat{\bfseries}
\settocsubsectionformat{\mdseries}
\settocsubsubsectionformat{\mdseries}
\setsectionpnumformat{\bfseries}
\setsubsectionpnumformat{\mdseries}
\setsubsubsectionpnumformat{\mdseries}

\let\oldtableofcontents\tableofcontents
\renewcommand{\tableofcontents}{%
  \vspace*{-\linespacing}
  \oldtableofcontents}

\usepackage[tikz]{mdframed}
\usepackage[all]{xy}
\usepackage{multirow}

\usepackage{leftidx}

\theoremstyle{plain}
\newtheorem{thm}{Theorem}[section] 
\newtheorem{lem}[thm]{Lemma}

\theoremstyle{plain}
\newtheorem*{thm*}{Theorem}
\newtheorem*{lem*}{Lemma}
\newtheorem*{prop*}{Proposition}
\newtheorem*{cor*}{Corollary}
\newtheorem*{remark*}{{\sc Remark}}
\newtheorem*{problem*}{{\sc Problem}}
\newtheorem*{example*}{{\sc Example}}
\newtheorem*{defn*}{{\sc Definition}}
\newtheorem*{conjecture*}{{Conjecture}}

\newcommand{\G}{\mathbb{G}}

\newcommand{\rB}{{\rm B}}

\newcommand{\Hom}{{\rm Hom}}
\newcommand{\bm}[1]{\mbox{\boldmath$#1$}}

\newcommand{\diag}{{\rm diag}}

\newcommand{\GL}{{\rm GL}}
\newcommand{\id}{{\rm id}}
\newcommand{\inn}{{\rm inn}}
\newcommand{\Inn}{{\rm Inn}}
\newcommand{\inv}{{\rm inv}}
\newcommand{\Mat}{{\rm Mat}}
\newcommand{\Mor}{{\rm Mor}}
\newcommand{\mf}{\mathfrak}

\newcommand{\SL}{{\rm SL}}

\newcommand{\wt}{\widetilde}

\newcommand{\cU}{\mathcal{U}}

\newcommand{\Z}{\mathbb{Z}}

\newcommand{\sfP}{{\sf P}}

\newcommand{\circled}[1]{\raise0.1ex\hbox{\textcircled{\scriptsize{\raise0.2ex\hbox{#1}}}}}
\newcommand{\ouparrow}{\raise0.1ex\hbox{\textcircled{\scriptsize{\raise0.1ex\hbox{$\uparrow$}}}}}
\newcommand{\orightarrow}{ \raise0.1ex\hbox{\textcircled{\scriptsize{\raise0.23ex\hbox{\hspace{0.1mm}$\rightarrow$}}}} }

\allowdisplaybreaks[4]

\makeatletter

\renewcommand\section{\@startsection{section}{1}%
  \z@{-.5\linespacing\@plus-.7\linespacing}{.5\linespacing}%
  {\bf \Large {\normalfont\scshape}}}

\renewcommand\subsection{\@startsection{subsection}{2}%
  \z@{-.5\linespacing\@plus-.7\linespacing}{.5\linespacing}%
  {\bf \large {\normalfont\scshape}}}

\renewcommand\subsubsection{\@startsection{subsubsection}{3}%
  \z@{.5\linespacing\@plus.7\linespacing}{.5\linespacing}
   {\bf  {\normalfont\scshape}}}
\makeatother

\makeatletter
\@namedef{subjclassname}{{\rm 2020} Mathematics Subject Classification}
\makeatother

\begin{document}

\title{Homomorphisms from $\SL(2, k)$ to $\SL(4, k)$ in positive characteristic}
\author[Ryuji Tanimoto]{Ryuji Tanimoto}
\address{Faculty of Education, Shizuoka University, 836 Ohya, Suruga-ku, Shizuoka 422-8529, Japan} 
\email{tanimoto.ryuji@shizuoka.ac.jp}
\maketitle

\begin{abstract}
Let $k$ be an algebraically closed field of positive characteistic $p$ 
and let $\SL(n, k)$ denote the special linear algebraic group of degree $n$ over $k$. 
In this paper, we describe homomorphisms from $\SL(2, k)$ to $\SL(4, k)$. 
As by-products of this description, 
we give a classification of homomorphisms from $\SL(2, k)$ to $\SL(4, k)$
and describe the indecomposable decompositions of homomorphisms 
from $\SL(2, k)$ to $\SL(4, k)$. 
\end{abstract}


\section*{Introduction}

As a continuation of our study of exponential matrices (see \cite{Tanimoto 2019, 
Tanimoto 2020}), 
we were concerned with fundamental representations of $\G_a$ into $\SL(n, k)$ 
in positive characteristic $p$. 
A representation $\varphi : \G_a \to \SL(n, k)$ is said to be {\it fundamental} if 
$\varphi$ factors through a representation of $\SL(2, k)$.  
For each $1 \leq n \leq 2$, any representation $\varphi : \G_a \to \SL(n, k)$ is fundamental. 
However, for $n = 3$, Fauntleroy found in 1977, with assuming $p \geq 3$, 
a non-fundamental representation $\varphi : \G_a \to \SL(3, k)$ 
(see \cite{Fauntleroy}). 
About $45$ years later from this example,
we classify three-dimensional fundamental representations 
$\varphi : \G_a \to \SL(3, k)$ in positive characteristic $p$ (see \cite{Tanimoto 2022}). 
And then we can find many non-fundamental representations $\varphi : \G_a \to \SL(3, k)$ 
in any characteristic $p \geq 2$. 
As a by-product of our classification of three-dimensional fundamental representations 
of $\G_a$, we can describe homomorphisms from $\SL(2, k)$ to $\SL(3, k)$ in positive characteristic $p$. 

Based on the above circumstances, we became interetsted in 
homomorphisms from $\SL(2, k)$ to $\SL(n, k)$ over an algebraically closed field $k$ 
in positive characteristic $p$. 
There is a one-to-one correspondence between the set of all homomorphisms from 
$\SL(2, k)$ to $\SL(n, k)$ and the set of all homomorphisms from $\SL(2, k)$ to $\GL(n, k)$
(cf. \cite[Lemma 2.7]{Tanimoto 2022}). 
But, in positive characteristic, 
there are indecomposable representations of $\SL(2, k)$  
which are not irreducible (cf. \cite{Nagata 1961}).

Here, we raise the following problem on which we can work and 
in which we are interested:   
\begin{problem*}
Describe the forms of homomorphisms from $\SL(2, k)$ to $\SL(4, k)$ 
in positive characteristic. 
\end{problem*}

The answer to the above problem is given in Theorem 5.26 and Subsection 5.1. 
We extract the following table from Theorem 5.26: 
\medskip 

{\scriptsize 
\begin{center} 
\renewcommand{\arraystretch}{1.1} 
\begin{tabular}{| p{5em} | p{5em} | p{5em} | p{5em} | p{5em} |}
\hline 
 $p = 2$ 
 & $ p = 3$ 
 & $p \geq 5$ 
 & $p \geq 2$ 
 & $d$ \\
\hline 
 & 
 & $\rm (I)^*$
 &  
 & $(0, 0)$ \\
\hline  
 & $\rm (II)^*$ 
 & 
 & 
 & $(0, 0)$ \\
\hline 
 $\rm (IV)^*$
 & $\rm (IV)^*$ 
 & $\rm (IV)^*$
 & $\rm (IV)^*$ 
 & $(0, 0)$ \\ 
\hline 
 $\rm (V)^*$
 & 
 & 
 & 
 & $(1, 1)$ \\
\hline 
 & $\rm (VII)^*$
 &  
 & 
 & $(0, 0)$ \\
\hline  
 & $\rm (IX)^*$
 & $\rm (IX)^*$ 
 & 
 & $(1, 1)$ \\
\hline 
 $\rm (XI)^*$
 & 
 & 
 & 
 & $(1, 2)$ \\
\hline
 $\rm (XV)^*$
 & $\rm (XV)^*$ 
 & $\rm (XV)^*$
 & $\rm (XV)^*$ 
 & $(0, 0)$ \\
\hline
 $\rm (XIX)^*$
 & 
 & 
 & 
 & $(2, 1)$ \\
\hline
 $\rm (XXIV)^*$
 & $\rm (XXIV)^*$ 
 & $\rm (XXIV)^*$ 
 & $\rm (XXIV)^*$ 
 & $(2, 2)$ \\
\hline
 $\rm (XXVI)^*$
 & $\rm (XXVI)^*$ 
 & $\rm (XXVI)^*$
 & $\rm (XXVI)^*$ 
 & $(4, 4)$ \\
\hline 
$7$ types 
 & $7$ types
 & $6$ types 
 & $4$ types 
 \\
\cline{1-4} 
\end{tabular}
\end{center} 
}
\medskip 

The number of types of homomorphisms from $\SL(2, k)$ to $\SL(4, k)$ 
can be stated below. 
If $p =2$, seven types appear. 
If $p = 3$, seven types appear. 
If $p \geq 5$, six types appear. 
And four types appear in common among all characteristics $p \geq 2$.

As known from the above table, 
it is too hard for us to answer intuitively to the above problem 
(the appearing types vary in accordance with the characteristics $p$). 
We go on working steadily. We solve the problem by separating many cases. 
So, it would be better to write the overall flow of the solution in this introduction.  
\\

\begin{center}
\bf  ------  The outline of the solution ------
\end{center}

In order to describe the forms of homomorphisms from $\SL(2, k)$ to $\SL(4, k)$, 
we focus on the following three subgroups $U^+$, $T$, $U^-$ of $\SL(2, k)$: 
\begin{align*}
 U^+ 
 & := 
\left\{
\left. 
\left(
\begin{array}{c c}
 1 & t \\
 0 & 1
\end{array} 
\right)
\; \right| \; 
 t \in k
\right\} , \\
 T & := 
\left\{ 
\left. 
\left(
\begin{array}{c c}
 u & 0 \\
 0 & u^{-1}
\end{array} 
\right) 
\; \right| \; 
 u \in k \backslash \{ \, 0 \, \} 
\right\} , \\
 U^- 
 & := 
\left\{
\left. 
\left(
\begin{array}{c c}
 1 & 0 \\
 s & 1
\end{array} 
\right)
\; \right| \; 
 s \in k
\right\} . 
\end{align*}
It is well known that $\SL(2, k)$ is generated by the above three subgroups $U^+$, $T$, $U^-$. 
Furthermore, the subgroup $U^+$ is isomorphic to the additive group $\G_a$ of $k$, 
the subgroup $T$ is isomorphic to the multiplicative group $\G_m$ of $k$, 
and the subgroup $U^-$ is isomorphic to $\G_a$.

Two homomorphisms $h_1$ and $h_2$ from an algebraic group $G$ 
to $\SL(n, k)$ is said to be {\it equivalent} if there exists a regular matrix 
$P$ of $\GL(n, k)$ such that $P^{-1} \, h_1(g) P = h_2(g) $ for all $g \in G$. 
If two homomorphisms $h_1 : G \to \SL(n, k)$ and $h_2 : G \to \SL(n, k)$ 
are equivalent, we write $h_1 \sim h_2$. 
%
\\

\noindent 
{\bf $\blacksquare$ \underline{Reducing the problem.}} 
\medskip 

Given a homomorphism $\sigma : \SL(2, k) \to \SL(n, k)$, 
we can replace $\sigma$ with a conjugation $\sigma^*$ of $\sigma$ so that 
$\sigma^*$ is antisymmetirc (see Lemma 1.20 (1)).    
\[
\xymatrix@R=36pt@C=36pt@M=6pt{
 \SL(2, k) \ar[r]^\sigma \ar@{-->}[rd]_{\sigma^* \, : \,  \text{antisymmetric} \quad  } 
 & \SL(4, k)  \ar[d]^{\Inn_P} \\
 & \SL(4, k) 
}
\]
Here, $\Inn_P : \SL(4, k) \to \SL(4, k)$ is the isomorphism defined by $\Inn_P(A) := P^{-1} \, A \, P$ for all $A \in \SL(4, k)$, 
where $P \in \GL(4, k)$. 
We denote by $\Hom^a(\SL(2, k), \SL(4, k))$ the set of all 
antisymmetric homomorphisms from $\SL(2, k)$ to $\SL(4, k)$. 
We can obtain the following natural one-to-one correspondence: 
\[
 \Hom(\SL(2, k), \SL(4, k)) / \sim 
\; \cong \; 
\Hom^a(\SL(2, k), \SL(4, k)) / \sim . 
\] 
Thus we can reduce the description of $\Hom(\SL(2, k), \SL(4, k)) / \sim$ 
to the description of 
\[
\Hom^a(\SL(2, k), \SL(4, k)) / \sim . 
\] 
\medskip

\noindent 
{\bf $\blacksquare$ \underline{In advance.}}
\medskip 

Let $\rB(2, k)$ denote the Borel subgroup generated by $U^+$ and $T$ 
and let $\Hom^a(\rB(2, k), \SL(4, k))$ denote the set of 
all antisymmetric homomorphisms from $\rB(2, k)$ to $\SL(4, k)$. 
The following diagram is commutative: 
{\small 
\[
\xymatrix@R=36pt@C=36pt@M=6pt{
 \Hom(\SL(2, k), \SL(4, k)) \ar[d] 
 & \Hom^a(\SL(2, k), \SL(4, k)) \ar@{_(->}[l] \ar[r] \ar[d]
 & \Hom^a(\rB(2, k), \SL(4, k))  \ar[d] \\
 \Hom(\SL(2, k), \SL(4, k))  / \sim 
 & \Hom^a(\SL(2, k), \SL(4, k)) / \sim  \ar[l] \ar[r] 
 & \Hom^a(\rB(2, k), \SL(4, k)) / \sim 
}
\]
}

In advance, we explain how to use the description of 
$\Hom^a(\rB(2, k), \SL(4, k)) / \sim$ 
for describing $\Hom^a(\SL(2, k), \SL(4, k)) / \sim$. 
Consider any class appearing in $\Hom^a(\rB(2, k), \SL(4, k)) / \sim$ 
and write $\psi$ for a representative of the class. 
So,  $\psi \in \Hom^a(\SL(2, k), \SL(4, k))$. 
We determine whether or not $\psi$ is extendable to an element $\sigma^*$ of 
$\Hom^a(\SL(2, k), \SL(4, k))$ (see the above commutative diagram). 
If $\psi$ is extendable, we can show that $\psi$ is uniquely extendable to 
the element $\sigma^*$ of $\Hom^a(\SL(2, k), \SL(4, k))$. 
We collect such elements $\sigma^*$ of $\Hom^a(\SL(2, k), \SL(4, k))$, and 
then we can describe $\Hom^a(\SL(2, k), \SL(4, k))/ \sim$. 
\\

\noindent 
{\bf $\blacksquare$ \underline{A morphism from $\Hom^a(\rB(2, k), \SL(4, k))$ 
to $\cU_4 \times \Omega(4)$.}}
\medskip

Since $\rB(2, k)$ is generated by $U^+$ and $T$, 
we have two natural closed immersions $\imath_1' : \G_a \to \rB(2, k)$ and 
$\imath_2' : \G_m \to \rB(2, k)$ whose images are $\imath_1'(\G_a) = U^+$ and 
$\imath_2'(\G_m) = T$. 
Thus we have the following morphism: 
\[
\xymatrix@R=3pt@C=36pt@M=6pt{
 \Hom^a(\rB(2, k), \SL(4, k)) 
 \ar[r] 
 & \Hom(\G_a, \SL(4, k)) \times \Hom(\G_m, \SL(4, k)) \\
\psi \ar[r] 
 & (\psi \circ \imath_1', \, \psi \circ \imath_2')
}
\]
We can shrink the target. This is because any antisymmetric homomorphism 
$\psi : \rB(2, k) \to \SL(4, k)$ has the following crucial 
property (u) (see Lemma 1.10): 
\[
\text{
(u) \quad 
For any $t \in k$, the regular matrix $(\psi \circ \imath_1')(t)$ is 
an upper triangular matrix. 
}
\]
So, we denote by $\cU_4$ the set of all homomorphisms $\varphi : \G_a \to \SL(4, k)$ 
such that $\varphi(t)$ are upper triangular for all $t \in \G_a$.  
We denote by $\Omega(4)$ the set of all homomorphisms $\omega : \G_m \to \SL(4, k)$ 
such that $\omega$ is antisymmetric. 
Then we have the following commutative diagram:
\[
\xymatrix@R=36pt@C=36pt@M=6pt{
 \Hom^a(\rB(2, k), \SL(4, k)) 
 \ar[r] \ar[rd] 
 & \Hom(\G_a, \SL(4, k)) \times \Hom(\G_m, \SL(4, k)) \\
 & \cU_4 \times \Omega(4) \ar@{^(->}[u]
} 
\]

\noindent 
{\bf $\blacksquare$ \underline{Describing $\Hom^a(\rB(2, k), \SL(4, k)) / \sim$.}}
\medskip

The set $\cU_4$ is the disjoint union of the following eight subsets (see Lemma 1.2): 
\[
 \cU_{[4]}, \quad 
\cU_{[3, 1]} , \quad 
\cU_{[2, 2]} , \quad 
\cU_{[1, 3]} , \quad 
\cU_{[2, 1, 1]} , \quad 
\cU_{[1, 2, 1]} , \quad 
\cU_{[1, 1, 2]} , \quad 
\cU_{[1, 1, 1, 1]} . 
\]
Choose any pair $(\varphi, \omega)$ from $\cU_4 \times \Omega(4)$ and 
assume that $\varphi$ lies in one of the above subsets. 
We have a morphism $\psi_{\varphi, \, \omega} : \G_a \rtimes \G_m \to \SL(4, k)$ 
obtained from the pair $(\varphi, \omega)$ 
(see Subsection 1.3.1 for the definition of $\psi_{\varphi, \, \omega}$). 
So, $\psi_{\varphi, \, \omega} \in \Mor(\G_a \rtimes \G_m, \SL(4, k))$, 
where for varieties $X$, $Y$, we denote by $\Mor(X, Y)$ the set of all morphisms 
from $X$ to $Y$. 
We have the following commutative diagram: 
\[
\xymatrix@R=36pt@C=36pt@M=6pt{
& \Mor(\G_a \rtimes \G_m, \SL(4, k))
\\
 \Hom^a(\rB(2, k), \SL(4, k))  \ar@{-->}[ru]^{\text{the composition morphism} \qquad \qquad }
 \ar[r] \ar[rd] 
 & \Hom(\G_a, \SL(4, k)) \times \Hom(\G_m, \SL(4, k)) 
\ar[u]_{(\varphi, \, \omega) \; \mapsto \; \psi_{(\varphi, \, \omega)}} \\
 & \cU_4 \times \Omega(4) \ar@{^(->}[u]
} 
\]
We apply the homomorphism criterion (see Lemma 1.9) to the $\psi_{\varphi, \, \omega}$ 
so that we have $\psi_{\varphi, \, \omega} \in \Hom(\G_a \rtimes \G_m, \SL(4, k))$. 
We then almost determine the form of $(\varphi, \omega)$.  
Successively, using an appropriate regular matrix $P$ of $\GL(4, k)$, 
we can arrange the form of the pair $(\varphi, \omega)$. 
Through the process,  
we can describe $\Hom^a(\rB(2, k), \SL(4, k)) / \sim$, 
where we identify $\G_a \rtimes \G_m$ with $\rB(2, k)$ using the natural isomorphism 
$\jmath : \G_a \rtimes \G_m \to \rB(2, k)$ (see Theorem 3.1).  
This description gives $26$ antisymmetric homomorphisms 
from $\rB(2, k)$ to $\SL(4, k)$. 
\\

\noindent
{\bf $\blacksquare$ \underline{A necessary condition 
under assuming $\psi_{\varphi^*, \, \omega^*} \circ \jmath^{-1}$ is extendable.}} 
\medskip

Let $(\varphi^*, \omega^*)$ be a pair of $\cU(4) \times \Omega(4)$ such that 
$\psi_{\varphi^*, \, \omega^*} \circ \jmath^{-1} \in \Hom^a(\rB(2, k), \SL(4, k))$. 
Assume that $\psi_{\varphi^*, \, \omega^*} \circ \jmath^{-1}$ is extendable to 
an element $\sigma^*$ of $\Hom^a(\SL(2, k), \SL(4, k))$. 
Then we can define a homomoprhism $\phi^- : \G_a \to \SL(4, k)$ as 
\[
 \phi^-(s) 
 := 
\sigma^*
\left(
\begin{array}{c c}
 1 & 0 \\
 s & 1
\end{array}
\right) . 
\]
Then the following conditions (i) and (ii) hold true: 
\begin{enumerate}[label = {\rm (\roman*)}]
\item For any $s \in \G_a$, 
the regular matrix $\phi^-(s)$ is a lower triangular matrix.  

\item 
$
 \varphi^*(t) \, \phi^-(s) 
 = 
\phi^-\left( \frac{s}{1 + t \, s} \right) 
\, 
\omega^*(1 + t \, s) 
\;
\varphi^* \left( \frac{t}{1 + t \, s} \right)
$ 
for all $t, s \in k$ with $1 + t \, s \ne 0$.    
\end{enumerate} 
Using the above conditions (i) and (ii), 
we give a necessary condition for the element $\psi$ 
of $\Hom^a(\rB(2, k), \SL(4, k))$ to be extendable to 
an element $\sigma^*$ of $\Hom^a(\SL(2, k), \SL(4, k))$. 
As the necessary condition, we obtain either the unique form of $\phi^-$ 
or a contradiction (see Subsection 4.1). 
\\

\noindent
{\bf $\blacksquare$ \underline{Describing $\Hom^a(\SL(2, k), \SL(4, k)) / \sim$.}} 
\medskip

For any element $\psi_{\varphi^*, \, \omega^*} \circ \jmath^{-1} 
\in \Hom^a(\rB(2, k), \SL(4, k))$ which is extendable to an element $\sigma^*$ of $\Hom^a(\SL(2, k), \SL(4, k))$, 
using the unique form of $\phi^-$ and the forms of $\varphi^*$ and $\omega^*$, 
we can obtain the form of $\sigma^*$ (see Subsection 4.2). 
Then we can show that each $\sigma^* : \SL(2, k) \to \SL(4, k)$ 
becomes a homomorphism from $\SL(2, k)$ to $\SL(4, k)$.  
Collecting such homomorphisms $\sigma^* : \SL(2, k) \to \SL(4, k)$, 
we can describe $\Hom^a(\SL(2, k), \SL(4, k)) / \sim$ (see Theorem 5.26).   
\\

\begin{center}
\bf ------  By-products ------
\end{center}

As by-products of the description of $\Hom^a(\SL(2, k), \SL(4, k)) / \sim$, 
we give a classification of homomorphisms from $\SL(2, k)$ to $\SL(4, k)$ 
in positive characteristic (see Theorem 6.26) 
and describe the indecomposable decompositions of homomorphisms 
from $\SL(2, k)$ to $\SL(4, k)$ in positvie characteristic (see Section 7). 
\\

\begin{center}
\bf ------  Notations ------
\end{center}

Let $k$ be an algebraically closed field of positive characteristic $p$.  
\begin{itemize}
\item Let $\G_a$ denote the additive group of $k$ 
and let $\G_m$ denote the multiplicative group of $k$. 

\item Let $k[T]$ denote the polynomial ring in one variable over $k$. 

\item A polynomial $f(T)$ is said to be a {\it $p$-polynomial} 
if $f(T)$ has the following form: 
\begin{align*}
 f(T) 
 = a_0 \, T + a_1 \, T^p + \cdots + a_r \, T^{p^r} 
\qquad (\, a_0, a_1, \ldots, a_r \in k \,) . 
\end{align*}

\item We denote by $\sfP$ the set of all $p$-polynomials of $k[T]$. 

\item For algebraic groups $G$, $G'$ over $k$, we denote by $\Hom(G, G')$ 
the set of all homomorphisms from $G$ to $G'$ over $k$. 
\end{itemize} 

Let $R$ be a commutative ring with unity.  
\begin{itemize}
\item 
We denote by $\Mat_{m, n}(R)$ the left $R$-module 
of all $m \times n$ matrices whose entries belong to $R$. 
We denote by $O_{\Mat_{m, n}(R)}$ the zero matrix of $\Mat_{m, n}(R)$. 
We write $\Mat_{n, n}(R)$ as $\Mat(n, R)$. 
The $R$-module $\Mat(n, R)$ becomes an $R$-algebra. 
We denote by $I_{\Mat(n, R)}$ the identity matrix of $\Mat(n, R)$. 
The zero matrix $O_{\Mat_{m, n}(R)}$ is frequently referred as $O_{m, n}$ or $O$, 
and the identity matrix $I_{\Mat(n, R)}$ as $I_n$ or $I$.

We say that a matrix $\Mat(n, R)$ is {\it regular} if there exists a matrix 
$X \in \Mat(n, R)$ so that $A X = X A = I_n$. 
We denote by $\GL(n, R)$ the group of all regular matrices of $\Mat(n , R)$. 
We denote by $\SL(n, R)$ the subgroup of all regular matrices of $\GL(n, R)$ 
whose determinants are $1$.

\item For $a_1, \ldots, a_n \in R$, we denote by $\diag(a_1, \ldots, a_n)$ 
the diagonal matrix of $\Mat(n, R)$ whose $(i, i)$-th entry is $a_i$ for any 
$1 \leq i \leq n$.

\item Assume that $R$ is an algebraically closed field. 
For any regular matrix $P$ of $\GL(n, R)$, 
we can define an isomorphism $\Inn_P : \SL(n, R) \to \SL(n, R)$ of algebraic groups over $R$ 
as 
\[
 \Inn_P(A) := P^{-1} \, A \, P . 
\]

\item For $1 \leq \lambda < \mu \leq n$, we denote by $P_{\lambda,\, \mu} = (p_{i, j})$ 
be the regular matrix of $\GL(n, R)$ defined by 
\[
 p_{i, j}
 := 
 \left\{
 \begin{array}{l @{\qquad} l}
  1 & \text{ if \quad $(i, j) = (\lambda, \mu)$}, \\
  1 & \text{ if \quad $(i, j) = (\mu, \lambda)$}, \\
  1 & \text{ if \quad $i = j$ and $i \not \in \{\, \lambda, \, \mu \, \}$}, \\
  0 & \text{ otherwise} . 
 \end{array} 
 \right.
\]

\end{itemize}

\section{Homomorphisms} 

We consider homomorphisms from $\SL(2, k)$ to $\SL(4, k)$. 
So, we focus on the following three subgroups $U^+$, $T$, $U^-$ of $\SL(2, k)$: 
\begin{align*}
 U^+ 
 & := 
\left\{
\left. 
\left(
\begin{array}{c c}
 1 & t \\
 0 & 1
\end{array} 
\right)
\; \right| \; 
 t \in k
\right\} , \\
 T & := 
\left\{ 
\left. 
\left(
\begin{array}{c c}
 u & 0 \\
 0 & u^{-1}
\end{array} 
\right) 
\; \right| \; 
 u \in k \backslash \{ \, 0 \, \} 
\right\} , \\
 U^- 
 & := 
\left\{
\left. 
\left(
\begin{array}{c c}
 1 & 0 \\
 s & 1
\end{array} 
\right)
\; \right| \; 
 s \in k
\right\} . 
\end{align*}
Let $\rB(2, k)$ be the Borel subgroup of $\SL(2, k)$ generated by 
$U^+$ and $T$, i.e., 
\[
 \rB(2, k) := 
\left\{
\left. 
\left(
\begin{array}{c c}
 a & b \\
 c & d
\end{array}
\right) 
\in \SL(2, k) 
\;  \right| \;  
 c = 0\, 
\right\} . 
\] 
We can identify the Borel subgroup $\rB(2, k)$ with a semi-direct product 
$\G_a \rtimes \G_m$, where 
the product of elements $(t, u)$, $(t', u')$ of $\G_a \rtimes \G_m$ 
is given by 
\[
 (t,\, u) \cdot (t', \, u') := (t + u^2 \, t', \, u \, u' ) . 
\]
The identification is given by the isomorphism $\jmath : \G_a \rtimes \G_m \to \rB(2, k)$ 
defined by 
\[
\jmath(t, u) 
 := 
\left(
\begin{array}{c c}
 u & t \, u^{-1} \\
 0 & u^{-1} 
\end{array}
\right) 
\, 
\left(
= 
\left(
\begin{array}{c c}
 1 & t \\
 0 & 1
\end{array} 
\right) 
\left(
\begin{array}{c c}
 u & 0 \\
 0 & u^{-1} 
\end{array}
\right)
\right) .  
\]
Let $\imath_1 : \G_a \to \G_a \rtimes \G_m$ and 
$\imath_2 : \G_m \to \G_a \rtimes \G_m$ be homomorphisms defined by 
\[
 \imath_1(t) := (t, 1)  \qquad \text{ and } \qquad  
 \imath_2(u) := (0, u) . 
\]
Let $\imath_1' : \G_a \to \rB(2, k)$ and 
$\imath_2' : \G_m \to \rB(2, k)$ be the homomorphisms defined by 
\[
 \imath_1' := \jmath \circ \imath_1 
 \qquad \text{ and } \qquad 
 \imath_2' := \jmath \circ \imath_2 .  
\] 
We have the following commutative diagram: 
\[
\xymatrix@R=36pt@C=24pt@M=6pt{
 \G_a \ar[rd]_{\imath_1} \ar@/^12pt/[rrrd]^{\imath_1'} & & \\
        & \G_a \rtimes \G_m \ar[rr]^\jmath_\cong & & \rB(2, k) \\
 \G_m \ar[ru]^{\imath_2} \ar@/_12pt/[rrru]_{\imath_2'} & &
}
\]
Clearly, the image of $\imath_1'$ is $U^+$ and 
the image of $\imath_2'$ is $T$.

\subsection{Homomorphisms from $\G_a$ to $\SL(n, k)$}

\subsubsection{$\cU_n$, $\cU_{[\ell_1, \ldots , \ell_ \nu ]}$}

We denote by $\cU_n$ the set of all homomorphisms $\varphi : \G_a \to \SL(n, k)$ 
such that the following conditions {\rm (i)} and {\rm (ii)} hold true: 
\begin{enumerate}[label = {\rm (\roman*)}]
\item $ (\Inn_P \circ \varphi)(t)$ are upper triangular for all $t \in \G_a$. 

\item The diagonal entries of $ (\Inn_P \circ \varphi)(t)$ are $1$ for all $t \in \G_a$. 
\end{enumerate}

\begin{lem}
Let $\varphi : \G_a \to \SL(n, k)$ be a homomorphism. 
Then there exists a regular matrix $P$ of $\GL(n, k)$ such that 
$\Inn_P \circ \varphi \in \cU_n$. 
\end{lem}

\begin{proof}
See \cite[Lemmas 1.8 and 1.9]{Tanimoto 2019}. 
\end{proof}

An ordered sequence $[\ell_1, \ldots, \ell_\nu]$ of positive integers $\ell_i$ $(1 \leq i \leq \nu)$ 
is said to be an {\it ordered partition} of $n$ if 
$[\ell_1, \ldots, \ell_\nu]$ satisfies $\sum_{i = 1}^\nu \ell_i = n$.

Let $\varphi : \G_a \to \SL(n, k)$ be a homomorphism. 
So, there exists a unique polynomial matrix $A(T)$ of $\Mat(n, k[T])$ so that 
\[
 \varphi(t) = A(t) \qquad \text{ for all \quad $t \in \G_a$} . 
\]
Write $A(T) = ( a_{i, j}(T) )$. 
We say that $\varphi$ {\it has an ordered partision} $[\ell_1, \ldots, \ell_\nu]$ if 
$\varphi$ satisfies 
\[
\{\, i \in \{1 , \ldots, n \} \mid a_{i, \, i + 1}(T) = 0 \,\} 
= 
\{\,
 \ell_1, \; \ell_1 + \ell_2, \; \ldots, \; \ell_1 + \cdots + \ell_\nu 
\,\} . 
\]

We denote by $\cU_{[\ell_1, \ldots, \ell_\nu]}$ the set of all homomorphisms 
$\varphi : \G_a \to \SL(n, k)$ such that 
$\varphi$ has an ordered partition $[\ell_1, \ldots, \ell_\nu]$.

\begin{lem}
The set $\cU_{n}$ is a disjoint union of $\cU_{[\ell_1, \ldots, \ell_\nu]}$, 
where $[\ell_1, \ldots, \ell_\nu]$ are ordered partitions of $n$, 
i.e.,  
\[
 \cU_{n}
  = 
  \bigsqcup_{\text{$[\ell_1, \ldots, \ell_\nu]$ are ordered partitions of $n$}} 
  \cU_{[\ell_1, \ldots, \ell_\nu]} .  
\]
\end{lem}

\begin{proof}
The proof is straightforward. 
\end{proof}

\subsubsection{${^\tau\!}A$, ${^\tau\!}\varphi$}

For any matrix $A = (a_{i, j})$ of $\Mat(n, k)$, we can define a matrix 
${^\tau\!}A = (\alpha_{i, j})$ of $\Mat(n, k)$ as 
\[
 \alpha_{i, j} = a_{n - j + 1, \; n - i + 1} 
 \qquad \text{ for all \quad $1 \leq i, j \leq n$} . 
\]

\begin{lem}
Let $A \in \Mat(n, k)$. 
Let $J = (\jmath_{i, j}) \in \Mat(n, k)$ be the matrix defined by 
\[
 \jmath_{i, j} := 
\left\{
\begin{array}{l @{\qquad} l}
 1 & \text{ if \quad $i + j = n + 1$}, \\
 0 & \text{ otherwise}. 
\end{array}
\right.
\]
Then we have 
\[
 {^\tau\!}A = J \cdot {^t\!}A \cdot J . 
\]
\end{lem}

\begin{proof}
Write ${^\tau\!}A = (\alpha_{i, j})$ and ${^t\!}A = (a'_{i, j})$. 
For all $1 \leq i, j \leq n$, we have 
\begin{align*}
 \alpha_{i, j}
  = \sum_{1 \leq \lambda, \mu \leq n } 
  \jmath_{i, \lambda} \cdot a'_{\lambda, \mu} \cdot \jmath_{\mu,  j} 
  = a'_{n + 1 - i, \; n + 1- j} 
  = a_{n - j + 1, \; n - i + 1} . 
\end{align*}

\end{proof}

\begin{lem}
The following assertions {\rm (1)} and {\rm (2)} hold true: 
\begin{enumerate}[label = {\rm (\arabic*)}]
\item If $A \in \GL(n, k)$, then ${^\tau\!}A \in \GL(n, k)$. 

\item If $A \in \SL(n, k)$, then ${^\tau\!}A \in \SL(n, k)$. 
\end{enumerate} 
\end{lem}

\begin{proof}
(1) See \cite[Lemma 1.6 (5)]{Tanimoto 2019}. 

(2) See Lemma 1.3. 
\end{proof}

Let $\varphi : \G_a \to \SL(n, k)$ be a morphism. 
We can define a morphism ${^\tau\!}\varphi : \G_a \to \SL(n, k)$ as 
\[
 ({^\tau\!}\varphi)(t) := {^\tau\!}( \, \varphi(t) \, ) . 
\]

\begin{lem}
Let $\varphi : \G_a \to \SL(n, k)$ be a homomorphism. 
Then the following assertions {\rm (1)} and {\rm (2)} hold true: 
\begin{enumerate}[label = {\rm (\arabic*)}]
\item The morphism ${^\tau\!}\varphi : \G_a \to \SL(n, k)$ is also a homomorphism. 

\item Let $\phi : \G_a \to \SL(n, k)$ be a homomorphism such that $\varphi$ and $\phi$ are equivalent, i.e., 
$\varphi \sim \phi$. Then ${^\tau\!}\varphi$ and ${^\tau\!}\phi$ are equivalent, i.e., 
${^\tau\!}\varphi \sim {^\tau\!}\phi$. 
\end{enumerate} 
\end{lem}

\begin{proof}
(1) For all $t, t' \in \G_a$, we have 
\begin{align*}
 ({^\tau\!}\varphi)(t + t') 
 & = {^\tau\!}( \, \varphi(t + t') \, )
  = {^\tau\!}( \, \varphi(t' + t) \, ) 
  = {^\tau\!}( \, \varphi(t') \cdot \varphi(t) \, ) 
  = {^\tau\!}( \, \varphi(t) \, ) \cdot {^\tau\!}( \, \varphi(t') \, ) \\
 & = ({^\tau\!}\varphi)(t) \cdot ({^\tau\!}\varphi)(t') . 
\end{align*}

(2) There exists a regular matrix $P$ of $\GL(n, k)$ such that 
$P^{-1} \, \varphi(t) \, P = \phi(t)$ for all $t \in \G_a$. 
So, we have ${^\tau\!} P \cdot {^\tau\!}(\, \varphi(t) \,) \cdot {^\tau\!}( \, P^{-1} \,)  = {^\tau\!}(\, \phi(t) \,)$ for all $t \in \G_a$. 
Letting $Q := {^\tau\!}( \, P^{-1} \,)$, we have $Q^{-1} \cdot {^\tau\!}(\, \varphi(t) \,) \cdot Q = {^\tau\!}(\, \phi(t) \,)$ for all $t \in \G_a$, 
which implies that ${^\tau\!}\varphi$ and ${^\tau\!}\phi$ are equivalent.

\end{proof}

\subsection{Homomorphisms from $\G_m$ to $\SL(n, k)$}

For all integers $d_i$ ($1 \leq i \leq n$) satisfying $\sum_{i = 1}^n d_i = 0$, 
we can define a homomorphism $\omega_{d_1, \ldots, d_n} : \G_m \to \SL(n, k)$ as 
\begin{align*}
\omega_{d_1, \ldots, d_n}(u) 
 := 
\diag( \, u^{d_1}, \, \ldots, \, u^{d_n} \, ) . 
\end{align*}

A homomorphism $\omega : \G_m \to \SL(n, k)$ is said to be {\it antisymmetric} 
if $\omega$ has the form
\[
 \omega = \omega_{d_1, \ldots, d_n} , 
\qquad 
 d_1 \geq  \cdots \geq d_n,  
 \qquad d_i = - d_{n - i + 1} \quad (\, 1 \leq i \leq n \,)  . 
\]

We denote by $\Omega(n)$ the set of all antisymmetric homomorphisms 
$\omega : \G_m \to \SL(n, k)$.

For any homomorphism $\omega : \G_m \to \SL(n, k)$, we can define a homomorphism 
$\omega^\star : \G_m \to \SL(n, k)$ as 
\[
 \omega^\star(u) :=  {^\tau\!}(\, \omega(u)^{-1} \,) 
\]
for all $u \in \G_m$. Clearly, $\omega^{\star \star} = \omega$.

\begin{lem}
Let $\omega \in \Hom(\G_m, \SL(n, k))$. 
Then the following assertions {\rm (1)} and {\rm (2)} hold true: 
\begin{enumerate}[label = {\rm (\arabic*)}]
\item $\omega^{\star\star} = \omega$. 

\item If $\omega \in \Omega(n)$, then $\omega^\star = \omega$. 
\end{enumerate} 
\end{lem} 

\begin{proof}
The proofs of assertions (1) and (2) are straightforward. 
\end{proof}

\begin{lem}
Let $\omega \in \Omega(n)$ and let $P$ be a regular matrix of $\GL(n, k)$ 
such that $\Inn_P \circ \omega \in \Omega(n)$. 
Let $Q := {^\tau\!}(P^{-1})$. 
Then the following assertions {\rm (1)} and {\rm (2)} hold true: 
\begin{enumerate}[label = {\rm (\arabic*)}]
\item $\Inn_P \circ \omega = \omega$. 

\item $\Inn_Q \circ \omega = \omega$. 
\end{enumerate} 
\end{lem}

\begin{proof} 
(1) See \cite[Lemma 2.2]{Tanimoto 2022}.

(2) Since $P^{-1} \cdot \omega(u) \cdot P = \omega(u)$ for all $u \in \G_m$, 
we have $P^{-1} \cdot \omega(u^{-1}) \cdot P = \omega(u^{-1})$ for all $u \in \G_m$, 
which implies ${^\tau\!}P \cdot {^\tau\!}( \omega(u^{-1}) )\cdot {^\tau\!}(P^{-1}) = {^\tau\!}(\omega(u^{-1}))$ for all $u \in \G_m$. 
Using Lemma 1.6 (2), we have $\Inn_Q \circ \omega = \omega$. 
\end{proof}

\subsection{Homomorphisms from $\rB(2, k)$ to $\SL(n, k)$}

\subsubsection{Morphisms $\varphi_\psi$, $\omega_\psi$, $\psi_{\varphi, \, \omega}$}

Given a morphism $\psi : \rB(2, k) \to \SL(n, k)$, 
we can define a morphism $\varphi_\psi : \G_a \to \SL(n, k)$ as 
\[
 \varphi_\psi(t) 
 := 
\psi 
\left(
\begin{array}{c c}
 1 & t \\
 0 & 1 
\end{array} 
\right) , 
\]
and can also define a morphism $\omega_\psi : \G_m \to \SL(n, k)$ as 
\[
 \omega_\psi(u) 
 := 
\psi
\left(
\begin{array}{c c}
 u & 0 \\
 0 & u^{-1} 
\end{array}
\right) . 
\]
Clearly, if $\psi$ is a homomorphism, then $\varphi_\psi$ and $\omega_\psi$ 
are homomorphisms.

Conversely, given morphisms $\varphi : \G_a \to \SL(n, k)$ and $\omega : \G_m \to \SL(n, k)$, 
we can define a morphism $\psi_{\varphi, \, \omega} : \G_a \rtimes \G_m \to \SL(n, k)$ as 
\[
 \psi_{\varphi, \, \omega}(t, u) 
 := 
\varphi(t) 
\cdot 
\omega(u)  . 
\]

%
%
%


\begin{lem}
Let $\psi : \rB(2, k) \to \SL(n, k)$ be a homomorphism. 
Then the following assertions {\rm (1)} and {\rm (2)} hold true: 
\begin{enumerate}[label = {\rm (\arabic*)}]
\item $\psi \circ \jmath = \psi_{\varphi_\psi, \, \omega_\psi}$. 

\item There exist unique homomorphisms $\varphi : \G_a \to \SL(n, k)$ and 
$\omega : \G_m \to \SL(n, k)$ such that $\psi \circ \jmath = \psi_{\varphi, \, \omega}$. 
\end{enumerate} 
\end{lem}

\begin{proof}
(1) Let $\varphi : \G_a \to \SL(n, k)$ be the homomorphism defined by $\varphi := \varphi_\psi$. 
Let $\omega : \G_m \to \SL(n, k)$ be the homomorphism defined by $\omega = \omega_\psi$. 
Then
\begin{align*}
 (\psi \circ \jmath)(t, u)
 & = 
\psi 
\left( 
\begin{array}{c c}
 1 & t \\
 0 & 1
\end{array}
\right)
\cdot 
\psi  
\left(
\begin{array}{c c}
 u & 0 \\
 0 & u^{-1} 
\end{array}
\right) 
 = 
\varphi(t) \cdot \omega(u) 
 = \psi_{\varphi, \, \omega}(t, u) 
\end{align*}
for all $(t, u) \in \G_a \rtimes \G_m$. 

(2) We have only to show the uniqueness of $\varphi$ and $\psi$, 
which follows from the fact that 
$\psi_{\varphi, \, \omega}(t, 1) = \varphi(t)$ for all $t \in \G_a$ 
and 
$\psi_{\varphi, \, \omega}(0, u) = \omega(u)$ for all $u \in \G_m$. 
\end{proof}

\begin{lem}
Let $\varphi : \G_a \to \SL(n, k)$ and $\omega : \G_m \to \SL(n, k)$ be morphisms. 
Then the morphism $\psi_{\varphi, \, \omega} : \G_a \rtimes \G_m \to \SL(n, k)$ is a homomorphism 
if and only if the following conditions {\rm (1)}, {\rm (2)}, {\rm (3)} hold true: 
\begin{enumerate}[label = {\rm (\arabic*)}]
\item $\varphi$ is a homomorphism. 

\item $\omega$ is a homomorphism. 

\item $\omega(u) \, \varphi(t) \, \omega(u)^{-1} = \varphi(u^2 \, t)$ for all $(t, u) \in \G_a \rtimes \G_m$. 
\end{enumerate} 
\end{lem}

\begin{proof}
We first prove the implication $(\Longrightarrow)$. 
Assume $\psi_{\varphi, \, \omega}$ is a homomorphism. 
Since $\varphi(t) = \psi_{\varphi, \, \omega}(t, 1)$ for all $t \in \G_a$, we have 
\begin{align*}
 \varphi(t + t') 
 & = \psi_{\varphi, \, \omega}(t + t', 1) 
 = \psi_{\varphi, \, \omega}((t, 1) \cdot (t', 1))
 = \psi_{\varphi, \, \omega}(t, 1) \cdot \psi_{\varphi, \, \omega}(t', 1)
 = \varphi(t) \cdot \varphi(t')
\end{align*}
for all $t, t' \in \G_a$. So, assertion (1) holds true. 
Since $\omega(u) = \psi_{\varphi, \, \omega}(0, u)$ for all $u \in \G_m$, we have 
\begin{align*}
 \omega(u \, u') 
 & = \psi_{\varphi, \, \omega}(0, \, u \, u') 
 = \psi_{\varphi, \, \omega}((0, u) \cdot (0, u')) 
 = \psi_{\varphi, \, \omega}(0, u) \cdot \psi_{\varphi, \, \omega}(0, u') 
 = \omega(u) \cdot \omega(u')  
\end{align*}
for all $u, u' \in \G_m$. So, assertion (2) holds true. 
Since 
\[
 \psi_{\varphi, \, \omega}((0, u) \cdot (t', 1)) = \psi_{\varphi, \, \omega}(0, u) \cdot \psi_{\varphi, \, \omega}(t', 1)  
\]
for all $u \in \G_m$ and $t' \in \G_a$, we have
\[
  \varphi(u^2 \, t') \cdot \omega(u) = \omega(u) \cdot \varphi(t') 
\]
for all $u \in \G_m$ and $t' \in \G_a$, which implies assertion (3) holds true.

We next prove the implication $(\Longleftarrow)$. 
Assume that conditions (1), (2), (3) hold true. 
We have 
\begin{align*}
 \psi_{\varphi, \, \omega}((t, u) \cdot (t', u')) 
 & = \psi_{\varphi, \, \omega}(t + u^2 \, t', \, u \, u') \\
 & =  \varphi(t + u^2 \, t') \cdot \omega(u \, u') \\
 & = \varphi(t) \cdot \varphi(u^2 \, t') \cdot \omega(u) \cdot \omega(u') \\
 & = \varphi(t) \cdot \omega(u) \cdot \varphi(t') \cdot \omega(u)^{-1} \cdot \omega(u) \cdot \omega(u') \\
 & = \varphi(t) \cdot \omega(u) \cdot \varphi(t') \cdot \omega(u') \\
 & = \psi_{\varphi, \, \omega}(t, u) \cdot \psi_{\varphi, \, \omega}(t', u') 
\end{align*}
for all $(t, u), (t', u') \in \G_a \rtimes \G_m$. 
\end{proof}

\subsubsection{Antisymmetric homomorphisms from $\rB(2, k)$ to $\SL(n, k)$}

Let $\psi : \rB(2, k) \to \SL(n, k)$ be a homomorphism. 
We say that $\psi$ is an {\it antisymmetric} if $\omega_\psi \in \Omega(n)$. 

We denote by $\Hom^a(\rB(2, k), \SL(n, k))$ the set of all antisymmetric homomorphisms 
$\psi : \rB(2, k) \to \SL(n, k)$. 

\begin{lem}
Let $\psi \in \Hom^a(\rB(2, k), \SL(n, k))$. 
Then we have $\varphi_\psi \in \cU_n$. 
\end{lem}

\begin{proof}
See \cite[Lemma 2.5]{Tanimoto 2022}. 
\end{proof}

\subsubsection{$({^\tau\!}\varphi, \, \omega^\star)$}

\begin{lem}
Let $\varphi : \G_a \to \SL(n, k)$ and $\omega : \G_m \to \SL(n, k)$ be morphisms 
such that $\psi_{\varphi, \, \omega} : \G_a \rtimes \G_m \to \SL(n, k)$ is a homomorphism. 
Then the morphism $\psi_{\,{^\tau\!}\varphi,\, \omega^\star} : \G_a \rtimes \G_m \to \SL(n, k)$ 
is also a homomorphism. 
In particular, if $\omega \in \Omega(n)$, the morphism $\psi_{\,{^\tau\!}\varphi,\, \omega} : \G_a \rtimes \G_m \to \SL(n, k)$ 
is also a homomorphism. 
\end{lem}


\begin{proof}
By Lemma 1.9, we have the following: 
\begin{enumerate}[label = {\rm (\arabic*)}]
\item $\varphi$ is a homomorphism. 

\item $\omega$ is a homomorphism. 

\item $\omega(u) \, \varphi(t) \, \omega(u)^{-1} = \varphi(u^2 \, t)$ for all $(t, u) \in \G_a \rtimes \G_m$. 
\end{enumerate} 
Note the following: 
\begin{enumerate}
\item [\rm $(1)'$] ${^\tau\!}\varphi$ is a homomorphism. 

\item[\rm $(2)'$] $\omega^\star$ is a homomorphism.

\item[\rm $(3)'$] $\omega^\star(u) \cdot ({^\tau\!}\varphi)(t) \cdot \omega^\star(u)^{-1} = ({^\tau\!}\varphi)(u^2 \, t)$ 
for all $(t, u) \in \G_a \rtimes \G_m$. 
\end{enumerate}
For $(1)'$, see the above (1) and Lemma 1.5. 
For $(2)'$, see the above (2). 
For $(3)'$, using the above (3), we have ${^\tau\!}( \, \omega(u)^{-1} \,)  \cdot {^\tau\!}(\, \varphi(t) \,) \cdot {^\tau}\!( \, \omega(u) \,)  = {^\tau\!}( \, \varphi(u^2 \, t) \, )$ 
for all $(t, u) \in \G_a \rtimes \G_m$. 
By $(1)'$, $(2)'$, $(3)'$, we have the conclusion (see Lemma 1.9). 

In particular if $\omega \in \Omega(n)$, then $\psi_{\,{^\tau\!}\varphi,\, \omega^\star} = 
\psi_{\,{^\tau\!}\varphi,\, \omega}$ (see Lemma 1.6 (2)). 
\end{proof}

\subsubsection{$\psi^\star$, $({^\tau\!}\varphi \circ \inv, \, \omega^\star)$} 


Let $\psi : \rB(2, k) \to \SL(n, k)$ be a morphism. 
We can define a morphism $\psi^\star: \rB(2, k) \to \SL(n, k)$ as 
\begin{align*}
 \psi^\star(A) := {^\tau\!}(\, \psi(A^{-1}) \,) . 
\end{align*}

\begin{lem}
The following assertions {\rm (1)}, {\rm (2)}, {\rm (3)} hold true:
\begin{enumerate}[label = {\rm (\arabic*)}]
\item $\omega_{\psi^\star} = (\omega_\psi)^\star$. 

\item $\psi^{\star\star} = \psi$. 

\item If $\psi$ is a homomorpism, then $\psi^\star$ is also a homomorphism. 
\end{enumerate} 
\end{lem}

\begin{proof} 
(1) For any $u \in \G_m$, we have 
\begin{align*}
 \omega_{\psi^\star}(u) 
 & = 
 \psi^\star
 \left(
 \begin{array}{c c}
  u & 0 \\
  0 & u^{-1} 
 \end{array}
 \right) 
 =
  \prescript{\tau\!\!}{}
  {
  \left(  \, \psi
 \left( 
 \begin{array}{c c}
  u^{-1} & 0 \\
  0 & u 
 \end{array}
 \right) 
 \, \right)
 }
 =  {^\tau\!} \bigl(\, \omega_\psi(u^{-1} ) \, \bigr)
 = (\omega_\psi )^\star(u) .  
\end{align*}

(2), (3) The proofs are straightforward. 
\end{proof}

We can define an isomorphism $\inv : \G_a \to \G_a$ as 
\[
 \inv(t) := - t . 
\]

\begin{lem}
Let $\psi : \rB(2, k) \to \SL(n, k)$ be a homomorphism and write 
$\psi \circ \jmath = \psi_{\varphi, \, \omega}$ for some 
$(\varphi, \omega) \in \Hom(\G_a, \SL(n, k)) \times \Hom(\G_m, \SL(n, k))$. 
Then 
\[
 \psi^\star \circ \jmath
  = \psi_{({^\tau\!}\varphi ) \, \circ \, \inv, \; \omega^\star} . 
\]
In particular, if $\omega \in \Omega(n)$, then 
\[
 \psi^\star \circ \jmath
  = \psi_{({^\tau\!}\varphi ) \, \circ \, \inv, \; \omega} . 
\]
\end{lem}

\begin{proof}
For all $t \in \G_a$ and $u \in \G_m$, we have 
\begin{align*}
 \psi^\star 
 \left( 
 \begin{array}{c c}
  1 & t \\
  0 & 1 
 \end{array}
 \right)
 & = 
 {^\tau} \left( 
 \psi \left( 
 \begin{array}{r r} 
  1 & - t \\
   0 & 1
 \end{array}
 \right) 
 \right) 
 = {^\tau\!}( \, \varphi( - t) \,) 
 = ({^\tau\!}\varphi)(-t) 
 = \bigl(\, ({^\tau\!}\varphi) \circ \inv \, \bigr)(t) , \\ 
 \psi^\star 
 \left(
 \begin{array}{c c}
  u & 0 \\
   0 & u^{-1}  
 \end{array}
 \right) 
 & = 
 {^\tau} \left( 
 \psi \left( 
 \begin{array}{c c}
  u^{-1} & 0 \\
   0 & u 
 \end{array}
 \right)
 \right) 
 = {^\tau\!}\bigl( \,  \omega(u^{-1})   \, \bigr) 
 = \omega^\star(u)  . 
\end{align*}
Thus, for all $(t, u) \in \G_a \rtimes \G_m$, we have  
\begin{align*}
( \psi^\star \circ \jmath ) (t , u)
& = 
  \psi^\star 
 \left( 
 \begin{array}{c c}
  1 & t \\
  0 & 1 
 \end{array}
 \right) 
 \cdot 
  \psi^\star 
 \left(
 \begin{array}{c c}
  u & 0 \\
   0 & u^{-1}  
 \end{array}
 \right) 
 = 
 \psi_{({^\tau\!}\varphi ) \, \circ \, \inv, \; \omega^\star}(t, u) . 
\end{align*} 
If $\omega \in \Omega(n)$, we have $\omega = \omega^\star$ (see Lemma 1.6 (2)) 
and thereby have 
\[
 \psi^\star \circ \jmath
  = \psi_{({^\tau\!}\varphi ) \, \circ \, \inv, \; \omega} . 
\]
\end{proof}

\subsubsection{$(\varphi \circ \inv, \, \omega)$}


We can define an isomorphism $r_{\rB(2, k)} : \rB(2, k) \to \rB(2, k)$ as 
\[
r_{\rB(2, k)}
\left(
\begin{array}{c c}
 a & b \\
 0 & d
\end{array}
\right) 
 :=
 \left(
\begin{array}{r r}
 a & - b \\
 0 & d
\end{array}
\right) . 
\]
Clearly, $r_{\rB(2, k)}^2 = \id_{\rB(2, k)}$.

\begin{lem}
Let $\psi : \rB(2, k) \to \SL(n, k)$ be a homomorphism and write 
$\psi \circ \jmath = \psi_{\varphi, \, \omega}$, where 
$(\varphi, \omega) \in \Hom(\G_a, \SL(n, k)) \times \Hom(\G_m, \SL(n, k))$. 
Then 
\[
 \psi \circ r_{\rB(2, k)} \circ \jmath = \psi_{\varphi \, \circ \, \inv, \; \omega} . 
\]
\end{lem}

\begin{proof}
For all $t \in \G_a$ and $u \in \G_m$, we have 
\begin{align*}
(\psi \circ r_{\rB(2, k)} )
\left(
\begin{array}{c c}
 1 & t \\
 0 & 1 
\end{array}
\right)
 & = 
\psi
\left(
\begin{array}{r r}
 1 & - t \\
 0 & 1 
\end{array}
\right)
= \varphi( - t )
= (\varphi \circ \inv)(t) , \\
 (\psi \circ r_{\rB(2, k)})
 \left( 
 \begin{array}{c c}
  u & 0 \\
   0 & u^{- 1}
 \end{array}
 \right) 
 & = 
 \psi 
 \left( 
 \begin{array}{c c}
  u & 0 \\
   0 & u^{- 1}
 \end{array}
 \right) 
 = 
 \omega(u) .  
\end{align*}
Thus, for all $(t, u) \in \G_a \rtimes \G_m$, we have 
\begin{align*}
 (\psi \circ r_{\rB(2, k)} \circ \jmath) (t, u)
 & = 
 (\psi \circ r_{\rB(2, k)} )
\left(
\begin{array}{c c}
 1 & t \\
 0 & 1 
\end{array}
\right) 
\cdot 
 (\psi \circ r_{\rB(2, k)})
 \left( 
 \begin{array}{c c}
  u & 0 \\
   0 & u^{- 1}
 \end{array}
 \right)  \\
 & = (\varphi \circ \inv)(t) \cdot \omega(u) 
 =  \psi_{\varphi \, \circ \, \inv, \; \omega} (t, u) . 
\end{align*}
\end{proof}

\subsubsection{$({^\tau\!}\varphi, \, \omega^\star)$}


\begin{lem}
Let $\psi : \rB(2, k) \to \SL(n, k)$ be a homomorphism and write 
$\psi \circ \jmath = \psi_{\varphi, \, \omega}$, where 
$(\varphi, \omega) \in \Hom(\G_a, \SL(n, k)) \times \Hom(\G_m, \SL(n, k))$. 
Then 
\[
 \psi^\star \circ r_{\rB(2, k)} \circ \jmath = \psi_{\,{^\tau\!}\varphi, \; \omega^\star} . 
\]
In particular, if $\omega \in \Omega(n)$, then 
\[
 \psi^\star \circ r_{\rB(2, k)} \circ \jmath = \psi_{\, {^\tau\!}\varphi,  \; \omega} . 
\]
\end{lem}

\begin{proof}
We have 
\[
 \psi^\star \circ \jmath
  = \psi_{({^\tau\!}\varphi ) \, \circ \, \inv, \; \omega^\star} . 
\]
So, 
\[
 \psi^\star \circ r_{\rB(2, k)} \circ \jmath 
 = \psi_{({^\tau\!}\varphi) \, \circ \, \inv \, \circ \, \inv , \; \omega^\star} 
 = \psi_{\, {^\tau\!}\varphi , \; \omega^\star} . 
\]
If $\omega \in \Omega(n)$, we have $\omega = \omega^\star$ (see Lemma 1.6 (2)) 
and thereby have 
\[
 \psi^\star \circ r_{\rB(2, k)} \circ \jmath 
 = \psi_{\, {^\tau\!}\varphi, \; \omega} . 
\]
\end{proof}

\subsubsection{Equivalence of pairs of 
$\Hom(\G_a, \SL(n, k)) \times \Hom(\G_m\, \SL(n, k))$}

Let $(\varphi, \omega), (\varphi^*, \omega^*) 
\in \Hom(\G_a, \SL(n, k) ) \times \Hom(\G_m, \SL(n, k))$. 
The pairs $(\varphi, \omega)$ and $(\varphi^*, \omega^*)$ are {\it equivalent} if 
there exists a reguar matrix $P$ of $\GL(n, k)$ such that 
\[
 \Inn_P \circ \psi_{\varphi, \, \omega} = \psi_{\varphi^*, \, \omega^*} . 
\]
If the pairs $(\varphi, \omega)$ and $(\varphi^*, \omega^*)$ are equivalent, 
we write $(\varphi, \omega) \sim (\varphi^*, \omega^*)$.

\begin{lem}
Let $(\varphi, \omega) \in \Hom(\G_a, \SL(n, k) ) \times \Hom(\G_m, \SL(n, k))$. 
Let $P$ be a regular matrix of $\GL(n, k)$. 
Then we have 
\[
 \Inn_P \circ \psi_{\varphi, \, \omega} 
 = 
 \psi_{\varphi^*, \, \omega^*} , 
\]
where 
\[
 \varphi^* := \Inn_P \circ \varphi, \qquad 
 \omega^* := \Inn_P \circ \omega . 
\]
\end{lem}

\begin{proof}
The proof is straightforward. 
\end{proof}

\begin{lem}
Let $(\varphi, \omega), (\varphi^*, \omega^*) \in \Hom(\G_a, \SL(n, k)) \times \Omega(n)$ such 
that $\psi_{\varphi, \, \omega}$ and $\psi_{\varphi^*, \, \omega^*}$ are homomorphisms. 
Let $P$ be a regular matrix of $\GL(n, k)$ such that $\Inn_P \circ \psi_{\varphi, \, \omega} = \psi_{\varphi^*, \, \omega^*}$. 
Let  $Q := {^\tau\!}(P^{-1})$. 
Then the following assertions {\rm (1)} and {\rm (2)} hold true: 
\begin{enumerate}[label = {\rm (\arabic*)}]
\item $\omega = \omega^*$. 

\item $\Inn_Q \circ \psi_{\, {^\tau\!}\varphi, \, \omega} =  \psi_{\, {^\tau\!}(\varphi^*), \; \omega}$. 
\end{enumerate} 
\end{lem}

\begin{proof}
Since $\omega^* = \Inn_P \circ \omega$ and $\omega, \omega^* \in \Omega(n)$, 
we have $\omega^* = \omega$ (see Lemma 1.7 (1)). 
Assertion (1) is proved. 
Since $\varphi^* = \Inn_P \circ \varphi$, we have $P^{-1} \cdot \varphi(t) \cdot P = \varphi^*(t)$ for all $t \in \G_a$. 
So, ${^\tau\!}P \cdot {^\tau\!}( \varphi(t) )\cdot {^\tau\!}(P^{-1}) = {^\tau\!}(\varphi^*(t))$ for all $t \in \G_a$, 
which implies $\Inn_Q \circ  {^\tau\!}\varphi = {^\tau\!}(\varphi^*)$. 
By Lemma 1.7 (2), we have $\Inn_Q \circ \psi_{\, {^\tau\!}\varphi, \, \omega} =  \psi_{\, {^\tau\!}(\varphi^*), \; \omega}$. 
Assertion (2) is proved. 
\end{proof}

\subsection{Homomorphisms from $\SL(2, k)$ to $\SL(n, k)$}

\subsubsection{Morphisms $\varphi_\sigma$, $\omega_\sigma$, $\varphi_\sigma^-$}

Let $\sigma : \SL(2, k) \to \SL(n, k)$ be a morphism. 
We can define morphisms 
$\varphi_\sigma : \G_a \to \SL(n, k)$, 
$\omega_\sigma : \G_m \to \SL(n, k)$, 
$\varphi_\sigma^- : \G_a \to \SL(n, k)$, as follows: 
\begin{align*}
\varphi_\sigma (t) 
  := 
 \sigma 
 \left(
 \begin{array}{c c}
  1 & t \\
  0 & 1
 \end{array}
 \right) , \qquad 
\omega_\sigma(u) 
  := 
 \sigma 
 \left(
 \begin{array}{c c}
  u & 0 \\
  0 & u^{-1}
 \end{array}
 \right) , \qquad 
\varphi_\sigma^-(s) 
 := 
  \sigma 
 \left(
 \begin{array}{c c}
  1 & 0 \\
  s & 1
 \end{array}
 \right) . 
\end{align*}
Clearly, if $\sigma : \SL(2, k) \to \SL(n, k)$ is a homomorphism, 
the morphisms $\varphi_\sigma$, $\omega_\sigma$, $\varphi^-_\sigma$ are homomorphisms.

\begin{lem}
The following assertions {\rm (1)} and {\rm (2)} hold true: 
\begin{enumerate}[label = {\rm (\arabic*)}]
\item We have 
\[
\left(
\begin{array}{c c}
 1 & t \\
 0 & 1
\end{array}
\right) 
\left(
\begin{array}{c c}
 1 & 0 \\
 s & 1 
\end{array}
\right) 
= 
\left(
\begin{array}{c c}
 1 & 0 \\
 \frac{s}{1 + t \, s} & 1 
\end{array}
\right) 
\left(
\begin{array}{c c}
 1 +  t \,s & 0 \\
 0 & \frac{1}{1 + t \, s}
\end{array} 
\right) 
\left(
\begin{array}{c c}
 1 & \frac{t}{ 1 + t \, s} \\
 0 & 1
\end{array}
\right)
\]
for all $t, s \in k$ with $1 + t \, s \ne 0$.

\item Let $\sigma : \SL(2, k) \to \SL(n, k)$ be a homomorphism. 
Let $\varphi_\sigma : \G_a \to \SL(n, k)$, 
$\varphi_\sigma^- : \G_a \to \SL(n, k)$ and 
$\omega_\sigma : \G_m \to \SL(n, k)$ be 
the induced homomorphsims from $\sigma$. 
Then we have  
\[
 \varphi_\sigma(t) \, \varphi_\sigma^-(s) 
 = 
\varphi_\sigma^-\left( \frac{s}{1 + t \, s} \right) 
\, 
\omega_\sigma(1 + t \, s) 
\;
\varphi_\sigma \left( \frac{t}{1 + t \, s} \right)
\]
for all $t, s \in k$ with $1 + t \, s \ne 0$. 
\end{enumerate} 
\end{lem}

\begin{proof}
The proofs of assertions (1) and (2) are straightforward. 
\end{proof}

\begin{lem}
Let $(\varphi, \omega, \varphi^-) \in 
\Hom(\G_a, \SL(n, k)) \times \Omega(n) \times \Hom(\G_a, \SL(n, k))$.  
Then the following conditions {\rm (1)} and {\rm (2)} are equivalent: 
\begin{enumerate}[label = {\rm (\arabic*)}]

\item 
$
 \varphi(t) \, \varphi^-(s) 
 = 
\varphi^-\left( \frac{s}{1 + t \, s} \right) 
\, 
\omega(1 + t \, s) 
\;
\varphi \left( \frac{t}{1 + t \, s} \right)
$
for all $t, s \in k$ with $1 + t \, s \ne 0$.

\item 
$
({^\tau\!}\varphi)(t) \, ({^\tau\!}\varphi^-)(s) 
 = 
 ({^\tau\!}\varphi^-)\left( \frac{s}{1 + t \, s} \right) 
\, 
\omega(1 + t \, s) 
\;
({^\tau\!}\varphi) \left( \frac{t}{1 + t \, s} \right)
$
for all $t, s \in k$ with $1 + t \, s \ne 0$.

\end{enumerate} 
\end{lem}

\begin{proof}
We first prove the implication (1) $\Longrightarrow$ (2). 
We have 
\[
{^\tau}
\bigl(\; 
(  \varphi(t) \, \varphi^-(s) )^{-1}
\; \bigr)
 = 
\prescript{\tau\!}{}{
\left(\; 
\left( 
\varphi^-\left( \frac{s}{1 + t \, s} \right) 
\, 
\omega(1 + t \, s) 
\;
\varphi \left( \frac{t}{1 + t \, s} \right) 
\right)^{-1}
\; \right)
}
\]
for all $t, s \in k$ with $1 + t \, s \ne 0$. 
Thus 
\[
({^\tau\!}\varphi)( -t ) \, ({^\tau\!}\varphi^-)(- s) 
 = 
({^\tau\!}\varphi^-) \left( - \frac{s}{1 + t \,s} \right) 
\, 
\omega(1 + t \, s) 
\;
({^\tau\!}\varphi) \left( - \frac{t}{1 + t \, s} \right)
\]
for all $t, s \in k$ with $1 + t \, s \ne 0$, which implies condition (2) holds true. 

Using the implication (1) $\Longrightarrow$ (2), we can show the implication 
(2) $\Longrightarrow$ (1).  
\end{proof}

\subsubsection{Antisymmetric homomorphisms from $\SL(2, k)$ to $\SL(n, k)$}

Let $\sigma : \SL(2, k) \to \SL(n, k)$ be a homomorphism. 
We say that $\sigma$ is {\it antisymmetric} if $\varphi_\sigma \in \Omega(n)$. 

We denote by $\Hom^a(\SL(2, k), \SL(n, k))$ denote the set of all  
antisymmetric homomorphisms from $\SL(2, k)$ to $\SL(n, k)$. 

\begin{lem}
The following assertions {\rm (1)} and {\rm (2)} hold true: 
\begin{enumerate}[label = {\rm (\arabic*)}]
\item For any homomorphism $ \sigma : \SL(2, k) \to \SL(n, k)$, 
there exists an antisymmetric homomorphism $\SL(2, k) \to \SL(n, k)$ so that
$\sigma$ and $\sigma^*$ are equivalent, i.e., $\sigma \sim \sigma^*$.

\item We have the following commutative diagram: 
\[
\xymatrix@R=36pt@C=36pt@M=6pt{
 \Hom(\SL(2, k), \SL(n, k)) \ar@{->>}[r] & \Hom(\SL(2, k), \SL(n, k))  / \sim \\
 \Hom^a(\SL(2, k), \SL(n, k)) \ar@{->>}[r] 
 \ar@{^(->}[u] & \Hom^a(\SL(2, k), \SL(n, k))  / \sim  \ar@{=}[u]
}
\]
\end{enumerate} 
\end{lem}

\begin{proof}
(1) See \cite[Lemma 2.6]{Tanimoto 2022}. 

(2) The proof is straightforward. 
\end{proof}

\subsubsection{$\sigma^\star$}

Let $\sigma : \SL(2, k) \to \SL(n, k)$ be a morphism. 
We can define a morphism $\sigma^\star : \SL(2, k) \to \SL(n, k)$ as 
\[
 \sigma^\star(A) 
  := 
  {^\tau\!} (\, \sigma(A^{-1}) \,) . 
\]

\begin{lem}
Let $\sigma : \SL(2, k) \to \SL(n, k)$ be a morphism. 
Then the following assertions {\rm (1)} and {\rm (2)} hold true: 
\begin{enumerate}[label = {\rm (\arabic*)}]
\item $\sigma^{\star \star} = \sigma$. 

\item If $\sigma$ is a homomorphism, then $\sigma^\star$ is a homomorphism. 
\end{enumerate} 
\end{lem}

\begin{proof}
The proofs of assertions (1) and (2) are straightforward. 
\end{proof}

\subsubsection{${^\tau\!} \sigma^\tau $}

Let $\sigma : \SL(2, k) \to \SL(n, k)$ be a morphism. 
We can define a morphism ${^\tau\!} \sigma^\tau : \SL(2, k) \to \SL(n, k)$ as 
\[
 ( {^\tau\!} \sigma^\tau )(A) = {^\tau\!} \bigl(\,  \sigma({^\tau\!} A)  \, \bigr) . 
\]

\begin{lem}
Let $\sigma : \SL(2, k) \to \SL(n, k)$ be a morphism. 
Then the following assertions {\rm (1)} and {\rm (2)} hold true: 
\begin{enumerate}[label = {\rm (\arabic*)}]
\item ${^\tau\!} ( {^\tau\!} \sigma^\tau )^\tau = \sigma$. 

\item If $\sigma$ is a homomorphism, then ${^\tau\!} \sigma^\tau$ is a homomorphism. 
\end{enumerate} 
\end{lem}

\begin{proof}
(1) ${^\tau\!} ( {^\tau\!} \sigma^\tau )^\tau(A)
 = {^\tau\!} \bigl(\,   ( {^\tau\!} \sigma^\tau )({^\tau\!} A)  \, \bigr)
 =  {^\tau\!} \Bigl(\,   {^\tau\!} \bigl( \sigma( {^\tau\!}({^\tau\!} A))  \bigr) \, \Bigr)
  = \sigma(A)$. 

(2) is clear. 
\end{proof}

\begin{lem}
Let $\sigma : \SL(2, k) \to \SL(n, k)$ be a homomorphism. 
Then we have 
\[
 \varphi_{\, {^\tau\!} \sigma^\tau} = {^\tau\!}( \varphi_\sigma ) , \qquad 
 \omega_{\, {^\tau\!} \sigma^\tau} = (\omega_\sigma)^\star, \qquad 
 \varphi_{\, {^\tau\!} \sigma^\tau}^- = {^\tau\!}( \varphi_\sigma^- ) . 
\]
In particular if $\omega_\sigma \in \Omega(n)$, we have 
\[
 \varphi_{\, {^\tau\!} \sigma^\tau} = {^\tau\!}( \varphi_\sigma ) , \qquad 
 \omega_{\, {^\tau\!} \sigma^\tau} = \omega_\sigma, \qquad 
 \varphi_{\, {^\tau\!} \sigma^\tau}^- = {^\tau\!}( \varphi_\sigma^- ) . 
\]
\end{lem}

\begin{proof} 
The proof is straightforward. 
\end{proof}

\begin{lem}
Let $\sigma_i : \SL(2, k) \to \SL(n, k)$ $(i = 1, 2)$ be homomorphisms. 
Assume that $\sigma_1$ and $\sigma_2$ are equivalent. 
Then ${^\tau}(\sigma_1)^\tau$ and ${^\tau}(\sigma_2)^\tau$ are equivalent. 
\end{lem}

\begin{proof}
There exists a regular matrix $P$ of $\GL(n, k)$ such that 
$P^{-1} \, \sigma_1(A) \, P = \sigma_2(A)$ for all $A \in \SL(2, k)$. 
Let $Q : = {^\tau\!}(P^{-1})$. 
For all $ B \in \SL(2, k)$, 
we have 
\[ Q^{-1} \cdot {^\tau\!}(\sigma_1)^\tau (B) \cdot Q
 = {^\tau\!}P \cdot {^\tau\!}(\sigma_1)({^\tau\!}B) \cdot {^\tau\!}(P^{-1})
 = {^\tau\!} (\, P^{-1} \cdot \sigma_1({^\tau\!}B) \cdot P \,) 
 = (\, {^\tau\!}(\sigma_2)^\tau \,)(B) . 
\] 
\end{proof}

\section{Extending homomorphisms}

\subsection{Extending homomorphisms $\G_a \to \SL(n, k)$ to $\rB(2, k) \to \SL(n, k)$}

%

A homomorphism $\varphi : \G_a \to \SL(n, k)$ is said to be {\it $\rB(2, k)$-fundamental} if 
there exists a homomorphism $\psi : \rB(2, k) \to \SL(n, k)$ such that $\psi \circ \imath_1'  =\varphi$, i.e., 
the following diagram is commutative: 
\[
\xymatrix@M=6pt{
 \G_a \ar[r]^(.4)\varphi  \ar@{^(->}[d]_{\imath_1'}
  & \SL(n, k) \\
  \rB(2, k) \ar[ru]_\psi  
 & 
}
\]

\begin{lem}
Let $\varphi : \G_a \to \SL(n, k)$ be a homomorphism. 
Then the following conditions {\rm (1)} and {\rm (2)} are equivalent: 
\begin{enumerate}[label = {\rm (\arabic*)}]
\item $\varphi$ is a $\rB(2, k)$-fundamental homomorphism. 

\item There exists a homomorphism $\omega : \G_m \to \SL(n, k)$ such that 
$\psi_{\varphi, \, \omega} : \rB(2, k) \to \SL(n, k)$ is a homomorphism. 
\end{enumerate}
\end{lem}

\begin{proof}
We first prove the implication (1) $\Longrightarrow$ (2). 
There exists a homomorphism $\psi : \rB(2, k) \to \SL(n, k)$ such that 
$\psi \circ \imath_1' = \varphi$. 
So, $\psi \circ \jmath \circ \imath_1 = \varphi$. 
Since $\psi \circ \jmath = \psi_{\varphi_\psi, \, \omega_\psi}$, 
we have $\varphi_\psi = \varphi$, which implies that condition (2) holds true. 

We next prove the impliation (2) $\Longrightarrow$ (1). 
Let $\psi : \rB(2, k) \to \SL(n, k)$ be the homomorphism defined by 
$\psi := \psi_{\varphi, \, \omega} \circ \jmath^{-1}$. 
Then we have $\psi \circ \imath_1' = \psi \circ \jmath \circ \imath_1 = \varphi$, which implies that $\varphi$ is a $\rB(2, k)$-fundamental homomorphism. 
\end{proof}

\begin{lem}
Let $\varphi : \G_a \to \SL(n, k)$ be a $\rB(2, k)$-fundamental homomorphism. 
Then the following assertions {\rm (1)} and {\rm (2)} hold true: 
\begin{enumerate}[label = {\rm (\arabic*)}]
\item 
Let $\phi : \G_a \to \SL(n, k)$ be a homomorphim such that $\varphi$ and $\phi$ are equivalent, i.e., $\varphi \sim \phi$. 
Then $\phi$ is also a $\rB(2, k)$-fundamental homomorphism. 

\item The morphism ${^\tau\!}\varphi : \G_a \to \SL(n, k)$ is also a $\rB(2, k)$-fundamental homomorphism. 
\end{enumerate} 
\end{lem}

\begin{proof}
(1) There exists a regular matrix $P$ of $\GL(n, k)$ such that 
$\Inn_P \circ \varphi = \phi$. 
Since $\varphi$ is $\rB(2, k)$-fundamental, 
there exists a homomorphism $\omega : \G_m \to \SL(n, k)$ 
such that $\psi_{\varphi, \, \omega} : \rB(2, k) \to \SL(n, k)$ is a homomorphism. 
By Lemma 1.16, we have 
$\Inn_P \circ \psi_{\varphi, \, \omega}
 = \psi_{\phi, \; \Inn_P \circ \omega}$. 
Clearly, $\psi_{\phi, \; \Inn_P \circ \omega}$ is a homomorphism. 
Thus, $\phi$ is a $\rB(2, k)$-fundamental homomorphism.

(2) There exists a homomorphism $\omega : \G_m \to \SL(n, k)$ such that 
$\psi_{\varphi, \, \omega} : \rB(2, k) \to \SL(n, k)$ is a homomorphism. 
We know from Lemma 1.11 that ${^\tau\!}\varphi : \G_a \to \SL(n, k)$ is also a $\rB(2, k)$-fundamental homomorphism. 
\end{proof}

\subsection{Extending homomorphisms $\rB(2, k) \to \SL(n, k)$ to 
$\SL(2, k) \to \SL(n, k)$}

Let $\imath_{\rB(2, k)} : \rB(2, k) \to \SL(2, k)$ be the inclusion homomorphism. 
A homomorphism $\psi : \rB(2, k) \to \SL(n, k)$ is said to be {\it extendable} if 
there exists a homomorphism $\sigma : \SL(2, k) \to \SL(n, k)$ such that $\sigma \circ \imath_{\rB(2, k)}  = \psi$, i.e., 
the following diagram is commutative: 
\[
\xymatrix@M=6pt{
 \rB(2, k) \ar[r]^(.45)\psi  \ar@{^(->}[d]_{\imath_{\rB(2, k)}}
  & \SL(n, k) \\
  \SL(2, k) \ar[ru]_\sigma
 & 
}
\]

\begin{lem}
Let $\psi : \rB(2, k) \to \SL(n, k)$ and $\psi^* : \rB(2, k) \to \SL(n, k)$ be 
homomorphisms such that $\psi$ and $\psi^*$ are equivalent. 
If $\psi$ is extendable, then $\psi^*$ is extendable. 
\end{lem}

\begin{proof}
Since $\psi$ is extendable, there exists a homomorphism $\sigma : \SL(2, k) \to \SL(n, k)$ 
such that $\sigma \circ \imath_{\rB(2, k)} = \psi$. 
There exists a regular matrix $P$ of $\GL(n, k)$ such that $\Inn_P \circ \psi = \psi^*$. 
Let $\sigma^* : \SL(2, k) \to \SL(n, k)$ be the homomorphism defined by 
$\sigma^* := \Inn_P \circ \sigma$. 
Thus $\sigma^* \circ \imath_{\rB(2, k)} = \psi^*$, which implies $\psi^*$ is extendable. 

\end{proof}

A homomorphism $\psi : \rB(2, k) \to \SL(n, k)$ is said to be 
{\it antisymmentric} if 
there exists an element $(\varphi, \omega)$ 
of $\Hom(\G_a, \SL(n, k)) \times \Hom(\G_m, \SL(n, k))$ such that 
$\psi \circ \jmath = \psi_{\varphi, \, \omega}$ and $\omega \in \Omega(n)$.

\begin{lem}
Let $\psi : \rB(2, k) \to \SL(n, k)$ be an extendable homomorphism. 
Then there exists a homomorphism $\psi^* : \rB(2, k) \to \SL(n, k)$ 
such that $\psi$ and $\psi^*$ are equivalent, i.e, $\psi \sim \psi^*$,  
and $\psi^*$ is antisymmetric. 
\end{lem}

\begin{proof}
Write $\psi = \psi_{\varphi, \, \omega}$, where 
$(\varphi, \omega) \in \Hom(\G_a, \SL(n, k)) \times \Hom(\G_m, \SL(n, k))$. 
There exist a regular matrix $P$ of $\GL(n, k)$ and 
integers $d_1, \ldots, d_n$ ($d_1 \geq d_2 \geq \cdots \geq d_n$) such that 
$\Inn_P \circ \omega = \omega_{d_1, \, \ldots, \, d_n }$. 
Since $\psi$ is extendable, we have 
$d_i = - d_{n - i + 1}$ for all $1 \leq i \leq n$. 
So, $\omega_{d_1, \, \ldots, \, d_n } \in \Omega(n)$. 
Let $\psi^* : \rB(2, k) \to \GL(n, k)$ be the homomorphism defined by 
$\psi^* :=  \inn_P \circ \psi$. 
Thus $\psi \sim \psi^*$ and 
$\psi^*$ is antisymmetric (see Lemma 1.16). 
\end{proof}

\begin{lem}
Let $\psi : \rB(2, k) \to \SL(n, k)$ be an antisymmetric homomorphism. 
Then $\psi^\star$ is also an antisymmetric homomorphism. 
\end{lem}

\begin{proof}
See Lemma 1.12. 
\end{proof}

\begin{lem}
Let $\psi : \rB(2, k) \to \SL(n, k)$ be a homomorphism. 
Then the following assertions {\rm (1)} and {\rm (2)} hold true: 
\begin{enumerate}[label = {\rm (\arabic*)}]
\item $\psi$ is extendable if and only if $\psi^\star : \rB(2, k) \to \SL(n, k)$ is extendable. 

\item $\psi$ is uniquely extendable if and only if 
$\psi^\star : \rB(2, k) \to \SL(n, k)$ is uniquely extendable. 
\end{enumerate} 
\end{lem}

\begin{proof}
(1) If $\psi$ is extendable, there exists a homomorphism 
$\sigma : \SL(2, k) \to \SL(n, k)$ such that $\psi = \sigma \circ \imath_{\rB(2, k)}$. 
We can show $\psi^\star = \sigma^\star \circ \imath_{\rB(2, k)}$. 
Thus $\psi^\star$ is extendable (see Lemma 1.21 (2)). 
Conversely, if $\psi^\star$ is extendable, then $\psi^{\star\star}$ is extendable. 
So, $\psi$ is extendable (see Lemma 1.12).

(2) If $\psi$ is uniquely extendable, 
there exists a unique homomorphism 
$\sigma : \SL(2, k) \to \SL(n, k)$ such that $\psi = \sigma \circ \imath_{\rB(2, k)}$. 
Let $\tau_i : \SL(2, k) \to \SL(n, k)$ ($i = 1, 2$) be homomorphisms 
such that $\psi^\star = \tau_i \circ \imath_{\rB(2, k)}$. 
Then $\psi = \tau_i^\star \circ \imath_{\rB(2, k)}$, which implies $\tau_1^\star = \tau_2^\star$. Thus $\tau_1 = \tau_2$. 
Conversely, if $\psi^\star$ is uniquely extendable, 
then $\psi^{\star\star}$ is uniquely extendable. 
So, $\psi$ is uniquely extendable.  

\end{proof}

We can define an isomomorphism $r_{\SL(2, k)} : \SL(2, k) \to \SL(2, k)$ as 
\[
r_{\SL(2, k)}
\left(
\begin{array}{c c}
 a & b \\
 c & d
\end{array}
\right) 
 :=
 \left(
\begin{array}{r r}
 a & - b \\
 - c & d
\end{array}
\right) . 
\]
Clearly, $r_{\SL(2, k)}^2 = \id_{\SL(2, k)}$ and 
the following diagram is commutative: 
\[
\xymatrix@R=36pt@C=36pt@M=6pt{
 \rB(2, k) \ar[r]^{r_{\rB(2, k)}} \ar@{^(->}[d]^{\imath_{\rB(2, k)}}
 & \rB(2, k) \ar@{^(->}[d]^{\imath_{\rB(2, k)}} \\
 \SL(2, k) \ar[r]^{r_{\SL(2, k)}}
 & \SL(2, k) 
}
\]

\begin{lem}
Let $\psi : \rB(2, k) \to \SL(n, k)$ be a homomorphism. 
Then the following assertions {\rm (1)} and {\rm (2)} hold true: 
\begin{enumerate}[label = {\rm (\arabic*)}]
\item $\psi$ is extendable if and only if 
$\psi \circ r_{\rB(2, k)}$ is extendable. 

\item $\psi$ is uniquely extendable if and only if 
$\psi \circ r_{\rB(2, k)}$ is uniquely extendable. 
\end{enumerate} 
\end{lem}

\begin{proof}
(1) The proof is straightforward. See the above commutative diagram. 

(2) Assume $\psi$ is uniquely extendable. 
Let $\tau_i : \SL(2, k) \to \SL(n, k)$ $(i = 1, 2)$ be homomorphisms 
such that $\psi \circ r_{\rB(2, k)} = \tau_i \circ \imath_{\rB(2, k)}$. 
So, $\psi = \tau_i \circ \imath_{\rB(2, k)} \circ r_{\rB(2, k)} 
 = \tau_i \circ r_{\SL(2, k)} \circ \imath_{\rB(2, k)}$. 
Thus $\tau_1 \circ r_{\SL(2, k)} = \tau_2 \circ r_{\SL(2, k)}$. 
Thereby, $\tau_1 = \tau_2$. 
Conversely, assume $\psi \circ r_{\rB(2, k)}$ is uniquely extendable. 
Then $\psi \circ r_{\rB(2, k)} \circ r_{\rB(2, k)}$ is uniquely extendable. 
Thus $\psi$ is uniquely extendable.

\end{proof}

\begin{lem}
Let $(\varphi, \omega)$ be a pair of $\Hom(\G_a, \SL(n, k)) \times \Omega(n)$ 
such that $\psi_{\varphi, \, \omega} : \G_a \rtimes \G_m \to \SL(n, k)$ is a homomorphism. 
Then the following assertions {\rm (1)}, {\rm (2)}, {\rm (3)} hold true: 
\begin{enumerate}[label = {\rm (\arabic*)}]
\item 
$\psi_{\varphi, \, \omega} \circ \jmath^{-1} : \rB(2, k) \to \SL(n, k)$ 
is extendable if and only if 
$\psi_{\, ^\tau\!\varphi, \; \omega} \circ \jmath^{-1} : \rB(2, k) \to \SL(n, k)$ 
is extendable. 

\item $\psi_{\varphi, \, \omega} \circ \jmath^{-1} : \rB(2, k) \to \SL(n, k)$ 
is uniquely extendable if and only if 
$\psi_{\, ^\tau\!\varphi, \; \omega} \circ \jmath^{-1} : \rB(2, k) \to \SL(n, k)$ 
is uniquely extendable. 

\item Let $\sigma : \SL(2, k) \to \SL(n, k)$ be a homomorphism such that 
\[
 \psi_{\varphi, \, \omega} \circ \jmath^{-1}
 = \sigma \circ \imath_{\rB(2, k)} . 
\]
Then the homomorphism ${^\tau\!} \sigma^\tau : \SL(2, k) \to \SL(n, k)$ satisfies 
\[
\psi_{\, ^\tau\!\varphi, \; \omega} \circ \jmath^{-1}
 = ({^\tau\!} \sigma^\tau ) \circ \imath_{\rB(2, k)} . 
\]
\end{enumerate} 

\end{lem}

\begin{proof}
(1) Let $ \psi := \psi_{\varphi, \, \omega} \circ \jmath^{-1}$. 
If $\psi$ is extendable, then $\psi^\star \circ r_{\rB(2, k)}$ is also extendable 
(see Lemmas 2.6 (1) and 2.7 (1)). 
Since $\psi^\star \circ r_{\rB(2, k)} \circ \jmath = \psi_{\, ^\tau\!\varphi, \; \omega}$ 
(see Lemma 1.15), the homomorphism 
$\psi_{\, ^\tau\!\varphi, \; \omega} \circ \jmath^{-1} $ is extendable. 
Conversely, if $\psi_{\, ^\tau\!\varphi, \; \omega} \circ \jmath^{-1} $ is extendable, 
then $\psi^\star \circ r_{\rB(2, k)}$ is extendable. 
Thus $\psi$ is extendable (see Lemmas 2.6 (1) and 2.7 (1)).

(2) See Lemmas 2.6 (2) and 2.7 (2). 

(3) See Lemma 1.23. 
\end{proof}

\section{Antisymmetric  homomorphisms from $\rB(2, k)$ to $\SL(4, k)$} 

In Section 3, we describe antisymmetric homomorphisms $\rB(2, k) \to \SL(4, k)$ 
(see Theorem 3.1 in Subsection 3.2). 
We prepare such homomorphisms in Subsection 3.1.

\subsection{Antisymmetric  morphisms from $\rB(2, k)$ to $\SL(4, k)$}

In the following (I) -- (XXVI), we focus on certain integers $d_1$, $d_2$ 
(which give rise to antisymmetric homomorphisms $\G_m \to \SL(4, k)$) and 
define homomorphisms $\varphi^* : \G_a \to \SL(4, k)$.

\subsubsection{\rm (I)} 

Assume $p \geq 5$. Let $e_1$ be an integer such that 
\[ 
e_1 \geq 0 . 
\]
Let $d_1$ and $d_2$ be integers such that 
\[
\left\{
\begin{array}{r @{\,}l }
 d_1 & = 3 \, p^{e_1},  \\
 d_2 & = p^{e_1} . 
\end{array}
\right. 
\]
Let $\varphi^* : \G_a \to \SL(4, k)$ be the homomorphism defined by  
\[
 \varphi^*(t) 
:= 
\left( 
\begin{array}{c c c c}
 1 & t^{p^{e_1}} & \frac{1}{2} \, t^{2 \, p^{e_1}}  & \frac{1}{6} \, t^{3 \, p^{e_1}}\\
 0 & 1 & t^{p^{e_1}} &  \frac{1}{2} \, t^{2 \, p^{e_1}}\\
 0 & 0 & 1 & t^{p^{e_1}} \\
 0 & 0 & 0 & 1 
\end{array}
\right) . 
\]

\subsubsection{\rm (II)} 

Assume $p = 3$. 
Let $e_1$ be an integer such that 
\[
 e_1 \geq 0 . 
\]
Let $d_1$ and $d_2$ be integers such that
\[
\left\{
\begin{array}{r @{\,} l}
 d_1 & = p^{e_1 + 1}, \\
 d_2 & = p^{e_1} . 
\end{array}
\right.  
\]
Let $\varphi^* : \G_a \to \SL(4, k)$ be the homomorphism defined by  
\[
 \varphi^*(t) 
:= 
\left( 
\begin{array}{c c c c}
 1 & t^{p^{e_1}} & \frac{1}{2} \, t^{2 \, p^{e_1}} & t^{p^{e_1 + 1}}\\
 0 & 1 & t^{p^{e_1}} & 0  \\
 0 & 0 & 1 & 0 \\
 0 & 0 & 0 & 1 
\end{array}
\right) . 
\]

\subsubsection{\rm (III)} 

Assume $p \geq 3$. 
Let $e_1$ be an integer such that 
\[
 e_1 \geq 0 . 
\]
Let $d_1$ and $d_2$ be integers such that
\[
\left\{
\begin{array}{r @{\,} l}
 d_1 & = 3 \, p^{e_1}, \\ 
 d_2 & = p^{e_1} . 
\end{array}
\right.    
\]
Let $\varphi^* : \G_a \to \SL(4, k)$ be the homomorphism defined by  
\[
 \varphi^*(t) 
:= 
\left( 
\begin{array}{c c c c}
 1 & t^{p^{e_1}} & \frac{1}{2} \, t^{2 \, p^{e_1}} & 0 \\
 0 & 1 & t^{p^{e_1}} & 0  \\
 0 & 0 & 1 & 0 \\
 0 & 0 & 0 & 1 
\end{array}
\right) . 
\]

\subsubsection{\rm (IV)} 

Let $e_1$ and $e_2$ be integers such that 
\[ 
 e_2 > e_1 \geq 0 . 
\]
Let $d_1$ and $d_2$ be integers such that
\[ 
\left\{
\begin{array}{r @{\,} l}
 d_1 = p^{e_1} + p^{e_2}, \\
 d_2 = p^{e_2} - p^{e_1} . 
\end{array}
\right. 
\]
Let $\varphi^* : \G_a \to \SL(4, k)$ be the homomorphism defined by  
\[
 \varphi^*(t) 
:= 
\left( 
\begin{array}{c c c c}
 1 & t^{p^{e_1}} & t^{p^{e_2}} &  t^{p^{e_1} + p^{e_2}} \\
 0 & 1 & 0 & t^{p^{e_2}} \\
 0 & 0 & 1 & t^{p^{e_1}} \\
 0 & 0 & 0 & 1 
\end{array}
\right) . 
\]

\subsubsection{\rm (V)} 

Let $e_1$ and $f$ be integers such that 
\[
 e_1 \geq 0, \qquad  f \geq e_1 + 1 . 
\] 
Let $d_1$ and $d_2$ be integers such that
\[
\left\{
\begin{array}{r @{\,} l}
 d_1 & = p^{f}, \\
 d_2 & = p^{f} - 2 \, p^{e_1} . 
\end{array}
\right.  
\]
Let $\varphi^* : \G_a \to \SL(4, k)$ be the homomorphism defined by  
\[
 \varphi^*(t) 
:= 
\left( 
\begin{array}{c c c c}
 1 & t^{p^{e_1}} & 0 & t^{p^{f}} \\
 0 & 1 & 0 & 0 \\
 0 & 0 & 1 & t^{p^{e_1}} \\
 0 & 0 & 0 & 1 
\end{array}
\right) . 
\]

\subsubsection{\rm (VI)} 

Let $e_1$ be an integer such that 
\[
e_1 \geq 0 . 
\]
Let $d_1$ and $d_2$ be integers such that 
\[
\left\{
\begin{array}{r @{\,} l}
 d_1 & = d_2 + 2 \, p^{e_1}, \\ 
 d_2 & \geq 0 . 
\end{array}
\right.  
\]
Let $\varphi^* : \G_a \to \SL(4, k)$ be the homomorphism defined by  
\[
 \varphi^*(t) 
:= 
\left( 
\begin{array}{c c c c}
 1 & t^{p^{e_1}} & 0 & 0 \\
 0 & 1 & 0 & 0 \\
 0 & 0 & 1 & t^{p^{e_1}} \\
 0 & 0 & 0 & 1 
\end{array}
\right) . 
\]

\subsubsection{\rm (VII)} 

Assume $p = 3$. 
Let $e_1$ be an integer such that 
\[
 e_1 \geq 0 . 
\]
Let $d_1$ and $d_2$ be integers such that
\[
\left\{
\begin{array}{r @{\,} l}
 d_1 & = p^{e_1 + 1}, \\
 d_2 & = p^{e_1} . 
\end{array}
\right. 
\]
Let $\varphi^* : \G_a \to \SL(4, k)$ be the homomorphism defined by  
\[
 \varphi^*(t) 
:= 
\left( 
\begin{array}{c c c c}
 1 & 0 & 0  & t^{p^{e_1 + 1}}\\
 0 & 1 & t^{p^{e_1}} &  \frac{1}{2} \, t^{2 \, p^{e_1}}\\
 0 & 0 & 1 & t^{p^{e_1}} \\
 0 & 0 & 0 & 1 
\end{array}
\right) . 
\]

\subsubsection{\rm (VIII)}

Assume $p \geq 3$. 
Let $e_1$ be an integer such that 
\[
 e_1 \geq 0 . 
\]
Let $d_1$ and $d_2$ be integers such that
\[
\left\{
\begin{array}{r @{\,} l}
 d_1 & = 3 \, p^{e_1}, \\
 d_2 & = p^{e_1} . 
\end{array}
\right. 
\]
Let $\varphi^* : \G_a \to \SL(4, k)$ be the homomorphism defined by  
\[
 \varphi^*(t) 
:= 
\left( 
\begin{array}{c c c c}
 1 & 0 & 0  & 0 \\
 0 & 1 & t^{p^{e_1}} &  \frac{1}{2} \, t^{2 \, p^{e_1}} \\
 0 & 0 & 1 & t^{p^{e_1}} \\
 0 & 0 & 0 & 1 
\end{array}
\right) . 
\]

\subsubsection{\rm (IX)} 

Assume $p \geq 3$. Let $e_1$ be an integer such that 
\[
 e_1 \geq 0 . 
\]
Let $d_1$ and $d_2$ be integers such that
\[
\left\{
\begin{array}{r @{\,} l}
 d_1 & = 2 \, p^{e_1}, \\ 
 d_2 & = 0 . 
\end{array}
\right.  
\]
Let $\varphi^* : \G_a \to \SL(4, k)$ be the homomorphism defined by  
\[
 \varphi^*(t) 
:= 
\left( 
\begin{array}{c c c c}
 1 & t^{p^{e_1}} & 0 & \frac{1}{2} \, t^{p^{2 \, e_1}} \\
 0 & 1 & 0 & t^{p^{e_1}} \\
 0 & 0 & 1 & 0 \\
 0 & 0 & 0 & 1 
\end{array}
\right) . 
\]

\subsubsection{\rm (X)} 

Let $e_1$ and $e_2$ be integers such that 
\[
 e_2 > e_1 \geq 0 . 
\]
Let $d_1$ and $d_2$ be integers such that 
\[
\left\{
\begin{array}{r @{\,} l}
 d_1 & = p^{e_1} + p^{e_2}, \\ 
 d_2 & = p^{e_2} - p^{e_1} . 
\end{array}  
\right. 
\]
Let $\varphi^* : \G_a \to \SL(4, k)$ be the homomorphism defined by  
\[
 \varphi^*(t) 
:= 
\left( 
\begin{array}{c c c c}
 1 & t^{p^{e_1}} & t^{p^{e_2}} & 0 \\
 0 & 1 & 0 & 0 \\
 0 & 0 & 1 & 0 \\
 0 & 0 & 0 & 1 
\end{array}
\right) . 
\]

\subsubsection{\rm (XI)} 

Let $e_1$ and $e_3$ be integers such that 
\[
\left\{
\begin{array}{r @{\,}l }
 e_1 & \geq 0, \\ 
 e_3 & \geq e_1 + 1 . 
\end{array}
\right. 
\]
Let $d_1$ and $d_2$ be integers such that 
\[
\left\{
\begin{array}{r @{\,} l}
 d_1 & = p^{e_3}, \\ 
 d_2 & = p^{e_3} - 2 \, p^{e_1} . 
\end{array}
\right.  
\]
Let $\varphi^* : \G_a \to \SL(4, k)$ be the homomorphism defined by  
\[
 \varphi^*(t) 
:= 
\left( 
\begin{array}{c c c c}
 1 & t^{p^{e_1}} & 0 & t^{p^{e_3}}  \\
 0 & 1 & 0 & 0 \\
 0 & 0 & 1 & 0 \\
 0 & 0 & 0 & 1 
\end{array}
\right) . 
\]

\subsubsection{\rm (XII)} 

Let $e_1$ be an integer such that 
\[
 e_1 \geq 0 . 
\]
Let $d_1$ and $d_2$ be integers such that 
\[
\left\{
\begin{array}{r @{\,} l}
 d_1 & = 2 \, p^{e_1} + d_2,  \\
 d_2 & \geq 0 . 
\end{array}
\right. 
\]
Let $\varphi^* : \G_a \to \SL(4, k)$ be the homomorphism defined by  
\[
 \varphi^*(t) 
:= 
\left( 
\begin{array}{c c c c}
 1 & t^{p^{e_1}} & 0 & 0 \\
 0 & 1 & 0 & 0 \\
 0 & 0 & 1 & 0 \\
 0 & 0 & 0 & 1 
\end{array}
\right) . 
\]

\subsubsection{\rm (XIII)} 

Let $e_1$ and $e_3$ be integers such that 
\[ 
 e_1 > e_3 \geq 0 . 
\]
Let $d_1$ and $d_2$ be integers such that
\[
\left\{
\begin{array}{r @{\,} l}
 d_1 & = 2 \, p^{e_1} - p^{e_3}, \\
 d_2 & = p^{e_3} . 
\end{array}
\right.  
\]
Let $\varphi^* : \G_a \to \SL(4, k)$ be the homomorphism defined by  
\[
 \varphi^*(t) 
:= 
\left( 
\begin{array}{c c c c}
 1 & 0 & t^{p^{e_1}} & 0 \\
 0 & 1 & t^{p^{e_3}} & t^{p^{e_1}}  \\
 0 & 0 & 1 & 0 \\
 0 & 0 & 0 & 1 
\end{array}
\right) . 
\]

\subsubsection{\rm (XIV)} 

Let $e_1$ and $e_3$ be integers such that 
\[ 
 e_1 > e_3 \geq 0 . 
\]
Let $d_1$ and $d_2$ be integers such that
\[
\left\{
\begin{array}{r @{\,} l}
 d_1 & = 2 \, p^{e_1} - p^{e_3}, \\
 d_2 & = p^{e_3} . 
\end{array}
\right.  
\]
Let $\varphi^* : \G_a \to \SL(4, k)$ be the homomorphism defined by  
\[
 \varphi^*(t) 
:= 
\left( 
\begin{array}{c c c c}
 1 & 0 & t^{p^{e_1}} & 0 \\
 0 & 1 & t^{p^{e_3}} & 0  \\
 0 & 0 & 1 & 0 \\
 0 & 0 & 0 & 1 
\end{array}
\right) . 
\]

\subsubsection{\rm (XV)} 

 Let $e_2$ and $e_3$ be integers such that 
\[
 e_2 \geq e_3 \geq 0 . 
\]
Let $d_1$ and $d_2$ be integers such that
\[
\left\{
\begin{array}{r @{\,} l}
 d_1 & = p^{e_2},  \\
 d_2 & = p^{e_3} . 
\end{array}
\right.  
\]
Let $\varphi^* : \G_a \to \SL(4, k)$ be the homomorphism defined by  
\[
 \varphi^*(t) 
:= 
\left( 
\begin{array}{c c c c}
 1 & 0 & 0 & t^{p^{e_2}}  \\
 0 & 1 & t^{p^{e_3}} & 0  \\
 0 & 0 & 1 & 0 \\
 0 & 0 & 0 & 1 
\end{array}
\right) . 
\]

\subsubsection{\rm (XVI)} 

Let $e_3$ and $e_4$ be integers such that 
\[ 
 e_4 > e_3 \geq 0 . 
\]
Let $d_1$ and $d_2$ be integers such that
\[
\left\{
\begin{array}{r @{\,} l}
 d_1 & = 2 \, p^{e_4} - p^{e_3}, \\ 
 d_2 & = p^{e_3} . 
\end{array}
\right. 
\]
Let $\varphi^* : \G_a \to \SL(4, k)$ be the homomorphism defined by  
\[
 \varphi^*(t) 
:= 
\left( 
\begin{array}{c c c c}
 1 & 0 & 0 & 0 \\
 0 & 1 & t^{p^{e_3}} & t^{p^{e_4}}  \\
 0 & 0 & 1 & 0 \\
 0 & 0 & 0 & 1 
\end{array}
\right) . 
\]

\subsubsection{\rm (XVII)} 

Let $e_3$ be an integer such that 
\[
 e_3 \geq 0. 
\]
Let $d_1$ and $d_2$ be integers such that
\[
\left\{
\begin{array}{r @{\,} l}
 d_1 & \geq  p^{e_3} ,  \\ 
 d_2 & = p^{e_3} . 
\end{array}
\right. 
\]
Let $\varphi^* : \G_a \to \SL(4, k)$ be the homomorphism defined by  
\[
 \varphi^*(t) 
:= 
\left( 
\begin{array}{c c c c}
 1 & 0 & 0 & 0 \\
 0 & 1 & t^{p^{e_3}} & 0  \\
 0 & 0 & 1 & 0 \\
 0 & 0 & 0 & 1 
\end{array}
\right) . 
\]

\subsubsection{\rm (XVIII)} 

 Let $e_1$ and $e_2$ be integers such that 
\[
 e_2 > e_1 \geq 0 . 
\]
Let $d_1$ and $d_2$ be integers such that 
\[
\left\{
\begin{array}{r @{\,} l}
 d_1 & = p^{e_1} + p^{e_2}, \\
 d_2 & = p^{e_2} - p^{e_1} . 
\end{array}
\right. 
\]
Let $\varphi^* : \G_a \to \SL(4, k)$ be the homomorphism defined by  
\[
 \varphi^*(t) 
:= 
\left( 
\begin{array}{c c c c}
 1 & 0 & 0 & 0 \\
 0 & 1 & 0 & t^{p^{e_2}} \\
 0 & 0 & 1 & t^{p^{e_1}} \\
 0 & 0 & 0 & 1 
\end{array}
\right) . 
\]

\subsubsection{\rm (XIX)} 

Let $e_1$ and $e_3$ be integers such that 
\[
\left\{
\begin{array}{r @{\,} l}
 e_1 & \geq 0, \\
 e_3 & \geq e_1 + 1 . 
\end{array}
\right. 
\]
Let $d_1$ and $d_2$ be integers such that 
\[
\left\{
\begin{array}{r @{\,} l}
 d_1 & = p^{e_3}, \\
 d_2 & = p^{e_3} - 2 \, p^{e_1} . 
\end{array}
\right.  
\]
Let $\varphi^* : \G_a \to \SL(4, k)$ be the homomorphism defined by  
\[
 \varphi^*(t) 
:= 
\left( 
\begin{array}{c c c c}
 1 & 0 & 0 & t^{p^{e_3}}  \\
 0 & 1 & 0 & 0 \\
 0 & 0 & 1 & t^{p^{e_1}} \\
 0 & 0 & 0 & 1 
\end{array}
\right) . 
\]

\subsubsection{\rm (XX)} 

Let $e_1$ be an integer such that 
\[
 e_1 \geq 0 . 
\]
Let $d_1$ and $d_2$ be integers such that 
\[
\left\{
\begin{array}{r @{\,} l}
 d_1 & = 2 \, p^{e_1} + d_2  , \\
 d_2 & \geq 0 . 
\end{array}
\right. 
\]
Let $\varphi^* : \G_a \to \SL(4, k)$ be the homomorphism defined by  
\[
 \varphi^*(t) 
:= 
\left( 
\begin{array}{c c c c}
 1 & 0 & 0 & 0 \\
 0 & 1 & 0 & 0 \\
 0 & 0 & 1 & t^{p^{e_1}} \\
 0 & 0 & 0 & 1 
\end{array}
\right) . 
\]

\subsubsection{\rm (XXI)}  

Assume $p = 2$. 
Let $e_1$ be an integer such that 
\[
 e_1 \geq 0 . 
\]
Let $d_1$ and $d_2$ be integers such tht 
\[
\left\{
\begin{array}{r @{\,} l}
 d_1 & = p^{e_1 + 1} , \\
 d_2 & = 0 .
\end{array}
\right.
\]
Let $\varphi^* : \G_a \to \SL(4, k)$ be the homomorphism defined by  
\[
 \varphi^*(t) 
:= 
\left( 
\begin{array}{c c c c}
 1 & 0 & t^{p^{e_1}} & t^{p^{e_1 + 1}} \\
 0 & 1 & 0 & t^{p^{e_1}}  \\
 0 & 0 & 1 & 0 \\
 0 & 0 & 0 & 1 
\end{array}
\right) . 
\]

\subsubsection{\rm (XXII)} 

Let $e_1$ be an integer such that 
\[
 e_1 \geq 0 . 
\] 
Let $d_1$ and $d_2$ be integers such that
\[
\left\{
\begin{array}{l}
 2 \, p^{e_1} \geq d_1 \geq p^{e_1}, \\
 d_2 = 2 \, p^{e_1} - d_1. 
\end{array}
\right. 
\]
Let $\varphi^* : \G_a \to \SL(4, k)$ be the homomorphism defined by  
\[
 \varphi^*(t) 
:= 
\left( 
\begin{array}{c c c c}
 1 & 0 & t^{p^{e_1}} & 0 \\
 0 & 1 & 0 & t^{p^{e_1}}  \\
 0 & 0 & 1 & 0 \\
 0 & 0 & 0 & 1 
\end{array}
\right) . 
\]

\subsubsection{\rm (XXIII)} 

Let $e_1$ be an integer such that 
\[
 e_1 \geq 0 . 
\]
Let $d_1$ and $d_2$ be integers such that 
\[
\left\{
\begin{array}{l}
 2 \, p^{e_1} \geq d_1 \geq p^{e_1}, \\ 
 d_2 = 2 \, p^{e_1} - d_1 . 
\end{array}
\right. 
\]
Let $\varphi^* : \G_a \to \SL(4, k)$ be the homomorphism defined by  
\[
 \varphi^*(t) 
:= 
\left( 
\begin{array}{c c c c}
 1 & 0 & t^{p^{e_1}} & 0 \\
 0 & 1 & 0 & 0  \\
 0 & 0 & 1 & 0 \\
 0 & 0 & 0 & 1 
\end{array}
\right) . 
\]

\subsubsection{\rm (XXIV)} 

Let $e_2$ be an integer such that 
\[
 e_2 \geq 0 . 
\]
Let $d_1$ and $d_2$ be integers such that
\[
\left\{
\begin{array}{l}
  d_1 = p^{e_2} , \\
  p^{e_2} \geq d_2 \geq 0 .
\end{array}
\right.
\]
Let $\varphi^* : \G_a \to \SL(4, k)$ be the homomorphism defined by  
\[
 \varphi^*(t) 
:= 
\left( 
\begin{array}{c c c c}
 1 & 0 & 0 & t^{p^{e_2}} \\
 0 & 1 & 0 & 0  \\
 0 & 0 & 1 & 0 \\
 0 & 0 & 0 & 1 
\end{array}
\right) . 
\]

\subsubsection{\rm (XXV)} 

Let $e_3$ be an integer such that 
\[
 e_3 \geq 0 . 
\]
Let $d_1$ and $d_2$ be integers such that
\[
\left\{
\begin{array}{l}
 2 \, p^{e_3} \geq d_1 \geq p^{e_3}, \\
 d_2 = 2 \, p^{e_3} - d_1 . 
\end{array}
\right. 
\]
Let $\varphi^* : \G_a \to \SL(4, k)$ be the homomorphism defined by  
\[
 \varphi^*(t) 
:= 
\left( 
\begin{array}{c c c c}
 1 & 0 & 0 & 0 \\
 0 & 1 & 0 & t^{p^{e_3}}  \\
 0 & 0 & 1 & 0 \\
 0 & 0 & 0 & 1 
\end{array}
\right) . 
\]

\subsubsection{\rm (XXVI)} 

Let $d_1$ and $d_2$ be integers such that
\[ 
 d_1 \geq d_2 \geq 0 . 
\]
Let $\varphi^* : \G_a \to \SL(4, k)$ be the homomorphism defined by  
\[
 \varphi^*(t) 
:= 
I_4 . 
\]

\subsection{Antisymmetric  homomorphisms from $\rB(2, k)$ to $\SL(4, k)$}

Let $d_1$, $d_2$ be integers and let $\varphi^* : \G_a \to \SL(4, k)$ be 
a homomorphism such that $d_1$, $d_2$ and $\varphi^*$ have 
one of the above forms 
\text{\rm (I) -- (XXVI)}. 
We can define a homomorphism $\omega^* : \G_m \to \SL(4, k)$ as 
\[
 \omega^*(u) = \diag(u^{d_1}, \, u^{d_2}, \, u^{d_3}, \, u^{d_4}), 
\qquad d_1 \geq d_2 \geq d_3 \geq d_4  , 
\qquad d_3 = - d_2 , 
\qquad d_4 = - d_1 . 
\]
We say that the pair $(\varphi^*, \omega^*)$ has the form ($\nu$) if 
the integers $d_1$, $d_2$ and the morphism $\varphi^*$ have the forms given in ($\nu$), 
where $\nu = {\rm I}, \ldots, {\rm XXVI}$.

\begin{thm}
The following assertions {\rm (1)} and {\rm (2)} hold true: 
\begin{enumerate}[label = {\rm (\arabic*)}] 
\item Let $(\varphi^*, \omega^*)$ be a pair of the form {\rm ($\nu$)}, 
where $\nu = {\rm I}, \ldots, {\rm XXVI}$. 
Then the morphism $\psi_{\varphi^*, \, \omega^*} : \G_a \rtimes \G_m \to \SL(4, k)$ is a homomorphism and $\omega^* \in \Omega(4)$.

\item Let $\psi : \rB(2, k) \to \SL(4, k)$ be an antisymmetric homomorphism. 
Express $\psi$ as $\psi \circ \jmath = \psi_{\varphi, \, \omega}$ 
for some $(\varphi, \omega)$ of $\cU_4 \times \Omega(4)$. 
Then there exists an element $(\varphi^*, \omega^*)$ of $\cU_4 \times \Omega(4)$ 
such that the following conditions 
{\rm (a)} and {\rm (b)} hold true: 
\begin{enumerate}[label = {\rm (\alph*)}]
\item $(\varphi, \omega) \sim (\varphi^*, \omega^*)$. 

\item $(\varphi^*, \omega^*)$ has one of the above forms 
\text{\rm (I) -- (XXVI)}. 
\end{enumerate} 

\end{enumerate} 
\end{thm}

\noindent
{\it Proof of {\rm (1)}.} 
For each $\nu = {\rm I}, \ldots , {\rm XXVI}$, 
we can directly prove 
$\omega^*(u) \, \varphi^*(t) \, \omega^*(u)^{-1} = \varphi^*(u^2 \, t)$ 
for all $(t, u) \in \G_a \rtimes \G_m$. 
Thus $\psi_{\varphi^*, \, \omega^*}$ is a homomorphism (see Lemma 1.9). 
Since $d_1 \geq d_2 \geq  0$, we have $\omega^* \in \Omega(4)$.

\subsection{Proof of (2)}

We prepare the following Lemmas 3.2 and 3.3.

\begin{lem}
Let $\Delta \geq 0$ and let $a(t) \in k[t] \backslash \{ 0 \}$ such that 
\[
 u^\Delta \cdot a(t) = a(u^2 \, t) . 
\] 
Then $a(t)$ is a monomial whose degree $\delta$ satisfies $\Delta = 2 \, \delta$. 
\end{lem}

\begin{proof}
Write $a(t) = \sum_{i = 0}^\delta \lambda_i \, t^i$, 
where $\lambda_i \in k$ for all $0 \leq i \leq \delta$ and $\lambda_\delta \ne 0$. 
We have 
\[
\sum_{i = 0}^\delta \lambda_i \,  u^\Delta  \, t^i 
 = 
\sum_{i = 0}^\delta\lambda_i \, u^{2 \, i} \, t^i . 
\]
Since $\lambda_\delta \ne 0$, we have $u^\Delta \, t^\delta = u^{2 \, \delta} \, t^\delta$, which implies $\Delta = 2 \, \delta$ 
and $\lambda_i = 0$ for all $0 \leq i \leq \delta - 1$. 
\end{proof}

\begin{lem}
Let $\varphi : \G_a \to \SL(n, k)$ and let $\omega : \G_m \to \SL(n, k)$ be 
morphisms such that $\psi_{\varphi, \, \omega} : \G_a \rtimes \G_m \to \SL(n, k)$ 
is a homomorphism. 
Write $\varphi(t) = (a_{i, j}(t))$. 
Assume that $\varphi(t)$ is an upper triangular matrix, i.e., $a_{i, j}(t) = 0$ 
for all $1 \leq j < i \leq n$, 
and that $\omega \in \Omega(n)$, i.e., 
\[
\omega(u) = (u^{d_1}, \ldots, u^{d_n}) , \qquad d_1 \geq \cdots \geq d_n , 
\qquad d_i = - d_{n - i + 1} \quad (\, 1 \leq i \leq n \,)  . 
\]
Then each entry $a_{i, j}(t)$ $(1 \leq i \leq j \leq n)$ 
is either zero or a monomial of degree $(d_i - d_j)/2$. 
\end{lem}

\begin{proof}
Since $\omega(u) \cdot \varphi(t) \cdot \omega(u)^{-1} = \varphi(u^2 \, t)$, 
we have $u^{d_i - d_j} \, a_{i, j}(t) = a_{i, j}(u^2 \, t)$ for all $1 \leq i,  j \leq n$.  
we have the desired result. 
\end{proof}

Let $\varphi \in \cU_4$ (see Subsubsection 1.1.1). Then we have 
\[
 \varphi \, \in \, \cU_{[4]} \, \sqcup \, \cU_{[3, 1]} \, \sqcup \, \cU_{[2, 2]} \, \sqcup \, \cU_{[1, 3]} 
 \, \sqcup \, \cU_{[2, 1, 1]} \, \sqcup \, \cU_{[1, 2, 1]} \, \sqcup \, \cU_{[1, 1, 2]} 
 \, \sqcup \, \cU_{[1, 1, 1, 1]} . 
\] 
So, we shall prove assertion (2) of Theorem 3.1 by separating the following caseses: 
\begin{enumerate}[label = {\rm 3.3.\arabic*.}]
\item $\varphi \in \cU_{[4]}$ and $\omega \in \Omega(4)$. 

\item $\varphi \in \cU_{[3, 1]}$ and $\omega \in \Omega(4)$. 

\item $\varphi \in \cU_{[2, 2]}$ and $\omega \in \Omega(4)$. 

\item $\varphi \in \cU_{[1, 3]}$ and $\omega \in \Omega(4)$. 

\item $\varphi \in \cU_{[2, 1, 1]}$ and $\omega \in \Omega(4)$. 

\item $\varphi \in \cU_{[1, 2, 1]}$ and $\omega \in \Omega(4)$. 

\item $\varphi \in \cU_{[1, 1, 2]}$ and $\omega \in \Omega(4)$. 

\item $\varphi \in \cU_{[1, 1, 1, 1]}$ and $\omega \in \Omega(4)$. 
\end{enumerate}

\subsubsection{$\varphi \in \cU_{[4]}$ and $\omega \in \Omega(4)$}

\begin{lem} 
Let $\varphi \in \cU_{[4]}$ and $\omega \in \Omega(4)$. 
Assume that $\psi_{\varphi, \, \omega}$ is a homomorphism. 
Then there exists an element $(\varphi^*, \omega^*)$ of $\Hom(\G_a, \SL(4, k)) \times \Omega(4)$ such that  the following conditions 
{\rm (1)} and {\rm (2)} hold true: 
\begin{enumerate}[label = {\rm (\arabic*)}]
\item $(\varphi, \, \omega) \sim (\varphi^*, \, \omega^*)$. 

\item $(\varphi^*, \, \omega^*)$ has the form {\rm (I)}. 
\end{enumerate} 
\end{lem}

\begin{proof}
Since $\varphi \in \cU_{[4]}$, we can express $\varphi$ as 
\[
\varphi(t) 
 = 
\left(
\begin{array}{c c c c}
 1 & a_1 & a_2 & a_3 \\
 0 & 1 & a_4 & a_5 \\
 0 & 0 & 1 & a_6 \\
 0 & 0 & 0 & 1
\end{array}
\right)
\qquad 
(\, a_1, a_4, a_6 \in \sfP \backslash \{ 0 \}, 
\quad a_2, a_3, a_5 \in k[T]  \,) . 
\]
Write $a_i = \lambda_i \, t^{p^{e_i}}$ $(i = 1, 4, 6)$, 
where $\lambda_1, \lambda_4, \lambda_6 \in k \backslash \{ 0 \}$ 
and $e_1, e_4, e_6 \geq 0$. 
Since $\varphi$ is a homomorphism, we have $e_1 = e_4 = e_6$. 
In fact, for all $t, t' \in k$, we have 
\begin{align*}
\left\{
\begin{array}{l}
a_2(t') + a_1(t) \, a_4(t') + a_2(t) = a_2(t + t') , \\
a_5(t') + a_4(t) \, a_6(t') + a_5(t) = a_5(t + t') . 
\end{array}
\right. 
\end{align*}
Thus $a_1(t) \, a_4(t') =  a_1(t') \, a_4(t)$ and $a_4(t) \, a_6(t') = a_4(t') \, a_6(t)$ 
for all $t, t' \in k$, which implies $e_1 = e_4 = e_6$. 
Since $\omega \in \Omega(4)$, we can express $\omega$ as 
\[
 \omega(u) 
 = \diag(u^{d_1}, \, u^{d_2}, \, u^{d_3}, \, u^{d_4} ) , 
\qquad d_1 \geq d_2 \geq d_3 \geq d_4 , 
\qquad d_3 = - d_2, 
\qquad d_4 = - d_1 . 
\]
Since $\psi_{\varphi, \, \omega}$ is a homomorphism, 
we have $d_1 - d_2 = 2 \, p^{e_1}$, $d_2 - d_3 = 2 \, p^{e_4}$, $d_3 - d_4 = 2 \, p^{e_6}$ (see Lemma 3.3). 
Now, we have 
\[
\left\{
\begin{array}{r @{\,} l @{\qquad} r}
 d_1 - d_2 & = 2 \, p^{e_1}  &  \textcircled{\scriptsize 1} \\ 
 2 \, d_2 & = 2 \, p^{e_1} &  \textcircled{\scriptsize 2}   
\end{array}
\right.
\]
Thus
\[
\left\{ 
\begin{array}{r @{\,} l } 
 d_1 & = 3 \, p^{e_1} ,  \\
 d_2 & = p^{e_1}  .  
\end{array}
\right.
\]
We can express $\varphi$ as 
\[
\varphi(t) 
 = 
\left( 
\begin{array}{c c c c}
 1 & \lambda_1 \, t^{p^e_1} & a_2 & a_3 \\
 0 & 1 & \lambda_4 \, t^{p^{e_1}} & a_5 \\
 0 & 0 & 1 & \lambda_6 \, t^{p^{e_1}} \\
 0 & 0 & 0 & 1 
\end{array}
\right) . 
\]
Let 
\[
P := \diag(\lambda_1 \, \lambda_4 \, \lambda_6, \;  \lambda_4 \, \lambda_6, \;  \lambda_6, \; 1)  \in \GL(4, k) .  
\]
We can deform $(\Inn_P \circ \varphi)(t)$ as 
\[
 (\Inn_P \circ \varphi)(t) 
 = 
\left(
\begin{array}{c c c c}
 1 & t^{p^{e_1}} & \wt{a}_2 & \wt{a}_3 \\
 0 & 1 & t^{p^{e_1}} & \wt{a}_5 \\
 0 & 0 & 1 & t^{p^{e_1}} \\
 0 & 0 & 0  & 1 
\end{array}
\right) , 
\]
where 
\[
 \wt{a}_2 = \frac{1}{\lambda_1 \, \lambda_4} \, a_2, \qquad \wt{a}_3 = \frac{1}{\lambda_1 \, \lambda_4 \, \lambda_6} \, a_3 , 
\qquad \wt{a}_5 = \frac{1}{\lambda_4 \, \lambda_6} \, a_5 . 
\]
Since $\Inn_P \circ \varphi$ is a homomorphism, 
we have
\begin{align*}
\left\{
\begin{array}{l}
 \wt{a}_2(t') + t^{p^{e_1}} \cdot {t'}^{p^{e_1}} + \wt{a}_2(t) = \wt{a}_2(t + t') , \\
 \wt{a}_5(t') + t^{p^{e_1}} \cdot {t'}^{p^{e_1}} + \wt{a}_5(t) = \wt{a}_5(t + t')
\end{array}
\right. 
\end{align*}
for all $t, t' \in k$. 
Let 
\[
\renewcommand{\arraystretch}{1.5}
\left\{
\begin{array}{r @{\,} l}
 \alpha_2 & : = \wt{a}_2 - \frac{1}{2} \, t^{2 \, p^{e_1}} , \\
 \alpha_5 & : = \wt{a}_5 - \frac{1}{2} \, t^{2 \, p^{e_1}} . 
\end{array}
\right.
\]
Thus $\alpha_2$ and $\alpha_5$ are $p$-polynomials. 
We must have $\alpha_2 = \alpha_5 = 0$ (see Lemma 3.3 and use $p \geq 5$). 
We can express $\Inn_P \circ \varphi$ as 
\[
\renewcommand{\arraystretch}{1.2}
 (\Inn_P \circ \varphi)(t) 
 = 
\left(
\begin{array}{c c c c}
 1 & t^{p^{e_1}} & \frac{1}{2} \, t^{2 \, p^{e_1}} & \wt{a}_3 \\
 0 & 1 & t^{p^{e_1}} & \frac{1}{2} \, t^{2 \,p^{e_1}} \\
 0 & 0 & 1 & t^{p^{e_1}} \\
 0 & 0 & 0  & 1 
\end{array}
\right) . 
\]
Let 
\[
 \alpha_3 := \wt{a}_3 - \frac{1}{6} \, t^{3 \, p^{e_1}} . 
\]
Note that $\alpha_3$ is a $p$-polynomial. 
We must have $\alpha_3 = 0$. 
Thus we can express $\Inn_P \circ \varphi$ as 
\[
\renewcommand{\arraystretch}{1.2}
 (\Inn_P \circ \varphi) (t) 
 = 
\left(
\begin{array}{c c c c}
 1 & t^{p^{e_1}} & \frac{1}{2} \, t^{2 \, p^{e_1}}  & \frac{1}{6} \, t^{3 \, p^{e_1}} \\
 0 & 1 & t^{p^{e_1}} & \frac{1}{2} \, t^{2 \, p^{e_1}} \\
 0 & 0 & 1 & t^{p^{e_1}} \\
 0 & 0 & 0 & 1 
\end{array}
\right) . 
\]
So, let $\varphi^* := \Inn_P \circ \omega$ and $\omega^* := \Inn_P \circ \omega$. 
Clearly, $\omega^* = \omega$. 
The pair $(\varphi^*, \omega^*)$ has the form (I). 
\end{proof}

\subsubsection{$\varphi \in \cU_{[3, 1]}$ and $\omega \in \Omega(4)$}

\begin{lem} 
Let $\varphi \in \cU_{[3, 1]}$ and $\omega \in \Omega(4)$. 
Assume that $\psi_{\varphi, \, \omega}$ is a homomorphism. 
Then there exists an element $(\varphi^*, \omega^*)$ of $\Hom(\G_a, \SL(4, k)) \times \Omega(4)$ such that  the following conditions 
{\rm (1)} and {\rm (2)} hold true: 
\begin{enumerate}[label = {\rm (\arabic*)}]
\item $(\varphi, \, \omega) \sim (\varphi^*, \, \omega^*)$. 

\item $(\varphi^*, \, \omega^*)$ has one of the forms {\rm (II)}, {\rm (III)}. 
\end{enumerate} 
\end{lem}

\begin{proof}
Since $\varphi \in \cU_{[3, 1]}$, we can express $\varphi$ as 
\[
\varphi(t) 
 = 
\left(
\begin{array}{c c c c}
 1 & a_1 & a_2 & a_3 \\
 0 & 1 & a_4 & a_5 \\
 0 & 0 & 1 & 0 \\
 0 & 0 & 0 & 1
\end{array}
\right)
\qquad 
(\, a_1, a_4 \in \sfP \backslash \{ 0 \}, 
\quad a_2, a_3, a_5 \in k[T]  \,) . 
\]
Write $a_i = \lambda_i \, t^{p^{e_i}}$ $(i = 1, 4, 5)$, 
where $\lambda_1, \lambda_4 \in k \backslash\{ 0 \}$, $\lambda_5 \in k$ and $e_1, e_4, e_5 \geq 0$. 
Since $\varphi$ is a homomorphism, we have $e_1 = e_4$ (cf. the proof of Lemma 3.4). 
Since $\omega \in \Omega(4)$, we can express $\omega$ as 
\[
 \omega(u) 
 = \diag(u^{d_1}, \, u^{d_2}, \, u^{d_3}, \, u^{d_4} ) , 
\qquad d_1 \geq d_2 \geq d_3 \geq d_4 ,
 \qquad d_3 = - d_2, 
\qquad d_4 = - d_1 . 
\]
Since $\psi_{\varphi, \, \omega}$ is a homomorphism, we have 
$d_1 - d_2 = 2 \, p^{e_1}$ and $d_2 - d_3 = 2^{p^{e_4}}$ (see Lemma 3.3). 
Now, we have 
\[
\left\{
\begin{array}{r @{\,} l @{\qquad} r}
 d_1 - d_2 & = 2 \, p^{e_1}  &  \textcircled{\scriptsize 1} \\
 2 \, d_2 & = 2 \, p^{e_1}   & \textcircled{\scriptsize 2} 
\end{array} 
\right.
\]
Thus
\[
\left\{ 
\begin{array}{r @{\,} l }
 d_1 & = 3 \, p^{e_1} , \\ 
 d_2 & = p^{e_1}  . 
\end{array}
\right.
\]

Suppose to the contrary that $\lambda_5 \ne 0$. 
We have  $d_2 - d_4 = 2 \, p^{e_5}$. 
Since $\varphi$ is a homomorphism, we have $e_1 = e_5$. So, $d_1 + d_2 = 2 \, p^{e_1}$ 
(since $d_4 =  - d_1$). 
This equality implies $4 \, p^{e_1} = 2 \, p^{e_1}$. This is a contradiction.
So, we can express $\varphi$ as 
\[
 \varphi(t) 
 = 
\left( 
\begin{array}{c c c c}
 1 & \lambda_1 \, t^{p^{e_1}} & a_2 & a_3 \\
 0 & 1 & \lambda_4 \, t^{p^{e_1}} & 0 \\
 0 & 0 & 1 & 0 \\
 0 & 0 & 0 & 1
\end{array}
\right) . 
\]
Write $a_3 = \lambda_3 \, t^{p^{e_3}}$, where $\lambda_3 \in k$ and $e_3 \geq 0$. 
We argue by separating the following two cases: 
\begin{enumerate}[label = {\rm (ii.\arabic*)}]
\item $\lambda_3 = 0$. 

\item $\lambda_3 \ne 0$. 
\end{enumerate} 

{\bf Case (ii.1).} 
Let
\[
P := \diag(\lambda_1 \, \lambda_4 , \;  \lambda_4, \; 1, \;  1)  \in \GL(4, k) .  
\]
We can deform $(\Inn_P \circ \varphi)(t)$ as  
\[
 (\Inn_P \circ \varphi)(t) 
 = 
\left(
\begin{array}{c c c c}
 1 & t^{p^{e_1}} & \wt{a}_2 & 0 \\
 0 & 1 & t^{p^{e_1}} & 0 \\
 0 & 0 & 1 & 0 \\
 0 & 0 & 0 & 1
\end{array}
\right) , 
\]
where 
\[
 \wt{a}_2 = \frac{1}{\lambda_1 \, \lambda_4} \, a_2 . 
\]
So, we can express $\Inn_P \circ \varphi$ as 
\[
 (\Inn_P \circ \varphi)(t) 
 = 
\left(
\begin{array}{c c c c}
 1 & t^{p^{e_1}} & \frac{1}{2} \, t^{2 \, p^{e_1}} & 0 \\
 0 & 1 & t^{p^{e_1}} & 0 \\
 0 & 0 & 1 & 0 \\
 0 & 0 & 0 & 1 
\end{array}
\right) . 
\]
So, let $\varphi^* := \Inn_P \circ \omega$ and $\omega^* := \Inn_P \circ \omega$. 
Clearly, $\omega^* = \omega$. 
The pair $(\varphi^*, \omega^*)$ has the form (III).

{\bf Case (ii.2).} We have $d_1 - d_4 = 2 \, p ^{e_3}$, which implies $d_1 = p^{e_3}$. 
So, $3 \, p^{e_1} = p^{e_3}$. Therefore $p = 3$ and $e_3 = e_1 + 1$. 
Let 
\[
P := \diag(\lambda_1 \, \lambda_4 , \;  \lambda_4, \;  1, \;  1/\lambda_3)  
\in \GL(4, k). 
\]
We can show 
\[
 (\Inn_P \circ \varphi)(t) 
 = 
\left(
\begin{array}{c c c c}
 1 & t^{p^{e_1}} & \frac{1}{2} \, t^{2 \, p^{e_1}} &t^{p^{e_1 + 1}} \\
 0 & 1 & t^{p^{e_1}} & 0 \\
 0 & 0 & 1 & 0 \\
 0 & 0 & 0 & 1 
\end{array}
\right) . 
\]
So, let $\varphi^* := \Inn_P \circ \omega$ and $\omega^* := \Inn_P \circ \omega$. 
Clearly, $\omega^* = \omega$. 
The pair $(\varphi^*, \omega^*)$ has the form (II). 
\end{proof}

\subsubsection{$\varphi \in \cU_{[2, 2]}$ and $\omega \in \Omega(4)$}

\begin{lem} 
Let $\varphi \in \cU_{[2, 2]}$ and $\omega \in \Omega(4)$. 
Assume that $\psi_{\varphi, \, \omega}$ is a homomorphism. 
Then there exists an element $(\varphi^*, \omega^*)$ of $\Hom(\G_a, \SL(4, k)) \times \Omega(4)$ such that  the following conditions 
{\rm (1)} and {\rm (2)} hold true: 
\begin{enumerate}[label = {\rm (\arabic*)}]
\item $(\varphi, \, \omega) \sim (\varphi^*, \, \omega^*)$. 

\item $(\varphi^*, \, \omega^*)$ has one of the forms
{\rm (IV)}, {\rm (V)}, {\rm (VI)}, {\rm (IX)}, {\rm (XXIII)}. 
\end{enumerate} 
\end{lem}

\begin{proof}
Since $\varphi \in \cU_{[2, 2]}$, we can express $\varphi$ as 
\[
\varphi(t) 
 = 
\left(
\begin{array}{c c c c}
 1 & a_1 & a_2 & b \\
 0 & 1 & 0 & a_3 \\
 0 & 0 & 1 & a_4 \\
 0 & 0 & 0 & 1
\end{array}
\right)
\qquad 
(\, a_1, a_2, a_3, a_4 \in \sfP, 
\quad a_1 \ne 0, \quad a_4 \ne 0, 
\quad b \in k[T]  \,) . 
\]
Since $\omega \in \Omega(4)$, we can express $\omega$ as 
\[
 \omega(u) 
 = \diag(u^{d_1}, \, u^{d_2}, \, u^{d_3}, \, u^{d_4} ) , 
\qquad d_1 \geq d_2 \geq d_3 \geq d_4  , 
\qquad d_3 = - d_2, 
\qquad d_4 = - d_1. 
\]
We argue by separating the following cases: 
\[
\begin{array}{c | @{\quad} r @{\qquad}  r }
 & a_2 & a_3  \\
\hline 
 {\rm (i)} &  0 & 0  \\
 {\rm (ii)} &  0 & \ne 0 \\
 {\rm (iii)} & \ne 0 & 0 \\
  {\rm (iv)} & \ne 0 & \ne 0 
\end{array}
\]
\begin{enumerate}[label = {\rm (\roman*)}]
\item Write $a_1 = \lambda_1 \, t^{p^{e_1}}$, $a_4 = \lambda_4 \, t^{p^{e_4}}$ 
and $b = \mu \, t^{p^f}$,  
where $\lambda_1, \lambda_4 , \mu \in k$ with $\lambda_i \ne 0$ $(i = 1, 4)$ 
and $e_1, e_4, f \geq 0$. 
So, 
\[
 \varphi(t) 
 = 
\left(
\begin{array}{c c c c}
 1 & \lambda_1 \, t^{p^{e_1}} & 0 &  \mu \, t^{p^f}  \\
 0 & 1 & 0 & 0 \\
 0 & 0 & 1 & \lambda_4 \, t^{p^{e_4}} \\
 0 & 0 & 0 & 1  
\end{array}
\right) . 
\]
We argue by separting the following cases: 
\begin{enumerate}[label = {\rm (i.\arabic*)}]
\item $\mu = 0$. 

\item $\mu \ne 0$. 
\end{enumerate} 

{\bf Case (i.1).}
So, 
\[
 \varphi(t) 
 = 
\left(
\begin{array}{c c c c}
 1 & \lambda_1 \, t^{p^{e_1}} & 0 & 0 \\
 0 & 1 & 0 & 0 \\
 0 & 0 & 1 & \lambda_4 \, t^{p^{e_4}} \\
 0 & 0 & 0 & 1  
\end{array}
\right) . 
\]
Since $\psi_{\varphi, \, \omega}$ is a homomorphism, we have 
$d_1 - d_2 = 2 \, p^{e_1}$ and $d_3 - d_4 = 2 \, p^{e_4}$. 
Now, we have 
\[
\left\{
\begin{array}{r @{\,} l }
 d_1 - d_2 & = 2 \, p^{e_1}  , \\
 d_1 - d_2 & = 2 \, p^{e_4} . 
\end{array}
\right.
\]
So, $e_1 = e_4$. 
Thus 
\[
\left\{
\begin{array}{r @{\,} l }
 d_1 & = d_2 + 2 \, p^{e_1} ,  \\
 d_2 & \geq 0 .
\end{array}
\right.
\]
Let 
\[
P := \diag(1, \; 1/\lambda_1, \; \lambda_4, \;1) 
\in \GL(4, k) . 
\]
We can deform $(\Inn_P \circ \varphi)(t)$ as 
\begin{align*}
(\Inn_P \circ \varphi)(t) 
 & = 
\left(
\begin{array}{c c c c}
 1 & t^{p^{e_1}} & 0 & 0 \\
 0 & 1 & 0 & 0 \\
 0 & 0 & 1 & t^{p^{e_1}} \\
 0 & 0 & 0 & 1  
\end{array}
\right) . 
\end{align*}
So, let $\varphi^* := \Inn_P \circ \varphi$ and $\omega^* := \Inn_P \circ \omega$. 
Clearly, $\omega^* = \omega$. 
The pair $(\varphi^*, \omega^*)$ has the form (VI).

{\bf Case (i.2).} 
Since $\psi_{\varphi, \, \omega}$ is a homomorphism, we have 
$d_1 - d_2 = 2 \, p^{e_1}$, $d_1 - d_4 = 2 \, p^f$, $d_3 - d_4 = 2 \, p^{e_4}$. 
Now, we have 
\[
\left\{
\begin{array}{r @{\,} l @{\qquad} r}
 d_1 - d_2 & = 2 \, p^{e_1}  & \textcircled{\scriptsize 1} \\
 2 \, d_1 & = 2 \, p^f  & \textcircled{\scriptsize 2} \\
 d_1 - d_2 & = 2 \, p^{e_4}  & \textcircled{\scriptsize 3}   
\end{array}
\right.
\]
Thus 
\[
\left\{
\begin{array}{r @{\,} l @{\qquad} l}
 d_1 & = p^f  & (\text{ see \textcircled{\scriptsize 2} }), \\
 d_2 & = p^f - 2 \, p^{e_1} & (\text{ see \textcircled{\scriptsize 1} and use $d_1 = p^f$ }), \\
 e_1 & = e_4  &   (\text{ see \textcircled{\scriptsize 1} and \textcircled{\scriptsize 3} }) . 
\end{array} 
\right. 
\]
Since $d_2 \geq 0$, we have 
\[
 f \geq e_1 + 1 . 
\]
Let 
\[
P := \diag(\mu, \; \mu/\lambda_1, \; \lambda_4, \;1) 
\in \GL(4, k) . 
\]
We can deform $(\Inn_P \circ \varphi)(t)$ as 
\begin{align*}
(\Inn_P \circ \varphi)(t) 
 & = 
\left(
\begin{array}{c c c c}
 1 & t^{p^{e_1}} & 0 & t^{p^f} \\
 0 & 1 & 0 & 0 \\
 0 & 0 & 1 & t^{p^{e_1}} \\
 0 & 0 & 0 & 1  
\end{array}
\right) . 
\end{align*}
So, let $\varphi^* := \Inn_P \circ \varphi$ and $\omega^* := \Inn_P \circ \omega$. 
Clearly, $\omega^* = \omega$. 
The pair $(\varphi^*, \omega^*)$ has the form (V).

\item Write $a_i = \lambda_i \, t^{p^{e_i}}$ ($i = 1, 3, 4$), 
where $\lambda_1, \lambda_3, \lambda_4\in k \backslash \{ 0 \}$ 
and $e_1, e_3, e_4 \geq 0$. 
So, 
\[
 \varphi(t) 
 = 
\left(
\begin{array}{c c c c}
 1 & \lambda_1 \, t^{p^{e_1}} & 0 & b \\
 0 & 1 & 0 & \lambda_3 \, t^{p^{e_3}} \\
 0 & 0 & 1 & \lambda_4 \, t^{p^{e_4}} \\
 0 & 0 & 0 & 1  
\end{array}
\right) . 
\]
Since $\varphi$ is a homomorphism, we have $e_1 = e_3$. 
Since $\psi_{\varphi, \, \omega}$ is a homomorphism, 
we have $d_1 - d_2 = 2 \, p^{e_1}$, $d_2 - d_4 = 2 \, p^{e_3}$, 
$d_3 - d_4 = 2 \, p^{e_4}$. 
Now, we have 
\[
\left\{
\begin{array}{r @{\,} l  @{\qquad} r}
 d_1 - d_2 & = 2 \, p^{e_1}  & \textcircled{\scriptsize 1} \\
 d_1 + d_2 & = 2 \, p^{e_1}  & \textcircled{\scriptsize 2} \\
 d_1 - d_2 & = 2 \, p^{e_4}  & \textcircled{\scriptsize 3}
\end{array}
\right.
\]
Thus 
\[
\left\{
\begin{array}{r @{\,} l  @{\qquad} l}
 d_1 & = 2 \, p^{e_1} & (\text{ see \textcircled{\scriptsize 1} and \textcircled{\scriptsize 2} }) , \\
 d_2 & = 0 & (\text{ see \textcircled{\scriptsize 1} and \textcircled{\scriptsize 2} }) , \\ 
 e_1 & = e_4 & (\text{ see \textcircled{\scriptsize 3} and 
use $d_1 = 2 \, p^{e_1}$ and $d_2 = 0$ }) . 
\end{array}
\right. 
\]
Let 
\begin{align*}
 P  & := 
\left(
\begin{array}{c c c c}
 1 & 0 & 0 & 0 \\
 0 & 1 & 0 & 0 \\
 0 & \lambda_4/\lambda_3 & 1 & 0 \\
 0 & 0 & 0 & 1
\end{array}
\right)
 \cdot 
\diag(\lambda_1 \, \lambda_3, \;  \lambda_3, \;  1, \;  1) \\
& = 
\left(
\begin{array}{c c c c}
 \lambda_1 \, \lambda_3 & 0 & 0 & 0 \\
 0 & \lambda_3 & 0 & 0 \\
 0 & \lambda_4 & 1 & 0 \\
 0 & 0 & 0 & 1 
\end{array}
\right) 
\in \GL(4, k) . 
\end{align*}
We can deform $(\Inn_P \circ \varphi)(t)$ as 
\begin{align*}
 (\Inn_P \circ \varphi)(t) 
 = 
\left(
\begin{array}{c c c c}
 1 & t^{p^{e_1}} & 0 & \wt{b} \\
 0 & 1 & 0 & t^{p^{e_1}} \\
 0 & 0 & 1 & 0 \\
 0 & 0 & 0 & 1 
\end{array}
\right) , 
\end{align*}
where 
\[
 \wt{b} = \frac{1}{\lambda_1 \, \lambda_3} \, b . 
\]
We must have 
\[
 \wt{b} = \frac{1}{2} \, t^{2 \, p^{e_1}} . 
\]
So, we can express $\Inn_P \circ \varphi$ as 
\begin{align*}
 (\Inn_P \circ \varphi)(t) 
 = 
\left(
\begin{array}{c c c c}
 1 & t^{p^{e_1}} & 0 & \frac{1}{2} \, t^{2 \, p^{e_1}}  \\
 0 & 1 & 0 & t^{p^{e_1}} \\
 0 & 0 & 1 & 0 \\
 0 & 0 & 0 & 1 
\end{array}
\right) , 
\end{align*}
So, let $\varphi^* := \Inn_P \circ \varphi$ and 
$\omega^* := \Inn_P \circ \omega$. 
Clearly, $\omega^* = \omega$. 
The pair $(\varphi^*, \omega^*)$ has the form (IX).

\item Write $a_i = \lambda_i \, t^{p^{e_i}}$ ($i = 1, 2, 4$), 
where $\lambda_1, \lambda_2, \lambda_4 \in k \backslash \{ 0 \}$ 
and $e_1, e_2, e_4 \geq 0$. 
So, 
\[
 \varphi(t) 
 = 
\left(
\begin{array}{c c c c}
 1 & \lambda_1 \, t^{p^{e_1}} & \lambda_2 \, t^{p^{e_2}} & b \\
 0 & 1 & 0 & 0 \\
 0 & 0 & 1 & \lambda_4 \, t^{p^{e_4}} \\
 0 & 0 & 0 & 1  
\end{array}
\right) . 
\]
Since $\varphi$ is a homomorphism, we have $e_2 = e_4$. 
Since $\psi_{\varphi, \, \omega}$ is a homomorphism, 
we have $d_1 - d_2 = 2 \, p^{e_1}$, $d_1 - d_3 = 2 \, p^{e_2}$, 
$d_3 - d_4 = 2 \, p^{e_4}$. 
Now, we have 
\[
\left\{
\begin{array}{r @{\,} l @{\qquad} r}
 d_1 - d_2 & = 2 \, p^{e_1}  & \textcircled{\scriptsize 1} \\
 d_1 + d_2 & = 2 \, p^{e_2}  & \textcircled{\scriptsize 2} \\
 d_1 - d_2 & = 2 \, p^{e_2}  & \textcircled{\scriptsize 3}   
\end{array}
\right.
\]
Thus 
\[
\left\{
\begin{array}{r @{\,} l @{\qquad} l}
 e_1 & = e_2  &   (\text{ see \textcircled{\scriptsize 1} and \textcircled{\scriptsize 3} }) , \\
 d_1 & = 2 \, p^{e_1}& (\text{ see \textcircled{\scriptsize 2} and 
\textcircled{\scriptsize 3} and use $e_1 = e_2$ }), \\
 d_2 & = 0 & (\text{ see \textcircled{\scriptsize 2} and 
\textcircled{\scriptsize 3} }) . 
\end{array} 
\right. 
\]
Let 
\begin{align*}
 P & := 
\left(
\begin{array}{c c c c}
 1 & 0 & 0 & 0 \\
 0 & 1 & 0 & 0 \\
 0 & - \lambda_1/\lambda_2 & 1 & 0 \\
 0 & 0 & 0 & 1
\end{array}
\right)
 \cdot 
\diag(\lambda_2 \, \lambda_4, \;  1, \; \lambda_4, \, 1)
 \cdot 
 \left(
\begin{array}{c c c c}
 1 & 0 & 0 & 0 \\
 0 & 0 & 1 & 0 \\
 0 & 1 & 0 & 0 \\
 0 & 0 & 0 & 1
\end{array}
\right) \\
& = 
\left(
\begin{array}{c c c c}
 \lambda_2 \, \lambda_4 & 0 & 0 & 0 \\
 0 & 0 & 1 & 0 \\
 0 & \lambda_4 &  - \lambda_1/\lambda_2 & 0 \\
 0 & 0 & 0 & 1
\end{array}
\right) 
 \in \GL(4, k) . 
\end{align*}
We can deform $(\Inn_P \circ \varphi)(t)$ as 
\[
 (\Inn_P \circ \varphi)(t) 
= 
\left(
\begin{array}{c c c c}
 1 & t^{p^{e_1}} & 0 & \wt{b} \\
 0 & 1 & 0 & t^{p^{e_1}} \\
 0 & 0 & 1 & 0 \\
 0 & 0 & 0 & 1
\end{array}
\right) , 
\]
where 
\[
 \wt{b} = \frac{1}{\lambda_2 \, \lambda_4} \, b . 
\]
We must have 
\[
 \wt{b} = \frac{1}{2} \, t^{2 \, p^{e_1}} . 
\]
So, we can express $\Inn_P \circ \varphi$ as 
\begin{align*}
 (\Inn_P \circ \varphi)(t) 
 = 
\left(
\begin{array}{c c c c}
 1 & t^{p^{e_1}} & 0 & \frac{1}{2} \, t^{2 \, p^{e_1}}  \\
 0 & 1 & 0 & t^{p^{e_1}} \\
 0 & 0 & 1 & 0 \\
 0 & 0 & 0 & 1 
\end{array}
\right) . 
\end{align*}
So, let $\varphi^* := \Inn_P \circ \varphi$ and 
$\omega^* := \Inn_P \circ \omega$. 
Clearly, $\omega^* = \omega$. 
The pair $(\varphi^*, \omega^*)$ has the form (IX).

\item Write $a_i = \lambda_i \, t^{p^{e_i}}$ ($i = 1, 2, 3, 4)$, 
where $\lambda_1, \lambda_2, \lambda_3, \lambda_4 \in k \backslash \{ 0 \}$ 
and $e_1, e_2, e_3, e_4 \geq 0$. 
So, 
\[
 \varphi(t) 
 = 
\left(
\begin{array}{c c c c}
 1 & \lambda_1 \, t^{p^{e_1}} & \lambda_2 \, t^{p^{e_2}} & b \\
 0 & 1 & 0 & \lambda_3 \, t^{p^{e_3}} \\
 0 & 0 & 1 & \lambda_4 \, t^{p^{e_4}} \\
 0 & 0 & 0 & 1 
\end{array}
\right) . 
\]
Since $\psi_{\varphi, \, \omega}$ is a homomorphism, we have 
$d_1 - d_2 = 2 \, p^{e_1}$, 
$d_1 - d_3 = 2 \, p^{e_2}$, 
$d_2 - d_4 = 2 \, p^{e_3}$, 
$d_3 - d_4 = 2 \, p^{e_4}$. 
Now, we have 
\[
\left\{
\begin{array}{r @{\,} l @{\qquad} r}
 d_1 - d_2 &= 2 \, p^{e_1} & \textcircled{\scriptsize 1} \\
 d_1 + d_2 &= 2 \, p^{e_2} & \textcircled{\scriptsize 2} \\
 d_1 + d_2 & = 2 \, p^{e_3}  & \textcircled{\scriptsize 3} \\
 d_1 - d_2 & = 2 \, p^{e_4}  & \textcircled{\scriptsize 4}   
\end{array}
\right.
\]
Thus 
\[
\left\{
\begin{array}{r @{\,} l @{\qquad} l}
 e_1 & = e_4  &   (\text{ see \textcircled{\scriptsize 1} and \textcircled{\scriptsize 4} }) , \\
 e_2 & = e_3  &   (\text{ see \textcircled{\scriptsize 2} and \textcircled{\scriptsize 3} }) , \\
 d_1 & = p^{e_1} + p^{e_2} & (\text{ see \textcircled{\scriptsize1} and 
\textcircled{\scriptsize 2} }), \\
 d_2 & = p^{e_2} - p^{e_1} & (\text{ see \textcircled{\scriptsize 1} and 
\textcircled{\scriptsize 2} }) . 
\end{array} 
\right. 
\]
Let 
\[
 \mu := \frac{\lambda_2 \, \lambda_4}{\lambda_1 \, \lambda_3} 
\]
and let 
\[
 P_1 := \diag(1, \; 1/\lambda_1, \; 1/\lambda_2, \; 1/(\lambda_1 \, \lambda_3)) . 
\]
We can deform $\Inn_{P_1} \circ \varphi$ as 
\[
(\Inn_{P_1} \circ \varphi)(t) 
 = 
\left(
\begin{array}{c c c c}
 1 & t^{p^{e_1}} & t^{p^{e_2}} & \wt{b} \\
 0 & 1 & 0 & t^{p^{e_2}} \\
 0 & 0 & 1 & \mu \, t^{p^{e_1}} \\
 0 & 0 & 0 & 1
\end{array}
\right) , 
\]
where 
\[
 \wt{b} = \frac{1}{\lambda_1 \, \lambda_3} \, b . 
\]
Since $d_2 \geq 0$, we have $e_2 \geq e_1$. 

We argue by separting the following two cases: 
\begin{enumerate}[label = {(iv.\arabic*)}]
\item $e_1 = e_2$. 

\item $e_2 > e_1$. 
\end{enumerate}
\medskip

{\bf Case (iv.1).} 
We have 
\begin{align*}
\left\{
\begin{array}{r @{\,} l}
 d_1 & = 2 \, p ^{e_1} , \\
 d_2 & = 0 . 
\end{array}
\right.
\end{align*}
Let 
\[
 P_2 := 
\left(
\begin{array}{c c c c}
 1 & 0 & 0 & 0 \\
 0 & 1 & - 1 & 0 \\
 0 & 0 & 1 & 0 \\
 0 & 0 & 0 & 1
\end{array}
\right) 
\cdot \diag(1, \; 1, \; \mu, \; 1) 
= 
\left(
\begin{array}{c c c c}
 1 & 0 & 0 & 0 \\
 0 & 1 & - \mu& 0 \\
 0 & 0 & \mu & 0 \\
 0 & 0 & 0 & 1
\end{array}
\right) 
\in \GL(4, k) . 
\]
We can deform $(\Inn_{P_2} \circ \Inn_{P_1} \circ \varphi)(t)$ as 
\[
 (\Inn_{P_2} \circ \Inn_{P_1} \circ \varphi)(t)
 = 
\left(
\begin{array}{c c c c}
 1 & t^{p^{e_1}} & 0 & \wt{b} \\
 0 & 1 & 0 & (\mu + 1) \, t^{p^{e_1}} \\
 0 & 0 & 1 & t^{p^{e_1}} \\
 0 & 0 & 0 & 1
\end{array}
\right) .  
\]
Now, we separte the following two cases: 
\begin{enumerate}[label = {\rm (iv.1.\arabic*)}]
\item $\mu + 1 = 0$.

\item $\mu + 1 \ne 0$. 
\end{enumerate}

{\bf Case (iv.1.1).}  We have 
\[
 (\Inn_{P_2} \circ \Inn_{P_1} \circ \varphi)(t)
 = 
\left(
\begin{array}{c c c c}
 1 & t^{p^{e_1}} & 0 & \wt{b} \\
 0 & 1 & 0 & 0 \\
 0 & 0 & 1 & t^{p^{e_1}} \\
 0 & 0 & 0 & 1
\end{array}
\right) .  
\]
Write $\wt{b} = \nu \, t^{p^f}$, where $\nu \in k$ and $f \geq 0$. 
If $\nu \ne 0$, we have $d_1 - d_4 = 2 \, p^f$, which implies $p = 2$ 
and $f = e_1 + 1$ (since $d_1 = 2 \, p^{e_1}$). 
So, let 
\[
 P := 
\left\{
\begin{array}{l @{\qquad} l}
 P_1 \cdot P_2 \cdot 
\left(
\begin{array}{c c c c}
 1 & 0 & 0 & 0 \\
 0 & 0 & 1 & 0 \\
 0 & 1 & 0 & 0 \\
 0 & 0 & 0 & 1
\end{array}
\right) 
& \text{ if \quad $\nu = 0$, } \\ [3em] 
 P_1 \cdot P_2 \cdot \diag(1, \; 1, \; 1/\nu, \; 1/\nu) 
 & \text{ if \quad $\nu \ne 0$}. 
\end{array}
\right.
\]
We can express $\Inn_P \circ \varphi$ as 
\[
 (\Inn_P \circ \varphi)(t) 
 = 
\left\{
\begin{array}{l @{\qquad} l}
\left(
\begin{array}{c c c c}
 1 & 0 & t^{p^{e_1}} & 0 \\
 0 & 1 & 0 & t^{p^{e_1}} \\
 0 & 0 & 1 & 0  \\
 0 & 0 & 0 & 1
\end{array}
\right) 
& 
\text{ if \quad $\nu = 0$} , 
\\ [3em] 
\left(
\begin{array}{c c c c}
 1 & t^{p^{e_1}} & 0 & t^{p^{e_1 + 1}} \\
 0 & 1 & 0 & 0 \\
 0 & 0 & 1 &  t^{p^{e_1}} \\
 0 & 0 & 0 & 1
\end{array}
\right) 
& 
\text{ if \quad $\nu \ne 0$} . 
\end{array}
\right. 
\]
So, let $\varphi^* := \Inn_P \circ \varphi$ and 
$\omega^* := \Inn_P \circ \omega$. 
Clearly, $\omega^* = \omega$. 
Thus if $\mu + 1 = 0$ and $\nu = 0$, 
the pair $(\varphi^*, \omega^*)$ has the form (XXIII); 
and if $\mu + 1 = 0$ and $\nu \ne 0$, the pair $(\varphi^*, \omega^*)$ 
has the form (V).

{\bf Case (iv.1.2).}  

Let 
\[
 P := P_1 \cdot P_2 \cdot 
\left(
\begin{array}{c c c c}
 1 & 0 & 0 & 0 \\
 0 & 1 & 0 & 0 \\
 0 & 1/(\mu + 1) & 1 & 0 \\
 0 & 0 & 0 & 1 
\end{array}
\right) 
 \cdot 
\diag(1, \; 1, \; 1, \; 1/(\mu + 1)) . 
\]
We can deform $(\Inn_P \circ \varphi)(t)$ as 
\[
 (\Inn_P \circ \varphi)(t) 
 = 
\left(
\begin{array}{c c c c}
 1 & t^{p^{e_1}} & 0 & \widehat{b} \\
 0 & 1 & 0 &  t^{p^{e_1}} \\
 0 & 0 & 1 & 0 \\
 0 & 0 & 0 & 1
\end{array}
\right) , 
\]
where 
\[
 \widehat{b} 
 = \frac{1}{\mu + 1} \, \wt{b} . 
\] 
So, we can express $\Inn_P \circ \varphi$ as 
\[
 (\Inn_P \circ \varphi)(t) 
 = 
\left(
\begin{array}{c c c c}
 1 & t^{p^{e_1}} & 0 & \frac{1}{2} \, t^{2 \, p^{e_1}} \\
 0 & 1 & 0 & t^{p^{e_1}} \\
 0 & 0 & 1 & 0 \\
 0 & 0 & 0 & 1
\end{array}
\right) .  
\]
So, let $\varphi^* := \Inn_P \circ \varphi$ and 
$\omega^* := \Inn_P \circ \omega$. 
Clearly, $\omega^* = \omega$. 
The pair $(\varphi^*, \omega^*)$ has the form (IX). 
\medskip

{\bf Case (iv.2). }
Let  
\[
 P
 := 
P_1 . 
\]
We can show that $\mu = 1$ and $\wt{b} = t^{p^{e_1} + p^{e_2}} + \wt{\beta}$ for some 
$p$-polynomial $\wt{\beta}$ (see \cite[Theorem 3.3]{Tanimoto 2019}). 
Since $\wt{b}$ is a monomial, we must have $\wt{b} = t^{p^{e_1} + p^{e_2}}$. 
So, let $\varphi^* := \Inn_P \circ \varphi$ and 
$\omega^* := \Inn_P \circ \omega$. 
Clearly, $\omega^* = \omega$. 
The pair $(\varphi^*, \omega^*)$ has the form (IV). 
\end{enumerate} 

\end{proof}

\subsubsection{$\varphi \in \cU_{[1, 3]}$ and $\omega \in \Omega(4)$}

\begin{lem} 
Let $\varphi \in \cU_{[1, 3]}$ and $\omega \in \Omega(4)$. 
Assume that $\psi_{\varphi, \, \omega}$ is a homomorphism. 
Then there exists an element $(\varphi^*, \omega^*)$ of $\Hom(\G_a, \SL(4, k)) \times \Omega(4)$ such that  the following conditions 
{\rm (1)} and {\rm (2)} hold true: 
\begin{enumerate}[label = {\rm (\arabic*)}]
\item $(\varphi, \, \omega) \sim (\varphi^*, \, \omega^*)$. 

\item $(\varphi^*, \, \omega^*)$ has one of the forms {\rm (VII)}, {\rm (VIII)}. 
\end{enumerate} 
\end{lem}

\begin{proof}
Let $\phi := {^\tau\!}\varphi$. 
So, $\phi \in \mf{U}_{[3, 1]}$ and 
$\psi_{\phi, \, \omega}$ is a homomorphism (see Lemma 1.11). 
By Lemma 3.5, there exists an element $(\phi', \omega')$ of 
$\Hom(\G_a, \SL(4, k)) \times \Omega(4)$ such that the following 
conditions (i) and (ii) hold true: 
\begin{enumerate}[label = {\rm (\roman*)}]
\item $(\phi, \, \omega) \sim (\phi', \, \omega')$. 

\item $(\phi', \, \omega')$ has one of the forms {\rm (II)}, {\rm (III)}. 
\end{enumerate} 
Let $Q$ be a regular matrix of $\GL(4, k)$ such that 
\[
 \Inn_Q \circ \psi_{\phi, \, \omega} = \psi_{\phi', \, \omega'} 
\]
and let $P := {^\tau\!}(Q^{-1})$. 
We know from Lemma 1.7 that $\omega' = \omega$ and 
\[
 \Inn_P \circ \psi_{\varphi, \, \omega} 
 = \psi_{{^\tau\!}(\phi'), \; \omega} . 
\]
We know from condition (ii) that the pair $({^\tau\!}(\phi'), \, \omega)$ 
has one of the forms (VII), (VIII). 
 
\end{proof}

\subsubsection{$\varphi \in \cU_{[2, 1, 1]}$ and $\omega \in \Omega(4)$}

\begin{lem} 
Let $\varphi \in \cU_{[2, 1, 1]}$ and $\omega \in \Omega(4)$. 
Assume that $\psi_{\varphi, \, \omega}$ is a homomorphism. 
Then there exists an element $(\varphi^*, \omega^*)$ of $\Hom(\G_a, \SL(4, k)) \times \Omega(4)$ such that  the following conditions 
{\rm (1)} and {\rm (2)} hold true: 
\begin{enumerate}[label = {\rm (\arabic*)}]
\item $(\varphi, \, \omega) \sim (\varphi^*, \, \omega^*)$. 

\item $(\varphi^*, \, \omega^*)$ has one of the forms 
{\rm (IX)}, {\rm (X)}, {\rm (XI)}, {\rm (XII)}. 
\end{enumerate} 
\end{lem}

\begin{proof}
Since $\varphi \in \cU_{[2, 1, 1]}$, we can express $\varphi$ as 
\[
\varphi(t) 
 = 
\left(
\begin{array}{c c c c}
 1 & a_1 & a_2 & a_3 \\
 0 & 1 & 0 & a_4 \\
 0 & 0 & 1 & 0 \\
 0 & 0 & 0 & 1
\end{array}
\right)
\qquad 
(\, a_1 \in \sfP \backslash \{ 0 \}, \quad a_2, a_3, a_4 \in k[T]  \,) . 
\]
Since $\omega \in \Omega(4)$, we can express $\omega$ as 
\[
 \omega(u) 
 = \diag(u^{d_1}, \, u^{d_2}, \, u^{d_3}, \, u^{d_4} ) , 
\qquad d_1 \geq d_2 \geq d_3 \geq d_4  , 
\qquad d_3 = - d_2, 
\qquad d_4 = - d_1 . 
\]
We argue by separting the following cases: 
\[
\begin{array}{c | @{\quad} r @{\qquad} r @{\qquad} r }
 & a_2 & a_3 & a_4 \\
\hline 
 {\rm (i)} &  0 & 0 & 0 \\
  {\rm (ii)} & 0 & 0 & \ne 0 \\
 {\rm (iii)} & 0 & \ne 0 & 0 \\
 {\rm (iv)} & \ne 0 & 0 & 0 \\
 {\rm (v)} &  0 & \ne 0 & \ne 0 \\
 {\rm (vi)} & \ne 0 & 0 & \ne 0 \\
 {\rm (vii)} & \ne 0 & \ne 0 & 0 \\
  {\rm (viii)} & \ne 0 & \ne 0 & \ne 0 \\
\end{array}
\]
\begin{enumerate}[label = {\rm (\roman*)}]
\item 
Write $a_1 = \lambda_1 \, t^{p^{e_1}}$, 
where $\lambda_1 \in k \backslash \{ 0 \}$ and $e_1 \geq 0$. 
So, 
\[
 \varphi(t) 
= 
\left(
\begin{array}{c c c c}
 1 & \lambda_1 \, t^{p^{e_1}} & 0 & 0 \\
 0 & 1 & 0 & 0 \\
 0 & 0 & 1 & 0 \\
 0 & 0 & 0 & 1
\end{array}
\right) . 
\] 
We have $d_1 - d_2 = 2 \, p^{e_1}$. 
Thus 
\[
\left\{
\begin{array}{r @{\,} l}
 d_1 & = d_2+ 2 \, p^{e_1} , \\
 d_2 & \geq 0 . 
\end{array}
\right. 
\]
Let 
\[
P := \diag(\lambda_1, \; 1, \; 1, \; 1) \in \GL(4, k) . 
\] 
So, let $\varphi^* := \Inn_P \circ \varphi$ 
and $\omega^* := \Inn_P \circ \omega$. 
Clearly, $\omega^* = \omega$. 
The pair $(\varphi^*, \omega^*)$ has the form (XII).

\item Suppose to the contrary that this case happens. 
Since $\varphi$ is a homomorphism, we must have 
\[
 a_3(t + t') - a_3(t) - a_3(t') = a_1(t) \, a_4(t') , 
\]
which implies $a_1(t) \cdot a_4(t') = 0$. 
This is a contradiction.

\item Write $a_i = \lambda_i \, t^{p^{e_i}}$ ($i  = 1, 3$), 
where $\lambda_1, \lambda_3 \in k \backslash \{ 0 \}$ and $e_1, e_3 \geq 0$. 
So, 
\[
 \varphi(t) 
= 
\left(
\begin{array}{c c c c}
 1 & \lambda_1 \, t^{p^{e_1}} & 0 & \lambda_3 \, t^{p^{e_3}} \\
 0 & 1 & 0 & 0 \\
 0 & 0 & 1 & 0 \\
 0 & 0 & 0 & 1
\end{array}
\right) . 
\] 
We have $d_1 - d_2 = 2 \, p^{e_1}$ and $d_1 - d_4 = 2 \,p^{e_3}$. 
Now, we have
\[
\left\{
\begin{array}{r @{\,} l @{\qquad} r}
 d_1 - d_2 &= 2 \, p^{e_1} &  \textcircled{\scriptsize 1} \\
 2 \, d_1 &= 2 \, p^{e_3} &  \textcircled{\scriptsize 2}    
\end{array}
\right.
\]
Thus 
\[
\left\{
\begin{array}{r @{\,} l}
 d_1 & = p^{e_3}  , \\
 d_2 & = p^{e_3} - 2 \, p^{e_1} . 
\end{array} 
\right. 
\]
Since $d_2 \geq 0$, we have $e_3 \geq e_1 + 1$. 
Let 
\[
 P := \diag(1, \; 1/\lambda_1, \; 1, \; 1/\lambda_3) \in \GL(4, k) . 
\]
We can deform $(\Inn_P \circ \varphi)(t)$ as 
\[
 (\Inn_P \circ \varphi)(t) 
 = 
\left(
\begin{array}{c c c c}
 1 & t^{p^{e_1}} & 0 & t^{p^{e_3}} \\
 0 & 1 & 0 & 0 \\
 0 & 0 & 1 & 0 \\
 0 & 0 & 0 & 1
\end{array}
\right) . 
\]
So, let $\varphi^* := \Inn_P \circ \varphi$ 
and $\omega^* := \Inn_P \circ \omega$. 
Clearly, $\omega^* = \omega$. 
The pair $(\varphi^*, \omega^*)$ has the form (XI).

\item Write $a_i = \lambda_i \, t^{p^{e_i}}$ ($i = 1, 2$), 
where $\lambda_1, \lambda_2 \in k \backslash \{ 0 \}$ and $e_1, e_2 \geq 0$. 
So, 
\[
 \varphi(t) 
= 
\left(
\begin{array}{c c c c}
 1 & \lambda_1 \, t^{p^{e_1}} & \lambda_2 \, t^{p^{e_2}} & 0  \\
 0 & 1 & 0 & 0 \\
 0 & 0 & 1 & 0 \\
 0 & 0 & 0 & 1
\end{array}
\right) . 
\] 
We have $d_1 - d_2 = 2 \, p^{e_1}$ and $d_1 - d_3 = 2 \, p^{e_2}$. 
Now, we have 
\[
\left\{
\begin{array}{r @{\,} l @{\qquad} r}
 d_1 - d_2 &= 2 \, p^{e_1} &  \textcircled{\scriptsize 1} \\
 d_1 + d_2 &= 2 \, p^{e_2} &  \textcircled{\scriptsize 2}    
\end{array}
\right.
\]
Thus 
\[
\left\{
\begin{array}{r @{\,} l}
 d_1 & = p^{e_1} + p^{e_2} , \\
 d_2 & = p^{e_2} - p^{e_1} . 
\end{array} 
\right. 
\]
Since $d_2 \geq 0$, we have $e_2 \geq e_1$. 
We argue by separating the following cases: 
\begin{enumerate}[label = {\rm (iv.\arabic*)}]
\item $e_1 = e_2$. 

\item $e_2 > e_1$. 
\end{enumerate} 

\noindent 
{\bf Case (iv.1).} We have 
\[
\left\{
\begin{array}{r @{\,} l }
 d_1 & = 2 \, p^{e_1} ,  \\
 d_2 & = 0 . 
\end{array} 
\right. 
\]
Let 
\begin{align*}
 P
 & := 
 \diag(1, \; 1/\lambda_1, \; 1/\lambda_2, \; 1)
 \cdot 
\left(
\begin{array}{c c c c}
 1 & 0 & 0 & 0 \\
 0 & 1 & -1 & 0 \\
 0 & 0 & 1 & 0 \\
 0 & 0 & 0 & 1
\end{array}
\right) \\
 & = 
\left(
\begin{array}{c c c c}
 1 & 0 & 0 & 0 \\
 0 & 1/\lambda_1 & - 1/\lambda_1 & 0 \\
 0 & 0 & 1/\lambda_2 & 0 \\
 0 & 0 & 0 & 1
\end{array}
\right) 
 \in \GL(4, k) . 
\end{align*}
We can deform $(\Inn_P \circ \varphi)(t)$ as 
\[
 (\Inn_P \circ \varphi)(t) 
 = 
\left(
\begin{array}{c c c c}
 1 & t^{p^{e_1}} & 0 & 0 \\
 0 & 1 & 0 & 0 \\
 0 & 0 & 1 & 0 \\
 0 & 0 & 0 & 1
\end{array}
\right) . 
\]
So, let $\varphi^* := \Inn_P \circ \varphi$ 
and $\omega^* := \Inn_P \circ \omega$. 
Clearly, $\omega^* = \omega$. 
The pair $(\varphi^*, \omega^*)$ has the form (XII).

\noindent
{\bf Case (iv.2).} Let 
\[
 P := \diag(1, \; 1/\lambda_1, \; 1/\lambda_2, \; 1) \in \GL(4, k) .  
\]
We can deform $(\Inn_P \circ \varphi)(t)$ as 
\[
 (\Inn_P \circ \varphi)(t) 
 = 
\left(
\begin{array}{c c c c}
 1 & t^{p^{e_1}} & t^{p^{e_2}} & 0 \\
 0 & 1 & 0 & 0 \\
 0 & 0 & 1 & 0 \\
 0 & 0 & 0 & 1
\end{array}
\right) . 
\]
So, let $\varphi^* := \Inn_P \circ \varphi$ 
and $\omega^* := \Inn_P \circ \omega$. 
Clearly, $\omega^* = \omega$. 
The pair $(\varphi^*, \omega^*)$ has the form (X).

\item Write $a_i = \lambda_i \, t^{p^{e_i}}$ ($i = 1, 4$), 
where $\lambda_1, \lambda_4 \in k \backslash \{ 0 \}$ and $e_1, e_4 \geq 0$. 
So, 
\[
 \varphi(t) 
= 
\left(
\begin{array}{c c c c}
 1 & \lambda_1 \, t^{p^{e_1}} & 0 & a_3  \\
 0 & 1 & 0 & \lambda_4 \, t^{p^{e_4}} \\
 0 & 0 & 1 & 0 \\
 0 & 0 & 0 & 1
\end{array}
\right) . 
\] 
We have $e_1 = e_4$, $d_1 - d_2 = 2 \, p^{e_1}$ and $d_2 - d_4 = 2 \, p^{e_4}$. 
Now, we have 
\[
\left\{
\begin{array}{r @{\,} l }
 d_1 - d_2 &= 2 \, p^{e_1} , \\
 d_1 + d_2 &= 2 \, p^{e_1} . 
\end{array}
\right.
\]
Thus 
\[
\left\{
\begin{array}{r @{\,} l }
 d_1 & = 2 \, p^{e_1}  , \\
 d_2 & = 0  . 
\end{array} 
\right. 
\]
Let
\[
 P := 
\diag(1, \; 1/\lambda_1, \; 1, \; 1/(\lambda_1 \, \lambda_4))  \in \GL(4, k) . 
\]
We can deform $(\Inn_P \circ \varphi)(t)$ as 
\[
 (\Inn_P \circ \varphi)(t) 
= 
\left(
\begin{array}{c c c c}
 1 & t^{p^{e_1}} & 0 & \wt{a}_3  \\
 0 & 1 & 0 & t^{p^{e_1}} \\
 0 & 0 & 1 & 0 \\
 0 & 0 & 0 & 1
\end{array}
\right) , 
\] 
where 
\[
 \wt{a}_3 = \frac{1}{\lambda_1 \, \lambda_4} \, a_3 . 
\]
We can show $\wt{a}_3 = (1/2) \, t^{2 \, p^{e_1}}$ and 
express $\Inn_P \circ \varphi$ as 
\[
 (\Inn_P \circ \varphi)(t) 
= 
\left(
\begin{array}{c c c c}
 1 & t^{p^{e_1}} & 0 &  \frac{1}{2} \, t^{2 \, p^{e_1}} \\
 0 & 1 & 0 & t^{p^{e_1}} \\
 0 & 0 & 1 & 0 \\
 0 & 0 & 0 & 1
\end{array}
\right) . 
\] 
So, let $\varphi^* := \Inn_P \circ \varphi$ 
and $\omega^* := \Inn_P \circ \omega$. 
Clearly, $\omega^* = \omega$. 
The pair $(\varphi^*, \omega^*)$ has the form (IX).

\item This case cannot happen (cf. the proof written in (ii)).

\item Write $a_i = \lambda_i \, t^{p^{e_i}}$ ($i = 1, 2, 3$), 
where $\lambda_1, \lambda_2, \lambda_3 \in k \backslash \{ 0 \}$ 
and $e_1, e_2, e_3 \geq 0$. 
So, 
\[
 \varphi(t) 
= 
\left(
\begin{array}{c c c c}
 1 & \lambda_1 \, t^{p^{e_1}} & \lambda_2 \, t^{p^{e_2}}  & \lambda_3 \, t^{p^{e_3}}  \\
 0 & 1 & 0 & 0 \\
 0 & 0 & 1 & 0 \\
 0 & 0 & 0 & 1
\end{array}
\right) . 
\] 
We have $d_1 - d_2 = 2 \, p^{e_1}$, $d_1 - d_3 = 2 \, p^{e_2}$, $d_1 - d_4 = 2 \, p ^{e_3}$. 
Now, we have 
\[
\left\{
\begin{array}{r @{\,} l @{\qquad} r}
 d_1 - d_2 &= 2 \, p^{e_1} &  \textcircled{\scriptsize 1} \\
 d_1 + d_2 &= 2 \, p^{e_2} &  \textcircled{\scriptsize 2} \\
 2 \, d_1 &= 2 \, p^{e_3} &  \textcircled{\scriptsize 3}     
\end{array}
\right.
\]
Thus 
\[
\left\{
\begin{array}{r @{\,} l @{\qquad} l}
 d_1 & = p^{e_1} + p^{e_2}  & (\text{ see \textcircled{\scriptsize 1} and 
\textcircled{\scriptsize 2}  }), \\
 d_2 & = p^{e_2} - p^{e_1} & (\text{ see \textcircled{\scriptsize 1} and 
\textcircled{\scriptsize 2} })  , \\
 d_1 & = p^{e_3} &  (\text{ see \textcircled{\scriptsize 3} } ) . 
\end{array} 
\right. 
\]
Since $d_2 \geq 0$, we have $e_2 \geq e_1$. 
Since $p^{e_1} + p^{e_2} = p^{e_3}$, we have $1 + p^{e_2 - e_1} = p^{e_3 - e_1}$, which implies $e_3 - e_1 > 0$. 
So, $e_2 = e_1$, $p = 2$ and $e_3 - e_1 = 1$. 
Thus 
\[
\left\{
\begin{array}{r @{\,} l }
 d_1 & = p^{e_1 + 1}  , \\
 d_2 & = 0  . 
\end{array} 
\right. 
\]
Let
\begin{align*}
 P 
& := 
\diag(1, \; 1/\lambda_1, \; 1/\lambda_2, \; 1/\lambda_3) 
 \cdot 
\left(
\begin{array}{c c c c}
 1 & 0 & 0 & 0 \\
 0 & 1  & - 1 & 0 \\
 0 & 0 & 1 & 0 \\
 0 & 0 & 0 & 1  
\end{array}
\right)  \\ 
& = 
\left(
\begin{array}{c c c c}
 1 & 0 & 0 & 0 \\
 0 & 1 / \lambda_1 & - 1 / \lambda_1 & 0 \\
 0 & 0 & 1 / \lambda_2 & 0 \\
 0 & 0 & 0 & 1 / \lambda_3 
\end{array}
\right)  
\in \GL(4, k) . 
\end{align*}
We can deform $(\Inn_P \circ \varphi)(t)$ as 
\begin{align*}
 (\Inn_P \circ \varphi)(t) 
 & = 
\left(
\begin{array}{c c c c}
 1 & t^{p^{e_1}} & 0 & t^{p^{e_1 + 1}}  \\
 0 & 1 & 0 & 0 \\
 0 & 0 & 1 & 0 \\
 0 & 0 & 0 & 1
\end{array}
\right) . 
\end{align*}
So, let $\varphi^* := \Inn_P \circ \varphi$ 
and $\omega^* := \Inn_P \circ \omega$. 
Clearly, $\omega^* = \omega$. 
The pair $(\varphi^*, \omega^*)$ has the form (XI).

\item Write $a_i = \lambda_i \, t^{p^{e_i}}$ $(i = 1, 2, 4)$, 
where $\lambda_1, \lambda_2, \lambda_4 \in k \backslash \{ 0 \}$ 
and $e_1, e_2, e_4 \geq 0$. 
So, 
\[
 \varphi(t) 
= 
\left(
\begin{array}{c c c c}
 1 & \lambda_1 \, t^{p^{e_1}} & \lambda_2 \, t^{p^{e_2}}  & a_3 \\
 0 & 1 & 0 & \lambda_4 \, t^{p^{e_4}}  \\
 0 & 0 & 1 & 0 \\
 0 & 0 & 0 & 1
\end{array}
\right) . 
\] 
We have $e_1 = e_4$, $d_1 - d_2 = 2 \, p^{e_1}$, $d_1 - d_3 = 2 \, p^{e_2}$, $d_2 - d_4 = 2 \, p ^{e_4}$. 
Now, we have 
\[
\left\{
\begin{array}{r @{\,} l @{\qquad} r}
 d_1 - d_2 &= 2 \, p^{e_1} &  \textcircled{\scriptsize 1} \\
 d_1 + d_2 &= 2 \, p^{e_2} &  \textcircled{\scriptsize 2} \\
 d_1 + d_2 &= 2 \, p^{e_1} &  \textcircled{\scriptsize 3}     
\end{array}
\right.
\]
Thus 
\[
\left\{
\begin{array}{r @{\,} l @{\qquad} l}
 d_1 & = 2 \, p^{e_1}  & (\text{ see \textcircled{\scriptsize 1} and 
\textcircled{\scriptsize 3}  }), \\
 d_2 & = 0 & (\text{ see \textcircled{\scriptsize 1} and 
\textcircled{\scriptsize 3} })  , \\
 d_1 & = 2 \, p^{e_2} &  (\text{ see \textcircled{\scriptsize 2} and 
use $d_2 = 0$ }) . 
\end{array} 
\right. 
\]
So, $e_1 = e_2$. 
Let
\begin{align*}
 P 
& := 
\diag(1, \; 1/\lambda_1, \; 1/\lambda_2, \; 1/(\lambda_1 \, \lambda_4)) 
 \cdot 
\left(
\begin{array}{c c c c}
 1 & 0 & 0 & 0 \\
 0 & 1  & - 1 & 0 \\
 0 & 0 & 1 & 0 \\
 0 & 0 & 0 & 1  
\end{array}
\right)  \\ 
& = 
\left(
\begin{array}{c c c c}
 1 & 0 & 0 & 0 \\
 0 & 1 / \lambda_1 & - 1 / \lambda_1 & 0 \\
 0 & 0 & 1 / \lambda_2 & 0 \\
 0 & 0 & 0 & 1 / ( \lambda_1 \, \lambda_4 ) 
\end{array}
\right)  
\in \GL(4, k) . 
\end{align*}
We can express $\Inn_P \circ \varphi$ as 
\[
 (\Inn_P \circ \varphi)(t) 
= 
\left(
\begin{array}{c c c c}
 1 & t^{p^{e_1}} & 0 &  \frac{1}{2} \, t^{2 \, p^{e_1}} \\
 0 & 1 & 0 & t^{p^{e_1}} \\
 0 & 0 & 1 & 0 \\
 0 & 0 & 0 & 1
\end{array}
\right) . 
\] 
So, let $\varphi^* := \Inn_P \circ \varphi$ 
and $\omega^* := \Inn_P \circ \omega$. 
Clearly, $\omega^* = \omega$. 
The pair $(\varphi^*, \omega^*)$ has the form (IX).

\end{enumerate}

\end{proof}

\subsubsection{$\varphi \in \cU_{[1, 2, 1]}$ and $\omega \in \Omega(4)$}

\begin{lem} 
Let $\varphi \in \cU_{[1, 2, 1]}$ and $\omega \in \Omega(4)$. 
Assume that $\psi_{\varphi, \, \omega}$ is a homomorphism. 
Then there exists an element $(\varphi^*, \omega^*)$ of $\Hom(\G_a, \SL(4, k)) \times \Omega(4)$ such that  the following conditions 
{\rm (1)} and {\rm (2)} hold true: 
\begin{enumerate}[label = {\rm (\arabic*)}]
\item $(\varphi, \, \omega) \sim (\varphi^*, \, \omega^*)$. 

\item $(\varphi^*, \, \omega^*)$ has one of the forms 
{\rm (XIII)}, {\rm (XIV)}, {\rm (XV)}, {\rm (XVI)}, {\rm (XVII)}, {\rm (XXII)}. 
\end{enumerate} 
\end{lem}

\begin{proof}

Since $\varphi \in \cU_{[1, 2, 1]}$, we can express $\varphi$ as 
\[
\varphi(t) 
 = 
\left(
\begin{array}{c c c c}
 1 & 0 & a_1 & a_2 \\
 0 & 1 & a_3 & a_4 \\
 0 & 0 & 1 & 0 \\
 0 & 0 & 0 & 1
\end{array}
\right)
\qquad 
(\, a_1, a_2, a_3, a_4 \in \sfP, \quad a_3 \ne 0  \,) . 
\]
Since $\omega \in \Omega(4)$, we can express $\omega$ as 
\[
 \omega(u) 
 = \diag(u^{d_1}, \, u^{d_2}, \, u^{d_3}, \, u^{d_4} ) , 
\qquad d_1 \geq d_2 \geq d_3 \geq d_4  , 
\qquad d_3 =  - d_2, 
\qquad d_4 =  - d_1 .  
\]
We argue by separting the following cases:  
\[
\begin{array}{c | @{\quad} r @{\qquad} r @{\qquad} r }
 & a_1 & a_2 & a_4 \\
\hline 
 {\rm (i)} &  0 & 0 & 0 \\
  {\rm (ii)} & 0 & 0 & \ne 0 \\
 {\rm (iii)} & 0 & \ne 0 & 0 \\
 {\rm (iv)} & \ne 0 & 0 & 0 \\
 {\rm (v)} &  0 & \ne 0 & \ne 0 \\
 {\rm (vi)} & \ne 0 & 0 & \ne 0 \\
 {\rm (vii)} & \ne 0 & \ne 0 & 0 \\
  {\rm (viii)} & \ne 0 & \ne 0 & \ne 0 \\
\end{array}
\]
\begin{enumerate}[label = {\rm (\roman*)}]
\item Write $a_3 = \lambda_3 \, t^{p^{e_3}}$, 
where $\lambda_3 \in k \backslash \{ 0 \}$ and $e_3 \geq 0$. 
So, 
\[
 \varphi(t) 
= 
\left(
\begin{array}{c c c c}
 1 & 0 & 0 & 0 \\
 0 & 1 & \lambda_3 \, t^{p^{e_3}} & 0 \\
 0 & 0 & 1 & 0 \\
 0 & 0 & 0 & 1
\end{array}
\right) . 
\] 
We have $d_2 - d_3 = 2 \, p^{e_3}$. 
Thus 
\[
 d_1 \geq d_2 = p^{e_3} . 
\]
Let 
\[
P := \diag(1, \; 1, \; 1 / \lambda_3, \; 1) \in \GL(4, k) . 
\] 
We can express $\Inn_P \circ \varphi$ as 
\[
 (\Inn_P \circ \varphi)(t) 
= 
\left(
\begin{array}{c c c c}
 1 & 0 & 0 &  0 \\
 0 & 1 & t^{p^{e_3}} & 0 \\
 0 & 0 & 1 & 0 \\
 0 & 0 & 0 & 1
\end{array}
\right) . 
\] 
So, let $\varphi^* := \Inn_P \circ \varphi$ 
and $\omega^* := \Inn_P \circ \omega$. 
Clearly, $\omega^* = \omega$. 
The pair $(\varphi^*, \omega^*)$ has the form (XVII).

\item Write $a_i = \lambda_i \, t^{p^{e_i}}$ ($i = 3, 4$), 
where $\lambda_3, \lambda_4 \in k \backslash \{ 0 \}$ 
and $e_3, e_4 \geq 0$. 
So, 
\[
 \varphi(t) 
= 
\left(
\begin{array}{c c c c}
 1 & 0 & 0 & 0 \\
 0 & 1 & \lambda_3 \, t^{p^{e_3}} & \lambda_4 \, t^{p^{e_4}} \\
 0 & 0 & 1 & 0 \\
 0 & 0 & 0 & 1
\end{array}
\right) . 
\] 
We have $d_2 - d_3 = 2 \, p^{e_3}$ and $d_2 - d_4 = 2 \, p^{e_4}$. 
Now, we have 
\[
\left\{
\begin{array}{r @{\,} l @{\qquad} r}
 2 \, d_2 &= 2 \, p^{e_3} &  \textcircled{\scriptsize 1} \\
 d_1 + d_2 &= 2 \, p^{e_4} &  \textcircled{\scriptsize 2}     
\end{array}
\right.
\]
Thus 
\[
\left\{
\begin{array}{r @{\,} l}
 d_2 & = p^{e_3} , \\
 d_1 & = 2 \, p^{e_4} - p^{e_3} . 
\end{array} 
\right. 
\]
Since $d_1 \geq d_2$, we have $e_4 \geq e_3$. 
We argue by separating the following cases: 
\begin{enumerate}[label = {\rm (ii.\arabic*)}]
\item $e_4 = e_3$. 

\item $e_4 > e_3$. 
\end{enumerate}

{\bf Case (ii.1).} Let 
\begin{align*}
 P & := \diag(1, \; 1, \; 1/\lambda_3, \; 1/\lambda_4) 
 \cdot 
\left(
\begin{array}{c c c c}
 1 & 0 & 0 & 0 \\
 0 & 1 & 0 & 0 \\
 0 & 0 & 1 & - 1 \\
 0 & 0 & 0 & 1 
\end{array}
\right) \\
 & = 
\left(
\begin{array}{c c c c}
 1 & 0 & 0 & 0 \\
 0 & 1 & 0 & 0 \\
 0 & 0 & 1 / \lambda_3 & - 1/ \lambda_3 \\
 0 & 0 & 0 & 1 / \lambda_4
\end{array}
\right)
 \in \GL(4, k) . 
\end{align*}
We can express $\Inn_P \circ \varphi$ as 
\[
 (\Inn_P \circ \varphi)(t) 
= 
\left(
\begin{array}{c c c c}
 1 & 0 & 0 &  0 \\
 0 & 1 & t^{p^{e_3}} & 0 \\
 0 & 0 & 1 & 0 \\
 0 & 0 & 0 & 1
\end{array}
\right) . 
\] 
So, let $\varphi^* := \Inn_P \circ \varphi$ 
and $\omega^* := \Inn_P \circ \omega$. 
Clearly, $\omega^* = \omega$. 
The pair $(\varphi^*, \omega^*)$ has the form (XVII).

{\bf Case (ii.2).} Let 
\[
P := \diag(1, \; 1, \; 1/\lambda_3, \; 1/\lambda_4) \in \GL(4, k) . 
\]
We can express $\Inn_P \circ \varphi$ as 
\[
 (\Inn_P \circ \varphi)(t) 
= 
\left(
\begin{array}{c c c c}
 1 & 0 & 0 &  0 \\
 0 & 1 & t^{p^{e_3}} & t^{p^{e_4}} \\
 0 & 0 & 1 & 0 \\
 0 & 0 & 0 & 1
\end{array}
\right) . 
\] 
So, let $\varphi^* := \Inn_P \circ \varphi$ 
and $\omega^* := \Inn_P \circ \omega$. 
Clearly, $\omega^* = \omega$. 
The pair $(\varphi^*, \omega^*)$ has the form (XVI).

\item Write $a_i = \lambda_i \, t^{p^{e_i}}$ $(i = 2, 3)$, 
where $\lambda_2, \lambda_3 \in k \backslash \{ 0 \}$ and $e_2, e_3 \geq 0$. 
So, 
\[
 \varphi(t) 
= 
\left(
\begin{array}{c c c c}
 1 & 0 & 0 & \lambda_2 \, t^{p^{e_2}}  \\
 0 & 1 & \lambda_3 \, t^{p^{e_3}} & 0 \\
 0 & 0 & 1 & 0 \\
 0 & 0 & 0 & 1
\end{array}
\right) . 
\]
We have $d_1 - d_4 = 2 \, p^{e_2}$ and $d_2 - d_3 = 2 \, p^{e_3}$. 
So, 
\[
\left\{
\begin{array}{r @{\,} l }
 d_1 & = p^{e_2} , \\
 d_2 & = p^{e_3} . 
\end{array}
\right. 
\]
Since $d_1 \geq d_2$, we have $e_2 \geq e_3$. 
Let 
\[
P := \diag(1, \; 1, \; 1/\lambda_3,  \; 1/\lambda_2) \in \GL(4, k) . 
\]
We can express $\Inn_P \circ \varphi$ as 
\[
 (\Inn_P \circ \varphi)(t) 
= 
\left(
\begin{array}{c c c c}
 1 & 0 & 0 &  t^{p^{e_2}} \\
 0 & 1 & t^{p^{e_3}} & 0 \\
 0 & 0 & 1 & 0 \\
 0 & 0 & 0 & 1
\end{array}
\right) . 
\] 
So, let $\varphi^* := \Inn_P \circ \varphi$ 
and $\omega^* := \Inn_P \circ \omega$. 
Clearly, $\omega^* = \omega$. 
The pair $(\varphi^*, \omega^*)$ has the form (XV).

\item Write $a_i = \lambda_i \, t^{p^{e_i}}$ $(i = 1, 3)$, where $\lambda_1, \lambda_3 \in k \backslash \{ 0 \}$ and $e_1, e_3 \geq 0$. 
So, 
\[
 \varphi(t) 
= 
\left(
\begin{array}{c c c c}
 1 & 0 & \lambda_1 \, t^{p^{e_1}} & 0 \\
 0 & 1 & \lambda_3 \, t^{p^{e_3}} & 0 \\
 0 & 0 & 1 & 0 \\
 0 & 0 & 0 & 1
\end{array}
\right) . 
\] 
We have $d_1 - d_3 = 2 \, p^{e_1}$ and $d_2 - d_3 = 2 \, p^{e_3}$. 
Now, we have 
\[
\left\{
\begin{array}{r @{\,} l @{\qquad} r}
 d_1 + d_2 &= 2 \, p^{e_1} &  \textcircled{\scriptsize 1} \\
 2 \, d_2 &= 2 \, p^{e_3} &  \textcircled{\scriptsize 2}     
\end{array}
\right.
\]
Thus 
\[
\left\{
\begin{array}{r @{\,} l}
 d_2 & = p^{e_3}  , \\
 d_1 & = 2 \, p^{e_1} - p^{e_3} . 
\end{array} 
\right. 
\]
Since $d_1 \geq d_2$, we have $e_1 \geq e_3$. 
We argue by separating the following cases: 
\begin{enumerate}[label = {\rm (iv.\arabic*)}]
\item $e_1 = e_3$. 

\item $e_1 > e_3$. 
\end{enumerate} 

{\bf (iv.1)} So, 
\begin{align*}
\left\{
\begin{array}{r @{\,} l}
 d_1 & = p^{e_3} , \\
 d_2 & = p^{e_3} . 
\end{array}
\right.
\end{align*}
Let 
\begin{align*}
  P &:= \diag(\lambda_1, \; \lambda_3, \; 1, \; 1) 
 \cdot 
\left(
\begin{array}{c c c c}
 1 & 1 & 0 & 0 \\
 0 & 1 & 0 & 0 \\
 0 & 0 & 1 & 0 \\
 0 & 0 & 0 & 1
\end{array}
\right) \\
 & = 
\left(
\begin{array}{c c c c}
 \lambda_1 & \lambda_1 & 0 & 0 \\
 0 & \lambda_3 & 0 & 0 \\
 0 & 0 & 1 & 0 \\
 0 & 0 & 0 & 1
\end{array}
\right) 
\in \GL(4, k) .  
\end{align*}
We can express $\Inn_P \circ \varphi$ as 
\[
 (\Inn_P \circ \varphi)(t) 
= 
\left(
\begin{array}{c c c c}
 1 & 0 & 0 & 0 \\
 0 & 1 & t^{p^{e_3}} & 0 \\
 0 & 0 & 1 & 0 \\
 0 & 0 & 0 & 1
\end{array}
\right) . 
\] 
So, let $\varphi^* := \Inn_P \circ \varphi$ 
and $\omega^* := \Inn_P \circ \omega$. 
Clearly, $\omega^* = \omega$. 
The pair $(\varphi^*, \omega^*)$ has the form (XVII).

{\bf (iv.2)} Let 
\[
  P := \diag(\lambda_1, \; \lambda_3, \; 1, \; 1) \in \GL(4, k) .  
\]
We can express $\Inn_P \circ \varphi$ as 
\[
 (\Inn_P \circ \varphi)(t) 
= 
\left(
\begin{array}{c c c c}
 1 & 0 & t^{p^{e_1}} & 0 \\
 0 & 1 & t^{p^{e_3}} & 0 \\
 0 & 0 & 1 & 0 \\
 0 & 0 & 0 & 1
\end{array}
\right) . 
\] 
So, let $\varphi^* := \Inn_P \circ \varphi$ 
and $\omega^* := \Inn_P \circ \omega$. 
Clearly, $\omega^* = \omega$. 
The pair $(\varphi^*, \omega^*)$ has the form (XIV).

\item Write $a_i = \lambda_i \, t^{p^{e_i}}$ $(i = 2, 3, 4)$, 
where $\lambda_2, \lambda_3, \lambda_4 \in k \backslash \{ 0 \}$ and $e_2, e_3, e_4 \geq 0$. 
So, 
\[
 \varphi(t) 
= 
\left(
\begin{array}{c c c c}
 1 & 0 & 0 & \lambda_2 \, t^{p^{e_2}} \\
 0 & 1 & \lambda_3 \, t^{p^{e_3}} & \lambda_4 \, t^{p^{e_4}} \\
 0 & 0 & 1 & 0 \\
 0 & 0 & 0 & 1
\end{array}
\right) . 
\] 
We have $d_1 - d_4 = 2 \, p^{e_2}$ and $d_2 - d_3 = 2 \, p^{e_3}$, 
$d_2 - d_4 = 2 \, p^{e_4}$.  
Now, we have 
\[
\left\{
\begin{array}{r @{\,} l @{\qquad} r}
 2 \, d_1 &= 2 \, p^{e_2} &  \textcircled{\scriptsize 1} \\
 2 \, d_2 &= 2 \, p^{e_3} &  \textcircled{\scriptsize 2} \\
 d_1 + d_2 & = 2 \, p^{e_4} & \textcircled{\scriptsize 3}
\end{array}
\right.
\]
Thus 
\[
\left\{
\begin{array}{r @{\,} l @{\qquad} l}
 d_1 & = p^{e_2} & (\text{ see \textcircled{\scriptsize 1} })  , \\
 d_2 & = p^{e_3} & (\text{ see \textcircled{\scriptsize 2} })  , \\
 p^{e_2} + p^{e_3} & = 2 \, p^{e_4}  &  (\text{ see \textcircled{\scriptsize 3} and 
use $d_1 = p^{e_2} $ and $d_2 = p^{e_3}$ }) . 
\end{array} 
\right. 
\]
Since $d_1 \geq d_2$, we have $e_2 \geq e_3$. 
So, $p^{e_2 - e_3} + 1 = 2 \, p^{e_4 - e_3}$. 
Thus $e_4 - e_3 = 0$ and $e_2 - e_3 = 0$. 
Let 
\begin{align*}
P & := 
\left( 
\begin{array}{c c c c}
 1 & 0 & 0 & 0 \\
 \lambda_4/\lambda_2 & 1 & 0 & 0 \\
 0 & 0 & 1 & 0 \\
 0 & 0 & 0 & 1
\end{array}
\right) 
\cdot 
\diag(\lambda_2, \;  \lambda_3, \;  1, \;  1) \\
& = 
\left( 
\begin{array}{c c c c}
 \lambda_2 & 0 & 0 & 0 \\
 \lambda_4 & \lambda_3 & 0 & 0 \\
 0 & 0 & 1 & 0 \\
 0 & 0 & 0 & 1
\end{array}
\right)  \in \GL(4, k) . 
\end{align*}
We can express $\Inn_P \circ \varphi$ as 
\[
 (\Inn_P \circ \varphi)(t) 
= 
\left(
\begin{array}{c c c c}
 1 & 0 & 0 & t^{p^{e_2}} \\
 0 & 1 & t^{p^{e_3}} & 0 \\
 0 & 0 & 1 & 0 \\
 0 & 0 & 0 & 1
\end{array}
\right) . 
\] 
So, let $\varphi^* := \Inn_P \circ \varphi$ 
and $\omega^* := \Inn_P \circ \omega$. 
Clearly, $\omega^* = \omega$. 
The pair $(\varphi^*, \omega^*)$ has the form (XV).

\item Write $a_i = \lambda_i \, t^{p^{e_i}}$ $(i = 1, 3, 4)$, 
where $\lambda_1, \lambda_3, \lambda_4 \in k \backslash \{ 0 \}$ and $e_1, e_3, e_4 \geq 0$. 
So, 
\[
 \varphi(t) 
= 
\left(
\begin{array}{c c c c}
 1 & 0 & \lambda_1 \, t^{p^{e_1}} & 0 \\
 0 & 1 & \lambda_3 \, t^{p^{e_3}} & \lambda_4 \, t^{p^{e_4}}  \\
 0 & 0 & 1 & 0 \\
 0 & 0 & 0 & 1
\end{array}
\right) . 
\] 
We have $d_1 - d_3 = 2 \, p^{e_1}$, 
$d_2 - d_3 = 2 \, p^{e_3}$, 
$d_2 - d_4 = 2 \, p^{e_4}$. 
Now, we have 
\[
\left\{
\begin{array}{r @{\,} l @{\qquad} r}
 d_1 + d_2 &= 2 \, p^{e_1} &  \textcircled{\scriptsize 1} \\
 2 \, d_2 &= 2 \, p^{e_3} &  \textcircled{\scriptsize 2} \\
 d_1 + d_2 & = 2 \, p^{e_4} & \textcircled{\scriptsize 3}
\end{array}
\right.
\]
Thus 
\[
\left\{
\begin{array}{r @{\,} l @{\qquad} l}
 d_2 & = p^{e_3} & (\text{ see \textcircled{\scriptsize 3} })  , \\
 d_1 & = 2 \, p^{e_1} - p^{e_3} & (\text{ see \textcircled{\scriptsize 1} and use $d_2 = p^{e_3}$ })  , \\
 e_1  & = e_4 &  (\text{ see \textcircled{\scriptsize 1} and 
\textcircled{\scriptsize 3} }) . 
\end{array} 
\right. 
\]
Since $d_1 \geq d_2$, we have $e_1 \geq e_3$. 
We argue by separating the following cases: 
\begin{enumerate}[label = {\rm (vi.\arabic*)}]
\item $e_1 = e_3$. 

\item $e_1 > e_3$. 
\end{enumerate} 

{\bf Case (vi.1).} 
So, 
\[
\left\{
\begin{array}{r @{\,} l}
 d_1 & = p^{e_1} , \\
 d_2 & = p^{e_1} . 
\end{array}
\right.
\]
Let 
\begin{align*}
P & := \diag(\lambda_1/\lambda_3, \; 1, \; 1/\lambda_3, \; 1/\lambda_4) 
\cdot 
\left(
\begin{array}{c c c c}
 1 & 0 & 0 & 0 \\
 1 & 1 & 0 & 0 \\
 0 & 0 & 1 & 0 \\
 0 & 0 & 0 & 1
\end{array} 
\right) \\
 & = 
\left(
\begin{array}{c c c c}
 \lambda_1/\lambda_3 & 0 & 0 & 0 \\
 1 & 1 & 0 & 0 \\
 0 & 0 & 1/\lambda_3 & 0 \\
 0 & 0 & 0 & 1/\lambda_4
\end{array} 
\right)
 \in \GL(4, k) . 
\end{align*}
We can express $\Inn_P \circ \varphi$ as 
\[
 (\Inn_P \circ \varphi)(t) 
= 
\left(
\begin{array}{c c c c}
 1 & 0 & t^{p^{e_1}} & 0 \\
 0 & 1 & 0 & t^{p^{e_1}}  \\
 0 & 0 & 1 & 0 \\
 0 & 0 & 0 & 1
\end{array}
\right) . 
\] 
So, let $\varphi^* := \Inn_P \circ \varphi$ 
and $\omega^* := \Inn_P \circ \omega$. 
Clearly, $\omega^* = \omega$. 
The pair $(\varphi^*, \omega^*)$ has the form (XXII).

{\bf Case (vi.2).} 
Let 
\[
P := \diag(\lambda_1/\lambda_3, \; 1, \; 1/\lambda_3, \; 1/\lambda_4) \in \GL(4, k) . 
\] 
We can express $\Inn_P \circ \varphi$ as 
\[
 (\Inn_P \circ \varphi)(t) 
= 
\left(
\begin{array}{c c c c}
 1 & 0 & t^{p^{e_1}} & 0 \\
 0 & 1 & t^{p^{e_3}} & t^{p^{e_1}}  \\
 0 & 0 & 1 & 0 \\
 0 & 0 & 0 & 1
\end{array}
\right) . 
\] 
So, let $\varphi^* := \Inn_P \circ \varphi$ 
and $\omega^* := \Inn_P \circ \omega$. 
Clearly, $\omega^* = \omega$. 
The pair $(\varphi^*, \omega^*)$ has the form (XIII).

\item Write $a_i = \lambda_i \, t^{p^{e_i}}$ $(i = 1, 2, 3)$, 
where $\lambda_1, \lambda_2, \lambda_3 \in k \backslash \{ 0 \}$ and $e_1, e_2, e_3 \geq 0$. 
So, 
\[
 \varphi(t) 
= 
\left(
\begin{array}{c c c c}
 1 & 0 & \lambda_1 \, t^{p^{e_1}} & \lambda_2 \, t^{p^{e_2}} \\
 0 & 1 & \lambda_3 \, t^{p^{e_3}} & 0 \\
 0 & 0 & 1 & 0 \\
 0 & 0 & 0 & 1
\end{array}
\right) . 
\] 
We have $d_1 - d_3 = 2 \, p^{e_1}$, 
$d_1 - d_4 = 2 \, p^{e_2}$, 
$d_2 - d_3 = 2 \, p^{e_3}$. 
Now, we have 
\[
\left\{
\begin{array}{r @{\,} l @{\qquad} r}
 d_1 + d_2 & = 2 \, p^{e_1} &  \textcircled{\scriptsize 1} \\
 2 \, d_1 & = 2 \, p^{e_2} &  \textcircled{\scriptsize 2} \\
 2 \, d_2 & = 2 \, p^{e_3} & \textcircled{\scriptsize 3}
\end{array}
\right.
\]
Thus 
\[
\left\{
\begin{array}{r @{\,} l @{\qquad} l}
 d_1 & = p^{e_2} & (\text{ see \textcircled{\scriptsize 2} })  , \\
 d_2 & = p^{e_3} & (\text{ see \textcircled{\scriptsize 3} })  , \\
 p^{e_2} + p^{e_3} & = 2 \, p^{e_1} &  (\text{ see \textcircled{\scriptsize 1} and use $d_1 = p^{e_2}$ and 
$d_2 = p^{e_3}$ }) . 
\end{array} 
\right. 
\]
Since $d_1 \geq d_2$, we have $e_2 \geq e_3$. 
Thus $e_1 = e_2 = e_3$ and $d_1 = d_2 = p^{e_1}$. 
Let 
\begin{align*}
P & := 
\left( 
\begin{array}{c c c c}
 1 & \lambda_1/\lambda_3 & 0 & 0 \\
 0 & 1 & 0 & 0 \\
 0 & 0 & 1 & 0 \\
 0 & 0 & 0 & 1
\end{array}
\right) 
\cdot 
\diag(\lambda_2, \;  \lambda_3, \;  1, \;  1) \\
& = 
\left( 
\begin{array}{c c c c}
 \lambda_2 & \lambda_1 & 0 & 0 \\
 0 & \lambda_3& 0 & 0 \\
 0 & 0 & 1 & 0 \\
 0 & 0 & 0 & 1
\end{array}
\right) \in \GL(4, k) . 
\end{align*}
We can express $\Inn_P \circ \varphi$ as 
\[
 (\Inn_P \circ \varphi)(t) 
= 
\left(
\begin{array}{c c c c}
 1 & 0 & 0 & t^{p^{e_2}} \\
 0 & 1 & t^{p^{e_2}} & 0  \\
 0 & 0 & 1 & 0 \\
 0 & 0 & 0 & 1
\end{array}
\right) . 
\] 
So, let $\varphi^* := \Inn_P \circ \varphi$ 
and $\omega^* := \Inn_P \circ \omega$. 
Clearly, $\omega^* = \omega$. 
The pair $(\varphi^*, \omega^*)$ has the form (XV).

\item Write $a_i = \lambda_i \, t^{p^{e_i}}$ ($i = 1, 2, 3, 4$), 
where $\lambda_1, \lambda_2, \lambda_3, \lambda_4 \in k \backslash \{ 0 \}$ 
and $e_1, e_2, e_3, e_4 \geq 0$. 
So, 
\[
 \varphi(t) 
= 
\left(
\begin{array}{c c c c}
 1 & 0 & \lambda_1 \, t^{p^{e_1}} & \lambda_2 \, t^{p^{e_2}} \\
 0 & 1 & \lambda_3 \, t^{p^{e_3}} & \lambda_4 \, t^{p^{e_4}}  \\
 0 & 0 & 1 & 0 \\
 0 & 0 & 0 & 1
\end{array}
\right) . 
\] 
We have $d_1 - d_3 = 2 \, p^{e_1}$, 
$d_1 - d_4 = 2 \, p^{e_2}$, 
$d_2 - d_3 = 2 \, p^{e_3}$, 
$d_2 - d_4 = 2 \, p^{e_4}$. 
Now, we have 
\[
\left\{
\begin{array}{r @{\,} l @{\qquad} r}
 d_1 + d_2 & = 2 \, p^{e_1} &  \textcircled{\scriptsize 1} \\
 2 \, d_1 & = 2 \, p^{e_2} &  \textcircled{\scriptsize 2} \\
 2 \, d_2 & = 2 \, p^{e_3} & \textcircled{\scriptsize 3} \\
 d_1 + d_2 & = 2 \, p^{e_4} &  \textcircled{\scriptsize 4} 
\end{array}
\right.
\]
Thus 
\[
\left\{
\begin{array}{r @{\,} l @{\qquad} l}
 d_1 & = p^{e_2} & (\text{ see \textcircled{\scriptsize 2} })  , \\
 d_2 & = p^{e_3} & (\text{ see \textcircled{\scriptsize 3} })  , \\
 p^{e_2} + p^{e_3} & = 2 \, p^{e_1} &  (\text{ see \textcircled{\scriptsize 1} and use $d_1 = p^{e_2}$ and 
$d_2 = p^{e_3}$ }) , \\
 e_1 & = e_4 & (\text{ see \textcircled{\scriptsize 1} and 
\textcircled{\scriptsize 4} }) . 
\end{array} 
\right. 
\]
Since $d_1 \geq d_2$, we have $e_2 \geq e_3$. 
Since $p^{e_2} + p^{e_3} = 2 \, p^{e_1}$, we have $e_1 = e_2 = e_3$ 
and $d_1 = d_2 = p^{e_2}$.   
Let $\lambda' := (\lambda_2 \, \lambda_3 - \lambda_1 \, \lambda_4) / \lambda_3$ and let 
\begin{align*}
P & := 
\left( 
\begin{array}{c c c c}
 1 & \lambda_1/\lambda_3 & 0 & 0 \\
 0 & 1 & 0 & 0 \\
 0 & 0 & 1 & - \lambda_4/ \lambda_3 \\
 0 & 0 & 0 & 1
\end{array}
\right) 
\cdot 
\diag(\lambda', \; \lambda_3, \; 1, \; 1) \\
& = 
\left( 
\begin{array}{c c c c}
 \lambda' & \lambda_1 & 0 & 0 \\
 0 & \lambda_3 & 0 & 0 \\
 0 & 0 & 1 & - \lambda_4/ \lambda_3 \\
 0 & 0 & 0 & 1
\end{array}
\right) 
\in \GL(4, k) . 
\end{align*}
We can express $\Inn_P \circ \varphi$ as 
\[
 (\Inn_P \circ \varphi)(t) 
= 
\left(
\begin{array}{c c c c}
 1 & 0 & 0 & t^{p^{e_2}} \\
 0 & 1 & t^{p^{e_2}} & 0  \\
 0 & 0 & 1 & 0 \\
 0 & 0 & 0 & 1
\end{array}
\right) . 
\] 
So, let $\varphi^* := \Inn_P \circ \varphi$ 
and $\omega^* := \Inn_P \circ \omega$. 
Clearly, $\omega^* = \omega$. 
The pair $(\varphi^*, \omega^*)$ has the form (XV). 

\end{enumerate} 
\end{proof}

\subsubsection{$\varphi \in \cU_{[1, 1, 2]}$ and $\omega \in \Omega(4)$}

\begin{lem} 
Let $\varphi \in \cU_{[1, 1, 2]}$ and $\omega \in \Omega(4)$. 
Assume that $\psi_{\varphi, \, \omega}$ is a homomorphism. 
Then there exists an element $(\varphi^*, \omega^*)$ of $\Hom(\G_a, \SL(4, k)) \times \Omega(4)$ such that  the following conditions 
{\rm (1)} and {\rm (2)} hold true: 
\begin{enumerate}[label = {\rm (\arabic*)}]
\item $(\varphi, \, \omega) \sim (\varphi^*, \, \omega^*)$. 

\item $(\varphi^*, \, \omega^*)$ has one of the forms 
{\rm (IX)}, {\rm (XVIII)}, {\rm (XIX)}, {\rm (XX)}. 
\end{enumerate} 
\end{lem}

\begin{proof}
Let $\phi := {^\tau\!}\varphi$. 
So, $\phi \in \mf{U}_{[2, 1, 1]}$ and 
$\psi_{\phi, \, \omega}$ is a homomorphism (see Lemma 1.11). 
By Lemma 3.8, there exists an element $(\phi', \omega')$ of 
$\Hom(\G_a, \SL(4, k)) \times \Omega(4)$ such that the following 
conditions (i) and (ii) hold true: 
\begin{enumerate}[label = {\rm (\roman*)}]
\item $(\phi, \, \omega) \sim (\phi', \, \omega')$. 

\item $(\phi', \, \omega')$ has one of the forms {\rm (IX)}, {\rm (X)}, {\rm (XI)}, {\rm (XII)}. 
\end{enumerate} 
Let $Q$ be a regular matrix of $\GL(4, k)$ such that 
\[
 \Inn_Q \circ \psi_{\phi, \, \omega} = \psi_{\phi', \, \omega'} 
\]
and let $P := {^\tau\!}(Q^{-1})$. 
We know from Lemma 1.7 that $\omega' = \omega$ and 
\[
 \Inn_P \circ \psi_{\varphi, \, \omega} 
 = \psi_{{^\tau\!}(\phi'), \; \omega} . 
\]
So, $(\varphi, \omega) \sim ({^\tau\!}(\phi'), \, \omega)$. 
If $(\phi', \omega)$ has one of the forms (X), (XI), (XII), then 
$({^\tau\!}(\phi'), \, \omega)$ has one of the forms (XVIII), (XIX), (XX), respectively. 
If $(\phi', \omega)$ has the form (IX), we let 
\[
 P' := 
 \left(
 \begin{array}{c c c c}
  1 & 0 & 0 & 0 \\
  0 & 0 & 1 & 0 \\
  0 & 1 & 0 & 0 \\
  0 & 0 & 0 & 1
 \end{array}
 \right) \in \GL(4, k) 
\]
and then have 
\[
 \Inn_{P'} \circ  \Inn_P \circ \psi_{\varphi, \, \omega} 
 = \psi_{\phi', \, \omega} . 
\]
Thus $(\varphi, \omega) \sim (\phi', \omega)$ and $(\phi', \omega)$ has the form (IX). 
\end{proof}

\subsubsection{$\varphi \in \cU_{[1, 1, 1, 1]}$ and $\omega \in \Omega(4)$}

\begin{lem} 
Let $\varphi \in \cU_{[1, 1, 1, 1]}$ and $\omega \in \Omega(4)$. 
Assume that $\psi_{\varphi, \, \omega}$ is a homomorphism. 
Then there exists an element $(\varphi^*, \omega^*)$ of $\Hom(\G_a, \SL(4, k)) \times \Omega(4)$ such that  the following conditions 
{\rm (1)} and {\rm (2)} hold true: 
\begin{enumerate}[label = {\rm (\arabic*)}]
\item $(\varphi, \, \omega) \sim (\varphi^*, \, \omega^*)$. 

\item $(\varphi^*, \, \omega^*)$ has one of the forms 
{\rm (XI)}, 
{\rm (XIX)}, {\rm (XXI)}, {\rm (XXII)}, {\rm (XXIII)}, {\rm (XXIV)}, 
{\rm (XXV)}, {\rm (XXVI)}. 
\end{enumerate} 
\end{lem}

\begin{proof}
Since $\varphi \in \cU_{[1, 1, 1, 1]}$, we can express $\varphi$ as 
\begin{align*}
\varphi(t) 
 = 
\left(
\begin{array}{c c c c}
 1 & 0 & a_1 & a_2 \\
 0 & 1 & 0 & a_3 \\
 0 & 0 & 1 & 0 \\
 0 & 0 & 0 & 1 
\end{array}
\right) 
\qquad 
(\, a_1, a_2, a_3 \in \sfP  \,) . 
\end{align*}
Since $\omega \in \Omega(4)$, we can express $\omega$ as 
\[
 \omega(u) 
 = \diag(u^{d_1}, \, u^{d_2}, \, u^{d_3}, \, u^{d_4} ) , 
\qquad d_1 \geq d_2 \geq d_3 \geq d_4  , 
\qquad d_3 = - d_2 , 
\qquad d_4 = - d_1 .  
\]
We argue by separting the following cases: 
\[
\begin{array}{c | @{\quad} r @{\qquad} r @{\qquad} r }
 & a_1 & a_2 & a_3 \\
\hline 
 {\rm (i)} &  0 & 0 & 0 \\
 {\rm (ii)} & 0 & 0 & \ne 0 \\
 {\rm (iii)} & 0 & \ne 0 & 0 \\
 {\rm (iv)} & \ne 0 & 0 & 0 \\
 {\rm (v)} & 0 & \ne 0 & \ne 0 \\
 {\rm (vi)} & \ne 0 & 0 & \ne 0 \\
 {\rm (vii)} &  \ne 0 & \ne 0 &  0 \\
 {\rm (viii)} & \ne 0 & \ne 0 & \ne 0 \\
\end{array}
\]

\begin{enumerate}[label = {\rm (\roman*)}]
\item So, $\varphi(t) = I_4$. 
Let $(\varphi^*, \omega^*) := (\varphi, \omega)$. 
Then the pair $(\varphi^*, \omega^*)$ has the form (XXVI).

\item Write $a_3 = \lambda_3 \, t^{p^{e_3}}$, where $\lambda_3 \in k \backslash \{ 0 \}$ and $e_3 \geq 0$. So, 
\[
 \varphi(t) 
= 
\left(
\begin{array}{c c c c}
 1 & 0 & 0 & 0 \\
 0 & 1 & 0 & \lambda_3 \, t^{p^{e_3}}  \\
 0 & 0 & 1 & 0 \\
 0 & 0 & 0 & 1
\end{array}
\right) . 
\] 
We have $d_2 - d_4 = 2 \, p^{e_3}$. 
So, $d_1 + d_2 = 2 \, p^{e_3}$. 
Since $d_1 \geq d_2 \geq 0$, we have 
\[
\left\{
\begin{array}{l}
 2 \, p^{e_3} \geq d_1 \geq p^{e_3} , \\
 d_2 = 2 \, p^{e_3} - d_1 .   
\end{array}
\right.
\]
Let 
\[
 P := \diag(1, \; \lambda_3, \;  1, \; 1) \in \GL(4, k) .
\]
We can express $\Inn_P \circ \varphi$ as 
\[
 (\Inn_P \circ \varphi)(t) 
= 
\left(
\begin{array}{c c c c}
 1 & 0 & 0 & 0 \\
 0 & 1 & 0 & t^{p^{e_3}} \\
 0 & 0 & 1 & 0 \\
 0 & 0 & 0 & 1
\end{array}
\right) . 
\] 
So, let $\varphi^* := \Inn_P \circ \varphi$ 
and $\omega^* := \Inn_P \circ \omega$. 
Clearly, $\omega^* = \omega$. 
The pair $(\varphi^*, \omega^*)$ has the form (XXV).

\item Write $a_2 = \lambda_2 \, t^{p^{e_2}}$, where $\lambda_2 \in k \backslash \{ 0 \}$ and $e_2 \geq 0$. So, 
\[
 \varphi(t) 
= 
\left(
\begin{array}{c c c c}
 1 & 0 & 0 & \lambda_2 \, t^{p^{e_2}} \\
 0 & 1 & 0 & 0  \\
 0 & 0 & 1 & 0 \\
 0 & 0 & 0 & 1
\end{array}
\right) . 
\] 
We have $d_1 - d_4 = 2 \, p^{e_2}$.  
So, $d_1 = p^{e_2}$. 
Let 
\[
 P := \diag(\lambda_2, \; 1, \;  1, \; 1) \in \GL(4, k) .
\]
We can express $\Inn_P \circ \varphi$ as 
\[
 (\Inn_P \circ \varphi)(t) 
= 
\left(
\begin{array}{c c c c}
 1 & 0 & 0 & t^{p^{e_2}} \\
 0 & 1 & 0 & 0 \\
 0 & 0 & 1 & 0 \\
 0 & 0 & 0 & 1
\end{array}
\right) . 
\] 
So, let $\varphi^* := \Inn_P \circ \varphi$ 
and $\omega^* := \Inn_P \circ \omega$. 
Clearly, $\omega^* = \omega$. 
The pair $(\varphi^*, \omega^*)$ has the form (XXIV).

\item Write $a_1 = \lambda_1 \, t^{p^{e_1}}$, where $\lambda_1 \in k \backslash \{ 0 \}$ and $e_1 \geq 0$. So, 
\[
 \varphi(t) 
= 
\left(
\begin{array}{c c c c}
 1 & 0 & \lambda_1 \, t^{p^{e_1}} & 0 \\
 0 & 1 & 0 & 0  \\
 0 & 0 & 1 & 0 \\
 0 & 0 & 0 & 1
\end{array}
\right) . 
\] 
We have $d_1 - d_3 = 2 \, p^{e_1}$.  
So, $d_1 + d_2 = 2 \, p^{e_1}$. Since $d_1 \geq d_2 \geq 0$, we have 
\[
\left\{
\begin{array}{l}
 2 \, p^{e_1} \geq d_1 \geq p^{e_1} , \\
 d_2 = 2 \, p^{e_1} - d_1 .   
\end{array}
\right.
\]
Let 
\[
 P := \diag(\lambda_1, \; 1, \;  1, \; 1) \in \GL(4, k) .
\]
We can express $\Inn_P \circ \varphi$ as 
\[
 (\Inn_P \circ \varphi)(t) 
= 
\left(
\begin{array}{c c c c}
 1 & 0 & t^{p^{e_1}} & 0 \\
 0 & 1 & 0 & 0 \\
 0 & 0 & 1 & 0 \\
 0 & 0 & 0 & 1
\end{array}
\right) . 
\] 
So, let $\varphi^* := \Inn_P \circ \varphi$ 
and $\omega^* := \Inn_P \circ \omega$. 
Clearly, $\omega^* = \omega$. 
The pair $(\varphi^*, \omega^*)$ has the form (XXIII).

\item Write $a_i = \lambda_i \, t^{p^{e_i}}$ ($i = 2, 3$), where $\lambda_2, \lambda_3 \in k \backslash \{ 0 \}$ and $e_2, e_3 \geq 0$. 
So, 
\[
 \varphi(t) 
= 
\left(
\begin{array}{c c c c}
 1 & 0 & 0 & \lambda_2 \, t^{p^{e_2}} \\
 0 & 1 & 0 & \lambda_3 \, t^{p^{e_3}}  \\
 0 & 0 & 1 & 0 \\
 0 & 0 & 0 & 1
\end{array}
\right) . 
\] 
We have $d_1 - d_4 = 2 \, p^{e_2}$ and $d_2 - d_4 = 2 \, p^{e_3}$.  
Now, we have 
\[
\left\{
\begin{array}{r @{\,} l}
 2 \, d_1 & = 2 \, p^{e_1}  , \\
 d_1 + d_2 & = 2 \, p^{e_3} .    
\end{array}
\right.
\]
Thus 
\[
\left\{
\begin{array}{r @{\,} l}
 d_1 & = p^{e_2} , \\
 d_2 & = 2 \, p^{e_3} - p^{e_2} .    
\end{array}
\right.
\]
Since $d_1 \geq d_2$, we have $e_2 \geq e_3$. 
Since $d_2 \geq 0$, we have $2 \, p^{e_3} \geq p^{e_2}$. 
So, $2 \geq p^{e_2 - e_3} \geq 1$, which implies that 
one of the following cases can occur: 
\begin{enumerate}[label = {\rm (v.\arabic*)}]
\item $e_2 = e_3 + 1$ and $p = 2$. 

\item $e_2 = e_3$ and $p \geq 2$. 
\end{enumerate}

{\bf Case (v.1).} We have 
\[
\left\{
\begin{array}{r @{\,} l}
 d_1 & = p^{e_3 + 1} , \\
 d_2 & = 0 .    
\end{array}
\right.
\]
Let 
\[
 P := \diag(\lambda_2, \; \lambda_3, \;  1, \; 1) 
\cdot 
\left(
\begin{array}{c c c c}
 1 & 0 & 0 & 0 \\
 0 & 0 & 1 & 0 \\
 0 & 1 & 0 & 0 \\
 0 & 0 & 0 & 1
\end{array}
\right) 
 = 
\left(
\begin{array}{c c c c}
 \lambda_2 & 0 & 0 & 0 \\
 0 & 0 & \lambda_3 & 0 \\
 0 & 1 & 0 & 0 \\
 0 & 0 & 0 & 1
\end{array}
\right) \in \GL(4, k) .
\]
We can express $\Inn_P \circ \varphi$ as 
\[
 (\Inn_P \circ \varphi)(t) 
= 
\left(
\begin{array}{c c c c}
 1 & 0 & 0 & t^{p^{e_3 + 1}} \\
 0 & 1 & 0 & 0 \\
 0 & 0 & 1 & t^{p^{e_3}} \\
 0 & 0 & 0 & 1
\end{array}
\right). 
\] 
So, let $\varphi^* := \Inn_P \circ \varphi$ 
and $\omega^* := \Inn_P \circ \omega$. 
Clearly, $\omega^* = \omega$. 
The pair $(\varphi^*, \omega^*)$ has the form (XIX),

{\bf Case (v.2).} 
We have 
\[
\left\{
\begin{array}{r @{\,} l}
 d_1 & = p^{e_3}, \\
 d_2 & = p^{e_3} .    
\end{array}
\right.
\]Let 
\begin{align*}
 P & : = 
\diag(\lambda_2, \; \lambda_3, \; 1, \; 1) 
\cdot 
\left(
\begin{array}{c c c c}
 1 & 1 & 0 & 0 \\
 0 & 1 & 0 & 0 \\
 0 & 0 & 1 & 0 \\
 0 & 0 & 0 & 1
\end{array}
\right) 
 = 
\left(
\begin{array}{c c c c}
 \lambda_2 & \lambda_2 & 0 & 0 \\
 0 & \lambda_3 & 0 & 0 \\
 0 & 0 & 1 & 0 \\
 0 & 0 & 0 & 1
\end{array}
\right) 
\in \GL(4, k) . 
\end{align*}
We can express $\Inn_P \circ \varphi$ as 
\[
 (\Inn_P \circ \varphi)(t) 
= 
\left(
\begin{array}{c c c c}
 1 & 0 & 0 & 0 \\
 0 & 1 & 0 & t^{p^{e_3}} \\
 0 & 0 & 1 & 0 \\
 0 & 0 & 0 & 1
\end{array}
\right) . 
\] 
So, let $\varphi^* := \Inn_P \circ \varphi$ 
and $\omega^* := \Inn_P \circ \omega$. 
Clearly, $\omega^* = \omega$. 
The pair $(\varphi^*, \omega^*)$ has the form (XXV).

\item Write $a_i = \lambda_i \, t^{p^{e_i}}$ $(i = 1, 3)$, where $\lambda_1, \lambda_3 \in k \backslash \{ 0 \}$ and $e_1, e_3 \geq 0$. 
So, 
\[
 \varphi(t) 
= 
\left(
\begin{array}{c c c c}
 1 & 0 & \lambda_1 \, t^{p^{e_1}} & 0 \\
 0 & 1 & 0 & \lambda_3 \, t^{p^{e_3}}  \\
 0 & 0 & 1 & 0 \\
 0 & 0 & 0 & 1
\end{array}
\right) . 
\] 
We have $d_1 - d_3 = 2 \, p^{e_1}$ and $d_2 - d_4 = 2 \, p^{e_3}$. 
Now, we have 
\[
\left\{
\begin{array}{r @{\,} l @{\qquad} r}
 d_1 + d_2 &= 2 \, p^{e_1} &  \textcircled{\scriptsize 1} \\
 d_1 + d_2 &= 2 \, p^{e_3} &  \textcircled{\scriptsize 2}     
\end{array}
\right.
\]
So, $e_1 = e_3$. Since $d_1 \geq d_2 \geq 0$, we have 
\[
\left\{
\begin{array}{l}
 2 \, p^{e_1} \geq d_1 \geq p^{e_1} , \\
 d_2 = 2 \, p^{e_1} - d_1 . 
\end{array}
\right. 
\]
Let 
\[
  P := \diag(\lambda_1, \; \lambda_3, \; 1, \; 1) \in \GL(4, k) .  
\]
We can express $\Inn_P \circ \varphi$ as 
\[
 (\Inn_P \circ \varphi)(t) 
= 
\left(
\begin{array}{c c c c}
 1 & 0 & t^{p^{e_1}} & 0 \\
 0 & 1 & 0 & t^{p^{e_1}} \\
 0 & 0 & 1 & 0 \\
 0 & 0 & 0 & 1
\end{array}
\right) . 
\] 
So, let $\varphi^* := \Inn_P \circ \varphi$ 
and $\omega^* := \Inn_P \circ \omega$. 
Clearly, $\omega^* = \omega$. 
The pair $(\varphi^*, \omega^*)$ has the form (XXII).

\item Write $a_i = \lambda_i \, t^{p^{e_i}}$ ($i = 1, 2$), 
where $\lambda_1, \lambda_2 \in k \backslash \{ 0 \}$ and $e_1, e_2 \geq 0$. 
So, 
\[
 \varphi(t) 
= 
\left(
\begin{array}{c c c c}
 1 & 0 & \lambda_1 \, t^{p^{e_1}} & \lambda_2 \, t^{p^{e_2}} \\
 0 & 1 & 0 & 0 \\
 0 & 0 & 1 & 0 \\
 0 & 0 & 0 & 1
\end{array}
\right) . 
\] 
We have $d_1 - d_3 = 2 \, p^{e_1}$ and $d_1 - d_4 = 2 \, p^{e_2}$. 
Now, we have 
\[
\left\{
\begin{array}{r @{\,} l @{\qquad} r}
 d_1 + d_2 &= 2 \, p^{e_1} &  \textcircled{\scriptsize 1} \\
 2 \, d_1 &= 2 \, p^{e_2} &  \textcircled{\scriptsize 2}     
\end{array}
\right.
\]
Thus 
\[
\left\{
\begin{array}{r @{\,} l}
 d_1 & = p^{e_2} , \\
 d_2 & = 2 \, p^{e_1} - p^{e_2} .    
\end{array}
\right.
\]
Since $d_1 \geq d_2$, we have $e_2 \geq e_1$. 
Since $d_2 \geq 0$, we have $2 \, p^{e_1} \geq p^{e_2}$, 
So, $2 \geq p^{e_2 - e_1} \geq 1$, which implies that one of the following 
cases can occur: 
\begin{enumerate}[label = {\rm (vii.\arabic*)}]
\item $e_2 = e_1 + 1$ and $p = 2$. 

\item $e_2 = e_1$ and $p \geq 2$. 
\end{enumerate}

{\bf Case (vii.1).} We have 
\[
\left\{
\begin{array}{r @{\,} l}
 d_1 & = p^{e_1 + 1} , \\
 d_2 & = 0 . 
\end{array}
\right.
\]
Let 
\begin{align*}
  P & := \diag(1, \; 1, \; 1 / \lambda_1, \; 1 / \lambda_2) 
\cdot 
\left(
\begin{array}{c c c c}
 1 & 0 & 0 & 0 \\
 0 & 0 & 1 & 0 \\
 0 & 1 & 0 & 0 \\
 0 & 0 & 0 & 1
\end{array}
\right) \\
 & = 
\left(
\begin{array}{c c c c}
 1 & 0 & 0 & 0 \\
 0 & 0 & 1 & 0 \\
 0 & 1 / \lambda_1 & 0 & 0 \\
 0 & 0 & 0 & 1 / \lambda_2 
\end{array}
\right)
\in \GL(4, k) .  
\end{align*}
We can express $\Inn_P \circ \varphi$ as 
\[
 (\Inn_P \circ \varphi)(t) 
= 
\left(
\begin{array}{c c c c}
 1 & t^{p^{e_1}} & 0 & t^{p^{e_1 + 1}} \\
 0 & 1 & 0 & 0 \\
 0 & 0 & 1 & 0 \\
 0 & 0 & 0 & 1
\end{array}
\right) . 
\] 
So, let $\varphi^* := \Inn_P \circ \varphi$ 
and $\omega^* := \Inn_P \circ \omega$. 
Clearly, $\omega^* = \omega$. 
The pair $(\varphi^*, \omega^*)$ has the form (XI).

{\bf Case (vii.2).} We have 
\[
\left\{
\begin{array}{r @{\,} l}
 d_1 & = p^{e_1} , \\
 d_2 & = p^{e_1} . 
\end{array}
\right.
\]
Let 
\begin{align*}
  P & := \diag(1, \; 1, \; 1 / \lambda_1, \; 1 / \lambda_2) 
 \cdot 
\left(
\begin{array}{c c c c}
 1 & 0 & 0 & 0 \\
 0 & 1 & 0 & 0 \\
 0 & 0 & 1 & -1 \\
 0 & 0 & 0 & 1
\end{array}
\right) \\
 & = 
\left(
\begin{array}{c c c c}
 1 & 0 & 0 & 0 \\
 0 & 1 & 0 & 0 \\
 0 & 0 & 1 / \lambda_1 & - 1 / \lambda_1 \\
 0 & 0 & 0 & 1 / \lambda_2 
\end{array}
\right)
\in \GL(4, k) .  
\end{align*}
We can express $\Inn_P \circ \varphi$ as 
\[
 (\Inn_P \circ \varphi)(t) 
= 
\left(
\begin{array}{c c c c}
 1 & 0 & t^{p^{e_1}} & 0  \\
 0 & 1 & 0 & 0 \\
 0 & 0 & 1 & 0 \\
 0 & 0 & 0 & 1
\end{array}
\right) . 
\] 
So, let $\varphi^* := \Inn_P \circ \varphi$ 
and $\omega^* := \Inn_P \circ \omega$. 
Clearly, $\omega^* = \omega$. 
The pair $(\varphi^*, \omega^*)$ has the form (XXIII).

\item Write $a_i = \lambda_i \, t^{p^{e_i}}$ ($i = 1, 2, 3$), 
where $\lambda_1, \lambda_2, \lambda_3 \in k \backslash \{ 0 \}$ and $e_1, e_2, e_3 \geq 0$. 
So, 
\[
 \varphi(t) 
= 
\left(
\begin{array}{c c c c}
 1 & 0 & \lambda_1 \, t^{p^{e_1}} & \lambda_2 \, t^{p^{e_2}} \\
 0 & 1 & 0 & \lambda_3 \, t^{p^{e_3}}  \\
 0 & 0 & 1 & 0 \\
 0 & 0 & 0 & 1
\end{array}
\right) . 
\] 
We have $d_1 - d_3 = 2 \, p^{e_1}$, 
$d_1 - d_4 = 2 \, p^{e_2}$, 
$d_2 - d_4 = 2 \, p^{e_3}$. 
Now, we have 
\[
\left\{
\begin{array}{r @{\,} l @{\qquad} r}
 d_1 + d_2 & = 2 \, p^{e_1} &  \textcircled{\scriptsize 1} \\
 2 \, d_1 & = 2 \, p^{e_2} &  \textcircled{\scriptsize 2} \\
 d_1 + d_2 & = 2 \, p^{e_3} & \textcircled{\scriptsize 3}
\end{array}
\right.
\]
Thus 
\[
\left\{
\begin{array}{r @{\,} l @{\qquad} l}
 d_1 & = p^{e_2} & (\text{ see \textcircled{\scriptsize 2} })  , \\
 d_2 & = 2 \, p^{e_1} - p^{e_2} & (\text{ see \textcircled{\scriptsize 1} and 
use $d_1 = p^{e_2}$ })  , \\
 e_1 & = e_3 &  (\text{ see \textcircled{\scriptsize 1} and \textcircled{\scriptsize 3} }) . 
\end{array} 
\right. 
\]
Since $d_1 \geq d_2$, we have $e_2 \geq e_1$. 
Since $d_2 \geq 0$, we have $2 \, p^{e_1} \geq p^{e_2}$. 
So, $2 \geq p^{e_2 - e_1} \geq 1$, which implies that 
one of the following cases can occur: 
\begin{enumerate}[label = {\rm (viii.\arabic*)}]
\item $e_2 = e_1 + 1$ and $p = 2$. 

\item $e_2 = e_1$ and $p \geq 2$. 
\end{enumerate}

{\bf Case (viii.1).} We have 
\[
\left\{
\begin{array}{r @{\,} l}
 d_1 & = p^{e_1 + 1} , \\
 d_2 & = 0 . 
\end{array}
\right.
\]
Let 
\begin{align*}
P & := 
\diag(\lambda_1, \;  (\lambda_1 \, \lambda_3) / \lambda_2, \;  1, 
\;  \lambda_1 / \lambda_2) \in \GL(4, k) . 
\end{align*}
We can express $\Inn_P \circ \varphi$ as 
\[
 (\Inn_P \circ \varphi)(t) 
= 
\left(
\begin{array}{c c c c}
 1 & 0 & t^{p^{e_1}} & t^{p^{e_1 + 1}} \\
 0 & 1 & 0 & t^{p^{e_1}}  \\
 0 & 0 & 1 & 0 \\
 0 & 0 & 0 & 1
\end{array}
\right) . 
\] 
So, let $\varphi^* := \Inn_P \circ \varphi$ 
and $\omega^* := \Inn_P \circ \omega$. 
Clearly, $\omega^* = \omega$. 
The pair $(\varphi^*, \omega^*)$ has the form (XXI).

{\bf Case (viii.2).} We have 
\[
\left\{
\begin{array}{r @{\,} l}
 d_1 & = p^{e_1} , \\
 d_2 & = p^{e_1} . 
\end{array}
\right.
\]
Let 
\begin{align*}
  P & := \diag(\lambda_1, \;  (\lambda_1 \, \lambda_3) / \lambda_2, \;  1, 
\;  \lambda_1 / \lambda_2) \cdot 
\left(
\begin{array}{c c c c}
 1 & 0 & 0 & 0 \\
 0 & 1 & 0 & 0 \\
 0 & 0 & 1 & -1 \\
 0 & 0 & 0 & 1
\end{array}
\right) \\
 & = 
\left(
\begin{array}{c c c c}
 \lambda_1 & 0 & 0 & 0 \\
 0 & (\lambda_1 \, \lambda_3) / \lambda_2 & 0& 0 \\
 0 & 0 & 1 & - 1 \\
 0 & 0 & 0 & \lambda_1 / \lambda_2 
\end{array}
\right)
\in \GL(4, k) .  
\end{align*}
We can express $\Inn_P \circ \varphi$ as 
\[
 (\Inn_P \circ \varphi)(t) 
= 
\left(
\begin{array}{c c c c}
 1 & 0 & t^{p^{e_1}} & 0  \\
 0 & 1 & 0 & t^{p^{e_1}} \\
 0 & 0 & 1 & 0 \\
 0 & 0 & 0 & 1
\end{array}
\right) . 
\] 
So, let $\varphi^* := \Inn_P \circ \varphi$ 
and $\omega^* := \Inn_P \circ \omega$. 
Clearly, $\omega^* = \omega$. 
The pair $(\varphi^*, \omega^*)$ has the form (XXII).

\end{enumerate}

\end{proof}

%
%

\section{On extending antisymmentric homomorphisms $\rB(2, k) \to \SL(4, k)$ 
to $\SL(2, k) \to \SL(4, k)$}

\subsection{The forms of homomorphisms $\phi^- : \G_a \to \SL(n, k)$}

Given an antisymmetric homomorphism $\psi : \rB(2, k) \to \SL(n, k)$, 
we can express $\psi$ as $\psi \circ \jmath = \psi_{\varphi, \, \omega}$ for some 
$(\varphi, \omega) \in \cU_n \times \Omega(n)$. 
If $\psi$ is extendable to a homomorphism $\sigma : \SL(2, k) \to \SL(n, k)$, then 
the following conditions {\rm (i)} and {\rm (ii)} hold true: 
\begin{enumerate}[label = {\rm (\roman*)}]
\item Let $\varphi^- : \G_a \to \SL(n, k)$ be the homomorphism defined by 
\[
 \varphi^-(s) 
  := 
  \sigma
  \left(
  \begin{array}{c c}
   1 & 0 \\
   s & 1
  \end{array}
  \right) . 
\]
Then, for any $s \in \G_a$, the regular matrix $\varphi^-(s)$ is 
a lower triangular matrix.

\item We have 
\[
 \varphi(t) \, \varphi^-(s) 
 = 
\varphi^-\left( \frac{s}{1 + t \, s} \right) 
\, 
\omega(1 + t \, s) 
\;
\varphi \left( \frac{t}{1 + t \, s} \right)
\]
for all $t, s \in k$ with $1 + t \, s \ne 0$ (see Lemma 1.18). 
\end{enumerate}

In this section, for any antisymmetric pair $(\varphi^*, \omega^*)$ of the form $(\nu)$, 
where $\nu = {\rm I}, {\rm II}, \ldots, {\rm XXVI}$, assuming that  
there exists a homomorphism $\phi^- : \G_a \to \SL(4, k)$ satisfying 
the following conditions {\rm (i)} and {\rm (ii)}: 
\begin{enumerate}[label = {\rm (\roman*)}]
\item For any $s \in \G_a$, the regular matrix $\phi^-(s)$ is a lower triangular matrix.  

\item 
$
 \varphi^*(t) \, \phi^-(s) 
 = 
\phi^-\left( \frac{s}{1 + t \, s} \right) 
\, 
\omega^*(1 + t \, s) 
\;
\varphi^* \left( \frac{t}{1 + t \, s} \right)
$ 
for all $t, s \in k$ with $1 + t \, s \ne 0$, 
\end{enumerate} 
we then express $\phi^-$ or conclude a contradition (i.e., there exists no homomorphism 
$\phi^-$ satisfying (i) and (ii)).

\subsubsection{$\rm (I)$}

\begin{lem}
Let $(\varphi^*, \omega^*)$ be of the form {\rm (I)}. 
Assume that there exists a homomorphism $\phi^- : \G_a \to \SL(4, k)$ 
satisfying the following conditions {\rm (i)} and {\rm (ii)}: 
\begin{enumerate}[label = {\rm (\roman*)}]
\item For any $s \in \G_a$, 
the regular matrix $\phi^-(s)$ is a lower triangular matrix.  

\item 
$
 \varphi^*(t) \, \phi^-(s) 
 = 
\phi^-\left( \frac{s}{1 + t \, s} \right) 
\, 
\omega^*(1 + t \, s) 
\;
\varphi^* \left( \frac{t}{1 + t \, s} \right)
$ 
for all $t, s \in k$ with $1 + t \, s \ne 0$.    
\end{enumerate} 
Then we can express $\phi^-$ as  
\[
\phi^-(s) 
= 
\left(
\begin{array}{c c c c}
 1 & 0 & 0 & 0 \\
 3 \, s^{p^{e_1}} & 1 & 0 & 0 \\
 6 \, s^{2 \, p^{e_1}} & 4 \, s^{p^{e_1}} & 1 & 0 \\
 6 \, s^{3 \, p^{e_1}} & 6 \, s^{2 \, p^{e_1}} & 3 \, s^{p^{e_1}} & 1  
\end{array}
\right) . 
\]
\end{lem}

\begin{proof} 
We can write $\phi^- : \G_a \to \SL(4, k)$ as 
\[
 \phi^-(s) = \bigl( \, b_{i, j}(s) \, \bigr)_{1 \leq i, j \leq 4} , 
\]
where the polynomials $b_{i, j}(s) \in k[s]$ ($1 \leq i, j \leq 4$) satisfy 
the following conditions (a) and (b): 
\begin{enumerate}[label = {\rm (\alph*)}]
\item $b_{i, i}(s) = 1$ for all $1 \leq i \leq 4$. 
\item $b_{i, j}(s) = 0$ for all $1 \leq i < j \leq 4$.  
\end{enumerate} 

By condition (ii), we have 
\begin{align*}
& \renewcommand{\arraystretch}{1.5} 
\left(
\begin{array}{r r r r}
 1 + b_{2, 1} + \frac{1}{2} \, b_{3, 1} + \frac{1}{6} \, b_{4, 1} 
 & 1 + \frac{1}{2} \, b_{3, 2} + \frac{1}{6} \, b_{4, 2} 
 & \frac{1}{2} + \frac{1}{6} \, b_{4, 3} 
 & \frac{1}{6} \\
 b_{2, 1} + b_{3, 1} + \frac{1}{2} \, b_{4, 1} 
  & 1 + b_{3, 2} + \frac{1}{2} \, b_{4, 2} 
 & 1 + \frac{1}{2} \, b_{4, 3} 
 & \frac{1}{2} \\ 
                                       b_{3, 1} + b_{4, 1} 
                                        & b_{3, 2} + b_{4, 2} 
                                        & 1 + b_{4, 3} 
                                        & 1 \\
                                                    b_{4, 1} 
                                                    & b_{4, 2} 
                                                    & b_{4, 3} 
                                                    & 1
\end{array}
\right) \\
&
= 
\left(
\begin{array}{l l}
 (1 + s)^{3 \, p^{e_1}} 
 & (1 + s)^{2 \, p^{e_1}} \\
 b_{2, 1}\left( \frac{s}{1 + s} \right) \, (1 + s)^{3 \, p^{e_1}} 
 & b_{2, 1}\left( \frac{s}{1 + s} \right)  \, (1 + s)^{2 \, p^{e_1}} + (1 + s)^{p^{e_1}} \\
  b_{3, 1}\left( \frac{s}{1 + s} \right) \, (1 + s)^{3 \, p^{e_1}} 
 & b_{3, 1}\left( \frac{s}{1 + s} \right)  \, (1 + s)^{2 \, p^{e_1}} 
   + b_{3, 2}\left( \frac{s}{1 + s} \right)  \, (1 + s)^{p^{e_1}}  \\
   b_{4, 1}\left( \frac{s}{1 + s} \right) \, (1 + s)^{3 \, p^{e_1}} 
 & b_{4, 1}\left( \frac{s}{1 + s} \right)  \, (1 + s)^{2 \, p^{e_1}} 
 + b_{4, 2}\left( \frac{s}{1 + s} \right)  \, (1 + s)^{p^{e_1}} \\
\end{array} 
\right.
\\
& \qquad 
\left. 
\begin{array}{l }
 \frac{1}{2} \, (1 + s)^{p^{e_1}} \\
 \frac{1}{2} \, b_{2, 1} \left( \frac{s}{1 + s} \right) \, (1 + s)^{p^{e_1}} + 1 \\
  \frac{1}{2} \, b_{3, 1} \left( \frac{s}{1 + s} \right) \, (1 + s)^{p^{e_1}} 
  + b_{3, 2} \left( \frac{s}{1 + s} \right) + \frac{1}{ (1 + s)^{p^{e_1}} } \\ 
   \frac{1}{2} \, b_{4, 1} \left( \frac{s}{1 + s} \right) \, (1 + s)^{p^{e_1}} 
  + b_{4, 2} \left( \frac{s}{1 + s} \right) 
  + b_{4, 3} \left( \frac{s}{1 + s} \right)  \,\frac{1}{ (1 + s)^{p^{e_1}} } \\
\end{array}
\right. \\
& \qquad 
\left. 
\begin{array}{l}
 \frac{1}{6}  \\
 \frac{1}{6} \, b_{2, 1} \left( \frac{s}{1 + s} \right) 
  + \frac{1}{2} \, \frac{1}{ (1 + s )^{p^{e_1} } } \\
  \frac{1}{6} \, b_{3, 1} \left( \frac{s}{1 + s} \right)  
  + \frac{1}{2} \, b_{3, 2} \left( \frac{s}{1 + s} \right) \, \frac{1}{ (1 + s)^{p^{e_1}} } 
  + \frac{1}{ (1 + s)^{2 \, p^{e_1}}}   \\ 
   \frac{1}{6} \, b_{4, 1} \left( \frac{s}{1 + s} \right) 
  + \frac{1}{2} \, b_{4, 2} \left( \frac{s}{1 + s} \right) \, \frac{1}{ (1 + s)^{p^{e_1}} } 
  + b_{4, 3} \left( \frac{s}{1 + s} \right)  \, \frac{1}{ (1 + s)^{2 \, p^{e_1}}}  
  + \frac{1}{ (1 + s)^{3 \, p^{e_1}} }  
\end{array}
\right) . 
\end{align*}

Comparing the$(1, 3)$-th entries of both sides of the equality, we have 
\[
 \frac{1}{2} + \frac{1}{6} \, b_{4, 3}  = \frac{1}{2} \, (1 + s)^{p^{e_1}} , 
\]
which implies 
\[
 b_{4, 3} = 3 \, s^{p^{e_1}} . 
\]

Comparing the$(2, 4)$-th entries of both sides of the equality, we have 
\[
 \frac{1}{2}  = \frac{1}{6} \, b_{2, 1} \left( \frac{s}{1 + s} \right) 
  + \frac{1}{2} \, \frac{1}{ (1 + s )^{p^{e_1} } }, 
\]
which implies 
\[
 b_{2, 1} = 3 \, s^{p^{e_1}} . 
\] 

Comparing the$(1, 2)$-th and the $(2, 2)$-th entries of 
both sides of the equality, we have 
\[
\left\{ 
\begin{array}{r @{\,} l}
 1 + \frac{1}{2} \, b_{3, 2} + \frac{1}{6} \, b_{4, 2} 
 & = (1 + s)^{2 \, p^{e_1}} , \\
1 + b_{3, 2} + \frac{1}{2} \, b_{4, 2} 
 & =   b_{2, 1}\left( \frac{s}{1 + s} \right)  \, (1 + s)^{2 \, p^{e_1}} + (1 + s)^{p^{e_1}} . 
\end{array}
\right. 
\]
So, 
\[
\left\{ 
\begin{array}{r @{\,} l}
 \frac{1}{2} \, b_{3, 2} + \frac{1}{6} \, b_{4, 2} 
 & = 2 \, s^{p^{e_1}} + s^{2 \, p^{e_1}} , \\
 b_{3, 2} + \frac{1}{2} \, b_{4, 2} 
 & =  4 \, s^{p^{e_1}} + 3 \, s^{2 \, p^{e_1}}  . 
\end{array}
\right. 
\]
Thus
\[
b_{4, 2} = 6 \, s^{2 \, p^{e_1}} , 
\qquad 
b_{3, 2} = 4 \, s^{p^{e_1}} . 
\]

Comparing the$(3, 2)$-th entries of both sides of the equality, we have 
\[
  b_{3, 2} + b_{4, 2} 
= 
b_{3, 1}\left( \frac{s}{1 + s} \right)  \, (1 + s)^{2 \, p^{e_1}} 
   + b_{3, 2}\left( \frac{s}{1 + s} \right)  \, (1 + s)^{p^{e_1}} 
\]
which implies 
\[
 b_{3, 1} = 6 \, s^{2 \, p^{e_1}} . 
\]

Comparing the$(4, 2)$-th entries of both sides of the equality, we have 
\[
  b_{4, 2} 
=    
 b_{4, 1}\left( \frac{s}{1 + s} \right)  \, (1 + s)^{2 \, p^{e_1}} 
 + b_{4, 2}\left( \frac{s}{1 + s} \right)  \, (1 + s)^{p^{e_1}} , 
\]
which implies 
\[
 b_{4, 1} = 6 \, s^{3 \, p^{e_1}} . 
\]

Thus 
\[
\phi^-(s) 
= 
\left(
\begin{array}{c c c c}
 1 & 0 & 0 & 0 \\
 3 \, s^{p^{e_1}} & 1 & 0 & 0 \\
  6 \, s^{2 \, p^{e_1}} & 4 \, s^{p^{e_1}} & 1 & 0 \\
   6 \, s^{3 \, p^{e_1}} & 6 \, s^{2 \, p^{e_1}} & 3 \, s^{p^{e_1}} & 1 
\end{array}
\right) . 
\]

\end{proof}

Let 
\[
 P : = \diag(1/36, \; 1/12, \; 1/3, \; 1 ) \in \GL(4, k) . 
\]
Then we have 
\[
 (\Inn_P \circ \phi^-)(s) 
  = 
  \left(
\begin{array}{c c c c}
 1 & 0 & 0 & 0 \\
 s^{p^{e_1}} & 1 & 0 & 0 \\
  \frac{1}{2} \, s^{2 \, p^{e_1}} & s^{p^{e_1}} & 1 & 0 \\
   \frac{1}{6} \, s^{3 \, p^{e_1}} & \frac{1}{2} \, s^{2 \, p^{e_1}} & s^{p^{e_1}} & 1 
\end{array}
\right) . 
\]
%
%
%

\subsubsection{$\rm (II)$}

\begin{lem}
Let $(\varphi^*, \omega^*)$ be of the form {\rm (II)}. 
Assume that there exists a homomorphism $\phi^- : \G_a \to \SL(4, k)$ 
satisfying the following conditions {\rm (i)} and {\rm (ii)}: 
\begin{enumerate}[label = {\rm (\roman*)}]
\item For any $s \in \G_a$, 
the regular matrix $\phi^-(s)$ is a lower triangular matrix.  

\item 
$
 \varphi^*(t) \, \phi^-(s) 
 = 
\phi^-\left( \frac{s}{1 + t \, s} \right) 
\, 
\omega^*(1 + t \, s) 
\;
\varphi^* \left( \frac{t}{1 + t \, s} \right)
$ 
for all $t, s \in k$ with $1 + t \, s \ne 0$.    
\end{enumerate} 
Then we can express $\phi^-$ as  
\[
\phi^-(s) 
= 
\left(
\begin{array}{c c c c}
 1 & 0 & 0 & 0 \\
 0 & 1 & 0 & 0 \\
 0 & s^{p^{e_1}} & 1 & 0 \\
 s^{p^{e_1 + 1}} & s^{2 \, p^{e_1}} & \frac{1}{2} \, s^{p^{e_1}} & 1  
\end{array}
\right) . 
\]
\end{lem}

\begin{proof} 
Write $\phi^-(s) = ( \, b_{i, j}(s) \, )_{1 \leq i, j \leq 4}$. 
By condition (ii), we have 
\begin{align*}
& \renewcommand{\arraystretch}{1.5} 
\left(
\begin{array}{r r r r}
 1 + b_{2, 1} + \frac{1}{2} \, b_{3, 1} + b_{4, 1} 
 & 1 + \frac{1}{2} \, b_{3, 2} + b_{4, 2} 
 & \frac{1}{2} + b_{4, 3} 
 & 1 \\
 b_{2, 1} + b_{3, 1} 
  & 1 + b_{3, 2}  
 & 1 
 & 0 \\ 
                                       b_{3, 1} 
                                        & b_{3, 2} 
                                        & 1 
                                        & 0 \\
                                                    b_{4, 1} 
                                                    & b_{4, 2} 
                                                    & b_{4, 3} 
                                                    & 1
\end{array}
\right) \\
&
= 
\left(
\begin{array}{l l}
 (1 + s)^{p^{e_1 + 1}} 
 & (1 + s)^{2 \, p^{e_1}} \\
 b_{2, 1}\left( \frac{s}{1 + s} \right) \, (1 + s)^{p^{e_1 + 1}} 
 & b_{2, 1}\left( \frac{s}{1 + s} \right)  \, (1 + s)^{2 \, p^{e_1}} + (1 + s)^{p^{e_1}} \\
  b_{3, 1}\left( \frac{s}{1 + s} \right) \, (1 + s)^{ p^{e_1 + 1}} 
 & b_{3, 1}\left( \frac{s}{1 + s} \right)  \, (1 + s)^{2 \, p^{e_1}}
   + b_{3, 2}\left( \frac{s}{1 + s} \right)  \, (1 + s)^{p^{e_1}}  \\
   b_{4, 1}\left( \frac{s}{1 + s} \right) \, (1 + s)^{ p^{e_1 + 1}} 
 & b_{4, 1}\left( \frac{s}{1 + s} \right)  \, (1 + s)^{2 \, p^{e_1}}  
 + b_{4, 2}\left( \frac{s}{1 + s} \right)  \, (1 + s)^{p^{e_1}} \\
\end{array} 
\right.
\\
& \qquad 
\left. 
\begin{array}{l l}
 \frac{1}{2} \, (1 + s)^{p^{e_1}}  
  & 1 \\
 \frac{1}{2} \, b_{2, 1} \left( \frac{s}{1 + s} \right) \, (1 + s)^{p^{e_1}} + 1 
 &  b_{2, 1} \left( \frac{s}{1 + s} \right) \\
  \frac{1}{2} \, b_{3, 1} \left( \frac{s}{1 + s} \right) \, (1 + s)^{p^{e_1}} 
  + b_{3, 2} \left( \frac{s}{1 + s} \right) + \frac{1}{ (1 + s)^{p^{e_1}} } 
   & b_{3, 1} \left( \frac{s}{1 + s} \right) \\ 
   \frac{1}{2} \, b_{4, 1} \left( \frac{s}{1 + s} \right) \, (1 + s)^{p^{e_1}} 
  + b_{4, 2} \left( \frac{s}{1 + s} \right) 
  + b_{4, 3} \left( \frac{s}{1 + s} \right)  \, \frac{1}{ (1 + s)^{p^{e_1}} } 
  &  b_{4, 1} \left( \frac{s}{1 + s} \right) + \frac{1}{ (1 + s)^{p^{e_1 + 1}} } 
\end{array}
\right) . 
\end{align*}

Comparing the $(2, 4)$-th entries of both sides of the equality, 
we have 
\[
 0 = b_{2, 1}\left( \frac{s}{1 + s} \right) , 
\]
which implies $b_{2, 1}(s) = 0$.

Comparing the $(3, 4)$-th entries of both sides of the equality, 
we have 
\[
 0 = b_{3, 1}\left( \frac{s}{1 + s} \right)  , 
\]
which implies $b_{3, 1}(s) = 0$.

Comparing the $(4, 4)$-th entries of both sides of the equality, 
we have 
\[
 1 = b_{4, 1} \left( \frac{s}{1 + s} \right) + \frac{1}{ (1 + s)^{p^{e_1 + 1}} } , 
\]
which implies $b_{4, 1}(s) = s^{p^{e_1} + 1}$.

Comparing the $(2, 2)$-th entries of both sides of the equality, 
we have 
\[
 1 + b_{3, 2} = b_{2, 1}\left( \frac{s}{1 + s} \right)  \, (1 + s) + (1 + s)^{p^{e_1}} , 
\]
which implies $b_{3, 2}(s) = s^{p^{e_1}}$.

Comparing the $(1, 3)$-th entries of both sides of the equality, 
we have 
\[
 \frac{1}{2} + b_{4, 3}  
 = \frac{1}{2} \, (1 + s)^{p^{e_1}}  , 
\]
which implies $b_{4, 3}(s) = \frac{1}{2} \, s^{p^{e_1}}$.

Comparing the $(1, 2)$-th entries of both sides of the equality, 
we have 
\[
 1 + \frac{1}{2} \, b_{3, 2} + b_{4, 2} 
 = (1 + s)^{2 \, p^{e_1}}   , 
\]
which implies $b_{4, 2}(s) = s^{2 \, p^{e_1}}$ (use $b_{3, 2} = s^{p^{e_1}}$).

Thus $\phi^-$ has the desired form. 

\end{proof}

Let 
\[
 P : = \diag(1, \; 2, \; 2, \; 1 ) \in \GL(4, k) . 
\]
Then we have 
\[
(\Inn_P \circ \phi^-)(s) 
= 
\left(
\begin{array}{c c c c}
 1 & 0 & 0 & 0 \\
 0 & 1 & 0 & 0 \\
 0 & s^{p^{e_1}} & 1 & 0 \\
 s^{p^{e_1 + 1}} & \frac{1}{2} \, s^{2 \, p^{e_1}} & s^{p^{e_1}} & 1  
\end{array}
\right) . 
\]
%

\subsubsection{$\rm (III)$}

\begin{lem}
Let $(\varphi^*, \omega^*)$ be of the form {\rm (III)}. 
Assume that there exists a homomorphism $\phi^- : \G_a \to \SL(4, k)$ 
satisfying the following conditions {\rm (i)} and {\rm (ii)}: 
\begin{enumerate}[label = {\rm (\roman*)}]
\item For any $s \in \G_a$, 
the regular matrix $\phi^-(s)$ is a lower triangular matrix.  

\item 
$
 \varphi^*(t) \, \phi^-(s) 
 = 
\phi^-\left( \frac{s}{1 + t \, s} \right) 
\, 
\omega^*(1 + t \, s) 
\;
\varphi^* \left( \frac{t}{1 + t \, s} \right)
$ 
for all $t, s \in k$ with $1 + t \, s \ne 0$.    
\end{enumerate} 
Then we have a contradiction. 
\end{lem}

\begin{proof}
Write $\phi^-(s) = ( \, b_{i, j}(s) \, )_{1 \leq i, j \leq 4}$. 
By condition (ii), we have 
\begin{align*}
& \renewcommand{\arraystretch}{1.5} 
\left(
\begin{array}{r r r r}
 1 + b_{2, 1} + \frac{1}{2} \, b_{3, 1} 
 & 1 + \frac{1}{2} \, b_{3, 2} 
 & \frac{1}{2} 
 & 0 \\
 b_{2, 1} + b_{3, 1} 
  & 1 + b_{3, 2}  
 & 1 
 & 0 \\ 
                                       b_{3, 1} 
                                        & b_{3, 2} 
                                        & 1 
                                        & 0 \\
                                                    b_{4, 1} 
                                                    & b_{4, 2} 
                                                    & b_{4, 3} 
                                                    & 1
\end{array}
\right) \\
&
= 
\left(
\begin{array}{l l}
 (1 + s)^{3 \, p^{e_1}} 
 & (1 + s)^{2 \, p^{e_1}} \\
 b_{2, 1}\left( \frac{s}{1 + s} \right) \, (1 + s)^{3 \, p^{e_1}} 
 & b_{2, 1}\left( \frac{s}{1 + s} \right)  \, (1 + s)^{2 \, p^{e_1}} + (1 + s)^{p^{e_1}} \\
  b_{3, 1}\left( \frac{s}{1 + s} \right) \, (1 + s)^{ 3 \, p^{e_1}} 
 & b_{3, 1}\left( \frac{s}{1 + s} \right)  \, (1 + s)^{2 \, p^{e_1}} 
   + b_{3, 2}\left( \frac{s}{1 + s} \right)  \, (1 + s)^{p^{e_1}}  \\
   b_{4, 1}\left( \frac{s}{1 + s} \right) \, (1 + s)^{ 3 \, p^{e_1}} 
 & b_{4, 1}\left( \frac{s}{1 + s} \right)  \, (1 + s)^{2 \, p^{e_1}}  
 + b_{4, 2}\left( \frac{s}{1 + s} \right)  \, (1 + s)^{p^{e_1}} \\
\end{array} 
\right.
\\
& \qquad 
\left. 
\begin{array}{l l}
 \frac{1}{2} \, (1 + s)^{p^{e_1}}  
  & 0 \\
 \frac{1}{2} \, b_{2, 1} \left( \frac{s}{1 + s} \right) \, (1 + s)^{p^{e_1}} + 1 
 &  0 \\
  \frac{1}{2} \, b_{3, 1} \left( \frac{s}{1 + s} \right) \, (1 + s)^{p^{e_1}} 
  + b_{3, 2} \left( \frac{s}{1 + s} \right) + \frac{1}{ (1 + s)^{p^{e_1}} } 
   & 0 \\ 
   \frac{1}{2} \, b_{4, 1} \left( \frac{s}{1 + s} \right) \, (1 + s)^{p^{e_1}} 
  + b_{4, 2} \left( \frac{s}{1 + s} \right) 
  + b_{4, 3} \left( \frac{s}{1 + s} \right)  \, \frac{1}{ (1 + s)^{p^{e_1}} } 
  &   \frac{1}{ (1 + s)^{3 \, p^{e_1}} } 
\end{array}
\right) . 
\end{align*}

Comparing the $(4, 4)$-th entries of both sides of the equality, 
we have 
\[
 1 = \frac{1}{ (1 + s)^{3 \, p^{e_1}} }   , 
\]
which implies a contradiction.

\end{proof}

\subsubsection{$\rm (IV)$}

\begin{lem}
Let $(\varphi^*, \omega^*)$ be of the form {\rm (IV)}. 
Assume that there exists a homomorphism $\phi^- : \G_a \to \SL(4, k)$ 
satisfying the following conditions {\rm (i)} and {\rm (ii)}: 
\begin{enumerate}[label = {\rm (\roman*)}]
\item For any $s \in \G_a$, 
the regular matrix $\phi^-(s)$ is a lower triangular matrix.  

\item 
$
 \varphi^*(t) \, \phi^-(s) 
 = 
\phi^-\left( \frac{s}{1 + t \, s} \right) 
\, 
\omega^*(1 + t \, s) 
\;
\varphi^* \left( \frac{t}{1 + t \, s} \right)
$ 
for all $t, s \in k$ with $1 + t \, s \ne 0$.    
\end{enumerate} 
Then we can express $\phi^-$ as  
\[
\phi^-(s) 
= 
\left(
\begin{array}{c c c c}
 1 & 0 & 0 & 0 \\
 s^{p^{e_1}} & 1 & 0 & 0 \\
 s^{p^{e_2}} & 0 & 1 & 0 \\
 s^{p^{e_1} + p^{e_2}} & s^{p^{e_2}} & s^{p^{e_1}} & 1  
\end{array}
\right) . 
\]
\end{lem}

\begin{proof}

Write $\phi^-(s) = ( \, b_{i, j}(s) \, )_{1 \leq i, j \leq 4}$. 
By condition (ii), we have 
\begin{align*}
& \renewcommand{\arraystretch}{1.5} 
\left(
\begin{array}{r r r r}
 1 + b_{2, 1} + b_{3, 1} + b_{4, 1} 
 & 1 + b_{3, 2} + b_{4, 2} 
 & 1 + b_{4, 3} 
 & 1 \\
 b_{2, 1} + b_{4, 1} 
  & 1 + b_{4, 2}  
 & b_{4, 3}
 & 1 \\ 
                                       b_{3, 1} + b_{4, 1}
                                        & b_{3, 2} + b_{4, 2} 
                                        & 1 + b_{4, 3}
                                        & 1 \\
                                                    b_{4, 1} 
                                                    & b_{4, 2} 
                                                    & b_{4, 3} 
                                                    & 1
\end{array}
\right) \\
&
= 
\left(
\begin{array}{l l}
 (1 + s)^{p^{e_1}+ p^{e_2}} 
 & (1 + s)^{p^{e_2}} \\
 b_{2, 1}\left( \frac{s}{1 + s} \right) \, (1 + s)^{p^{e_1}+ p^{e_2}} 
 & b_{2, 1}\left( \frac{s}{1 + s} \right)  \, (1 + s)^{p^{e_2}}  + (1 + s)^{p^{e_2} - p^{e_1}} \\
  b_{3, 1}\left( \frac{s}{1 + s} \right) \, (1 + s)^{p^{e_1}+ p^{e_2}} 
 & b_{3, 1}\left( \frac{s}{1 + s} \right)  \, (1 + s)^{p^{e_2}} 
   + b_{3, 2}\left( \frac{s}{1 + s} \right)  \, (1 + s)^{p^{e_2} - p^{e_1}}  \\
   b_{4, 1}\left( \frac{s}{1 + s} \right) \, (1 + s)^{p^{e_1}+ p^{e_2}} 
 & b_{4, 1}\left( \frac{s}{1 + s} \right)  \, (1 + s)^{p^{e_2}}   
 + b_{4, 2}\left( \frac{s}{1 + s} \right)  \, (1 + s)^{p^{e_2} - p^{e_1}}  \\
\end{array} 
\right.
\\
& \qquad 
\left. 
\begin{array}{l }
  (1 + s)^{p^{e_1}}  \\
 b_{2, 1} \left( \frac{s}{1 + s} \right) \, (1 + s)^{p^{e_1}} 
 \\
  b_{3, 1} \left( \frac{s}{1 + s} \right) \, (1 + s)^{p^{e_1}} 
  + \frac{1}{(1 + s)^{p^{e_2} - p^{e_1}}}  
\\ 
   b_{4, 1} \left( \frac{s}{1 + s} \right) \, (1 + s)^{p^{e_1}}  
  + b_{4, 3} \left( \frac{s}{1 + s} \right)  \, \frac{1}{(1 + s)^{p^{e_2} - p^{e_1}}}   
\end{array}
\right. \\
& \qquad 
\left. 
\begin{array}{l}
 1 \\
  b_{2, 1} \left( \frac{s}{1 + s} \right) + \frac{1}{ (1 + s)^{p^{e_1}} }  \\
  b_{3, 1} \left( \frac{s}{1 + s} \right)  
   +  b_{3, 2} \left( \frac{s}{1 + s} \right)  \, \frac{1}{ (1 + s)^{p^{e_1}} } 
   + \frac{1}{ (1 + s)^{p^{e_2}} } \\
  b_{4, 1} \left( \frac{s}{1 + s} \right)
   + b_{4, 2} \left( \frac{s}{1 + s} \right) \, \frac{1}{ (1 + s)^{p^{e_1}} } 
   + b_{4, 3} \left( \frac{s}{1 + s} \right) \, \frac{1}{ (1 + s)^{p^{e_2}} } 
   + \frac{1}{ (1 + s)^{p^{e_1} + p^{e_2}} } 
\end{array} 
\right) . 
\end{align*}

Comparing the $(1, 3)$-th entries of both sides of the equality, 
we have 
\[
  1 + b_{4, 3}
 =  (1 + s)^{p^{e_1}}  , 
\]
which implies $b_{4, 3}(s) = s^{p^{e_1}}$.

Comparing the $(2, 3)$-th entries of both sides of the equality, 
we have 
\[
  b_{4, 3}
 =   b_{2, 1} \left( \frac{s}{1 + s} \right) \, (1 + s)^{p^{e_1}}  , 
\]
which implies $b_{2, 1}(s) = s^{p^{e_1}}$ (use $b_{4, 3}(s) = s^{p^{e_1}}$).

Comparing the $(3, 3)$-th entries of both sides of the equality, 
we have 
\[
 1 + b_{4, 3}
 =     b_{3, 1} \left( \frac{s}{1 + s} \right) \, (1 + s)^{p^{e_1}} 
  +\frac{1}{(1 + s)^{p^{e_2} - p^{e_1}}}     , 
\]
which implies 
\[
 1 + s^{p^{e_1}}
  = 
  b_{3, 1} \left( \frac{s}{1 + s} \right) \, (1 + s)^{p^{e_1}} 
  + \frac{ (1 + s)^{p^{e_1}} }{(1 + s)^{p^{e_2}}} .
\]
Thus 
\[
 1 =  b_{3, 1} \left( \frac{s}{1 + s} \right) 
  + \frac{ 1 }{(1 + s)^{p^{e_2}}} .
\]
So, $b_{3, 1}(s) = s^{p^{e_2}}$.

Comparing the $(4, 3)$-th entries of both sides of the equality, 
we have 
\[
 b_{4, 3}
 =     
   b_{4, 1} \left( \frac{s}{1 + s} \right) \, (1 + s)^{p^{e_1}}  
  + b_{4, 3} \left( \frac{s}{1 + s} \right)  \, \frac{1}{(1 + s)^{p^{e_2} - p^{e_1}}}     ,  
\]
which implies 
\[
 s^{p^{e_1}}
  = 
   b_{4, 1} \left( \frac{s}{1 + s} \right) \, (1 + s)^{p^{e_1}}  
  + \frac{ s^{p^{e_1}} }{ (1 + s)^{p^{e_1}} }   \cdot \frac{1}{(1 + s)^{p^{e_2} - p^{e_1}}}    . 
\]
Thus $b_{4, 1}(s) = s^{p^{e_1} + p^{e_2}}$.

Comparing the $(2, 2)$-th entries of both sides of the equality, 
we have 
\[
 1 + b_{4, 2} 
 = b_{2, 1} \left( \frac{s}{1 + s} \right) \, (1 + s)^{p^{e_2}} + (1 + s)^{p^{e_2} - p^{e_1}}  , 
\]
which implies 
\begin{align*}
 1 + b_{4, 2}
 & = 
  \left(
  \frac{s}{1 + s}
  \right)^{p^{e_1}}
  \, (1 + s)^{p^{e_2}}
  + (1 + s)^{p^{e_2} - p^{e_1}} \\
  & = 1 + s^{p^{e_2}}. 
\end{align*}
Thus $b_{4, 2} = s^{p^{e_2}}$.

Comparing the $(1, 2)$-th entries of both sides of the equality, 
we have 
\[
 1 + b_{3, 2} + b_{4, 2} 
 = (1 + s)^{p^{e_2}}   
\]
which implies $b_{3, 2} = 0$ (use $b_{4, 2} = s^{p^{e_2}}$).

Hence $\phi^-$ has the desired form.

\end{proof}

\subsubsection{$\rm (V)$}

\begin{lem}
Let $(\varphi^*, \omega^*)$ be of the form {\rm (V)}. 
Assume that there exists a homomorphism $\phi^- : \G_a \to \SL(4, k)$ 
satisfying the following conditions {\rm (i)} and {\rm (ii)}: 
\begin{enumerate}[label = {\rm (\roman*)}]
\item For any $s \in \G_a$, 
the regular matrix $\phi^-(s)$ is a lower triangular matrix.  

\item 
$
 \varphi^*(t) \, \phi^-(s) 
 = 
\phi^-\left( \frac{s}{1 + t \, s} \right) 
\, 
\omega^*(1 + t \, s) 
\;
\varphi^* \left( \frac{t}{1 + t \, s} \right)
$ 
for all $t, s \in k$ with $1 + t \, s \ne 0$.    
\end{enumerate} 
Then we have $p = 2$ and $f = e_1 + 1$, and 
we can express $\phi^-$ as 
\[
\phi^-(s) 
= 
\left(
\begin{array}{c c c c}
 1 & 0 & 0 & 0 \\
 0 & 1 & 0 & 0 \\
 s^{p^{e_1}} & 0 & 1 & 0 \\
 s^{p^{e_1 + 1}} & s^{p^{e_1}} & 0 & 1  
\end{array}
\right) . 
\]
\end{lem}

\begin{proof}

Write $\phi^-(s) = ( \, b_{i, j}(s) \, )_{1 \leq i, j \leq 4}$. 
By condition (ii), we have 
\begin{align*}
& \renewcommand{\arraystretch}{1.5} 
\left(
\begin{array}{r r r r}
 1 + b_{2, 1} + b_{4, 1} 
 & 1 + b_{4, 2} 
 & b_{4, 3} 
 & 1 \\
 b_{2, 1}  
  & 1 
 & 0
 & 0 \\ 
                                       b_{3, 1} + b_{4, 1}
                                        & b_{3, 2} + b_{4, 2} 
                                        & 1 + b_{4, 3}
                                        & 1 \\
                                                    b_{4, 1} 
                                                    & b_{4, 2} 
                                                    & b_{4, 3} 
                                                    & 1
\end{array}
\right) \\
&
= 
\left(
\begin{array}{l l}
 (1 + s)^{p^f} 
 & (1 + s)^{p^f - p^{e_1}} \\
 b_{2, 1}\left( \frac{s}{1 + s} \right) \, (1 + s)^{p^f} 
 & b_{2, 1}\left( \frac{s}{1 + s} \right)  \, (1 + s)^{p^f - p^{e_1}}  + (1 + s)^{p^f - 2 \, p^{e_1}} \\
  b_{3, 1}\left( \frac{s}{1 + s} \right) \, (1 + s)^{p^f} 
 & b_{3, 1}\left( \frac{s}{1 + s} \right)  \, (1 + s)^{p^f - p^{e_1}}  
   + b_{3, 2}\left( \frac{s}{1 + s} \right)  \, (1 + s)^{p^f - 2 \, p^{e_1}}   \\
   b_{4, 1}\left( \frac{s}{1 + s} \right) \, (1 + s)^{p^f} 
 & b_{4, 1}\left( \frac{s}{1 + s} \right)  \, (1 + s)^{p^f - p^{e_1}}    
 + b_{4, 2}\left( \frac{s}{1 + s} \right)  \, (1 + s)^{p^f - 2 \, p^{e_1}}   \\
\end{array} 
\right.
\\
& \qquad 
\left. 
\begin{array}{l l}
  0 
  & 1\\
  0 
  &  b_{2, 1} \left( \frac{s}{1 + s} \right)  \\
\frac{1}{ (1 + s)^{p^f - 2 \, p^{e_1}} } 
& 
b_{3, 1} \left( \frac{s}{1 + s} \right)  
   +  \frac{1}{ (1 + s)^{p^f - p^{e_1}} }
\\ 
    b_{4, 3} \left( \frac{s}{1 + s} \right)  \, \frac{1}{ (1 + s)^{p^f - 2 \, p^{e_1}} }  
    & 
      b_{4, 1} \left( \frac{s}{1 + s} \right)
   + b_{4, 3} \left( \frac{s}{1 + s} \right) \, \frac{1}{ (1 + s)^{p^f - p^{e_1}} } 
   + \frac{1}{ (1 + s)^{p^f} }
\end{array}
\right) . 
\end{align*}

Comparing the $(1, 3)$-th entries of both sides of the equality, 
we have 
\[
 b_{4, 3} = 0 . 
\]

Comparing the $(3, 3)$-th entries of both sides of the equality, 
we have 
\[
1 + b_{4, 3}
 =
(1 + s)^{2 \, p^{e_1} - p^f}  , 
\]
which implies 
\[
1 = (1 + s)^{2 \, p^{e_1} - p^f}  . 
\]
So, 
\[
 2 \, p^{e_1} = p^f . 
\]
Therefore, $p = 2$ and $f = e_1 + 1$.

Comparing the $(2, 4)$-th entries of both sides of the equality, 
we have 
\[
 0 =  b_{2, 1} \left( \frac{s}{1 + s} \right)  , 
\]
which implies $b_{2, 1} = 0$.

Comparing the $(1, 2)$-th entries of both sides of the equality, 
we have 
\[
  1 + b_{4, 2}  = (1 + s)^{p^{e_1}}
\]
which implies $b_{4, 2} = s^{p^{e_1}}$.

Comparing the $(4, 2)$-th entries of both sides of the equality, 
we have 
\[
 b_{4, 2} = 
 b_{4, 1}\left( \frac{s}{1 + s} \right)  \, (1 + s)^{p^f - p^{e_1}}    
 + b_{4, 2}\left( \frac{s}{1 + s} \right)  \, (1 + s)^{p^f - 2 \, p^{e_1}}  , 
\]
which implies 
\[
 s^{p^{e_1}}
  = 
   b_{4, 1}\left( \frac{s}{1 + s} \right)  \, (1 + s)^{ p^{e_1}}    
 + \left( \frac{s}{1 + s} \right)^{p^{e_1}}  . 
\]
Thus 
\[
 b_{4, 1}\left( \frac{s}{1 + s} \right) 
 =
 \left( \frac{s}{1 + s} \right)^{2 \, p^{e_1}}  . 
\]
So, $b_{4, 1} = s^{2 \, p^{e_1}} = s^{p^f}$.

Comparing the $(3, 4)$-th entries of both sides of the equality, 
we have 
\[
 1 = b_{3, 1} \left( \frac{s}{1 + s} \right)  +  \frac{1}{ (1 + s)^{p^f - p^{e_1}} } , 
\]
which implies 
\[
  1 = b_{3, 1} \left( \frac{s}{1 + s} \right)  +  \frac{1}{ (1 + s)^{p^{e_1}} } . 
\]
So, $b_{3, 1} = s^{p^{e_1}}$.

Comparing the $(3, 2)$-th entries of both sides of the equality, 
we have 
\[
 b_{3, 2} + b_{4, 2} 
 =
b_{3, 1}\left( \frac{s}{1 + s} \right)  \, (1 + s)^{p^f - p^{e_1}}  
   + b_{3, 2}\left( \frac{s}{1 + s} \right)  \, (1 + s)^{p^f - 2 \, p^{e_1}}   ,  
\]
which implies 
\[
 b_{3, 2} + s^{p^{e_1}} 
 = 
 \left( \frac{s}{1 + s} \right)^{p^{e_1}}  \, (1 + s)^{p^f - p^{e_1}}  
   + b_{3, 2}\left( \frac{s}{1 + s} \right) .   
\]
So, 
\[
 b_{3, 2} 
 = 
 b_{3, 2}\left( \frac{s}{1 + s} \right) .   
\]
Thereby, $b_{3, 2} \in k$. Since $b_{3, 2} \in \sfP$, we have $b_{3, 2} = 0$. 

Hence $\phi^-$ has the desired form. 
\end{proof}

\subsubsection{$\rm (VI)$}

\begin{lem}
Let $(\varphi^*, \omega^*)$ be of the form {\rm (VI)}. 
Assume that there exists a homomorphism $\phi^- : \G_a \to \SL(4, k)$ 
satisfying the following conditions {\rm (i)} and {\rm (ii)}: 
\begin{enumerate}[label = {\rm (\roman*)}]
\item For any $s \in \G_a$, 
the regular matrix $\phi^-(s)$ is a lower triangular matrix.  

\item 
$
 \varphi^*(t) \, \phi^-(s) 
 = 
\phi^-\left( \frac{s}{1 + t \, s} \right) 
\, 
\omega^*(1 + t \, s) 
\;
\varphi^* \left( \frac{t}{1 + t \, s} \right)
$ 
for all $t, s \in k$ with $1 + t \, s \ne 0$.    
\end{enumerate} 
Then we have a contradiction. 
\end{lem}

\begin{proof}

Write $\phi^-(s) = ( \, b_{i, j}(s) \, )_{1 \leq i, j \leq 4}$. 
By condition (ii), we have 
\begin{align*}
& \renewcommand{\arraystretch}{1.5} 
\left(
\begin{array}{r r r r}
 1 + b_{2, 1} 
 & 1  
 & 0
 & 0 \\
 b_{3, 1}  
  & 1 
 & 0
 & 0 \\ 
                                       b_{3, 1} + b_{4, 1}
                                        & b_{3, 2} + b_{4, 2} 
                                        & 1 + b_{4, 3}
                                        & 1 \\
                                                    b_{4, 1} 
                                                    & b_{4, 2} 
                                                    & b_{4, 3} 
                                                    & 1
\end{array}
\right) \\
&
= 
\left(
\begin{array}{l l}
 (1 + s)^{d_2 + 2 \, p^{e_1}} 
 & (1 + s)^{d_2 + p^{e_1}} \\
 b_{2, 1}\left( \frac{s}{1 + s} \right) \, (1 + s)^{d_2 + 2 \, p^{e_1}} 
 & b_{2, 1}\left( \frac{s}{1 + s} \right)  \, (1 + s)^{d_2 + p^{e_1}}   + (1 + s)^{d_2} \\
  b_{3, 1}\left( \frac{s}{1 + s} \right) \,  (1 + s)^{d_2 + 2 \, p^{e_1}} 
 & b_{3, 1}\left( \frac{s}{1 + s} \right)  \, (1 + s)^{d_2 + p^{e_1}} 
   + b_{3, 2}\left( \frac{s}{1 + s} \right)  \, (1 + s)^{d_2}    \\
   b_{4, 1}\left( \frac{s}{1 + s} \right) \,  (1 + s)^{d_2 + 2 \, p^{e_1}} 
 & b_{4, 1}\left( \frac{s}{1 + s} \right)  \, (1 + s)^{d_2 + p^{e_1}}   
 + b_{4, 2}\left( \frac{s}{1 + s} \right)  \, (1 + s)^{d_2}    \\
\end{array} 
\right.
\\
& \qquad 
\left. 
\begin{array}{l l}
  0 
  & 0 \\
  0 
  &  0  \\
\frac{1}{(1 + s)^{d_2}}   
& 
\frac{1}{(1 + s)^{d_2 + p^{e_1}}}   
\\ 
    b_{4, 3} \left( \frac{s}{1 + s} \right)  \, \frac{1}{(1 + s)^{d_2}}   
    & 
   b_{4, 3} \left( \frac{s}{1 + s} \right) \, \frac{1}{(1 + s)^{d_2 + p^{e_1}}}   
   + \frac{1}{ (1 + s)^{d_2 + 2 \, p^{e_1}} }
\end{array}
\right) . 
\end{align*}

Comparing the $(1, 2)$-th entries of both sides of the equality, 
we have 
\[
 1 = (1 + s)^{d_2 + p^{e_1}} , 
\]
which implies a contradiction (since $e_1 \geq 0$). 

\end{proof}

\subsubsection{$\rm (VII)$}

\begin{lem}
Let $(\varphi^*, \omega^*)$ be of the form {\rm (VII)}. 
Assume that there exists a homomorphism $\phi^- : \G_a \to \SL(4, k)$ 
satisfying the following conditions {\rm (i)} and {\rm (ii)}: 
\begin{enumerate}[label = {\rm (\roman*)}]
\item For any $s \in \G_a$, 
the regular matrix $\phi^-(s)$ is a lower triangular matrix.  

\item 
$
 \varphi^*(t) \, \phi^-(s) 
 = 
\phi^-\left( \frac{s}{1 + t \, s} \right) 
\, 
\omega^*(1 + t \, s) 
\;
\varphi^* \left( \frac{t}{1 + t \, s} \right)
$ 
for all $t, s \in k$ with $1 + t \, s \ne 0$.    
\end{enumerate} 
Then we can express $\phi^-$ as  
\[
\phi^-(s) 
= 
\left(
\begin{array}{c c c c}
 1 & 0 & 0 & 0 \\
 \frac{1}{2} \, s^{p^{e_1}} & 1 & 0 & 0 \\
 s^{2 \, p^{e_1}} & s^{p^{e_1}} & 1 & 0 \\
 s^{p^{e_1 + 1}} & 0 & 0 & 1  
\end{array}
\right) . 
\]
\end{lem}

\begin{proof}
The pair $({^\tau\!}\varphi^*, \omega^*)$ has the form (II). 
Let $f^- := {^\tau\!}\phi^-$. 
By Lemma 1.19, 
the following conditions {\rm (i)} and {\rm (ii)} hold true: 
\begin{enumerate}[label = {\rm (\roman*)}]
\item For any $s \in \G_a$, 
the regular matrix $f^-(s)$ is a lower triangular matrix.  

\item 
$
 {^\tau\!}\varphi^*(t) \, f^-(s) 
 = 
f^-\left( \frac{s}{1 + t \, s} \right) 
\, 
\omega^*(1 + t \, s) 
\;
{^\tau\!}\varphi^* \left( \frac{t}{1 + t \, s} \right)
$ 
for all $t, s \in k$ with $1 + t \, s \ne 0$.  
\end{enumerate} 
By Lemma 4.2, we must have 
\[
f^-(s) 
= 
\left(
\begin{array}{c c c c}
 1 & 0 & 0 & 0 \\
 0 & 1 & 0 & 0 \\
 0 & s^{p^{e_1}} & 1 & 0 \\
 s^{p^{e_1 + 1}} & s^{2 \, p^{e_1}} & \frac{1}{2} \, s^{p^{e_1}} & 1  
\end{array}
\right) . 
\] 
Hence $\phi^-$ has the desired form. 

\end{proof}

\subsubsection{$\rm (VIII)$}

\begin{lem}
Let $(\varphi^*, \omega^*)$ be of the form {\rm (VIII)}. 
Assume that there exists a homomorphism $\phi^- : \G_a \to \SL(4, k)$ 
satisfying the following conditions {\rm (i)} and {\rm (ii)}: 
\begin{enumerate}[label = {\rm (\roman*)}]
\item For any $s \in \G_a$, 
the regular matrix $\phi^-(s)$ is a lower triangular matrix.  

\item 
$
 \varphi^*(t) \, \phi^-(s) 
 = 
\phi^-\left( \frac{s}{1 + t \, s} \right) 
\, 
\omega^*(1 + t \, s) 
\;
\varphi^* \left( \frac{t}{1 + t \, s} \right)
$ 
for all $t, s \in k$ with $1 + t \, s \ne 0$.    
\end{enumerate} 
Then we have a contradiction. 
\end{lem}

\begin{proof}
The pair $({^\tau}\varphi^*, \omega^*)$ has the form (III). 
Since $\psi_{\, {^\tau}\varphi^*,\; \omega^*} \circ \jmath^{-1}$ is not extendable, 
$\psi_{\varphi^*, \, \omega^*} \circ \jmath^{-1}$ is not extendable (see 
Lemma 2.8 (1)). 

\end{proof}

\subsubsection{$\rm (IX)$}

\begin{lem}
Let $(\varphi^*, \omega^*)$ be of the form {\rm (IX)}. 
Assume that there exists a homomorphism $\phi^- : \G_a \to \SL(4, k)$ 
satisfying the following conditions {\rm (i)} and {\rm (ii)}: 
\begin{enumerate}[label = {\rm (\roman*)}]
\item For any $s \in \G_a$, 
the regular matrix $\phi^-(s)$ is a lower triangular matrix.  

\item 
$
 \varphi^*(t) \, \phi^-(s) 
 = 
\phi^-\left( \frac{s}{1 + t \, s} \right) 
\, 
\omega^*(1 + t \, s) 
\;
\varphi^* \left( \frac{t}{1 + t \, s} \right)
$ 
for all $t, s \in k$ with $1 + t \, s \ne 0$.    
\end{enumerate} 
Then we can express $\phi^-$ as  
\[
\phi^-(s) 
= 
\left(
\begin{array}{c c c c}
 1 & 0 & 0 & 0 \\
 2 \, s^{p^{e_1}} & 1 & 0 & 0 \\
 0 & 0 & 1 & 0 \\
 2 \, s^{2 \, p^{e_1}} & 2 \, s^{p^{e_1}} & 0 & 1  
\end{array}
\right) . 
\]
\end{lem}

\begin{proof}

Write $\phi^-(s) = ( \, b_{i, j}(s) \, )_{1 \leq i, j \leq 4}$. 
By condition (ii), we have 
\begin{align*}
& \renewcommand{\arraystretch}{1.5} 
\left(
\begin{array}{r r r r}
 1 + b_{2, 1} + \frac{1}{2} \, b_{4, 1}
 & 1 + \frac{1}{2} \, b_{4, 2} 
 & \frac{1}{2} \, b_{4, 3}
 & \frac{1}{2} \\
 b_{2, 1} + b_{4, 1}  
  & 1 + b_{4, 2}
 & b_{4, 3}
 & 1 \\ 
                                       b_{3, 1} 
                                        & b_{3, 2} 
                                        & 1 
                                        & 0 \\
                                                    b_{4, 1} 
                                                    & b_{4, 2} 
                                                    & b_{4, 3} 
                                                    & 1
\end{array}
\right) \\
&
= 
\left(
\begin{array}{l l}
 (1 + s)^{2 \, p^{e_1}} 
 & (1 + s)^{p^{e_1}} \\
 b_{2, 1}\left( \frac{s}{1 + s} \right) \, (1 + s)^{2 \, p^{e_1}} 
 & b_{2, 1}\left( \frac{s}{1 + s} \right)  \, (1 + s)^{p^{e_1}}  + 1 \\
  b_{3, 1}\left( \frac{s}{1 + s} \right) \,  (1 + s)^{ 2 \, p^{e_1}} 
 & b_{3, 1}\left( \frac{s}{1 + s} \right)  \, (1 + s)^{ p^{e_1}} 
   + b_{3, 2}\left( \frac{s}{1 + s} \right)   \\
   b_{4, 1}\left( \frac{s}{1 + s} \right) \,  (1 + s)^{ 2 \, p^{e_1}} 
 & b_{4, 1}\left( \frac{s}{1 + s} \right)  \, (1 + s)^{ p^{e_1}}   
 + b_{4, 2}\left( \frac{s}{1 + s} \right)    \\
\end{array} 
\right.
\\
& \qquad 
\left. 
\begin{array}{l l}
  0 
  & \frac{1}{2} \\
  0 
  &  \frac{1}{2} \, b_{2, 1} \left( \frac{s}{1 + s} \right) + \frac{1}{(1 + s)^{p^{e_1}}}  \\
1  
& 
\frac{1}{2} \, b_{3, 1}\left( \frac{s}{1 + s} \right) 
+ b_{3, 2}\left( \frac{s}{1 + s}\right) \, \frac{1}{(1 + s)^{p^{e_1}}} 
\\ 
    b_{4, 3} \left( \frac{s}{1 + s} \right)  
    & 
   \frac{1}{2} \, b_{4, 1} \left( \frac{s}{1 + s} \right) 
   + b_{4, 2}\left( \frac{s}{1 + s} \right) \, \frac{1}{(1 + s)^{p^{e_1}}}   
   + \frac{1}{ (1 + s)^{2 \, p^{e_1}}}
\end{array}
\right) . 
\end{align*}

Comparing the $(1, 3)$-th entries of both sides of the equality, 
we have 
\[
 \frac{1}{2} \, b_{4, 3} = 0 , 
\]
which implies $b_{4, 3} = 0$.

Comparing the $(2, 4)$-th entries of both sides of the equality, 
we have 
\[
 1 = \frac{1}{2} \, b_{2, 1} \left( \frac{s}{1 + s} \right) + \frac{1}{(1 + s)^{p^{e_1}}} , 
\]
which implies $b_{2, 1} = 2 \, s^{p^{e_1}}$.

Comparing the $(2, 1)$-th entries of both sides of the equality, 
we have 
\[
 b_{2, 1} + b_{4, 1}
 =
 b_{2, 1}\left( \frac{s}{1 + s} \right) \, (1 + s)^{2 \, p^{e_1}}  , 
\]
which implies 
\[
 2 \, s^{p^{e_1}} + b_{4, 1}
 =
b_{3, 1}\left( \frac{s}{1 + s} \right) \,  (1 + s)^{ 2 \, p^{e_1}}  . 
\]
Thus $b_{4, 1} = 2 \, s^{2 \, p^{e_1}}$.

Comparing the $(3, 2)$-th and $(3, 4)$-th entries of both sides of the equality, 
we have 
\begin{align*}
\left\{
\begin{array}{r @{\,} l}
 b_{3, 2} & =  b_{3, 1}\left( \frac{s}{1 + s} \right)  \, (1 + s)^{ p^{e_1}} 
   + b_{3, 2}\left( \frac{s}{1 + s} \right)  , \\
 0 & =
 \frac{1}{2} \, b_{3, 1}\left( \frac{s}{1 + s} \right) 
+ b_{3, 2}\left( \frac{s}{1 + s}\right) \, \frac{1}{(1 + s)^{p^{e_1}}} .
\end{array}
\right.
\end{align*} 
Thus 
\[
 b_{3, 2}\left( \frac{s}{1 + s} \right) = - b_{3, 2} , 
\]
which implies $b_{3, 2} \in k$. Since $b_{3, 2} \in \sfP$, we have $b_{3, 2} = 0$. 
Therefore, $b_{3, 1} = 0$.

Comparing the $(1, 2)$-th entries of both sides of the equality, 
we have 
\[
 1 + \frac{1}{2} \, b_{4, 2} = (1 + s)^{p^{e_1}} , 
\]
which implies $b_{4, 2} = 2 \, s^{p^{e_1}}$.

Hence $\phi^-$ has the desired form. 
\end{proof}

Let 
\[
 P : = \diag(1 / 4, \; 1 / 2, \; 1, \; 1 ) \in \GL(4, k) . 
\]
Then we have 
\[
(\Inn_P \circ \phi^-)(s) 
= 
\left(
\begin{array}{c c c c}
 1 & 0 & 0 & 0 \\
 s^{p^{e_1}}  & 1 & 0 & 0 \\
 0 & 0 & 1 & 0 \\
 \frac{1}{2} \, s^{2 \, p^{e_1}} &  s^{p^{e_1}} & 0 & 1  
\end{array}
\right) . 
\]

\subsubsection{$\rm (X)$}

\begin{lem}
Let $(\varphi^*, \omega^*)$ be of the form {\rm (X)}. 
Assume that there exists a homomorphism $\phi^- : \G_a \to \SL(4, k)$ 
satisfying the following conditions {\rm (i)} and {\rm (ii)}: 
\begin{enumerate}[label = {\rm (\roman*)}]
\item For any $s \in \G_a$, 
the regular matrix $\phi^-(s)$ is a lower triangular matrix.  

\item 
$
 \varphi^*(t) \, \phi^-(s) 
 = 
\phi^-\left( \frac{s}{1 + t \, s} \right) 
\, 
\omega^*(1 + t \, s) 
\;
\varphi^* \left( \frac{t}{1 + t \, s} \right)
$ 
for all $t, s \in k$ with $1 + t \, s \ne 0$.    
\end{enumerate} 
Then we have a contradiction. 
\end{lem}

\begin{proof}

Write $\phi^-(s) = ( \, b_{i, j}(s) \, )_{1 \leq i, j \leq 4}$. 
By condition (ii), we have 
\begin{align*}
& \renewcommand{\arraystretch}{1.5} 
\left(
\begin{array}{r r r r}
 1 + b_{2, 1} + b_{3, 1}
 & 1 + b_{3, 2} 
 & 1
 & 0 \\
 b_{2, 1} 
  & 1 
 & 0
 & 0 \\ 
                                       b_{3, 1} 
                                        & b_{3, 2} 
                                        & 1 
                                        & 0 \\
                                                    b_{4, 1} 
                                                    & b_{4, 2} 
                                                    & b_{4, 3} 
                                                    & 1
\end{array}
\right) \\
&
= 
\left(
\begin{array}{l l}
 (1 + s)^{p^{e_1} + p^{e_2}} 
 & (1 + s)^{p^{e_2}} \\
 b_{2, 1}\left( \frac{s}{1 + s} \right) \,  (1 + s)^{p^{e_1} + p^{e_2}} 
 & b_{2, 1}\left( \frac{s}{1 + s} \right)  \,  (1 + s)^{ p^{e_2}} + (1 + s)^{p^{e_2} - p^{e_1}}  \\
  b_{3, 1}\left( \frac{s}{1 + s} \right) \,   (1 + s)^{p^{e_1} + p^{e_2}} 
 & b_{3, 1}\left( \frac{s}{1 + s} \right)  \, (1 + s)^{p^{e_2}}  
   + b_{3, 2}\left( \frac{s}{1 + s} \right) \, (1 + s)^{p^{e_2} - p^{e_1}}   \\
   b_{4, 1}\left( \frac{s}{1 + s} \right) \,  (1 + s)^{p^{e_1} + p^{e_2}} 
 & b_{4, 1}\left( \frac{s}{1 + s} \right)  \, (1 + s)^{ p^{e_2}}   
 + b_{4, 2}\left( \frac{s}{1 + s} \right)  \, (1 + s)^{p^{e_2} - p^{e_1}} 
\end{array} 
\right.
\\
& \qquad 
\left. 
\begin{array}{l l}
   (1 + s)^{p^{e_1}} 
  & 0 \\
    b_{2, 1}\left( \frac{s}{1 + s} \right) \,  (1 + s)^{p^{e_1} }  
  &  0  \\
b_{3, 1}\left( \frac{s}{1 + s} \right) \,  (1 + s)^{p^{e_1} } 
 + \frac{1}{(1 + s)^{p^{e_2} - p^{e_1}}}
& 
0
\\ 
    b_{4, 1} \left( \frac{s}{1 + s} \right)  \,  (1 + s)^{p^{e_1} }  
    + 
    b_{4, 3}\left( \frac{s}{1 + s} \right) \, \frac{1}{ (1 + s)^{p^{e_2} - p^{e_1} }  }
    & 
   \frac{1}{(1 + s)^{p^{e_1} + p^{e_2}}}
\end{array}
\right) . 
\end{align*}

Comparing the $(4, 4)$-th entries of both sides of the equality, 
we have 
\[
 1 = \frac{1}{(1 + s)^{p^{e_1} + p^{e_2}}} , 
\]
which implies a contradiction.

\end{proof}

\subsubsection{$\rm (XI)$}

\begin{lem}
Let $(\varphi^*, \omega^*)$ be of the form {\rm (XI)}. 
Assume that there exists a homomorphism $\phi^- : \G_a \to \SL(4, k)$ 
satisfying the following conditions {\rm (i)} and {\rm (ii)}: 
\begin{enumerate}[label = {\rm (\roman*)}]
\item For any $s \in \G_a$, 
the regular matrix $\phi^-(s)$ is a lower triangular matrix.  

\item 
$
 \varphi^*(t) \, \phi^-(s) 
 = 
\phi^-\left( \frac{s}{1 + t \, s} \right) 
\, 
\omega^*(1 + t \, s) 
\;
\varphi^* \left( \frac{t}{1 + t \, s} \right)
$ 
for all $t, s \in k$ with $1 + t \, s \ne 0$.    
\end{enumerate} 
Then we have $p = 2$ and $e_3 = e_1 + 1$, and we can express $\phi^-$ as 
\[
\phi^-(s) 
= 
\left(
\begin{array}{c c c c}
 1 & 0 & 0 & 0 \\
 0 & 1 & 0 & 0 \\
 0 & 0 & 1 & 0 \\
 s^{p^{e_3}} & s^{p^{e_1}} & 0 & 1  
\end{array}
\right) . 
\]
\end{lem}

\begin{proof}
Write $\phi^-(s) = ( \, b_{i, j}(s) \, )_{1 \leq i, j \leq 4}$. 
By condition (ii), we have 
\begin{align*}
& \renewcommand{\arraystretch}{1.5} 
\left(
\begin{array}{r r r r}
 1 + b_{2, 1} + b_{4, 1}
 & 1 + b_{4, 2} 
 & b_{4, 3} 
 & 1 \\
 b_{2, 1}
  & 1 
 & 0
 & 0 \\ 
                                       b_{3, 1} 
                                        & b_{3, 2} 
                                        & 1 
                                        & 0 \\
                                                    b_{4, 1} 
                                                    & b_{4, 2} 
                                                    & b_{4, 3} 
                                                    & 1
\end{array}
\right) \\
&
= 
\left(
\begin{array}{l l}
 (1 + s)^{p^{e_3}} 
 & (1 + s)^{p^{e_3} - p^{e_1}} \\
 b_{2, 1}\left( \frac{s}{1 + s} \right) \,  (1 + s)^{p^{e_3}} 
 & b_{2, 1}\left( \frac{s}{1 + s} \right)  \,  (1 + s)^{ p^{e_3} - p^{e_1}} + (1 + s)^{p^{e_3} - 2 \, p^{e_1}}  \\
  b_{3, 1}\left( \frac{s}{1 + s} \right) \,   (1 + s)^{p^{e_3}} 
 & b_{3, 1}\left( \frac{s}{1 + s} \right)  \, (1 + s)^{p^{e_3} - p^{e_1}}  
   + b_{3, 2}\left( \frac{s}{1 + s} \right) \, (1 + s)^{p^{e_3} - 2 \, p^{e_1}}   \\
   b_{4, 1}\left( \frac{s}{1 + s} \right) \,  (1 + s)^{p^{e_3}} 
 & b_{4, 1}\left( \frac{s}{1 + s} \right)  \, (1 + s)^{ p^{e_3} - p^{e_1}}   
 + b_{4, 2}\left( \frac{s}{1 + s} \right)  \, (1 + s)^{p^{e_3} - 2 \, p^{e_1}} 
\end{array} 
\right.
\\
& \qquad 
\left. 
\begin{array}{l l}
   0
  & 1 \\
    0
  &  b_{2, 1}\left( \frac{s}{1 + s} \right)   \\
\frac{1}{(1 + s)^{p^{e_3} - 2 \, p^{e_1}}}
& 
b_{3, 1}\left( \frac{s}{1 + s} \right) 
\\ 
    b_{4, 3} \left( \frac{s}{1 + s} \right) \, \frac{1}{ (1 + s)^{p^{e_3} - 2 \, p^{e_1} }  }
    & 
   b_{4, 1} \left(\frac{s}{1 + s} \right) + \frac{1}{(1 + s)^{p^{e_3}}}
\end{array}
\right) . 
\end{align*}

Comparing the $(3, 3)$-th entries of both sides of the equality, 
we have 
\[
 1 = \frac{1}{(1 + s)^{p^{e_3} - 2 \, p^{e_1}}} , 
\]
which implies $p^{e_3} = 2 \, p^{e_1}$. So, $p = 2$ and $e_3 = e_1 + 1$.

Comparing the $(2, 4)$-th entries of both sides of the equality, 
we have 
\[
 0 = b_{2, 1}\left( \frac{1}{1 + s} \right) , 
\]
which implies $b_{2, 1} = 0$.

Comparing the $(3, 4)$-th entries of both sides of the equality, 
we have 
\[
 0 = b_{3, 1}\left( \frac{1}{1 + s} \right) , 
\]
which implies $b_{3, 1} = 0$.

Comparing the $(1, 3)$-th entries of both sides of the equality, 
we have 
\[
 b_{4, 3} = 0.  
\]

Comparing the $(4, 4)$-th entries of both sides of the equality, 
we have 
\[
 1 = b_{4, 1} \left(\frac{s}{1 + s} \right) + \frac{1}{(1 + s)^{p^{e_3}}} , 
\]
which implies $b_{4, 1} = s^{p^{e_3}}$.

Comparing the $(1, 2)$-th entries of both sides of the equality, 
we have 
\[
 1 + b_{4, 2} 
 =  (1 + s)^{p^{e_3} - p^{e_1}}  , 
\]
which implies $b_{4, 2} = s^{p^{e_1}}$.

Comparing the $(3, 2)$-th entries of both sides of the equality, 
we have 
\[
  b_{3, 2} 
 = b_{3, 1}\left( \frac{s}{1 + s} \right)  \, (1 + s)^{p^{e_3} - p^{e_1}}  
   + b_{3, 2}\left( \frac{s}{1 + s} \right) \, (1 + s)^{p^{e_3} - 2 \, p^{e_1}} ,  
\]
which implies 
\[
 b_{3, 2} = b_{3, 2}\left( \frac{s}{1 + s} \right) . 
\]
So, $b_{3, 2} \in k$. Since $b_{3, 2} \in \sfP$, we have $b_{3, 2} = 0$. 

Hence $\phi^-$ has the desired form. 

\end{proof}

\subsubsection{$\rm (XII)$}

\begin{lem}
Let $(\varphi^*, \omega^*)$ be of the form {\rm (XII)}. 
Assume that there exists a homomorphism $\phi^- : \G_a \to \SL(4, k)$ 
satisfying the following conditions {\rm (i)} and {\rm (ii)}: 
\begin{enumerate}[label = {\rm (\roman*)}]
\item For any $s \in \G_a$, 
the regular matrix $\phi^-(s)$ is a lower triangular matrix.  

\item 
$
 \varphi^*(t) \, \phi^-(s) 
 = 
\phi^-\left( \frac{s}{1 + t \, s} \right) 
\, 
\omega^*(1 + t \, s) 
\;
\varphi^* \left( \frac{t}{1 + t \, s} \right)
$ 
for all $t, s \in k$ with $1 + t \, s \ne 0$.    
\end{enumerate} 
Then we have a contradiction. 
\end{lem}

\begin{proof}

Write $\phi^-(s) = ( \, b_{i, j}(s) \, )_{1 \leq i, j \leq 4}$. 
By condition (ii), we have 
\begin{align*}
& \renewcommand{\arraystretch}{1.5} 
\left(
\begin{array}{r r r r}
 1 + b_{2, 1} 
 & 1
 & 0 
 & 0 \\
 b_{2, 1}
 & 1
 & 0
 & 0 \\ 
                                       b_{3, 1} 
                                        & b_{3, 2} 
                                        & 1 
                                        & 0 \\
                                                    b_{4, 1} 
                                                    & b_{4, 2} 
                                                    & b_{4, 3} 
                                                    & 1
\end{array}
\right) \\
&
= 
\left(
\begin{array}{l l}
 (1 + s)^{2 \, p^{e_1} + d_2} 
 & (1 + s)^{p^{e_1} + d_2 } \\
 b_{2, 1}\left( \frac{s}{1 + s} \right) \,  (1 + s)^{2 \, p^{e_1} + d_2} 
 & b_{2, 1}\left( \frac{s}{1 + s} \right)  \, (1 + s)^{p^{e_1} + d_2 } + (1 + s)^{d_2}  \\
  b_{3, 1}\left( \frac{s}{1 + s} \right) \,   (1 + s)^{2 \, p^{e_1} + d_2} 
 & b_{3, 1}\left( \frac{s}{1 + s} \right)  \, (1 + s)^{p^{e_1} + d_2 } 
   + b_{3, 2}\left( \frac{s}{1 + s} \right) \, (1 + s)^{d_2}   \\
   b_{4, 1}\left( \frac{s}{1 + s} \right) \,  (1 + s)^{2 \, p^{e_1} + d_2} 
 & b_{4, 1}\left( \frac{s}{1 + s} \right)  \, (1 + s)^{p^{e_1} + d_2 }   
 + b_{4, 2}\left( \frac{s}{1 + s} \right)  \, (1 + s)^{d_2} 
\end{array} 
\right.
\\
& \qquad 
\left. 
\begin{array}{l l}
   0
  & 0 \\
   0 
  & 0 \\
\frac{1}{(1 + s)^{d_2}}
& 
0
\\ 
   b_{4, 3}\left( \frac{s}{1 + s} \right)  \, \frac{1}{(1 + s)^{d_2}} 
    & 
    \frac{1}{(1 + s)^{2 \, p^{e_1} + d_2}}
\end{array}
\right) . 
\end{align*}

Comparing the $(4, 4)$-th entries of both sides of the equality, 
we have 
\[
  1 
 =  \frac{1}{(1 + s)^{2 \, p^{e_1} + d_2}} , 
\]
which implies a contradiction (since $e_1 \geq 0$).

\end{proof}

\subsubsection{$\rm (XIII)$}

\begin{lem}
Let $(\varphi^*, \omega^*)$ be of the form {\rm (XIII)}. 
Assume that there exists a homomorphism $\phi^- : \G_a \to \SL(4, k)$ 
satisfying the following conditions {\rm (i)} and {\rm (ii)}: 
\begin{enumerate}[label = {\rm (\roman*)}]
\item For any $s \in \G_a$, 
the regular matrix $\phi^-(s)$ is a lower triangular matrix.  

\item 
$
 \varphi^*(t) \, \phi^-(s) 
 = 
\phi^-\left( \frac{s}{1 + t \, s} \right) 
\, 
\omega^*(1 + t \, s) 
\;
\varphi^* \left( \frac{t}{1 + t \, s} \right)
$ 
for all $t, s \in k$ with $1 + t \, s \ne 0$.    
\end{enumerate} 
Then we have a contradiction. 
\end{lem}

\begin{proof}

Write $\phi^-(s) = ( \, b_{i, j}(s) \, )_{1 \leq i, j \leq 4}$. 
By condition (ii), we have 
\begin{align*}
& \renewcommand{\arraystretch}{1.5} 
\left(
\begin{array}{r r r r}
 1 + b_{3, 1} 
 & b_{3, 2}
 & 1 
 & 0 \\
 b_{2, 1} + b_{3, 1} + b_{4, 1}
 & 1 + b_{3, 2} + b_{4, 2}
 & 1 + b_{4, 3}
 & 1 \\ 
                                       b_{3, 1} 
                                        & b_{3, 2} 
                                        & 1 
                                        & 0 \\
                                                    b_{4, 1} 
                                                    & b_{4, 2} 
                                                    & b_{4, 3} 
                                                    & 1
\end{array}
\right) \\
&
= 
\left(
\begin{array}{l l}
 (1 + s)^{2 \, p^{e_1} - p^{e_3}} 
 & 0 \\
 b_{2, 1}\left( \frac{s}{1 + s} \right) \,  (1 + s)^{2 \, p^{e_1} - p^{e_3}} 
 & (1 + s)^{p^{e_3}}  \\
  b_{3, 1}\left( \frac{s}{1 + s} \right) \,   (1 + s)^{2 \, p^{e_1} - p^{e_3}} 
 &  b_{3, 2}\left( \frac{s}{1 + s} \right) \, (1 + s)^{p^{e_3}}    \\
   b_{4, 1}\left( \frac{s}{1 + s} \right) \,  (1 + s)^{2 \, p^{e_1} - p^{e_3}} 
 &  b_{4, 2}\left( \frac{s}{1 + s} \right)  \, (1 + s)^{p^{e_3}} 
\end{array} 
\right.
\\
& \qquad 
\left. 
\begin{array}{l l}
    (1 + s)^{p^{e_1} - p^{e_3}} 
  & 0 \\
   b_{2, 1}\left( \frac{s}{1 + s} \right) \,  (1 + s)^{p^{e_1} - p^{e_3}} + 1 
  & (1 + s)^{p^{e_3} - p^{e_1}}  \\
 b_{3, 1}\left( \frac{s}{1 + s} \right) \,  (1 + s)^{p^{e_1} - p^{e_3}} 
  + b_{3, 2}\left( \frac{s}{1 + s} \right) 
  + \frac{1}{(1 + s)^{p^{e_3}}} 
& 
 b_{3, 2}\left( \frac{s}{1 + s} \right) \, (1 + s)^{p^{e_3} - p^{e_1}} 
\\ 
   b_{4, 1}\left( \frac{s}{1 + s} \right) \,  (1 + s)^{p^{e_1} - p^{e_3}} 
   + b_{4, 2}\left( \frac{s}{1 + s} \right) 
   + b_{4, 3} \left( \frac{s}{1 + s} \right) \,  \frac{1}{(1 + s)^{p^{e_3}}} 
    & 
    b_{4, 2}\left( \frac{s}{1 + s} \right) \, (1 + s)^{p^{e_3} - p^{e_1}} 
    +
    \frac{1}{(1 + s)^{2 \, p^{e_1} - p^{e_3}}}
\end{array}
\right) . 
\end{align*}

Comparing the $(2, 4)$-th entries of both sides of the equality, 
we have 
\[
  1 
 =  (1 + s)^{p^{e_3} - p^{e_1}}  , 
\]
which implies $e_3 = e_1$, which contradicts $e_1 > e_3$.

\end{proof}

\subsubsection{$\rm (XIV)$}

\begin{lem}
Let $(\varphi^*, \omega^*)$ be of the form {\rm (XIV)}. 
Assume that there exists a homomorphism $\phi^- : \G_a \to \SL(4, k)$ 
satisfying the following conditions {\rm (i)} and {\rm (ii)}: 
\begin{enumerate}[label = {\rm (\roman*)}]
\item For any $s \in \G_a$, 
the regular matrix $\phi^-(s)$ is a lower triangular matrix.  

\item 
$
 \varphi^*(t) \, \phi^-(s) 
 = 
\phi^-\left( \frac{s}{1 + t \, s} \right) 
\, 
\omega^*(1 + t \, s) 
\;
\varphi^* \left( \frac{t}{1 + t \, s} \right)
$ 
for all $t, s \in k$ with $1 + t \, s \ne 0$.    
\end{enumerate} 
Then we have a contradition. 
\end{lem}

\begin{proof}

Write $\phi^-(s) = ( \, b_{i, j}(s) \, )_{1 \leq i, j \leq 4}$. 
By condition (ii), we have 
\begin{align*}
& \renewcommand{\arraystretch}{1.5} 
\left(
\begin{array}{r r r r}
 1 + b_{3, 1} 
 & b_{3, 2}
 & 1 
 & 0 \\
 b_{2, 1} + b_{3, 1} 
 & 1 + b_{3, 2} 
 & 1 
 & 0 \\ 
                                       b_{3, 1} 
                                        & b_{3, 2} 
                                        & 1 
                                        & 0 \\
                                                    b_{4, 1} 
                                                    & b_{4, 2} 
                                                    & b_{4, 3} 
                                                    & 1
\end{array}
\right) \\
&
= 
\left(
\begin{array}{l l}
 (1 + s)^{2 \, p^{e_1} - p^{e_3}} 
 & 0 \\
 b_{2, 1}\left( \frac{s}{1 + s} \right) \,  (1 + s)^{2 \, p^{e_1} - p^{e_3}} 
 & (1 + s)^{p^{e_3}}  \\
  b_{3, 1}\left( \frac{s}{1 + s} \right) \,   (1 + s)^{2 \, p^{e_1} - p^{e_3}} 
 &  b_{3, 2}\left( \frac{s}{1 + s} \right) \, (1 + s)^{p^{e_3}}    \\
   b_{4, 1}\left( \frac{s}{1 + s} \right) \,  (1 + s)^{2 \, p^{e_1} - p^{e_3}} 
 &  b_{4, 2}\left( \frac{s}{1 + s} \right)  \, (1 + s)^{p^{e_3}} 
\end{array} 
\right.
\\
& \qquad 
\left. 
\begin{array}{l l}
    (1 + s)^{p^{e_1} - p^{e_3}} 
  & 0 \\
   b_{2, 1}\left( \frac{s}{1 + s} \right) \,  (1 + s)^{p^{e_1} - p^{e_3}} + 1 
  & 0 \\
 b_{3, 1}\left( \frac{s}{1 + s} \right) \,  (1 + s)^{p^{e_1} - p^{e_3}} 
  + b_{3, 2} \left( \frac{s}{1 + s} \right) 
  + \frac{1}{(1 + s)^{p^{e_3}}} 
& 
 0 
\\ 
   b_{4, 1}\left( \frac{s}{1 + s} \right) \,  (1 + s)^{p^{e_1} - p^{e_3}} 
   + b_{4, 2}\left( \frac{s}{1 + s} \right) 
   + b_{4, 3} \left( \frac{s}{1 + s} \right) \,  \frac{1}{(1 + s)^{p^{e_3}}} 
    & 
    \frac{1}{(1 + s)^{2 \, p^{e_1} - p^{e_3}}}
\end{array}
\right) . 
\end{align*}

Comparing the $(4, 4)$-th entries of both sides of the equality, 
we have 
\[
  1 
 =  (1 + s)^{2 \, p^{e_1} - p^{e_3}}  , 
\]
which implies $2 \, p^{e_1} = p^{e_3}$. 
So, $p = 2$ and $e_3 = e_1 +1$, which contradicts $e_1 > e_3$.

\end{proof}

\subsubsection{$\rm (XV)$}

\begin{lem}
Let $(\varphi^*, \omega^*)$ be of the form {\rm (XV)}. 
Assume that there exists a homomorphism $\phi^- : \G_a \to \SL(4, k)$ 
satisfying the following conditions {\rm (i)} and {\rm (ii)}: 
\begin{enumerate}[label = {\rm (\roman*)}]
\item For any $s \in \G_a$, 
the regular matrix $\phi^-(s)$ is a lower triangular matrix.  

\item 
$
 \varphi^*(t) \, \phi^-(s) 
 = 
\phi^-\left( \frac{s}{1 + t \, s} \right) 
\, 
\omega^*(1 + t \, s) 
\;
\varphi^* \left( \frac{t}{1 + t \, s} \right)
$ 
for all $t, s \in k$ with $1 + t \, s \ne 0$.    
\end{enumerate} 
Then we can express $\phi^-$ as 
\[
\phi^-(s) 
= 
\left(
\begin{array}{c c c c}
 1 & 0 & 0 & 0 \\
 0 & 1 & 0 & 0 \\
 0 & s^{p^{e_3}} & 1 & 0 \\
 s^{p^{e_2}} & 0 & 0 & 1  
\end{array}
\right) . 
\]
\end{lem}

\begin{proof}

Write $\phi^-(s) = ( \, b_{i, j}(s) \, )_{1 \leq i, j \leq 4}$. 
By condition (ii), we have 
\begin{align*}
& \renewcommand{\arraystretch}{1.5} 
\left(
\begin{array}{r r r r}
 1 + b_{4, 1} 
 & b_{4, 2}
 & b_{4, 3}
 & 1 \\
 b_{2, 1} + b_{3, 1} 
 & 1 + b_{3, 2} 
 & 1 
 & 0 \\ 
                                       b_{3, 1} 
                                        & b_{3, 2} 
                                        & 1 
                                        & 0 \\
                                                    b_{4, 1} 
                                                    & b_{4, 2} 
                                                    & b_{4, 3} 
                                                    & 1
\end{array}
\right) \\
&
= 
\left(
\begin{array}{l l}
 (1 + s)^{p^{e_2}} 
 & 0 \\
 b_{2, 1}\left( \frac{s}{1 + s} \right) \,  (1 + s)^{p^{e_2}} 
 & (1 + s)^{p^{e_3}}  \\
  b_{3, 1}\left( \frac{s}{1 + s} \right) \,   (1 + s)^{p^{e_2}} 
 &  b_{3, 2}\left( \frac{s}{1 + s} \right) \, (1 + s)^{p^{e_3}}    \\
   b_{4, 1}\left( \frac{s}{1 + s} \right) \,  (1 + s)^{p^{e_2}} 
 &  b_{4, 2}\left( \frac{s}{1 + s} \right)  \, (1 + s)^{p^{e_3}} 
\end{array} 
\right.
\\
& \qquad 
\left. 
\begin{array}{l l}
    0
  & 1 \\
  1  
  & b_{2, 1}\left( \frac{s}{1 + s} \right)  \\
  b_{3, 2}\left( \frac{s}{1 + s} \right) + \frac{1}{(1 + s)^{p^{e_3}}} 
& 
 b_{3, 1}\left( \frac{s}{1 + s} \right) 
\\ 
    b_{4, 2}\left( \frac{s}{1 + s} \right) 
   + b_{4, 3} \left( \frac{s}{1 + s} \right) \,  \frac{1}{(1 + s)^{p^{e_3}}} 
    & 
    b_{4, 1}\left( \frac{s}{1 + s} \right)  + \frac{1}{(1 + s)^{ p^{e_2}}}
\end{array}
\right) . 
\end{align*}

Comparing the $(2, 4)$-th entries of both sides of the equality, 
we have 
\[
  0 
 = b_{2, 1}\left( \frac{s}{1 + s} \right)   , 
\]
which implies $b_{2, 1} = 0$.

Comparing the $(3, 4)$-th entries of both sides of the equality, 
we have 
\[
  0 
 = b_{3, 1}\left( \frac{s}{1 + s} \right)   , 
\]
which implies $b_{3, 1} = 0$.

Comparing the $(4, 4)$-th entries of both sides of the equality, 
we have 
\[
  1 
 =   b_{4, 1}\left( \frac{s}{1 + s} \right)  + \frac{1}{(1 + s)^{ p^{e_2}}} , 
\]
which implies $b_{4, 1} = s^{p^{e_2}}$.

Comparing the $(1, 3)$-th entries of both sides of the equality, 
we have 
\[
  b_{4, 3}
 = 0 . 
\]

Comparing the $(1, 2)$-th entries of both sides of the equality, 
we have 
\[
  b_{4, 2} = 0 . 
\]

Comparing the $(2, 2)$-th entries of both sides of the equality, 
we have 
\[
 1 + b_{3, 2} = (1 + s)^{p^{e_3}} , 
\]
which implies $b_{3, 2} = s^{p^{e_3}}$.

Hence $\phi^-$ has the desired form.

\end{proof}

\subsubsection{$\rm (XVI)$}

\begin{lem}
Let $(\varphi^*, \omega^*)$ be of the form {\rm (XVI)}. 
Assume that there exists a homomorphism $\phi^- : \G_a \to \SL(4, k)$ 
satisfying the following conditions {\rm (i)} and {\rm (ii)}: 
\begin{enumerate}[label = {\rm (\roman*)}]
\item For any $s \in \G_a$, 
the regular matrix $\phi^-(s)$ is a lower triangular matrix.  

\item 
$
 \varphi^*(t) \, \phi^-(s) 
 = 
\phi^-\left( \frac{s}{1 + t \, s} \right) 
\, 
\omega^*(1 + t \, s) 
\;
\varphi^* \left( \frac{t}{1 + t \, s} \right)
$ 
for all $t, s \in k$ with $1 + t \, s \ne 0$.    
\end{enumerate} 
Then we have a contradiction. 
\end{lem}

\begin{proof} 

The pair $({^\tau}\varphi^*, \omega^*)$ has the form (XIV). 
Since $\psi_{\, {^\tau}\varphi^*,\; \omega^*} \circ \jmath^{-1}$ is not extendable, 
$\psi_{\varphi^*, \, \omega^*} \circ \jmath^{-1}$ is not extendable (see 
Lemma 2.8 (1)). 

\end{proof}

\subsubsection{$\rm (XVII)$}

\begin{lem}
Let $(\varphi^*, \omega^*)$ be of the form {\rm (XVII)}. 
Assume that there exists a homomorphism $\phi^- : \G_a \to \SL(4, k)$ 
satisfying the following conditions {\rm (i)} and {\rm (ii)}: 
\begin{enumerate}[label = {\rm (\roman*)}]
\item For any $s \in \G_a$, 
the regular matrix $\phi^-(s)$ is a lower triangular matrix.  

\item 
$
 \varphi^*(t) \, \phi^-(s) 
 = 
\phi^-\left( \frac{s}{1 + t \, s} \right) 
\, 
\omega^*(1 + t \, s) 
\;
\varphi^* \left( \frac{t}{1 + t \, s} \right)
$ 
for all $t, s \in k$ with $1 + t \, s \ne 0$.    
\end{enumerate} 
Then we have a contradiction. 
\end{lem}

\begin{proof}

Write $\phi^-(s) = ( \, b_{i, j}(s) \, )_{1 \leq i, j \leq 4}$. 
By condition (ii), we have 
\begin{align*}
& \renewcommand{\arraystretch}{1.5} 
\left(
\begin{array}{r r r r}
 1 
 & 0
 & 0 
 & 0 \\
 b_{2, 1} + b_{3, 1}
 & 1 + b_{3, 2}
 & 1 
 & 0 \\ 
                                       b_{3, 1} 
                                        & b_{3, 2} 
                                        & 1 
                                        & 0 \\
                                                    b_{4, 1} 
                                                    & b_{4, 2} 
                                                    & b_{4, 3} 
                                                    & 1
\end{array}
\right) \\
&
= 
\left(
\begin{array}{l l}
 (1 + s)^{d_1} 
 & 0 \\
 b_{2, 1}\left( \frac{s}{1 + s} \right) \,  (1 + s)^{d_1} 
 &  (1 + s)^{p^{e_3}}  \\
  b_{3, 1}\left( \frac{s}{1 + s} \right) \,  (1 + s)^{d_1} 
 &  b_{3, 2}\left( \frac{s}{1 + s} \right) \, (1 + s)^{p^{e_3}}   \\
   b_{4, 1}\left( \frac{s}{1 + s} \right) \,  (1 + s)^{d_1} 
 & b_{4, 2}\left( \frac{s}{1 + s} \right)  \, (1 + s)^{p^{e_3}} 
\end{array} 
\right.
\\
& \qquad 
\left. 
\begin{array}{l l}
   0
  & 0 \\
   1
  & 0 \\
b_{3, 2}\left( \frac{s}{1 + s} \right) +  \frac{1}{(1 + s)^{p^{e_3}}}
& 
0
\\ 
   b_{4, 2} \left( \frac{s}{1 + s} \right)  
   + 
   b_{4, 3} \left( \frac{s}{1 + s} \right)  \, \frac{1}{(1 + s)^{p^{e_3}}} 
    & 
    \frac{1}{(1 + s)^{d_1}}
\end{array}
\right) . 
\end{align*}

Comparing the $(4, 4)$-th entries of both sides of the equality, 
we have 
\[
 1 = \frac{1}{(1 + s)^{d_1}} , 
\]
which implies $d_1 = 0$. So, $d_2 = p^{e_3} = 0$, which contradicts $e_3 \geq 1$.

\end{proof}

\subsubsection{$\rm (XVIII)$}

\begin{lem}
Let $(\varphi^*, \omega^*)$ be of the form {\rm (XVIII)}. 
Assume that there exists a homomorphism $\phi^- : \G_a \to \SL(4, k)$ 
satisfying the following conditions {\rm (i)} and {\rm (ii)}: 
\begin{enumerate}[label = {\rm (\roman*)}]
\item For any $s \in \G_a$, 
the regular matrix $\phi^-(s)$ is a lower triangular matrix.  

\item 
$
 \varphi^*(t) \, \phi^-(s) 
 = 
\phi^-\left( \frac{s}{1 + t \, s} \right) 
\, 
\omega^*(1 + t \, s) 
\;
\varphi^* \left( \frac{t}{1 + t \, s} \right)
$ 
for all $t, s \in k$ with $1 + t \, s \ne 0$.    
\end{enumerate} 
Then we have a contradiction. 

\end{lem}

\begin{proof}

The pair $({^\tau}\varphi^*, \omega^*)$ has the form (X). 
Since $\psi_{\, {^\tau}\varphi^*,\; \omega^*} \circ \jmath^{-1}$ is not extendable, 
$\psi_{\varphi^*, \, \omega^*} \circ \jmath^{-1}$ is not extendable (see 
Lemma 2.8 (1)). 

\end{proof}

\subsubsection{$\rm (XIX)$}

\begin{lem}
Let $(\varphi^*, \omega^*)$ be of the form {\rm (XIX)}. 
Assume that there exists a homomorphism $\phi^- : \G_a \to \SL(4, k)$ 
satisfying the following conditions {\rm (i)} and {\rm (ii)}: 
\begin{enumerate}[label = {\rm (\roman*)}]
\item For any $s \in \G_a$, 
the regular matrix $\phi^-(s)$ is a lower triangular matrix.  

\item 
$
 \varphi^*(t) \, \phi^-(s) 
 = 
\phi^-\left( \frac{s}{1 + t \, s} \right) 
\, 
\omega^*(1 + t \, s) 
\;
\varphi^* \left( \frac{t}{1 + t \, s} \right)
$ 
for all $t, s \in k$ with $1 + t \, s \ne 0$.    
\end{enumerate} 
Then we have $p = 2$ and $e_3 = e_1 + 1$, and we can express $\phi^-$ as 
\[
\phi^-(s) 
= 
\left(
\begin{array}{c c c c}
 1 & 0 & 0 & 0 \\
 0 & 1 & 0 & 0 \\
 s^{p^{e_1}} & 0 & 1 & 0 \\
 s^{p^{e_3}} & 0 & 0 & 1  
\end{array}
\right) . 
\]
\end{lem}

\begin{proof}
The pair $({^\tau\!}\varphi^*, \omega^*)$ has the form (XI). 
Let $f^- := {^\tau\!}\phi^-$. 
By Lemma 1.19, the following conditions {\rm (i)} and {\rm (ii)} hold true: 
\begin{enumerate}[label = {\rm (\roman*)}]
\item For any $s \in \G_a$, 
the regular matrix $f^-(s)$ is a lower triangular matrix.  

\item 
$
 {^\tau\!}\varphi^*(t) \, f^-(s) 
 = 
f^-\left( \frac{s}{1 + t \, s} \right) 
\, 
\omega^*(1 + t \, s) 
\;
{^\tau\!}\varphi^* \left( \frac{t}{1 + t \, s} \right)
$ 
for all $t, s \in k$ with $1 + t \, s \ne 0$.  
\end{enumerate} 
By Lemma 4.11, we must have 
\[
f^-(s) 
= 
\left(
\begin{array}{c c c c}
 1 & 0 & 0 & 0 \\
 0 & 1 & 0 & 0 \\
 0 & 0 & 1 & 0 \\
 s^{p^{e_3}} & s^{p^{e_1}} & 0 & 1  
\end{array}
\right) . 
\]
Hence $\phi^-$ has the desired form. 
\end{proof}

\subsubsection{$\rm (XX)$}

\begin{lem}
Let $(\varphi^*, \omega^*)$ be of the form {\rm (XX)}. 
Assume that there exists a homomorphism $\phi^- : \G_a \to \SL(4, k)$ 
satisfying the following conditions {\rm (i)} and {\rm (ii)}: 
\begin{enumerate}[label = {\rm (\roman*)}]
\item For any $s \in \G_a$, 
the regular matrix $\phi^-(s)$ is a lower triangular matrix.  

\item 
$
 \varphi^*(t) \, \phi^-(s) 
 = 
\phi^-\left( \frac{s}{1 + t \, s} \right) 
\, 
\omega^*(1 + t \, s) 
\;
\varphi^* \left( \frac{t}{1 + t \, s} \right)
$ 
for all $t, s \in k$ with $1 + t \, s \ne 0$.    
\end{enumerate} 
Then we have a contradiction. 
\end{lem}

\begin{proof}

The pair $({^\tau}\varphi^*, \omega^*)$ has the form (XII). 
Since $\psi_{\, {^\tau}\varphi^*,\; \omega^*} \circ \jmath^{-1}$ is not extendable, 
$\psi_{\varphi^*, \, \omega^*} \circ \jmath^{-1}$ is not extendable (see 
Lemma 2.8 (1)). 

\end{proof}

\subsubsection{$\rm (XXI)$}

\begin{lem}
Let $(\varphi^*, \omega^*)$ be of the form {\rm (XXI)}. 
Assume that there exists a homomorphism $\phi^- : \G_a \to \SL(4, k)$ 
satisfying the following conditions {\rm (i)} and {\rm (ii)}: 
\begin{enumerate}[label = {\rm (\roman*)}]
\item For any $s \in \G_a$, 
the regular matrix $\phi^-(s)$ is a lower triangular matrix.  

\item 
$
 \varphi^*(t) \, \phi^-(s) 
 = 
\phi^-\left( \frac{s}{1 + t \, s} \right) 
\, 
\omega^*(1 + t \, s) 
\;
\varphi^* \left( \frac{t}{1 + t \, s} \right)
$ 
for all $t, s \in k$ with $1 + t \, s \ne 0$.    
\end{enumerate} 
Then we can express $\phi^-$ as 
\[
\phi^-(s) 
= 
\left(
\begin{array}{c c c c}
 1 & 0 & 0 & 0 \\
 s^{p^{e_1}} & 1 & 0 & 0 \\
 0 & 0 & 1 & 0 \\
 s^{2 \, p^{e_1}} & 0 & s^{p^{e_1}} & 1  
\end{array}
\right) . 
\]
\end{lem}

\begin{proof}

Write $\phi^-(s) = ( \, b_{i, j}(s) \, )_{1 \leq i, j \leq 4}$. 
By condition (ii), we have 
\begin{align*}
& \renewcommand{\arraystretch}{1.5} 
\left(
\begin{array}{r r r r}
 1 + b_{2, 1} + b_{4, 1} 
 & b_{3, 2} + b_{4, 2}
 & 1 + b_{4, 3}
 & 1 \\
 b_{2, 1} + b_{4, 1} 
 & 1 + b_{4, 2} 
 & b_{4, 3}
 & 1 \\ 
                                       b_{3, 1} 
                                        & b_{3, 2} 
                                        & 1 
                                        & 0 \\
                                                    b_{4, 1} 
                                                    & b_{4, 2} 
                                                    & b_{4, 3} 
                                                    & 1
\end{array}
\right) \\
&
= 
\left(
\begin{array}{l l}
 (1 + s)^{p^{e_1 + 1}} 
 & 0 \\
 b_{2, 1}\left( \frac{s}{1 + s} \right) \,  (1 + s)^{p^{e_1 + 1}} 
 & 1 \\
  b_{3, 1}\left( \frac{s}{1 + s} \right) \,   (1 + s)^{p^{e_1 + 1}} 
 &  b_{3, 2}\left( \frac{s}{1 + s} \right)     \\
   b_{4, 1}\left( \frac{s}{1 + s} \right) \,   (1 + s)^{p^{e_1 + 1}} 
 &  b_{4, 2}\left( \frac{s}{1 + s} \right)  
\end{array} 
\right.
\\
& \qquad 
\left. 
\begin{array}{l l}
    (1 + s)^{p^{e_1}} 
  & 1 \\
  b_{2, 1}\left( \frac{s}{1 + s} \right) \, (1 + s)^{p^{e_1}} 
  & b_{2, 1}\left( \frac{s}{1 + s} \right) + \frac{1}{(1 + s)^{p^{e_1}}}  \\
 b_{3, 1}\left( \frac{s}{1 + s} \right) \, (1 + s)^{p^{e_1}} + 1 
& 
 b_{3, 1}\left( \frac{s}{1 + s} \right) + b_{3, 2}\left( \frac{s}{1 + s} \right) \, \frac{1}{(1 + s)^{p^{e_1}}}  
\\ 
   b_{4, 1}\left( \frac{s}{1 + s} \right) \, (1 + s)^{p^{e_1}}
   + b_{4, 3} \left( \frac{s}{1 + s} \right) 
    & 
    b_{4, 1}\left( \frac{s}{1 + s} \right) + b_{4, 2}\left( \frac{s}{1 + s} \right) \, \frac{1}{(1 + s)^{p^{e_1}}}   + \frac{1}{(1 + s)^{ p^{e_1 + 1}}}
\end{array}
\right) . 
\end{align*}

Comparing the $(1, 3)$-th entries of both sides of the equality, 
we have 
\[
 1 + b_{4, 3} = (1 + s )^{p^{e_1}} , 
\]
which implies $b_{4, 3} = s^{p^{e_1}}$.

Comparing the $(2, 3)$-th entries of both sides of the equality, 
we have 
\[
b_{4, 3} 
= 
 b_{2, 1}\left( \frac{s}{1 + s} \right) \, (1 + s)^{p^{e_1}}  ,  
\]
which implies $b_{2, 1} = s^{p^{e_1}}$.

Comparing the $(3, 3)$-th entries of both sides of the equality, 
we have 
\[
1 
= 
b_{3, 1}\left( \frac{s}{1 + s} \right) \, (1 + s)^{p^{e_1}} + 1  ,  
\]
which implies $b_{3, 1} = 0$.

Comparing the $(2, 2)$-th entries of both sides of the equality, 
we have 
\[
1 + b_{4, 2}
= 
1  ,  
\]
which implies $b_{4, 2} = 0$.

Comparing the $(3, 2)$-th entries of both sides of the equality, 
we have 
\[
b_{3, 2}
= 
b_{3, 2} \left( \frac{s}{1 + s} \right) , 
\]
which implies $b_{3, 2} \in k$. Since $b_{3, 2} \in \sfP$, we have $b_{3, 2} = 0$.

Comparing the $(4, 3)$-th entries of both sides of the equality, 
we have 
\[
b_{4, 3}
= 
b_{4, 1}\left( \frac{s}{1 + s} \right) \, (1 + s)^{p^{e_1}}
   + b_{4, 3} \left( \frac{s}{1 + s} \right) , 
\]
which implies 
\[
 s^{p^{e_1}} 
  = 
  b_{4, 1}\left( \frac{s}{1 + s} \right) \, (1 + s)^{p^{e_1}} 
   + \left( \frac{s}{1 + s} \right)^{p^{e_1}} . 
\]
Thus $b_{4, 1} = s^{2 \, p^{e_1}}$.

Hence $\phi^-$ has the desired form.

\end{proof}

\subsubsection{$\rm (XXII)$}

\begin{lem}
Let $(\varphi^*, \omega^*)$ be of the form {\rm (XXII)}. 
Assume that there exists a homomorphism $\phi^- : \G_a \to \SL(4, k)$ 
satisfying the following conditions {\rm (i)} and {\rm (ii)}: 
\begin{enumerate}[label = {\rm (\roman*)}]
\item For any $s \in \G_a$, 
the regular matrix $\phi^-(s)$ is a lower triangular matrix.  

\item 
$
 \varphi^*(t) \, \phi^-(s) 
 = 
\phi^-\left( \frac{s}{1 + t \, s} \right) 
\, 
\omega^*(1 + t \, s) 
\;
\varphi^* \left( \frac{t}{1 + t \, s} \right)
$ 
for all $t, s \in k$ with $1 + t \, s \ne 0$.    
\end{enumerate} 
Then we have $d_1 = d_2 = p^{e_1}$, and we can express $\phi^-$ as 
\[
\phi^-(s) 
= 
\left(
\begin{array}{c c c c}
 1 & 0 & 0 & 0 \\
 0 & 1 & 0 & 0 \\
 s^{p^{e_1}} & 0 & 1 & 0 \\
 0 & s^{p^{e_1}} & 0 & 1  
\end{array}
\right) . 
\]
\end{lem}

\begin{proof} 

Write $\phi^-(s) = ( \, b_{i, j}(s) \, )_{1 \leq i, j \leq 4}$. 
By condition (ii), we have 
\begin{align*}
& \renewcommand{\arraystretch}{1.5} 
\left(
\begin{array}{r r r r}
 1 + b_{3, 1} 
 & b_{3, 2} 
 & 1 
 & 0 \\
 b_{2, 1} + b_{4, 1} 
 & 1 + b_{4, 2} 
 & b_{4, 3}
 & 1 \\ 
                                       b_{3, 1} 
                                        & b_{3, 2} 
                                        & 1 
                                        & 0 \\
                                                    b_{4, 1} 
                                                    & b_{4, 2} 
                                                    & b_{4, 3} 
                                                    & 1
\end{array}
\right) \\
&
= 
\left(
\begin{array}{l l}
 (1 + s)^{d_1} 
 & 0 \\
 b_{2, 1}\left( \frac{s}{1 + s} \right) \,  (1 + s)^{d_1} 
 & (1 + s)^{2 \, p^{e_1} - d_1}  \\
  b_{3, 1}\left( \frac{s}{1 + s} \right) \,  (1 + s)^{d_1} 
 &  b_{3, 2}\left( \frac{s}{1 + s} \right) \, (1 + s)^{2 \, p^{e_1} - d_1}     \\
   b_{4, 1}\left( \frac{s}{1 + s} \right) \, (1 + s)^{d_1} 
 &  b_{4, 2}\left( \frac{s}{1 + s} \right)  \, (1 + s)^{2 \, p^{e_1} - d_1} 
\end{array} 
\right.
\\
& \qquad 
\left. 
\begin{array}{l l}
    (1 + s)^{d_1 - p^{e_1}} 
  & 0 \\
  b_{2, 1}\left( \frac{s}{1 + s} \right) \, (1 + s)^{d_1 - p^{e_1}} 
  & \frac{1}{(1 + s)^{d_1 - p^{e_1}}}  \\
 b_{3, 1}\left( \frac{s}{1 + s} \right) \, (1 + s)^{d_1 - p^{e_1}} 
 + \frac{1}{ (1 + s)^{2 \, p^{e_1} - d_1} }
& 
  b_{3, 2}\left( \frac{s}{1 + s} \right) \, \frac{1}{(1 + s)^{d_1 - p^{e_1}}}  
\\ 
   b_{4, 1}\left( \frac{s}{1 + s} \right) \,  (1 + s)^{d_1 - p^{e_1}} 
   + b_{4, 3} \left( \frac{s}{1 + s} \right) \, \frac{1}{ (1 + s)^{2 \, p^{e_1} - d_1} }
    & 
  b_{4, 2}\left( \frac{s}{1 + s} \right) \, \frac{1}{(1 + s)^{d_1 - p^{e_1}}}   
  + \frac{1}{(1 + s)^{d_1}}
\end{array}
\right) . 
\end{align*}

Comparing the $(1, 3)$-th entries of both sides of the equality, 
we have 
\[
1 
= 
 (1 + s)^{d_1 - p^{e_1}} , 
\]
which implies 
\[
 d_1 = p^{e_1} . 
\]

Comparing the $(1, 1)$-th entries of both sides of the equality, 
we have 
\[
  1 + b_{3, 1}  
   = 
   (1 + s)^{d_1} , 
\]
which implies $b_{3, 1} = s^{p^{e_1}}$.

Comparing the $(1, 2)$-th entries of both sides of the equality, 
we have 
\[
b_{3, 2} = 0 . 
\]

Comparing the $(2, 2)$-th entries of both sides of the equality, 
we have 
\[
1 + b_{4, 2}
= 
 (1 + s)^{2 \, p^{e_1} - d_1}  , 
\]
which implies $ 1 + b_{4, 2} = (1 + s)^{p^{e_1}}$. 
So, $b_{4, 2} = s^{p^{e_1}}$.

Comparing the $(3, 4)$-th entries of both sides of the equality, 
we have 
\[
0 
 = 
 b_{3, 2}\left( \frac{s}{1 + s} \right) \, \frac{1}{(1 + s)^{d_1 - p^{e_1}}}  , 
\]
which implies $b_{3, 2} = 0$.

Comparing the $(2, 3)$-th, $(4, 3)$-th, $(2, 1)$-th entries of both sides of 
the equality, we have 
\begin{align*}
\left\{
\begin{array}{r @{\,} l @{\qquad} l}
b_{4, 3} 
& = 
b_{2, 1}\left( \frac{s}{1 + s} \right) 
& \textcircled{\scriptsize 1}  \\
b_{4, 3} 
& = b_{4, 1}\left( \frac{s}{1 + s} \right) \,  (1 + s)^{d_1 - p^{e_1}} 
   + b_{4, 3} \left( \frac{s}{1 + s} \right) \, \frac{1}{ (1 + s)^{2 \, p^{e_1} - d_1} } 
& \textcircled{\scriptsize 2}  
\\
b_{2, 1} + b_{4, 1} 
& = 
b_{2, 1}\left( \frac{s}{1 + s} \right) \, (1 + s)^{d_1} 
& \textcircled{\scriptsize 3} 
\end{array} 
\right. 
\end{align*}
From $\textcircled{\scriptsize 1}$ and $\textcircled{\scriptsize 2}$, 
we have 
\[
 b_{2, 1}\left( \frac{s}{1 + s} \right)  
  = 
  b_{4, 1}\left( \frac{s}{1 + s} \right) 
   + b_{4, 3} \left( \frac{s}{1 + s} \right) \, \frac{1}{ (1 + s)^{p^{e_1} } } , 
\]
which implies 
\[
 b_{2, 1}(x) = b_{4, 1}(x) + b_{4, 3}(x) \, (1 - x)^{p^{e_1}} . 
\]
Thus, 
\[
 b_{2, 1}(s) - b_{4, 1}(s) = b_{4, 3}(s) \, (1 - s)^{p^{e_1}} \qquad  
 \textcircled{\scriptsize 4} 
\]

From $\textcircled{\scriptsize 1}$, $\textcircled{\scriptsize 3}$ and $d_1 = p^{e_1}$, 
we have 
\[
 b_{2, 1} + b_{4, 1}
  =  b_{4, 3} \, (1 + s)^{p^{e_1}} \qquad  
 \textcircled{\scriptsize 5} 
\]

From $\textcircled{\scriptsize 4}$ and $\textcircled{\scriptsize 5}$, 
we have 
\[
 b_{2, 1} = b_{4, 3} . 
\]
By $\textcircled{\scriptsize 1}$, we have $b_{2, 1} \in k$. 
Since $b_{2, 1} \in \sfP$, we have $b_{2, 1} = 0$. 
Therefore, $b_{4, 3} = 0$. 
By $\textcircled{\scriptsize 3}$, we have $b_{4, 1} = 0$. 

Hence $\phi^-$ has the desired form.

\end{proof}

\subsubsection{$\rm (XXIII)$}

\begin{lem}
Let $(\varphi^*, \omega^*)$ be of the form {\rm (XXIII)}. 
Assume that there exists a homomorphism $\phi^- : \G_a \to \SL(4, k)$ 
satisfying the following conditions {\rm (i)} and {\rm (ii)}: 
\begin{enumerate}[label = {\rm (\roman*)}]
\item For any $s \in \G_a$, 
the regular matrix $\phi^-(s)$ is a lower triangular matrix.  

\item 
$
 \varphi^*(t) \, \phi^-(s) 
 = 
\phi^-\left( \frac{s}{1 + t \, s} \right) 
\, 
\omega^*(1 + t \, s) 
\;
\varphi^* \left( \frac{t}{1 + t \, s} \right)
$ 
for all $t, s \in k$ with $1 + t \, s \ne 0$.    
\end{enumerate} 
Then we have a contradiction. 
\end{lem}

\begin{proof}
Write $\phi^-(s) = ( \, b_{i, j}(s) \, )_{1 \leq i, j \leq 4}$. 
By condition (ii), we have 
\begin{align*}
& \renewcommand{\arraystretch}{1.5} 
\left(
\begin{array}{r r r r}
 1 + b_{3, 1} 
 & b_{3, 2} 
 & 1 
 & 0 \\
 b_{2, 1} 
 & 1 
 & 0 
 & 0 \\ 
                                       b_{3, 1} 
                                        & b_{3, 2} 
                                        & 1 
                                        & 0 \\
                                                    b_{4, 1} 
                                                    & b_{4, 2} 
                                                    & b_{4, 3} 
                                                    & 1
\end{array}
\right) \\
&
= 
\left(
\begin{array}{l l}
 (1 + s)^{d_1} 
 & 0 \\
 b_{2, 1}\left( \frac{s}{1 + s} \right) \,  (1 + s)^{d_1} 
 & (1 + s)^{2 \, p^{e_1} - d_1}  \\
  b_{3, 1}\left( \frac{s}{1 + s} \right) \,  (1 + s)^{d_1} 
 &  b_{3, 2}\left( \frac{s}{1 + s} \right) \, (1 + s)^{2 \, p^{e_1} - d_1}     \\
   b_{4, 1}\left( \frac{s}{1 + s} \right) \, (1 + s)^{d_1} 
 &  b_{4, 2}\left( \frac{s}{1 + s} \right)  \, (1 + s)^{2 \, p^{e_1} - d_1} 
\end{array} 
\right.
\\
& \qquad 
\left. 
\begin{array}{l l}
    (1 + s)^{d_1 - p^{e_1}} 
  & 0 \\
  b_{2, 1}\left( \frac{s}{1 + s} \right) \, (1 + s)^{d_1 - p^{e_1}} 
  & 0   \\
 b_{3, 1}\left( \frac{s}{1 + s} \right) \, (1 + s)^{d_1 - p^{e_1}} 
 + \frac{1}{ (1 + s)^{2 \, p^{e_1} - d_1} }
& 
  0
\\ 
   b_{4, 1}\left( \frac{s}{1 + s} \right) \,  (1 + s)^{d_1 - p^{e_1}} 
   + b_{4, 3} \left( \frac{s}{1 + s} \right) \, \frac{1}{ (1 + s)^{2 \, p^{e_1} - d_1} }
    & 
\frac{1}{(1 + s)^{d_1}}
\end{array}
\right) . 
\end{align*}

Comparing the $(4, 4)$-th entries of both sides of the equality, 
we have 
\[
1 
= 
 \frac{1}{(1 + s)^{d_1}}  , 
\]
which implies $d_1 = 0$. 
Since $d_1 = p^{e_1} \geq 1$, we have a contradiction. 

\end{proof}

\subsubsection{$\rm (XXIV)$}

\begin{lem}
Let $(\varphi^*, \omega^*)$ be of the form {\rm (XXIV)}. 
Assume that there exists a homomorphism $\phi^- : \G_a \to \SL(4, k)$ 
satisfying the following conditions {\rm (i)} and {\rm (ii)}: 
\begin{enumerate}[label = {\rm (\roman*)}]
\item For any $s \in \G_a$, 
the regular matrix $\phi^-(s)$ is a lower triangular matrix.  

\item 
$
 \varphi^*(t) \, \phi^-(s) 
 = 
\phi^-\left( \frac{s}{1 + t \, s} \right) 
\, 
\omega^*(1 + t \, s) 
\;
\varphi^* \left( \frac{t}{1 + t \, s} \right)
$ 
for all $t, s \in k$ with $1 + t \, s \ne 0$.    
\end{enumerate} 
Then we have $d_2 = 0$, and we can express $\phi^-$ as 
\[
\phi^-(s) 
= 
\left(
\begin{array}{c c c c}
 1 & 0 & 0 & 0 \\
 0 & 1 & 0 & 0 \\
 0 & 0 & 1 & 0 \\
 s^{p^{e_2}} & 0 & 0 & 1  
\end{array}
\right) . 
\]
\end{lem}

\begin{proof}

Write $\phi^-(s) = ( \, b_{i, j}(s) \, )_{1 \leq i, j \leq 4}$. 
By condition (ii), we have 
\begin{align*}
& \renewcommand{\arraystretch}{1.5} 
\left(
\begin{array}{r r r r}
 1 + b_{4, 1} 
 & b_{4, 2} 
 & b_{4, 3}
 & 1 \\
 b_{2, 1} 
 & 1 
 & 0 
 & 0 \\ 
                                       b_{3, 1} 
                                        & b_{3, 2} 
                                        & 1 
                                        & 0 \\
                                                    b_{4, 1} 
                                                    & b_{4, 2} 
                                                    & b_{4, 3} 
                                                    & 1
\end{array}
\right) \\
&
= 
\left(
\begin{array}{l l}
 (1 + s)^{p^{e_2}} 
 & 0 \\
 b_{2, 1}\left( \frac{s}{1 + s} \right) \,  (1 + s)^{p^{e_2}} 
 & (1 + s)^{d_2}  \\
  b_{3, 1}\left( \frac{s}{1 + s} \right) \,  (1 + s)^{p^{e_2}} 
 &  b_{3, 2}\left( \frac{s}{1 + s} \right) \, (1 + s)^{d_2}     \\
   b_{4, 1}\left( \frac{s}{1 + s} \right) \, (1 + s)^{p^{e_2}} 
 &  b_{4, 2}\left( \frac{s}{1 + s} \right)  \, (1 + s)^{d_2} 
\end{array} 
\right.
\\
& \qquad 
\left. 
\begin{array}{l l}
   0
  & 1 \\
  0 
  &  b_{2, 1}\left( \frac{s}{1 + s} \right)   \\
 \frac{1}{ (1 + s)^{d_2} }
& 
   b_{3, 1}\left( \frac{s}{1 + s} \right) 
\\ 
 b_{4, 3} \left( \frac{s}{1 + s} \right) \, \frac{1}{ (1 + s)^{d_2} }
    &  b_{4, 1}\left( \frac{s}{1 + s} \right) \, \frac{1}{(1 + s)^{p^{e_2}}}
\end{array}
\right) . 
\end{align*}

Comparing the $(3, 3)$-th entries of both sides of the equality, 
we have 
\[
1
 = 
 \frac{1}{(1 + s)^{d_2}}  ,  
\]
which implies $d_2 = 0$.

Comparing the $(1, 3)$-th entries of both sides of the equality, 
we have 
\[
 b_{4, 3} = 0 . 
\]

Comparing the $(2, 4)$-th entries of both sides of the equality, 
we have 
\[
 0 =  b_{2, 1}\left( \frac{s}{1 + s} \right)   , 
\]
which implies $b_{2, 1} = 0$.

Comparing the $(3, 4)$-th entries of both sides of the equality, 
we have 
\[
 0 = b_{3, 1}\left( \frac{s}{1 + s} \right)  , 
\]
which implies $b_{3, 1} = 0$.

Comparing the $(1, 2)$-th entries of both sides of the equality, 
we have 
\[
 b_{4, 2} = 0 . 
\]

Comparing the $(3, 2)$-th entries of both sides of the equality, 
we have 
\[
 b_{3, 2} 
  = 
  b_{3, 2}\left( \frac{s}{1 + s} \right) , 
\]
which implies $b_{3, 2} \in k$. Since $b_{3, 2} \in \sfP$, we have $b_{3, 2} = 0$.

Comparing the $(1, 1)$-th entries of both sides of the equality, 
we have 
\[
 1 + b_{4, 1} 
  = (1 + s)^{p^{e_2}} , 
\]
which implies $b_{4, 1} = s^{p^{e_2}}$.

Hence $\phi^-$ has the desired form. 
\end{proof}

\subsubsection{$\rm (XXV)$}

\begin{lem}
Let $(\varphi^*, \omega^*)$ be of the form {\rm (XXV)}. 
Assume that there exists a homomorphism $\phi^- : \G_a \to \SL(4, k)$ 
satisfying the following conditions {\rm (i)} and {\rm (ii)}: 
\begin{enumerate}[label = {\rm (\roman*)}]
\item For any $s \in \G_a$, 
the regular matrix $\phi^-(s)$ is a lower triangular matrix.  

\item 
$
 \varphi^*(t) \, \phi^-(s) 
 = 
\phi^-\left( \frac{s}{1 + t \, s} \right) 
\, 
\omega^*(1 + t \, s) 
\;
\varphi^* \left( \frac{t}{1 + t \, s} \right)
$ 
for all $t, s \in k$ with $1 + t \, s \ne 0$.    
\end{enumerate} 
Then we have a contradiction. 
\end{lem}

\begin{proof}

The pair $({^\tau}\varphi^*, \omega^*)$ has the form (XXIII). 
Since $\psi_{\, {^\tau}\varphi^*,\; \omega^*} \circ \jmath^{-1}$ is not extendable, 
$\psi_{\varphi^*, \, \omega^*} \circ \jmath^{-1}$ is not extendable (see 
Lemma 2.8 (1)). 

\end{proof}

\subsubsection{$\rm (XXVI)$}

\begin{lem}
Let $(\varphi^*, \omega^*)$ be of the form {\rm (XXVI)}. 
Assume that there exists a homomorphism $\phi^- : \G_a \to \SL(4, k)$ 
satisfying the following conditions {\rm (i)} and {\rm (ii)}: 
\begin{enumerate}[label = {\rm (\roman*)}]
\item For any $s \in \G_a$, 
the regular matrix $\phi^-(s)$ is a lower triangular matrix.  

\item 
$
 \varphi^*(t) \, \phi^-(s) 
 = 
\phi^-\left( \frac{s}{1 + t \, s} \right) 
\, 
\omega^*(1 + t \, s) 
\;
\varphi^* \left( \frac{t}{1 + t \, s} \right)
$ 
for all $t, s \in k$ with $1 + t \, s \ne 0$.    
\end{enumerate} 
Then we can express $\phi^-$ as 
\[
\phi^-(s) 
= 
I_4 . 
\]
\end{lem}

\begin{proof}
See \cite[Lemma 2.9]{Tanimoto 2022}. 
\end{proof}

\subsection{The forms of homomorphisms $\sigma^* : \SL(2, k) \to \SL(4, k)$}

In this subsection, for any pair $(\varphi^*, \omega^*)$ of the form 
$(\nu)$, where $\nu = \rm (I)$, $\rm (II)$, $\rm (IV)$, $\rm (V)$, $\rm (VII)$, 
$\rm (IX)$, $\rm (XI)$, 
$\rm (XV)$, $\rm (XIX)$, $\rm (XXI)$, $\rm (XXII)$, 
$\rm (XXIV)$, $\rm (XXVI)$, 
assuming that there exists a homomorphism 
$\sigma^* : \SL(2, k) \to \SL(4, k)$ such that 
$\sigma^* \circ \imath_{\rB(2, k)} = \psi_{\varphi^*, \, \omega^*}$, 
we find the form of $\sigma^*$. 
While obtaing the form, we use the equality 
\begin{align*}
\left(
\begin{array}{c c}
 a & b \\
 c & d
\end{array}
\right) 
= 
\left(
\begin{array}{c c}
 1 & 0 \\
 \frac{c}{a} & 1
\end{array}
\right) 
\left(
\begin{array}{c c}
 a & 0 \\
 0 & \frac{1}{a} 
\end{array}
\right) 
\left(
\begin{array}{c c}
 1 & \frac{b}{a}  \\
 0 & 1 
\end{array}
\right) 
\end{align*}
for any regular matrix 
\[
\left(
\begin{array}{c c}
 a & b \\
 c & d
\end{array}
\right) 
\in \SL(2, k) , \qquad a \ne0 . 
\]

\subsubsection{$\rm (I)$}

\begin{lem}
Let $(\varphi^*, \omega^*)$ be of the form {\rm (I)}. 
If there exists a homomorphsim $\sigma^* : \SL(2, k) \to \SL(4, k)$ such that 
$\sigma^* \circ \imath_{\rB(2, k)} = \psi_{\varphi^*, \, \omega^*}$. 
Then 
\[
\renewcommand{\arraystretch}{1.5} 
\sigma^*
\left(
\begin{array}{c c}
 a & b \\
 c & d
\end{array}
\right) 
 := 
\left(
\begin{array}{l l l l}
 a^{3 \, p^{e_1}} 
 & a^{2 \, p^{e_1}} \, b^{p^{e_1}} 
 & \frac{1}{2} \, a^{p^{e_1}} \, b^{2 \, p^{e_1}} 
 & \frac{1}{6} \, b^{3 \, p^{e_1}} \\
 3 \, a^{2 \, p^{e_1}} \, c^{p^{e_1}} 
 & a^{p^{e_1}} \cdot (a \, d + 2 \, b \, c)^{p^{e_1}} 
 & b^{p^{e_1}} \cdot (a \, d + \frac{1}{2} \, b \, c)^{p^{e_1}} 
 & \frac{1}{2} \, b^{2 \, p^{e_1}} \, d^{p^{e_1}} \\
 6 \, a^{p^{e_1}} \, c^{2 \, p^{e_1}} 
 &  4 \, c^{p^{e_1}} \cdot (a \, d + \frac{1}{2} \, b \, c)^{p^{e_1}} 
 & d^{p^{e_1}} \cdot (a \, d + 2 \, b \, c)^{p^{e_1}} 
 & b^{p^{e_1}} \, d^{2 \, p^{e_1}} \\
 6 \, c^{3 \, p^{e_1}} 
 & 6 \, c^{2 \, p^{e_1}} \, d^{p^{e_1}} 
 & 3 \, c^{p^{e_1}} \, d^{2 \, p^{e_1}} 
 & d^{3 \, p^{e_1}} 
\end{array}
\right) . 
\]
\end{lem}

\begin{proof}
We have 
\begin{align*}
\renewcommand{\arraystretch}{1.5} 
\sigma^*
\left(
\begin{array}{c c}
 a & b \\
 c & d
\end{array}
\right) 
 & = \phi^-\left(\frac{c}{a} \right) \, \omega^*(a) \, \varphi^*\left(\frac{b}{a}\right) \\
 & = 
 \renewcommand{\arraystretch}{1.5} 
\left(
\begin{array}{l l l}
 a^{3 \,p^e_1} 
 & a^{2 \,  p^{e_1}} \, b^{p^{e_1}} 
 & \frac{1}{2} \, a^{p^{e_1}} \, b^{2 \, p^{e_1}} 
\\
 3 \, a^{2 \, p^{e_1}} \, c^{p^{e_1}} 
 & 3 \, (a \, b \, c)^{p^{e_1}} + a^{p^{e_1}} 
 & \frac{3}{2} \, b^{2 \, p^{e_1}} \, c^{p^{e_1}} + b^{p^{e_1}} \\
 6 \, a^{p^{e_1}} \, c^{2 \, p^{e_1}} 
 & 6 \, b^{p^{e_1}} \, c^{2 \, p^{e_1}} + 4 \, c^{p^{e_1}} 
 & \frac{ 3 \, b^{2 \, p^{e_1}} \, c^{2 \, p^{e_1}} + 4 \, b^{p^{e_1}} \, c^{p^{e_1}} + 1}{a^{p^{e_1}}}\\
 6 \, c^{3 \, p^{e_1}}
 & 6 \, c^{2 \, p^{e_1}} \, \left( \frac{b \, c + 1}{a} \right)^{p^{e_1}}  
 & \frac{3 \, b^{2 \, p^{e_1}} \, c^{3 \, p^{e_1}} + 6 \, b^{p^{e_1}} \, c^{2 \, p^{e_1}} + 3 \, c^{p^{e_1}} }{a^{2 \, p^{e_1}}}
\end{array}
\right. \\
& 
 \renewcommand{\arraystretch}{1.5} 
\qquad 
\left. 
\begin{array}{l}
 \frac{1}{6} \, b^{3 \, p^{e_1}} \\
 \frac{1}{2} \, b^{2 \, p^{e_1}} \,
  \left( \frac{ b^{p^{e_1}} \, c^{p^{e_1}} + 1 }{ a^{ p^{e_1} }  } \right) \\
 \frac{ b^{3 \, p^{e_1}} \, c^{2 \, p^{e_1}} + 2 \, b^{2 \, p^{e_1}} \, c^{p^{e_1}} + b^{p^{e_1}} }{ a^{2 \, p^{e_1}} } \\
 \frac{ b^{3 \, p^{e_1}} \, c^{3 \, p^{e_1}} + 3 \, b^{2 \, p^{e_1}} \, c^{2 \, p^{e_1}} 
 + 3 \, b^{p^{e_1}} \, c^{p^{e_1}} + 1}{ a^{3 \, p^{e_1}} }
\end{array}
\right) . 
\end{align*}
Since $d = (b \, c + 1) / a$, 
we have 
\begin{align*}
\bigl( \text{ the $(2, 2)$-th entry } \bigr) 
 & = a^{p^{e_1}} \, (3 \, (b \, c)^{p^{e_1}} + 1) \\
& = a^{p^{e_1}} \, (3 \, (b \, c)^{p^{e_1}} + (a \, d - b \, c)^{p^{e_1}}) \\
 &= a^{p^{e_1}} \cdot (a \, d + 2 \, b \, c)^{p^{e_1}}  , \\
\bigl( \text{ the $(3, 2)$-th entry } \bigr) 
 & = c^{p^{e_1}} \cdot ( 6 \, b^{p^{e_1}} \, c^{p^{e_1}} + 4 ) \\ 
 & = c^{p^{e_1}} \cdot \bigl( 2 \, b^{p^{e_1}} \, c^{p^{e_1}} + 4  \, (a^{p^{e_1}} \, d^{p^{e_1}} - 1 ) + 4 \bigr) \\
 & = 4 \, c^{p^{e_1}} \cdot \left( a \, d + \frac{1}{2} \, b \, c \right)^{p^{e_1}} , \\ 
\bigl( \text{ the $(4, 2)$-th entry } \bigr) 
 & = 6 \, c^{2 \, p^{e_1}} \, d^{p^{e_1}} , \\
\bigl( \text{ the $(2, 3)$-th entry } \bigr) 
 & = \frac{1}{2} \, b^{2 \, p^{e_1}} \, c^{p^{e_1}} + b^{p^{e_1}} \cdot (a^{p^{e_1}} \, d^{p^{e_1}} - 1) 
 + b^{p^{e_1}} \\
 & =  b^{p^{e_1}} \cdot \left(  a^{p^{e_1}} \, d^{p^{e_1}}  + 
 \frac{1}{2} \, b^{p^{e_1}} \, c^{p^{e_1}} \right)  , 
 \\
\bigl( \text{ the $(3, 3)$-th entry } \bigr) 
 & = \frac{ 3 \, (a \, d - 1)^{2 \,p^{e_1}} + 4 \, ( a \, d - 1)^{p^{e_1}} + 1}{ a^{p^{e_1}} } \\ 
 & = 3 \, a^{p^{e_1}} \, d^{2 \, p^{e_1}} - 2 \, d^{p^{e_1}} \\
 & = d^{p^{e_1}} \cdot \left( 3 \, a^{p^{e_1}} \, d^{p^{e_1}} - 2 \right) \\
 & = d^{p^{e_1}} \cdot \left( a^{p^{e_1}} \, d^{p^{e_1}} + 2 \, b^{p^{e_1}} \, c^{p^{e_1}}  \right) , \\
\bigl( \text{ the $(4, 3)$-th entry } \bigr) 
 & = \frac{3 \, c^{p^{e_1}} \, 
 ( b^{2 \, p^{e_1}} \, c^{2 \, p^{e_1}} + 2 \, b^{p^{e_1}} \, c^{p^{e_1}} + 1) }{a^{2 \, p^{e_1}}} \\ 
 & = 3 \, c^{p^{e_1}} \, d^{2 \, p^{e_1}} , \\ 
\bigl( \text{ the $(2, 4)$-th entry } \bigr) 
 & = \frac{1}{2} \, b^{2 \, p^{e_1}} \cdot 
  \frac{ (b \, c + 1)^{p^{e_1}} }{ a^{ p^{e_1} }  } 
  = \frac{1}{2} \, b^{2 \, p^{e_1}} \, d^{p^{e_1}} ,  \\
\bigl( \text{ the $(3, 4)$-th entry } \bigr) 
 & = 
  \frac{ b^{p^{e_1}} \, 
  ( b^{2 \, p^{e_1}} \, c^{2 \, p^{e_1}} + 2 \, b^{ p^{e_1}} \, c^{p^{e_1}} + 1 ) }{ a^{2 \, p^{e_1}} }  
  = 
  b^{p^{e_1}} \, d^{2 \, p^{e_1}} , \\
\bigl( \text{ the $(4, 4)$-th entry } \bigr) 
 & =   \frac{  (b \, c + 1)^{3 \, p^{e_1}} }{ a^{3 \, p^{e_1}} } 
  = d^{3 \, p^{e_1}} . 
\end{align*}

\end{proof}

\subsubsection{${\rm (II)}$}

\begin{lem}
Let $(\varphi^*, \omega^*)$ be of the form {\rm (II)}. 
If there exists a homomorphsim $\sigma^* : \SL(2, k) \to \SL(4, k)$ such that 
$\sigma^* \circ \imath_{\rB(2, k)} = \psi_{\varphi^*, \, \omega^*}$. 
Then 
\begin{align*}
\renewcommand{\arraystretch}{1.5} 
\sigma^*
\left(
\begin{array}{c c}
 a & b \\
 c & d
\end{array}
\right) 
 = 
\left(
\begin{array}{l l l l}
 a^{p^{e_1 + 1}} 
 & a^{2 \, p^{e_1}} \, b^{p^{e_1}} 
 & \frac{1}{2} \, a^{p^{e_1}} \, b^{2 \, p^{e_1}} 
 & b^{p^{e_1 + 1}} \\
 0
 &  a^{p^{e_1}} 
 &  b^{p^{e_1}}
 &  0  \\
 0
 & c^{p^{e_1}} 
 & d^{p^{e_1}} 
 & 0 \\
 c^{p^{e_1 + 1}} 
 & c^{2 \, p^{e_1}} \, d^{p^{e_1}} 
 & \frac{1}{2} \, c^{p^{e_1}} \, d^{2 \, p^{e_1}} 
 & d^{p^{e_1 + 1}} 
\end{array}
\right) . 
\end{align*}
\end{lem}

\begin{proof}
We have 
\begin{align*}
\renewcommand{\arraystretch}{1.5} 
\sigma^*
\left(
\begin{array}{c c}
 a & b \\
 c & d
\end{array}
\right) 
 & = \phi^-\left(\frac{c}{a} \right) \, \omega^*(a) \, \varphi^*\left(\frac{b}{a}\right) \\
 & = 
 \renewcommand{\arraystretch}{1.5} 
\left(
\begin{array}{l l l l}
a^{p^{e_1 + 1}} 
 & a^{2 \, p^{e_1}} \, b^{p^{e_1} } 
 & \frac{1}{2} \, a^{p^{e_1}} \, b^{2 \, p^{e_1}} 
 & b^{p^{e_1 + 1}} 
\\ 
0
 & a^{p^{e_1}} 
 & b^{p^{e_1}} 
 & 0 \\
0
 & c^{p^{e_1}}  
 & \frac{1 + b^{p^{e_1}} \, c^{p^{e_1}} }{a^{p^{e_1}}}
 & 0 \\
c^{p^{e_1 + 1}} 
 & \frac{b^{p^{e_1}} \, c^{p^{e_1 + 1}} + c^{2 \, p^{e_1}}}{a^{p^{e_1}}}
 & \frac{b^{2 \, p^{e_1}} \, c^{3 \, p^{e_1}} + 2 \, b^{p^{e_1}} \, c^{2 \, p^{e_1}} + c^{p^{e_1}}}{ 2 \,  a^{2 \, p^{e_1}} }
  & 
   \frac{ b^{p^{e_1 + 1}} \, c^{p^{e_1 + 1}} + 1 }{ a^{p^{e_1 + 1}} } 
\end{array}
\right) . 
\end{align*}
Since $d = (b \, c + 1) / a$, 
we have 
\begin{align*}
\bigl( \text{ the $(4, 2)$-th entry } \bigr) 
 & =  \frac{b^{p^{e_1}} \, c^{p^{e_1 + 1}} + c^{2 \, p^{e_1}}}{a^{p^{e_1}}} 
 = c^{2 \, p^{e_1}} \, d^{p^{e_1}} , \\
\bigl( \text{ the $(3, 3)$-th entry } \bigr) 
 & =  d^{p^{e_1}} , \\
\bigl( \text{ the $(4, 3)$-th entry } \bigr) 
 & =  \frac{1}{2} \, c^{p^{e_1}} \cdot 
 \frac{b^{2 \, p^{e_1}} \, c^{2 \, p^{e_1}} + 2 \, b^{p^{e_1}} \, c^{p^{e_1}} + 1}{ a^{2 \, p^{e_1}} } \\
 & = \frac{1}{2} \, c^{p^{e_1}} \left( \frac{b \, c + 1}{a} \right)^{2 \, p^{e_1}} \\
 & = \frac{1}{2} \, c^{p^{e_1}} \, d^{2\, p^{e_1}} , \\
\bigl( \text{ the $(4, 4)$-th entry } \bigr) 
 & = d^{p^{e_1 + 1}} . 
\end{align*}

\end{proof}

\subsubsection{$\rm (IV)$}

\begin{lem}
Let $(\varphi^*, \omega^*)$ be of the form {\rm (IV)}. 
If there exists a homomorphsim $\sigma^* : \SL(2, k) \to \SL(4, k)$ such that 
$\sigma^* \circ \imath_{\rB(2, k)} = \psi_{\varphi^*, \, \omega^*}$. 
Then 
\begin{align*}
\renewcommand{\arraystretch}{1.5} 
\sigma^*
\left(
\begin{array}{c c}
 a & b \\
 c & d
\end{array}
\right) 
& = 
\renewcommand{\arraystretch}{1.5} 
\left(
\begin{array}{l l l l}
 a^{p^{e_2}} \cdot a^{p^{e_1}} 
 & a^{p^{e_2}} \cdot b^{p^{e_1}} 
 & b^{p^{e_2}} \cdot a^{p^{e_1}} 
 & b^{p^{e_2}} \cdot b^{p^{e_1}} \\
  a^{p^{e_2}} \cdot c^{p^{e_1}} 
 & a^{p^{e_2}} \cdot d^{p^{e_1}} 
 & b^{p^{e_2}} \cdot c^{p^{e_1}} 
 & b^{p^{e_2}} \cdot d^{p^{e_1}} \\
 c^{p^{e_2}} \cdot a^{p^{e_1}} 
 & c^{p^{e_2}} \cdot b^{p^{e_1}} 
 & d^{p^{e_2}} \cdot a^{p^{e_1}} 
 & d^{p^{e_2}} \cdot b^{p^{e_1}} \\
 c^{p^{e_2}} \cdot c^{p^{e_1}} 
 & c^{p^{e_2}} \cdot d^{p^{e_1}} 
 & d^{p^{e_2}} \cdot c^{p^{e_1}} 
 & d^{p^{e_2}} \cdot d^{p^{e_1}}
\end{array}
\right) \\
& = 
\renewcommand{\arraystretch}{1.5} 
\left(
\begin{array}{c c}
a^{p^{e_2}}
\left(
\begin{array}{c c}
 a^{p^{e_1}} & b^{p^{e_1}} \\
 c^{p^{e_1}} & d^{p^{e_1}} 
\end{array}
\right)
& 
b^{p^{e_2}}
\left(
\begin{array}{c c}
 a^{p^{e_1}} & b^{p^{e_1}} \\
 c^{p^{e_1}} & d^{p^{e_1}} 
\end{array}
\right) \\ [1.5em] 
c^{p^{e_2}}
\left(
\begin{array}{c c}
 a^{p^{e_1}} & b^{p^{e_1}} \\
 c^{p^{e_1}} & d^{p^{e_1}} 
\end{array}
\right)
& 
d^{p^{e_2}}
\left(
\begin{array}{c c}
 a^{p^{e_1}} & b^{p^{e_1}} \\
 c^{p^{e_1}} & d^{p^{e_1}} 
\end{array}
\right)
\end{array}
\right) . 
\end{align*}
\end{lem}

\begin{proof}
We have 
\begin{align*}
\renewcommand{\arraystretch}{1.5} 
\sigma^*
\left(
\begin{array}{c c}
 a & b \\
 c & d
\end{array}
\right) 
 & = \phi^-\left(\frac{c}{a} \right) \, \omega^*(a) \, \varphi^*\left(\frac{b}{a}\right) \\
 & = 
 \renewcommand{\arraystretch}{1.5} 
\left(
\begin{array}{l l l}
a^{p^{e_1} + p^{e_2}} 
 & a^{p^{e_2}} \, b^{p^{e_1}} 
 & a^{p^{e_1}} \, b^{p^{e_2}} 
\\
a^{p^{e_2}} \, c^{p^{e_1}} 
 & a^{p^{e_2} - p^{e_1}} \, b^{p^{e_1}} \, c^{p^{e_1}} +  a^{p^{e_2} - p^{e_1}} 
 & b^{p^{e_2}} \, c^{p^{e_1}} \\
a^{p^{e_1}} \, c^{p^{e_2}} 
 & b^{p^{e_1}} \, c^{p^{e_2}} 
 & a^{p^{e_1} - p^{e_2}} \, b^{p^{e_2}} \, c^{p^{e_2}} + a^{p^{e_1} - p^{e_2}} \\
 c^{p^{e_1} + p^{e_2}} 
 & \frac{ b^{p^{e_1}} \, c^{p^{e_1} + p^{e_2}} + c^{p^{e_2}} }{ a^{p^{e_1}} }
 & \frac{ b^{p^{e_2}} \, c^{p^{e_1} + p^{e_2}} + c^{p^{e_1}} }{a^{p^{e_2}}}
\end{array}
\right. \\
& 
\qquad 
\renewcommand{\arraystretch}{1.5} 
\left. 
\begin{array}{l}
 b^{p^{e_1} + p^{e_2}} \\
 \frac{ b^{p^{e_1} + p^{e_2}} \, c^{p^{e_1}} + b^{p^{e_2}} }{a^{p^{e_1}}} \\
 \frac{ b^{p^{e_1} + p^{e_2}} \, c^{p^{e_2}} + b^{p^{e_1}} }{a^{p^{e_2}}}  \\
\frac{ b^{p^{e_1} + p^{e_2}} \, c^{p^{e_1} + p^{e_2}} 
+ b^{p^{e_2}} \, c^{p^{e_2}} 
+ b^{p^{e_1}} \, c^{p^{e_1}} 
+ 1 }{a^{p^{e_1} + p^{e_2}}} 
\end{array}
\right) . 
\end{align*}
Since $d = (b \, c + 1) / a$, 
we have 
\begin{align*}
\bigl( \text{ the $(2, 2)$-th entry } \bigr) 
 & =  a^{p^{e_2}} \, \left( \frac{ b^{p^{e_1}} \, c^{p^{e_1}} + 1 }{ a^{p^{e_1}} } \right) 
 = a^{p^{e_2}} \, d^{p^{e_1}} , \\
\bigl( \text{ the $(4, 2)$-th entry } \bigr) 
 & =  c^{p^{e_2}} \, \left( \frac{ b^{p^{e_1}} \, c^{p^{e_1}} + 1 }{ a^{p^{e_1}} } \right) 
 = c^{p^{e_2}} \, d^{p^{e_1}} , \\
\bigl( \text{ the $(3, 3)$-th entry } \bigr) 
 & = a^{p^{e_1}} \, \left( \frac{ b^{p^{e_2}} \, c^{p^{e_2}} + 1 }{ a^{p^{e_2}} } \right) 
 = a^{p^{e_1}} \, d^{p^{e_2}} , \\ 
\bigl( \text{ the $(4, 3)$-th entry } \bigr) 
 & = c^{p^{e_1}} \, \left( \frac{ b^{p^{e_2}} \, c^{p^{e_2}} + 1 }{ a^{p^{e_2}} } \right) 
 = c^{p^{e_1}} \, d^{p^{e_2}} , \\
\bigl( \text{ the $(2, 4)$-th entry } \bigr) 
 & =   b^{p^{e_2}} \, \left( \frac{ b^{p^{e_1}} \, c^{p^{e_1}} + 1 }{ a^{p^{e_1}} } \right) 
 = a^{p^{e_2}} \, d^{p^{e_1}} , 
 \\
\bigl( \text{ the $(3, 4)$-th entry } \bigr) 
 & = b^{p^{e_1}} \, \left( \frac{ b^{p^{e_2}} \, c^{p^{e_2}} + 1 }{ a^{p^{e_2}} } \right) 
 = b^{p^{e_1}} \, d^{p^{e_2}} , 
 \\
\bigl( \text{ the $(4, 4)$-th entry } \bigr) 
 & = \left( \frac{ b^{p^{e_1}} \, c^{p^{e_1}} + 1 }{ a^{p^{e_1}} } \right) 
 \left( \frac{ b^{p^{e_2}} \, c^{p^{e_2}} + 1 }{ a^{p^{e_2}} } \right) 
 = d^{p^{e_1} + p^{e_2}} . 
\end{align*}

\end{proof}

\subsubsection{$\rm (V)$}

\begin{lem}
Let $(\varphi^*, \omega^*)$ be of the form {\rm (V)}. 
If there exists a homomorphsim $\sigma^* : \SL(2, k) \to \SL(4, k)$ such that 
$\sigma^* \circ \imath_{\rB(2, k)} = \psi_{\varphi^*, \, \omega^*}$. 
Then $p = 2$ and 
\begin{align*}
\renewcommand{\arraystretch}{1.5} 
\sigma^*
\left(
\begin{array}{c c}
 a & b \\
 c & d
\end{array}
\right) 
 = 
\left(
\begin{array}{l l l l}
 a^{p^{e_1 + 1}} 
 & a^{p^{e_1}} \, b^{p^{e_1}} 
 & 0
 & b^{p^{e_1 + 1}} \\
 0
 &  1
 &  0 
 &  0  \\
  a^{p^{e_1}} \, c^{p^{e_1}} 
 &  b^{p^{e_1}} \, c^{p^{e_1}} 
 & 1 
 &  b^{p^{e_1}} \, d^{p^{e_1}}  \\
 c^{p^{e_1 + 1}} 
 & c^{p^{e_1}} \, d^{p^{e_1}} 
 & 0
 & d^{p^{e_1 + 1}} 
\end{array}
\right) . 
\end{align*}
\end{lem}

\begin{proof}
We have 
\begin{align*}
\renewcommand{\arraystretch}{1.5} 
\sigma^*
\left(
\begin{array}{c c}
 a & b \\
 c & d
\end{array}
\right) 
 & = \phi^-\left(\frac{c}{a} \right) \, \omega^*(a) \, \varphi^*\left(\frac{b}{a}\right) \\
 & = 
 \renewcommand{\arraystretch}{1.5} 
 \left(
\begin{array}{l l l l}
 a^{p^{e_1 + 1}} 
 & a^{p^{e_1}} \, b^{p^{e_1}} 
 & 0
 & b^{p^{e_1 + 1}} \\
 0
 &  1
 &  0 
 &  0  \\
  a^{p^{e_1}} \, c^{p^{e_1}} 
 &  b^{p^{e_1}} \, c^{p^{e_1}} 
 & 1 
 &  \frac{ b^{p^{e_1 + 1} } \, c^{p^{e_1}} + b^{p^{e_1}} }{ a^{p^{e_1}} } \\
 c^{p^{e_1 + 1}} 
 & \frac{ b^{p^{e_1}} \, c^{p^{e_1 + 1}} + c^{p^{e_1}} }{ a^{p^{e_1}} }
 & 0
 & \frac{ b^{p^{e_1 + 1}} \, c^{p^{e_1 + 1}} + 1 }{ a^{p^{e_1 + 1}} }
\end{array}
\right) . 
\end{align*}
Since $d = (b \, c + 1) / a$, 
we have 
\begin{align*}
\bigl( \text{ the $(4, 2)$-th entry } \bigr) 
 & =   c^{p^{e_1}} \, \left( \frac{ b^{p^{e_1}} \, c^{p^{e_1}} + 1 }{ a^{p^{e_1}} } \right) 
 = c^{p^{e_1}} \, d^{p^{e_1}} , \\
\bigl( \text{ the $(3, 4)$-th entry } \bigr) 
 & =  b^{p^{e_1}} \, \left( \frac{ b^{p^{e_1}} \, c^{p^{e_1}} + 1 }{ a^{p^{e_1}} } \right) 
 = b^{p^{e_1}} \, d^{p^{e_1}} , \\
\bigl( \text{ the $(4, 4)$-th entry } \bigr) 
 & =  \frac{ b^{p^{e_1 + 1}} \, c^{p^{e_1 + 1}} + 1 }{ a^{p^{e_1 + 1}} } 
 = d^{p^{e_1 + 1}} . 
\end{align*}

\end{proof}

\subsubsection{$\rm (VII)$}

\begin{lem}
Let $(\varphi^*, \omega^*)$ be of the form {\rm (VII)}. 
If there exists a homomorphsim $\sigma^* : \SL(2, k) \to \SL(4, k)$ such that 
$\sigma^* \circ \imath_{\rB(2, k)} = \psi_{\varphi^*, \, \omega^*}$. 
Then 
\begin{align*}
\renewcommand{\arraystretch}{1.5} 
\sigma^*
\left(
\begin{array}{c c}
 a & b \\
 c & d
\end{array}
\right) 
 = 
\left(
\begin{array}{l l l l}
 a^{p^{e_1 + 1}} 
 & 0
 & 0
 & b^{p^{e_1 + 1}} \\
 \frac{1}{2} \, a^{2 \, p^{e_1}} \, c^{p^{e_1}} 
 &  a^{p^{e_1}} 
 &  b^{p^{e_1}} 
 &  \frac{1}{2} \, b^{2 \, p^{e_1}} \, d^{p^{e_1}}  \\
 a^{p^{e_1}} \, c^{2 \, p^{e_1}} 
 &  c^{p^{e_1}} 
 &  d^{p^{e_1}}
 &  b^{p^{e_1}} \, d^{2 \, p^{e_1}}  \\
 c^{p^{e_1 + 1}} 
 & 0 
 & 0
 & d^{p^{e_1 + 1}} 
\end{array}
\right) . 
\end{align*}
\end{lem}

\begin{proof}
We have 
\begin{align*}
\renewcommand{\arraystretch}{1.5} 
\sigma^*
\left(
\begin{array}{c c}
 a & b \\
 c & d
\end{array}
\right) 
  & = \phi^-\left(\frac{c}{a} \right) \, \omega^*(a) \, \varphi^*\left(\frac{b}{a}\right) \\
 & = 
 \left(
\begin{array}{l l l l}
 a^{p^{e_1 + 1}} 
 & 0
 & 0
 & b^{p^{e_1 + 1}} \\
 \frac{1}{2} \, a^{2 \, p^{e_1}} \, c^{p^{e_1}} 
 &  a^{p^{e_1}} 
 &  b^{p^{e_1}} 
 &  \frac{1}{2} \, \frac{b^{3 \, p^{e_1}} \, c^{p^{e_1}} + b^{2 \, p^{e_1}} }{ a^{p^{e_1}} }  \\
 a^{p^{e_1}} \, c^{2 \, p^{e_1}} 
 &  c^{p^{e_1}} 
 &  \frac{b^{p^{e_1}} \, c^{p^{e_1}} + 1}{a^{p^{e_1}}} 
 &  \frac{b^{3 \, p^{e_1}} \, c^{2 \, p^{e_1}} + 2 \, b^{2 \, p^{e_1}} \, c^{p^{e_1}} + b^{p^{e_1}} }{ a^{2 \, p^{e_1}} }   \\
 c^{p^{e_1 + 1}} 
 & 0 
 & 0
 & \frac{b^{p^{e_1 + 1}} \, c^{p^{e_1 + 1}} + 1}{a^{p^{e_1 + 1}}} 
\end{array}
\right) . 
\end{align*}
Since $d = (b \, c + 1) / a$, 
we have 
\begin{align*}
\bigl( \text{ the $(3, 3)$-th entry } \bigr) 
 & =  d^{p^{e_1}} , \\
\bigl( \text{ the $(2, 4)$-th entry } \bigr) 
 & =  \frac{1}{2} \, b^{2 \, p^{e_1}} \left(   \frac{b^{p^{e_1}} \, c^{p^{e_1}} + 1}{a^{p^{e_1}}}   \right) 
  = \frac{1}{2} \, b^{2 \, p^{e_1}} \, d^{p^{e_1}} , 
 \\
\bigl( \text{ the $(3, 4)$-th entry } \bigr) 
 & = b^{p^{e_1}} 
 \cdot 
 \frac{b^{2 \, p^{e_1}} \, c^{2 \, p^{e_1}} + 2 \, b^{p^{e_1}} \, c^{p^{e_1}} + 1 }{ a^{2 \, p^{e_1}} }  
 = b^{p^{e_1}} \, d^{2 \, p^{e_1}} , \\ 
\bigl( \text{ the $(4, 4)$-th entry } \bigr) 
 & = d^{p^{e_1 + 1}} . 
\end{align*}

\end{proof}

\subsubsection{$\rm (IX)$}

\begin{lem}
Let $(\varphi^*, \omega^*)$ be of the form {\rm (IX)}. 
If there exists a homomorphsim $\sigma^* : \SL(2, k) \to \SL(4, k)$ such that 
$\sigma^* \circ \imath_{\rB(2, k)} = \psi_{\varphi^*, \, \omega^*}$. 
Then 
\begin{align*}
\renewcommand{\arraystretch}{1.5} 
\sigma^*
\left(
\begin{array}{c c}
 a & b \\
 c & d
\end{array}
\right) 
 & = 
 \renewcommand{\arraystretch}{1.5} 
\left(
\begin{array}{l l l l}
 a^{2 \,p^{e_1}} 
 & a^{p^{e_1}} \, b^{p^{e_1}} 
 & 0
 & \frac{1}{2} \, b^{2 \, p^{e_1}} \\
 2 \, a^{p^{e_1}} \, c^{p^{e_1}} 
 &  a^{p^{e_1}} \, d^{p^{e_1}}  +  b^{p^{e_1}} \, c^{p^{e_1}} 
 &  0 
 &  b^{p^{e_1}} \, d^{p^{e_1}}  \\
 0
 & 0 
 & 1 
 & 0 \\
 2 \, c^{2 \, p^{e_1}} 
 & 2 \, c^{p^{e_1}} \, d^{p^{e_1}} 
 & 0
 & d^{2 \, p^{e_1}} 
\end{array}
\right) . 
\end{align*}
\end{lem}

\begin{proof}
We have 
\begin{align*}
\renewcommand{\arraystretch}{1.5} 
\sigma^*
\left(
\begin{array}{c c}
 a & b \\
 c & d
\end{array}
\right) 
 & = \phi^-\left(\frac{c}{a} \right) \, \omega^*(a) \, \varphi^*\left(\frac{b}{a}\right) \\
& = 
\renewcommand{\arraystretch}{1.5} 
\left(
\begin{array}{l l l l}
 a^{2 \,p^{e_1}} 
 & a^{p^{e_1}} \, b^{p^{e_1}} 
 & 0
 & \frac{1}{2} \, b^{2 \, p^{e_1}} \\
 2 \, a^{p^{e_1}} \, c^{p^{e_1}} 
 &  2 \, b^{p^{e_1}} \, c^{p^{e_1}} + 1  
 &  0 
 &   \frac{ b^{2 \, p^{e_1}} \, c^{p^{e_1}} + b^{p^{e_1}} }{a^{p^{e_1}}}  \\
 0
 & 0 
 & 1 
 & 0 \\
 2 \, c^{2 \, p^{e_1}} 
 & \frac{ 2 \, b^{p^{e_1}} \, c^{2 \, p^{e_1}} + 2 \, c^{p^{e_1}} }{a^{p^{e_1}}}
 & 0
 & \frac{b^{2 \, p^{e_1}} \, c^{2 \, p^{e_1}} + 2 \, b^{p^{e_1}} \, c^{p^{e_1}} + 1 }{ a^{2 \, p^{e_1}} } 
\end{array}
\right) . 
\end{align*}
Since $d = (b \, c + 1) / a$, 
we have 
\begin{align*}
\bigl( \text{ the $(2, 2)$-th entry } \bigr) 
 & =  a^{p^{e_1}} \, d^{p^{e_1}}  +  b^{p^{e_1}} \, c^{p^{e_1}}  , \\
\bigl( \text{ the $(4, 2)$-th entry } \bigr) 
 & =  2 \, c^{p^{e_1}} \, d^{p^{e_1}} , \\
\bigl( \text{ the $(2, 4)$-th entry } \bigr) 
 & =  b^{p^{e_1}} \, d^{p^{e_1}} , \\
\bigl( \text{ the $(4, 4)$-th entry } \bigr) 
 & = d^{2 \, p^{e_1}} . 
\end{align*}

\end{proof}

\subsubsection{$\rm (XI)$}

\begin{lem}
Let $(\varphi^*, \omega^*)$ be of the form {\rm (XI)}. 
If there exists a homomorphsim $\sigma^* : \SL(2, k) \to \SL(4, k)$ such that 
$\sigma^* \circ \imath_{\rB(2, k)} = \psi_{\varphi^*, \, \omega^*}$. 
Then $p = 2$ and 
\begin{align*}
\renewcommand{\arraystretch}{1.5} 
\sigma^*
\left(
\begin{array}{c c}
 a & b \\
 c & d
\end{array}
\right) 
 = 
\left(
\begin{array}{l l l l}
 a^{2 \, p^{e_1}} 
 & a^{p^{e_1}} \, b^{p^{e_1}} 
 & 0
 & b^{2 \, p^{e_1}} \\
 0 
 & 1 
 & 0 
 & 0  \\
 0
 & 0 
 & 1 
 & 0 \\
c^{2 \, p^{e_1}} 
 & c^{p^{e_1}} \, d^{p^{e_1}} 
 & 0
 & d^{2 \, p^{e_1}} 
\end{array}
\right) . 
\end{align*}
\end{lem}

\begin{proof}
We have 
\begin{align*}
\renewcommand{\arraystretch}{1.5} 
\sigma^*
\left(
\begin{array}{c c}
 a & b \\
 c & d
\end{array}
\right) 
 & = \phi^-\left(\frac{c}{a} \right) \, \omega^*(a) \, \varphi^*\left(\frac{b}{a}\right) \\
 & = 
\left(
\begin{array}{l l l l}
 a^{2 \, p^{e_1}} 
 & a^{p^{e_1}} \, b^{p^{e_1}} 
 & 0
 & b^{2 \, p^{e_1}} \\
 0 
 & 1 
 & 0 
 & 0  \\
 0
 & 0 
 & 1 
 & 0 \\
c^{2 \, p^{e_1}} 
 & \frac{ b^{p^{e_1}} \, c^{2 \, p^{e_1}} + c^{p^{e_1}} }{a^{p^{e_1}}}
 & 0
 & \frac{ (b \, c)^{2 \, p^{e_1}} + 1 }{a^{2 \, p^{e_1}}}
\end{array}
\right) . 
\end{align*}
Since $d = (b \, c + 1) / a$, 
we have 
\begin{align*}
\bigl( \text{ the $(4, 2)$-th entry } \bigr) 
 & =  c^{p^{e_1}} \, d^{p^{e_1}} , \\
\bigl( \text{ the $(4, 4)$-th entry } \bigr) 
 & = d^{2 \, p^{e_1}} . 
\end{align*}

\end{proof}

\subsubsection{$\rm (XV)$}

\begin{lem}
Let $(\varphi^*, \omega^*)$ be of the form {\rm (XV)}. 
If there exists a homomorphsim $\sigma^* : \SL(2, k) \to \SL(4, k)$ such that 
$\sigma^* \circ \imath_{\rB(2, k)} = \psi_{\varphi^*, \, \omega^*}$. 
Then 
\begin{align*}
\renewcommand{\arraystretch}{1.5} 
\sigma^*
\left(
\begin{array}{c c}
 a & b \\
 c & d
\end{array}
\right) 
 = 
\left(
\begin{array}{l l l l}
 a^{p^{e_2}} 
 & 0 
 & 0
 & b^{p^{e_2}} \\
 0 
 & a^{p^{e_3}}
 & b^{p^{e_3}} 
 & 0  \\
 0
 & c^{p^{e_3}}
 & d^{p^{e_3}} 
 & 0 \\
c^{p^{e_2}} 
 & 0 
 & 0
 & d^{p^{e_2}} 
\end{array}
\right) . 
\end{align*}
\end{lem}

\begin{proof}
We have 
\begin{align*}
\renewcommand{\arraystretch}{1.5} 
\sigma^*
\left(
\begin{array}{c c}
 a & b \\
 c & d
\end{array}
\right)
 & = \phi^-\left(\frac{c}{a} \right) \, \omega^*(a) \, \varphi^*\left(\frac{b}{a}\right) \\  
 & = 
\left(
\begin{array}{l l l l}
 a^{p^{e_2}} 
 & 0 
 & 0
 & b^{p^{e_2}} \\
 0 
 & a^{p^{e_3}}
 & b^{p^{e_3}} 
 & 0  \\
 0
 & c^{p^{e_3}}
 & \frac{ b^{p^{e_3}} \, c^{p^{e_3}} + 1 }{a^{p^{e_3}}} 
 & 0 \\
c^{p^{e_2}} 
 & 0 
 & 0
 &  \frac{ b^{p^{e_2}} \, c^{p^{e_2}} + 1 }{a^{p^{e_2}}} 
\end{array}
\right) . 
\end{align*}
Since $d = (b \, c + 1) / a$, 
we have 
\begin{align*}
\bigl( \text{ the $(3, 3)$-th entry } \bigr) 
 & =  d^{p^{e_3}} , \\
\bigl( \text{ the $(4, 4)$-th entry } \bigr) 
 & = d^{p^{e_2}} . 
\end{align*}

\end{proof}

\subsubsection{$\rm (XIX)$}

\begin{lem}
Let $(\varphi^*, \omega^*)$ be of the form {\rm (XIX)}. 
If there exists a homomorphsim $\sigma^* : \SL(2, k) \to \SL(4, k)$ such that 
$\sigma^* \circ \imath_{\rB(2, k)} = \psi_{\varphi^*, \, \omega^*}$. 
Then $p = 2$ and 
\begin{align*}
\renewcommand{\arraystretch}{1.5} 
\sigma^*
\left(
\begin{array}{c c}
 a & b \\
 c & d
\end{array}
\right) 
 = 
\left(
\begin{array}{l l l l}
 a^{2 \, p^{e_1}} 
 & 0 
 & 0
 & b^{2 \, p^{e_1}} \\
 0 
 & 1 
 & 0  
 & 0  \\
 a^{p^{e_1}} \, c^{p^{e_1}} 
 & 0 
 & 1 
 & b^{p^{e_1}} \, d^{p^{e_1}}  \\
c^{2 \, p^{e_1}} 
 & 0 
 & 0
 & d^{2 \, p^{e_1}} 
\end{array}
\right) . 
\end{align*}
\end{lem}

\begin{proof}
We have 
\begin{align*}
\renewcommand{\arraystretch}{1.5} 
\sigma^*
\left(
\begin{array}{c c}
 a & b \\
 c & d
\end{array}
\right) 
 & = \phi^-\left(\frac{c}{a} \right) \, \omega^*(a) \, \varphi^*\left(\frac{b}{a}\right) \\ 
 & = 
\left(
\begin{array}{l l l l}
 a^{2 \, p^{e_1}} 
 & 0 
 & 0
 & b^{2 \, p^{e_1}} \\
 0 
 & 1 
 & 0  
 & 0  \\
 a^{p^{e_1}} \, c^{p^{e_1}} 
 & 0 
 & 1 
 &  \frac{ b^{p^{e_1 + 1}} \, c^{p^{e_1}} + b^{p^{e_1}} }{a^{p^{e_1}}}  \\
c^{2 \, p^{e_1}} 
 & 0 
 & 0
 &  \frac{ b^{p^{e_1 + 1}} \, c^{p^{e_1 + 1}} + 1 }{a^{p^{e_1 + 1}}} 
\end{array}
\right) . 
\end{align*}
Since $d = (b \, c + 1) / a$, 
we have 
\begin{align*}
\bigl( \text{ the $(3, 4)$-th entry } \bigr) 
 & = b^{p^{e_1}} \, d^{p^{e_1}} , \\
\bigl( \text{ the $(4, 4)$-th entry } \bigr) 
 & = d^{2 \, p^{e_1}} . 
\end{align*}

\end{proof}

\subsubsection{$\rm (XXI)$}

\begin{lem}
Let $(\varphi^*, \omega^*)$ be of the form {\rm (XXI)}. 
If there exists a homomorphsim $\sigma^* : \SL(2, k) \to \SL(4, k)$ such that 
$\sigma^* \circ \imath_{\rB(2, k)} = \psi_{\varphi^*, \, \omega^*}$. 
Then 
\begin{align*}
\renewcommand{\arraystretch}{1.5} 
\sigma^*
\left(
\begin{array}{c c}
 a & b \\
 c & d
\end{array}
\right) 
 = 
\left(
\begin{array}{l l l l}
 a^{2 \, p^{e_1}} 
 & 0 
 & a^{p^{e_1}} \, b^{p^{e_1}}
 & b^{2 \, p^{e_1}} \\
 a^{p^{e_1}} \, c^{p^{e_1}} 
 &  1
 & b^{p^{e_1}} \, c^{p^{e_1}}
 & b^{p^{e_1}} \, d^{p^{e_1}} \\
 0
 & 0 
 & 1
 & 0 \\
 c^{2 \, p^{e_1}} 
 & 0 
 & c^{p^{e_1}} \, d^{p^{e_1}} 
 & d^{2 \, p^{e_1}} 
\end{array}
\right) . 
\end{align*}
\end{lem}

\begin{proof}
We have 
\begin{align*}
\renewcommand{\arraystretch}{1.5} 
\sigma^*
\left(
\begin{array}{c c}
 a & b \\
 c & d
\end{array}
\right) 
 & = \phi^-\left(\frac{c}{a} \right) \, \omega^*(a) \, \varphi^*\left(\frac{b}{a}\right) \\
 & = 
\left(
\begin{array}{l l l l}
 a^{2 \, p^{e_1}} 
 & 0 
 & a^{p^{e_1}} \, b^{p^{e_1}}
 & b^{2 \, p^{e_1}} \\
 a^{p^{e_1}} \, c^{p^{e_1}} 
 &  1
 & b^{p^{e_1}} \, c^{p^{e_1}}
 & b^{p^{e_1}} \cdot \frac{b^{p^{e_1}} \, c^{p^{e_1}} + 1}{a^{p^{e_1}}}   \\
 0
 & 0 
 & 1
 & 0 \\
 c^{2 \, p^{e_1}} 
 & 0 
 & c^{p^{e_1}} \cdot  \frac{b^{p^{e_1}} \, c^{p^{e_1}} + 1}{a^{p^{e_1}}}   
 &  \frac{b^{2 \, p^{e_1}} \, c^{2 \, p^{e_1}} + 1}{a^{2 \, p^{e_1}}}  
\end{array}
\right) . 
\end{align*}
Since $d = (b \, c + 1) / a$, 
we have 
\begin{align*}
\bigl( \text{ the $(4, 3)$-th entry } \bigr) 
 & =  c^{p^{e_1}} \, d^{p^{e_1}} , \\
\bigl( \text{ the $(2, 4)$-th entry } \bigr) 
 & =  b^{p^{e_1}} \, d^{p^{e_1}} , \\
\bigl( \text{ the $(4, 4)$-th entry } \bigr) 
 & = d^{2 \, p^{e_1}} .  
\end{align*}

\end{proof}

\subsubsection{$\rm (XXII)$}

\begin{lem}
Let $(\varphi^*, \omega^*)$ be of the form {\rm (XXII)}. 
If there exists a homomorphsim $\sigma^* : \SL(2, k) \to \SL(4, k)$ such that 
$\sigma^* \circ \imath_{\rB(2, k)} = \psi_{\varphi^*, \, \omega^*}$. 
Then $d_1 = d_2 = p^{e_1}$ and 
\begin{align*}
\renewcommand{\arraystretch}{1.5} 
\sigma^*
\left(
\begin{array}{c c}
 a & b \\
 c & d
\end{array}
\right) 
 = 
\left(
\begin{array}{l l l l}
 a^{p^{e_1}} 
 & 0 
 & b^{p^{e_1}}
 & 0 \\
 0 
 &  a^{p^{e_1}} 
 & 0 
 & b^{p^{e_1}} \\
 c^{p^{e_1}} 
 & 0 
 & d^{p^{e_1}} 
 & 0 \\
0
 & c^{p^{e_1}} 
 & 0 
 & d^{p^{e_1}} 
\end{array}
\right) . 
\end{align*}
\end{lem}

\begin{proof}
We have 
\begin{align*}
\renewcommand{\arraystretch}{1.5} 
\sigma^*
\left(
\begin{array}{c c}
 a & b \\
 c & d
\end{array}
\right) 
 & = \phi^-\left(\frac{c}{a} \right) \, \omega^*(a) \, \varphi^*\left(\frac{b}{a}\right) \\
 & = 
\left(
\begin{array}{l l l l}
 a^{p^{e_1}} 
 & 0 
 & b^{p^{e_1}}
 & 0 \\
 0 
 &  a^{p^{e_1}} 
 & 0 
 & b^{p^{e_1}} \\
 c^{p^{e_1}} 
 & 0 
 & \frac{b^{p^{e_1}} \, c^{p^{e_1}} + 1}{a^{p^{e_1}}}   
 & 0 \\
0
 & c^{p^{e_1}} 
 & 0 
 & \frac{b^{p^{e_1}} \, c^{p^{e_1}} + 1}{a^{p^{e_1}}}   
\end{array}
\right) . 
\end{align*}
Since $d = (b \, c + 1) / a$, 
we have 
\begin{align*}
\bigl( \text{ the $(3, 3)$-th entry } \bigr) 
 & = d^{p^{e_1}} ,  \\
\bigl( \text{ the $(4, 4)$-th entry } \bigr) 
 & = d^{p^{e_1}} . 
\end{align*}

\end{proof}

\subsubsection{$\rm (XXIV)$}

\begin{lem}
Let $(\varphi^*, \omega^*)$ be of the form {\rm (XXIV)}. 
If there exists a homomorphsim $\sigma^* : \SL(2, k) \to \SL(4, k)$ such that 
$\sigma^* \circ \imath_{\rB(2, k)} = \psi_{\varphi^*, \, \omega^*}$. 
Then $d_2 = 0$ and 
\begin{align*}
\renewcommand{\arraystretch}{1.5} 
\sigma^*
\left(
\begin{array}{c c}
 a & b \\
 c & d
\end{array}
\right) 
 = 
\left(
\begin{array}{l l l l}
 a^{p^{e_2}} 
 & 0 
 & 0
 & b^{p^{e_2}} \\
 0 
 & 1 
 & 0
 & 0  \\
 0
 & 0 
 & 1
 & 0 \\
c^{p^{e_2}} 
 & 0 
 & 0
 & d^{p^{e_2}} 
\end{array}
\right) . 
\end{align*}
\end{lem}

\begin{proof}
We have 
\begin{align*}
\renewcommand{\arraystretch}{1.5} 
\sigma^*
\left(
\begin{array}{c c}
 a & b \\
 c & d
\end{array}
\right) 
 & = \phi^-\left(\frac{c}{a} \right) \, \omega^*(a) \, \varphi^*\left(\frac{b}{a}\right) \\
 & = 
\left(
\begin{array}{l l l l}
 a^{p^{e_2}} 
 & 0 
 & 0
 & b^{p^{e_2}} \\
 0 
 & 1 
 & 0
 & 0  \\
 0
 & 0 
 & 1
 & 0 \\
c^{p^{e_2}} 
 & 0 
 & 0
 & \frac{b^{p^{e_2}} \, c^{p^{e_2}} + 1}{a^{p^{e_2}}}  
\end{array}
\right) . 
\end{align*}
Since $d = (b \, c + 1) / a$, 
we have 
\begin{align*}
\bigl( \text{ the $(4, 4)$-th entry } \bigr) 
 & = d^{p^{e_2}} . 
\end{align*}

\end{proof}

\subsubsection{$\rm (XXVI)$}

\begin{align*}
\renewcommand{\arraystretch}{1.5} 
\sigma^*
\left(
\begin{array}{c c}
 a & b \\
 c & d
\end{array}
\right) 
 = 
I_4 . 
\end{align*}

\section{An overlapping classification of homomorphisms from $\SL(2, k)$ to $\SL(4, k)$}

Let $V(n)$ denote the $n$-dimensional column vector space over $k$ and 
let $W(n)$ denote the $n$-dimensional row vector space over $k$. 
For any homomorphism $\sigma : \SL(2, k) \to \SL(n, k)$, we let $\SL(2, k)$ act linearly on $V(n)$ from the left 
and $\SL(2, k)$ act linearly on $W(n)$ from the right. 
We denote by $V(n)^\sigma$ the subspace consisting of all $\sigma$-fixed column vectors and 
by $W(n)^\sigma$ the subspace consisting of all $\sigma$-fixed row vectors, i.e., 
\begin{align*}
 V(n)^\sigma & := \{\, v \in V(n) \mid \sigma(A) \, v = v \; \text{ for all $A \in \SL(2, k)$}  \, \} 
\intertext{and}
 W(n)^\sigma & := \{ \, w \in W(n) \mid w \, \sigma(A)  = w \; \text{ for all $A \in \SL(2, k)$}  \, \} . 
\end{align*}
Let 
\[
d(\sigma) 
:= 
\bigl(\,
 \dim_k V(n)^\sigma, \; \dim_k W(n)^\sigma 
\,\bigr) . 
\]

\begin{lem}
Two homomorphisms $\sigma : \SL(2, k) \to \SL(n, k)$ and 
$\sigma^* : \SL(2, k) \to \SL(n, k)$ are equivalent. Then 
\[
 d(\sigma) = d(\sigma^*) . 
\]
\end{lem}

\begin{proof}
The proof is straightforward. 
\end{proof}

\subsection{Homomorphisms $\sigma^* : \SL(2, k) \to \SL(4, k)$ 
and $\sigma^+ : \SL(2, k) \to \SL(4, k)$}

\subsubsection{${\rm (I)}^*$}

Assume $p \geq 5$. 
Let $e_1$ be an integer such that 
\[ 
e_1 \geq 0 . 
\]
We can define a morphism $\sigma^* : \SL(2, k) \to \SL(4, k)$ as 
\begin{align*}
\renewcommand{\arraystretch}{1.5} 
\sigma^*
\left(
\begin{array}{c c}
 a & b \\
 c & d
\end{array}
\right) 
 := 
\left(
\begin{array}{l l l l}
 a^{3 \, p^{e_1}} 
 & a^{2 \, p^{e_1}} \, b^{p^{e_1}} 
 & \frac{1}{2} \, a^{p^{e_1}} \, b^{2 \, p^{e_1}} 
 & \frac{1}{6} \, b^{3 \, p^{e_1}} \\
 3 \, a^{2 \, p^{e_1}} \, c^{p^{e_1}} 
 & a^{p^{e_1}} \cdot (a \, d + 2 \, b \, c)^{p^{e_1}} 
 & b^{p^{e_1}} \cdot (a \, d + \frac{1}{2} \, b \, c)^{p^{e_1}} 
 & \frac{1}{2} \, b^{2 \, p^{e_1}} \, d^{p^{e_1}} \\
 6 \, a^{p^{e_1}} \, c^{2 \, p^{e_1}} 
 &  4 \, c^{p^{e_1}} \cdot (a \, d + \frac{1}{2} \, b \, c)^{p^{e_1}} 
 & d^{p^{e_1}} \cdot (a \, d + 2 \, b \, c)^{p^{e_1}} 
 & b^{p^{e_1}} \, d^{2 \, p^{e_1}} \\
 6 \, c^{3 \, p^{e_1}} 
 & 6 \, c^{2 \, p^{e_1}} \, d^{p^{e_1}} 
 & 3 \, c^{p^{e_1}} \, d^{2 \, p^{e_1}} 
 & d^{3 \, p^{e_1}} 
\end{array}
\right) . 
\end{align*}
We can define a morphism $\sigma^+ : \SL(2, k) \to \SL(4, k)$ as 
\begin{align*}
\renewcommand{\arraystretch}{1.5} 
\sigma^+ 
\left(
\begin{array}{c c}
 A & B \\
 C & D
\end{array}
\right) 
 := 
\left(
\begin{array}{l l l l}
 A^3 
 & A^2 \, B 
 & \frac{1}{2} \, A \, B^2 
 & \frac{1}{6} \, B^3 \\
3 \, A^2 \, C 
 & A \, (A \, D + 2 \, B \, C) 
 & B \, (A \, D + \frac{1}{2} \, B \, C) 
 & \frac{1}{2} \, B^2 \, D \\
6 \, A \, C^2 
 & 4 \, C \,  (A \, D + \frac{1}{2} \, B \, C) 
 & D \, (A \, D + 2 \, B \, C) 
 & B \, D^2 \\
6 \, C^3 
 & 6 \, C^2 \, D 
 & 3 \, C \, D^2 
 & D^3  
\end{array}
\right) . 
\end{align*}

\begin{lem}
Let $\sigma^* : \SL(2, k) \to \SL(4, k)$ and  
$\sigma^+ : \SL(2, k) \to \SL(4, k)$ be as above. 
Then the following assertions {\rm (1)}, {\rm (2)}, {\rm (3)} hold true: 
\begin{enumerate}[label = {\rm (\arabic*)}]
\item $\sigma^* =\sigma^+ \circ  F^{e_1}$. 

\item $\sigma^+$ is a homomorphism. 

\item $\sigma^*$ is a homomorphism. 
\end{enumerate} 
\end{lem}

\begin{proof}
The proof is straightforward. 
\end{proof}

\begin{lem}
The following assertions {\rm (1)}, {\rm (2)}, {\rm (3)}, {\rm (4)} hold true: 
\begin{enumerate}[label = {\rm (\arabic*)}]
\item 
$
 V(4)^{\sigma^+} = 0
 $. 

\item 
$
 W(4)^{\sigma^+} = 0 
 $. 
 
\item $d(\sigma^+) = (0, 0)$. 

\item $d(\sigma^*) = (0, 0)$. 
\end{enumerate} 
\end{lem}

\begin{proof}
Use the fact that the homomorphism $\sigma^+ : \SL(2, k) \to \SL(4, k)$ is irreducible.

\end{proof}

\subsubsection{${\rm (II)}^*$}

Assume $p = 3$. 
Let $e_1$ be an integer such that 
\[ 
e_1 \geq 0 . 
\]
We can define a morphism $\sigma^* : \SL(2, k) \to \SL(4, k)$ as 
\begin{align*}
\renewcommand{\arraystretch}{1.5} 
\sigma^*
\left(
\begin{array}{c c}
 a & b \\
 c & d
\end{array}
\right) 
 := 
\left(
\begin{array}{l l l l}
 a^{p^{e_1 + 1}} 
 & a^{2 \, p^{e_1}} \, b^{p^{e_1}} 
 & \frac{1}{2} \, a^{p^{e_1}} \, b^{2 \, p^{e_1}} 
 & b^{p^{e_1 + 1}} \\
 0
 &  a^{p^{e_1}} 
 &  b^{p^{e_1}}
 &  0  \\
 0
 & c^{p^{e_1}} 
 & d^{p^{e_1}} 
 & 0 \\
 c^{p^{e_1 + 1}} 
 & c^{2 \, p^{e_1}} \, d^{p^{e_1}} 
 & \frac{1}{2} \, c^{p^{e_1}} \, d^{2 \, p^{e_1}} 
 & d^{p^{e_1 + 1}} 
\end{array}
\right) . 
\end{align*}
We can define a morphism $\sigma^+ : \SL(2, k) \to \SL(4, k)$ as 
\begin{align*}
\renewcommand{\arraystretch}{1.5} 
\sigma^+ 
\left(
\begin{array}{c c}
 A & B \\
 C & D
\end{array}
\right) 
 := 
\left(
\begin{array}{l l | l l}
 A^3 
 & B^3 
 & A^2 \, B 
 & \frac{1}{2} \, A \, B^2 \\
C^3 
 & D^3  
 & C^2 \, D
 & \frac{1}{2} \, C \, D^2  \\ 
\hline 
0 
 & 0  
 & A 
 & B \\
0 
 & 0 
 & C 
 & D 
\end{array}
\right) . 
\end{align*}

\begin{lem}
Let $\sigma^* : \SL(2, k) \to \SL(4, k)$ and  
$\sigma^+ : \SL(2, k) \to \SL(4, k)$ be as above. 
Let $P := P_{3,\, 4} \, P_{2,\, 3} \in \GL(4, k)$. 
Then the following assertions {\rm (1)}, {\rm (2)}, {\rm (3)} hold true: 
\begin{enumerate}[label = {\rm (\arabic*)}]
\item $\Inn_P \circ \sigma^* =\sigma^+ \circ  F^{e_1}$. 

\item $\sigma^+$ is a homomorphism. 

\item $\sigma^*$ is a homomorphism. 
\end{enumerate} 
\end{lem}

\begin{proof}
The proof is straightforward. 
\end{proof}

\begin{lem}
The following assertions {\rm (1)}, {\rm (2)}, {\rm (3)}, {\rm (4)} hold true: 
\begin{enumerate}[label = {\rm (\arabic*)}]
\item 
$
 V(4)^{\sigma^+} = 0 
 $. 

\item 
$
 W(4)^{\sigma^+} = 0 
 $. 
 
\item $d(\sigma^+) = (0, 0)$. 

\item $d(\sigma^*) = (0, 0)$. 
\end{enumerate} 
\end{lem}

\begin{proof}
(1) 
Let 
\[
\left(
\begin{array}{c}
\bm{v}_1 \\
\bm{v}_2
\end{array}
\right)
 \in 
  V^{\sigma^+} , 
\qquad 
\bm{v}_1 \in k^{\oplus 2}, \quad \bm{v}_2 \in k^{\oplus 2} . 
\]
For any 
\[
\renewcommand{\arraystretch}{1.5} 
\left(
\begin{array}{c c}
 A & B \\
 C & D 
\end{array}
\right)
\in \SL(2, k) , 
\]
we have 
\begin{align*}
\renewcommand{\arraystretch}{1.5} 
\left(
\begin{array}{c c}
 A^3 & B^3 \\
 C^3 & D^3
\end{array}
\right) 
\bm{v}_1
+ 
\left(
\begin{array}{c c}
 A^2 \, B  & \frac{1}{2} \, A \, B^2 \\
 C^2 \, D & \frac{1}{2} \, C \, D^2
\end{array}
\right) 
\bm{v}_2 
& = \bm{0}, \\
\renewcommand{\arraystretch}{1.5} 
\left(
\begin{array}{c c}
 A & B \\
 C & D
\end{array}
\right) 
\bm{v}_2
& = 
\bm{0} . 
\end{align*}
From the latter equality, we have $\bm{v}_2 = \bm{0}$, and then 
using the former equality, we have $\bm{v}_1 = \bm{0}$. 

(2) The proof is similar to the proof of the above assertion (1). 

(3) The proof is straightforward. 

(4) Use the above assertion (3). 
\end{proof}

\subsubsection{${\rm (IV)}^*$}

Let $e_1$ and $e_2$ be integers such that 
\[ 
 e_2 > e_1 \geq 0 . 
\]
We can define a morphism $\sigma^* : \SL(2, k) \to \SL(4, k)$ as 
\begin{align*}
\renewcommand{\arraystretch}{1.5} 
\sigma^*
\left(
\begin{array}{c c}
 a & b \\
 c & d
\end{array}
\right) 
& := 
\renewcommand{\arraystretch}{1.5} 
\left(
\begin{array}{l l l l}
 a^{p^{e_2}} \cdot a^{p^{e_1}} 
 & a^{p^{e_2}} \cdot b^{p^{e_1}} 
 & b^{p^{e_2}} \cdot a^{p^{e_1}} 
 & b^{p^{e_2}} \cdot b^{p^{e_1}} \\
  a^{p^{e_2}} \cdot c^{p^{e_1}} 
 & a^{p^{e_2}} \cdot d^{p^{e_1}} 
 & b^{p^{e_2}} \cdot c^{p^{e_1}} 
 & b^{p^{e_2}} \cdot d^{p^{e_1}} \\
 c^{p^{e_2}} \cdot a^{p^{e_1}} 
 & c^{p^{e_2}} \cdot b^{p^{e_1}} 
 & d^{p^{e_2}} \cdot a^{p^{e_1}} 
 & d^{p^{e_2}} \cdot b^{p^{e_1}} \\
 c^{p^{e_2}} \cdot c^{p^{e_1}} 
 & c^{p^{e_2}} \cdot d^{p^{e_1}} 
 & d^{p^{e_2}} \cdot c^{p^{e_1}} 
 & d^{p^{e_2}} \cdot d^{p^{e_1}}
\end{array}
\right) \\
& = 
\renewcommand{\arraystretch}{1.5} 
\left(
\begin{array}{c c}
a^{p^{e_2}}
\left(
\begin{array}{c c}
 a^{p^{e_1}} & b^{p^{e_1}} \\
 c^{p^{e_1}} & d^{p^{e_1}} 
\end{array}
\right)
& 
b^{p^{e_2}}
\left(
\begin{array}{c c}
 a^{p^{e_1}} & b^{p^{e_1}} \\
 c^{p^{e_1}} & d^{p^{e_1}} 
\end{array}
\right) \\ [1.5em] 
c^{p^{e_2}}
\left(
\begin{array}{c c}
 a^{p^{e_1}} & b^{p^{e_1}} \\
 c^{p^{e_1}} & d^{p^{e_1}} 
\end{array}
\right)
& 
d^{p^{e_2}}
\left(
\begin{array}{c c}
 a^{p^{e_1}} & b^{p^{e_1}} \\
 c^{p^{e_1}} & d^{p^{e_1}} 
\end{array}
\right)
\end{array}
\right) . 
\end{align*}

%

Let $\Delta : \SL(2, k) \to \SL(2, k) \times \SL(2, k)$ be the homomorphism 
defined by 
\[
 \Delta(X) := (X, X) . 
\]

Let $F^{e_2} \times F^{e_1} : \SL(2, k) \times \SL(2, k) \to \SL(2, k) \times \SL(2, k)$ 
be the homomorphism defined by 
\begin{align*}
&
\renewcommand{\arraystretch}{1.5} 
 ( F^{e_2} \times F^{e_1}) 
\left( \, 
\left(
\begin{array}{c c}
 A & B \\
 C & D
\end{array}
\right) , 
\; 
\left(
\begin{array}{l l }
 A' & B' \\
 C' & D'
\end{array}
\right) 
\, \right) \\
& 
\renewcommand{\arraystretch}{1.5} 
\qquad  := 
\left( \,  
\left(
\begin{array}{l l }
 A^{p^{e_2}} & B^{p^{e_2}} \\
 C^{p^{e_2}} & D^{p^{e_2}}
\end{array}
\right)  , 
\; 
\left(
\begin{array}{l l }
 {A'}^{p^{e_1}} & {B'}^{p^{e_1}} \\
 {C'}^{p^{e_1}} & {D'}^{p^{e_1}}
\end{array}
\right)  
\, \right) . 
\end{align*}

Let $\vartheta : \SL(2, k) \times \SL(2, k) \to \SL(4, k)$ be the homomorphism 
defined by 
\[
 \vartheta(X_1, X_2) := X_1 \otimes X_2 . 
\]

\begin{lem}
Let $\sigma^* : \SL(2, k) \to \SL(4, k)$ be as above. 
Then the following assertions {\rm (1)} and {\rm (2)} hold true: 
\begin{enumerate}[label = {\rm (\arabic*)}]
\item $\sigma^* = \vartheta \circ (F^{e_2} \times F^{e_1}) \circ \Delta$, i.e., 
\[
\xymatrix@R=36pt@C=36pt@M=6pt{
 \SL(2, k) \ar[r]^(.38)\Delta \ar@/_24pt/[rrr]_{\sigma^*}
 &  \SL(2, k) \times \SL(2, k) \ar[r]^{F^{e_2} \times F^{e_1}} 
 &   \SL(2, k) \times \SL(2, k) \ar[r]^(.62)\vartheta 
 & \SL(4, k) 
}
\]

\item $\sigma^*$ is a homomorphism. 
\end{enumerate} 
\end{lem}

\begin{proof}
The proof is straightforward. 
\end{proof}

We can define a homomorphism $\sigma^+ : \SL(2, k) \to \SL(4, k)$ as 
\[
\renewcommand{\arraystretch}{1.5} 
\sigma^+
\left(
\begin{array}{c c}
 A & B \\
 C & D
\end{array}
\right) 
 := 
\left(
\begin{array}{l l }
 A & B \\
 C & D
\end{array}
\right) 
\otimes 
\left(
\begin{array}{l l }
 A & B \\
 C & D
\end{array}
\right)  . 
\]
Clearly, $\sigma^+ = \theta \circ \Delta$.

\begin{lem}
The following assertions {\rm (1)}, {\rm (2)}, {\rm (3)}, {\rm (4)} hold true: 
\begin{enumerate}[label = {\rm (\arabic*)}]
\item 
$ V(4)^{\sigma^+} = 0 $. 

\item 
$
 W(4)^{\sigma^+} = 0 
 $. 
 
\item $d(\sigma^+) = (0, 0)$. 

\item $d(\sigma^*) = (0, 0)$. 
\end{enumerate} 
\end{lem}

\begin{proof}
(1) Let 
\[
\left(
\begin{array}{c}
\bm{v}_1 \\
\bm{v}_2
\end{array}
\right)
 \in 
 V^{\sigma^+} , 
\qquad 
\bm{v}_1 \in k^{\oplus 2}, \quad \bm{v}_2 \in k^{\oplus 2} . 
\]
For any 
\[
\renewcommand{\arraystretch}{1.5} 
\left(
\begin{array}{c c}
 A & B \\
 C & D 
\end{array}
\right)
\in \SL(2, k) , 
\]
we have 
\begin{align*}
\renewcommand{\arraystretch}{1.5} 
A 
\left(
\begin{array}{c c}
 A & B \\
 C & D
\end{array}
\right) 
\bm{v}_1
+ 
B
\left(
\begin{array}{c c}
 A & B \\
 C & D
\end{array}
\right) 
\bm{v}_2 
& = \bm{0}, \\
\renewcommand{\arraystretch}{1.5} 
C 
\left(
\begin{array}{c c}
 A & B \\
 C & D
\end{array}
\right) 
\bm{v}_1 
+ 
D 
\left(
\begin{array}{c c}
 A & B \\
 C & D
\end{array}
\right) 
\bm{v}_2 
& = 
\bm{0} . 
\end{align*}
We can deform these equalities as 
\begin{align*}
\renewcommand{\arraystretch}{1.5} 
\left(
\begin{array}{c c}
 A & B \\
 C & D
\end{array}
\right) 
(\, A \, \bm{v}_1 + B \, \bm{v}_2 \, ) 
& = \bm{0}, \\
\renewcommand{\arraystretch}{1.5} 
\left(
\begin{array}{c c}
 A & B \\
 C & D
\end{array}
\right) 
(\, C \, \bm{v}_1 + D \, \bm{v}_2 \, )  
& = 
\bm{0} . 
\end{align*}
So, 
\begin{align*}
\renewcommand{\arraystretch}{1.5} 
\left(
\begin{array}{c c}
 \bm{v}_1 & \bm{v}_2
\end{array}
\right)
\left(
\begin{array}{c c}
 A & C \\
 B & D
\end{array}
\right)
= 
\left(
\begin{array}{c c}
 \bm{0} & \bm{0}
\end{array}
\right) . 
\end{align*}
Thus $\bm{v}_1 = \bm{v}_2 = \bm{0}$.

(2) The proof is similar to the proof of the above assertion (1).

(3) The proof is straightforward. 

(4) Use the above assertion (3). 
\end{proof}

\subsubsection{${\rm (V)}^*$}

Assume $p = 2$. 
Let $e_1$ be an integer such that 
\[ 
e_1 \geq 0 . 
\]
We can define a morphism $\sigma^* : \SL(2, k) \to \SL(4, k)$ as 
\begin{align*}
\renewcommand{\arraystretch}{1.5} 
\sigma^*
\left(
\begin{array}{c c}
 a & b \\
 c & d
\end{array}
\right) 
 := 
\left(
\begin{array}{l l l l}
 a^{p^{e_1 + 1}} 
 & a^{p^{e_1}} \, b^{p^{e_1}} 
 & 0
 & b^{p^{e_1 + 1}} \\
 0
 &  1
 &  0 
 &  0  \\
  a^{p^{e_1}} \, c^{p^{e_1}} 
 &  b^{p^{e_1}} \, c^{p^{e_1}} 
 & 1 
 &  b^{p^{e_1}} \, d^{p^{e_1}}  \\
 c^{p^{e_1 + 1}} 
 & c^{p^{e_1}} \, d^{p^{e_1}} 
 & 0
 & d^{p^{e_1 + 1}} 
\end{array}
\right) . 
\end{align*}
We can define a morphism $\sigma^+ : \SL(2, k) \to \SL(4, k)$ as 
\begin{align*}
\renewcommand{\arraystretch}{1.5} 
\sigma^+ 
\left(
\begin{array}{c c}
 A & B \\
 C & D
\end{array}
\right) 
 := 
\left(
\begin{array}{l l l l}
1 
 & 0 
 & 0 
 & 0 \\
A \, B  
 & A^2   
 & B^2 
 & 0  \\ 
C \, D 
 & C^2 
 & D^2  
 & 0  \\
B \, C 
 & A \, C
 & B \, D
 & 1 
\end{array}
\right) . 
\end{align*}

\begin{lem}
Let $\sigma^* : \SL(2, k) \to \SL(4, k)$ and  
$\sigma^+ : \SL(2, k) \to \SL(4, k)$ be as above. 
Let $P := P_{3, \, 4} \, P_{1, \, 2} \in \GL(4, k)$. 
Then the following assertions {\rm (1)}, {\rm (2)}, {\rm (3)} hold true: 
\begin{enumerate}[label = {\rm (\arabic*)}]
\item $\Inn_P \circ \sigma^* =\sigma^+ \circ  F^{e_1}$. 

\item $\sigma^+$ is a homomorphism. 

\item $\sigma^*$ is a homomorphism. 
\end{enumerate} 
\end{lem}

\begin{proof}
The proof is straightforward. 
\end{proof}

\begin{lem}
The following assertions {\rm (1)}, {\rm (2)}, {\rm (3)}, {\rm (4)} hold true: 
\begin{enumerate}[label = {\rm (\arabic*)}]
\item 
$
 V(4)^{\sigma^+} = k 
 \left(
 \begin{array}{c}
  0 \\
  0 \\
  0 \\
  1
 \end{array}
 \right) 
 $. 

\item 
$
 W(4)^{\sigma^+} = k  \left( \begin{array}{c c c c} 1 & 0 & 0 & 0  \end{array} \right) 
 $. 

\item $d(\sigma^+) = (1, 1)$. 

\item $d(\sigma^*) = (1, 1)$. 
\end{enumerate} 
\end{lem}

\begin{proof}
(1) Let 
\[
\bm{v} 
= 
\left(
\begin{array}{c}
 v_1 \\
 v_2 \\
 v_3 \\
 v_4
\end{array}
\right) 
\in 
V^{\sigma^+} . 
\]
Then 
\[
\left(
\begin{array}{c}
 v_1 \\
 v_2 \\
 v_3
\end{array}
\right)
 \in 
 (\, k^{\oplus 3} \,)^{\tau^+} , 
\]
where $\tau^+ : \SL(2, k) \to \SL(3, k)$ be the homomorphism defined by 
\[
\renewcommand{\arraystretch}{1.5} 
 \tau^+
 \left(
 \begin{array}{c c}
  A & B \\
  C & D
 \end{array}
 \right)
 := 
 \left(
 \begin{array}{c c c}
  1 & 0 & 0 \\
  A \, B & A^2 & B^2 \\
  C \, D & C^2 & D^2
 \end{array}
 \right) . 
\]
Since $(\, k^{\oplus 3} \,)^{\tau^+} = \bm{0}$, we have 
\[
\bm{v}
= 
\left(
\begin{array}{c}
 0 \\
 0 \\
 0 \\
 v_4
\end{array}
\right) . 
\]
Thus we have the desired equality.

(2) The proof is similar to the proof of the above assertion (1).

(3) The proof is straightforward.

(4) Use the above assertion (3). 
\end{proof}

\subsubsection{${\rm (VII)}^*$}

Assume $p = 3$. Let $e_1$ be an integer such that 
\[
 e_1 \geq 0. 
\]
We can define a homomorphism $\sigma^* : \SL(2, k) \to \SL(4, k)$ as 
\begin{align*}
\renewcommand{\arraystretch}{1.5} 
\sigma^*
\left(
\begin{array}{c c}
 a & b \\
 c & d
\end{array}
\right) 
 := 
\left(
\begin{array}{l l l l}
 a^{p^{e_1 + 1}} 
 & 0
 & 0
 & b^{p^{e_1 + 1}} \\
 \frac{1}{2} \, a^{2 \, p^{e_1}} \, c^{p^{e_1}} 
 &  a^{p^{e_1}} 
 &  b^{p^{e_1}} 
 &  \frac{1}{2} \, b^{2 \, p^{e_1}} \, d^{p^{e_1}}  \\
 a^{p^{e_1}} \, c^{2 \, p^{e_1}} 
 &  c^{p^{e_1}} 
 &  d^{p^{e_1}}
 &  b^{p^{e_1}} \, d^{2 \, p^{e_1}}  \\
 c^{p^{e_1 + 1}} 
 & 0 
 & 0
 & d^{p^{e_1 + 1}} 
\end{array}
\right) . 
\end{align*}
We can define a homomorphism $\sigma^+ : \SL(2, k) \to \SL(4, k)$ as 
\begin{align*}
\renewcommand{\arraystretch}{1.5} 
\sigma^+ 
\left(
\begin{array}{c c}
 A & B \\
 C & D
\end{array}
\right) 
 := 
\left(
\begin{array}{l l | l l}
 A 
 & B
 & \frac{1}{2} \, A^2 \, C
 & \frac{1}{2} \, B^2 \, D  \\
C
 & D  
 & A \, C^2 
 & B \, D^2  \\ 
\hline 
0 
 & 0  
 & A^3 
 & B^3 \\
0 
 & 0 
 & C^3 
 & D^3 
\end{array}
\right) . 
\end{align*}

\begin{lem}
Let $\sigma^* : \SL(2, k) \to \SL(4, k)$ and  
$\sigma^+ : \SL(2, k) \to \SL(4, k)$ be as above. 
Let $P := P_{1, \, 2} \,  P_{2, \, 3} \in \GL(4, k)$. 
Then the following assertions {\rm (1)}, {\rm (2)}, {\rm (3)} hold true: 
\begin{enumerate}[label = {\rm (\arabic*)}]
\item $\Inn_P \circ \sigma^* =\sigma^+ \circ  F^{e_1}$. 

\item $\sigma^+$ is a homomorphism. 

\item $\sigma^*$ is a homomorphism. 
\end{enumerate} 
\end{lem}

\begin{proof}
The proof is straightforward. 
\end{proof}

\begin{lem}
The following assertions {\rm (1)}, {\rm (2)}, {\rm (3)}, {\rm (4)} hold true: 
\begin{enumerate}[label = {\rm (\arabic*)}]
\item 
$
 V(4)^{\sigma^+} = 0 
 $. 

\item 
$
 W(4)^{\sigma^+} = 0 
 $. 
 
\item $d(\sigma^+) = (0, 0)$. 

\item $d(\sigma^*) = (0, 0)$. 
\end{enumerate} 
\end{lem}

\begin{proof}
Refer to the proof of Lemma 5.5. 
\end{proof}

\subsubsection{${\rm (IX)}^*$}

Assume $p \geq 3$. 
Let $e_1$ be an integer such that 
\[ 
e_1 \geq 0 . 
\]
We can define a morphism $\sigma^* : \SL(2, k) \to \SL(4, k)$ as 
\begin{align*}
\renewcommand{\arraystretch}{1.5} 
\sigma^*
\left(
\begin{array}{c c}
 a & b \\
 c & d
\end{array}
\right) 
 & := 
 \renewcommand{\arraystretch}{1.5} 
\left(
\begin{array}{l l l l}
 a^{2 \,p^{e_1}} 
 & a^{p^{e_1}} \, b^{p^{e_1}} 
 & 0
 & \frac{1}{2} \, b^{2 \, p^{e_1}} \\
 2 \, a^{p^{e_1}} \, c^{p^{e_1}} 
 &  a^{p^{e_1}} \, d^{p^{e_1}}  +  b^{p^{e_1}} \, c^{p^{e_1}} 
 &  0 
 &  b^{p^{e_1}} \, d^{p^{e_1}}  \\
 0
 & 0 
 & 1 
 & 0 \\
 2 \, c^{2 \, p^{e_1}} 
 & 2 \, c^{p^{e_1}} \, d^{p^{e_1}} 
 & 0
 & d^{2 \, p^{e_1}} 
\end{array}
\right) .  
\end{align*}
We can define a morphism $\sigma^+ : \SL(2, k) \to \SL(4, k)$ as 
\begin{align*}
\renewcommand{\arraystretch}{1.5} 
\sigma^+ 
\left(
\begin{array}{c c}
 A & B \\
 C & D
\end{array}
\right) 
 := 
\left(
\begin{array}{l l l | l}
 A^2 
 & A \, B
 & \frac{1}{2} \, B^2 
 & 0  \\
2 \, A \, C
 & A \, D + B \, C   
 & B \, D  
 & 0  \\ 
2 \, C^2 
 & 2 \, C \, D
 & D^2 
 & 0 \\ 
\hline 
0 
 & 0 
 & 0 
 & 1 
\end{array}
\right) . 
\end{align*}

\begin{lem}
Let $\sigma^* : \SL(2, k) \to \SL(4, k)$ and  
$\sigma^+ : \SL(2, k) \to \SL(4, k)$ be as above. 
Let $P := P_{3, \, 4} \in \GL(4, k)$. 
Then the following assertions {\rm (1)}, {\rm (2)}, {\rm (3)} hold true: 
\begin{enumerate}[label = {\rm (\arabic*)}]
\item $\Inn_P \circ \sigma^* =\sigma^+ \circ  F^{e_1}$. 

\item $\sigma^+$ is a homomorphism. 

\item $\sigma^*$ is a homomorphism. 
\end{enumerate} 
\end{lem}

\begin{proof}
The proof is straightforward. 
\end{proof}

\begin{lem}
The following assertions {\rm (1)}, {\rm (2)}, {\rm (3)}, {\rm (4)} hold true: 
\begin{enumerate}[label = {\rm (\arabic*)}]
\item 
$
 V(4)^{\sigma^+} = 
 k 
 \left(
 \begin{array}{c}
  0 \\
  0 \\
  0 \\
  1
 \end{array}
 \right)
 $. 

\item 
$
 W(4)^{\sigma^+} 
 = 
 k \left(
 \begin{array}{c c c c}
  0 & 0 & 0 & 1 
 \end{array}
 \right) 
 $. 
 
\item $d(\sigma^+) = (1, 1)$. 

\item $d(\sigma^*) = (1, 1)$. 
\end{enumerate} 
\end{lem}

\begin{proof}
Consider the two regular matrices 
\[
\left(
\begin{array}{c c}
 A & B \\
 C & D 
\end{array}
\right) 
= 
\left(
\begin{array}{c c}
 1 & 1 \\
 0 & 1 
\end{array}
\right) , \;
\left(
\begin{array}{c c}
 1 & 0 \\
 1 & 1 
\end{array}
\right)  
\in \SL(2, k) . 
\]
\end{proof}

\subsubsection{${\rm (XI)}^*$}

Assume $p = 2$. 
Let $e_1$ be an integer such that 
\[ 
e_1 \geq 0 . 
\]
We can define a morphism $\sigma^* : \SL(2, k) \to \SL(4, k)$ as 
\begin{align*}
\renewcommand{\arraystretch}{1.5} 
\sigma^*
\left(
\begin{array}{c c}
 a & b \\
 c & d
\end{array}
\right) 
 := 
\left(
\begin{array}{l l l l}
 a^{2 \, p^{e_1}} 
 & a^{p^{e_1}} \, b^{p^{e_1}} 
 & 0
 & b^{2 \, p^{e_1}} \\
 0 
 & 1 
 & 0 
 & 0  \\
 0
 & 0 
 & 1 
 & 0 \\
c^{2 \, p^{e_1}} 
 & c^{p^{e_1}} \, d^{p^{e_1}} 
 & 0
 & d^{2 \, p^{e_1}} 
\end{array}
\right) . 
\end{align*}
We can define a morphism $\sigma^+ : \SL(2, k) \to \SL(4, k)$ as 
\begin{align*}
\renewcommand{\arraystretch}{1.5} 
\sigma^+ 
\left(
\begin{array}{c c}
 A & B \\
 C & D
\end{array}
\right) 
 := 
\left(
\begin{array}{l l l | l}
 A^2 
 & B^2
 & A \, B
 & 0  \\
C^2 
 & D^2 
 & C \, D 
 & 0 \\ 
0 
 & 0  
 & 1 
 & 0 \\
\hline 
0 
 & 0 
 & 0
 & 1
\end{array}
\right) . 
\end{align*}

\begin{lem}
Let $\sigma^* : \SL(2, k) \to \SL(4, k)$ and  
$\sigma^+ : \SL(2, k) \to \SL(4, k)$ be as above. 
Let $P := P_{3, \, 4} \, P_{2, \, 3} \in \GL(4, k)$. 
Then the following assertions {\rm (1)}, {\rm (2)}, {\rm (3)} hold true: 
\begin{enumerate}[label = {\rm (\arabic*)}]
\item $\Inn_P \circ \sigma^* =\sigma^+ \circ  F^{e_1}$. 

\item $\sigma^+$ is a homomorphism. 

\item $\sigma^*$ is a homomorphism. 
\end{enumerate} 
\end{lem}

\begin{proof}
The proof is straightforward. 
\end{proof}

\begin{lem}
The following assertions {\rm (1)}, {\rm (2)}, {\rm (3)}, {\rm (4)} hold true: 
\begin{enumerate}[label = {\rm (\arabic*)}]
\item 
$
 V(4)^{\sigma^+} = 
 k 
 \left(
 \begin{array}{c}
  0 \\
  0 \\
  0 \\
  1
 \end{array}
 \right) 
 $. 

\item 
$
 W(4)^{\sigma^+} 
 = 
 k \left(
 \begin{array}{c c c c}
  0 & 0 & 1 & 0 
 \end{array}
 \right)
 \oplus 
 k \left(
 \begin{array}{c c c c}
  0 & 0 & 0 & 1 
 \end{array}
 \right) 
 $. 
 
\item $d(\sigma^+) = (1, 2)$. 

\item $d(\sigma^*) = (1, 2)$. 
\end{enumerate} 
\end{lem}

\begin{proof}
(1) Consider the regular matrices 
\[
\left(
\begin{array}{c c}
 A & B \\
 C & D 
\end{array}
\right) 
= 
\left(
\begin{array}{c c}
 1 & 1 \\
 0 & 1 
\end{array}
\right) , \;
\left(
\begin{array}{c c}
 1 & 0 \\
 1 & 1 
\end{array}
\right)  ,  \; 
\left(
\begin{array}{c c}
 u & 0 \\
 0 & u^{-1} 
\end{array}
\right)  
\in \SL(2, k) 
\qquad (\, u \in k \backslash \{\, 0, \, 1 \,\}\,) . 
\]

(2) Consider the two regular matrices 
\[
\left(
\begin{array}{c c}
 A & B \\
 C & D 
\end{array}
\right) 
= 
\left(
\begin{array}{c c}
 1 & 1 \\
 0 & 1 
\end{array}
\right) , \;
\left(
\begin{array}{c c}
 1 & 0 \\
 1 & 1 
\end{array}
\right)  
\in \SL(2, k) . 
\]

(3) The proof is straightforward.

(4) Use the above assertion (3). 
\end{proof}

\subsubsection{${\rm (XV)}^*$}

Let $e_2$ and $e_3$ be integers such that 
\[
 e_2 \geq e_3 \geq 0 . 
\]
We can define a morphism $\sigma^* : \SL(2, k) \to \SL(4, k)$ as 
\begin{align*}
\renewcommand{\arraystretch}{1.5} 
\sigma^*
\left(
\begin{array}{c c}
 a & b \\
 c & d
\end{array}
\right) 
 := 
\left(
\begin{array}{l l l l}
 a^{p^{e_2}} 
 & 0 
 & 0
 & b^{p^{e_2}} \\
 0 
 & a^{p^{e_3}}
 & b^{p^{e_3}} 
 & 0  \\
 0
 & c^{p^{e_3}}
 & d^{p^{e_3}} 
 & 0 \\
c^{p^{e_2}} 
 & 0 
 & 0
 & d^{p^{e_2}} 
\end{array}
\right) . 
\end{align*}

Let $\Delta : \SL(2, k) \to \SL(2, k) \times \SL(2, k)$ be the homomorphism 
defined by 
\[
 \Delta(X) := (X, X) . 
\]

Let $F^{e_2} \times F^{e_3} : \SL(2, k) \times \SL(2, k) \to \SL(2, k) \times \SL(2, k)$ 
be the homomorphism defined by 
\begin{align*}
&
\renewcommand{\arraystretch}{1.5} 
 ( F^{e_2} \times F^{e_3}) 
\left( \, 
\left(
\begin{array}{c c}
 A & B \\
 C & D
\end{array}
\right) , 
\; 
\left(
\begin{array}{l l }
 A' & B' \\
 C' & D'
\end{array}
\right) 
\, \right) \\
& 
\renewcommand{\arraystretch}{1.5} 
\qquad  := 
\left( \,  
\left(
\begin{array}{l l }
 A^{p^{e_2}} & B^{p^{e_2}} \\
 C^{p^{e_2}} & D^{p^{e_2}}
\end{array}
\right)  , 
\; 
\left(
\begin{array}{l l }
 {A'}^{p^{e_3}} & {B'}^{p^{e_3}} \\
 {C'}^{p^{e_3}} & {D'}^{p^{e_3}}
\end{array}
\right)  
\, \right) . 
\end{align*}

Let $i : \SL(2, k) \times \SL(2, k) \to \SL(4, k)$ be the homomorphism 
defined by 
\[
 i (X_1, X_2) 
:= 
\left(
\begin{array}{c c}
 X_1 & O_2 \\
 O_2 & X_2
\end{array}
\right) . 
\]

%
%

\begin{lem}
Let $\sigma^* : \SL(2, k) \to \SL(4, k)$ and  
$\sigma^+ : \SL(2, k) \to \SL(4, k)$ be as above. 
Let $P := P_{3, \, 4} \, P_{2, \, 3} \in \GL(4, k)$. 
Then the following assertions {\rm (1)} and {\rm (2)} hold true: 
\begin{enumerate}[label = {\rm (\arabic*)}]
\item $\Inn_P \circ \sigma^* 
= i \circ (F^{e_2} \times F^{e_3}) \circ \Delta$, i.e., 
\[
\xymatrix@R=36pt@C=36pt@M=6pt{
 \SL(2, k) \ar[r]^(.38)\Delta \ar@/_23pt/[rrr]_{\Inn_P \circ \sigma^*}
 & \SL(2, k) \times \SL(2, k) \ar[r]^{F^{e_2} \times F^{e_3}} 
 & \SL(2, k) \times \SL(2, k) \ar[r]^(.62)i
 & \SL(4, k)
}
\]

\item $\sigma^*$ is a homomorphism. 
\end{enumerate} 
\end{lem}

\begin{proof}
The proof is straightforward. 
\end{proof}

We can define a morphism $\sigma^+ : \SL(2, k) \to \SL(4, k)$ as 
\begin{align*}
\renewcommand{\arraystretch}{1.5} 
\sigma^+ 
\left(
\begin{array}{c c}
 A & B \\
 C & D
\end{array}
\right) 
  := 
\left( 
\begin{array}{c c}
 1 & 0 \\
 0 & 1
\end{array} 
\right) 
\otimes 
\left(
\begin{array}{c c}
 A & B \\
 C & D
\end{array}
\right) 
=
\left( 
\begin{array}{c c | c c}
A & B & 0 & 0\\
C & D & 0 & 0 \\
\hline 
 0 & 0 & A & B \\
 0 & 0 & C & D 
\end{array}
\right) . 
\end{align*}
Clearly, $\sigma^+ = i \circ \Delta$.

\begin{lem}
The following assertions {\rm (1)}, {\rm (2)}, {\rm (3)}, {\rm (4)} hold true: 
\begin{enumerate}[label = {\rm (\arabic*)}]
\item 
$
 V(4)^{\sigma^+} = 0 
 $. 

\item 
$
 W(4)^{\sigma^+} = 0
 $. 
 
\item $d(\sigma^+) = (0, 0)$. 

\item $d(\sigma^*) = (0, 0)$. 
\end{enumerate} 
\end{lem}

\begin{proof}
Consider the two regular matrices 
\[
\left(
\begin{array}{c c}
 A & B \\
 C & D 
\end{array}
\right) 
= 
\left(
\begin{array}{c c}
 1 & 1 \\
 0 & 1 
\end{array}
\right) , \;
\left(
\begin{array}{c c}
 1 & 0 \\
 1 & 1 
\end{array}
\right)  
\in \SL(2, k) . 
\]
\end{proof}

\subsubsection{${\rm (XIX)}^*$}

Assume $p = 2$. 
Let $e_1$ be an integer such that 
\[ 
e_1 \geq 0 . 
\]
We can define a morphism $\sigma^* : \SL(2, k) \to \SL(4, k)$ as 
\begin{align*}
\renewcommand{\arraystretch}{1.5} 
\sigma^*
\left(
\begin{array}{c c}
 a & b \\
 c & d
\end{array}
\right) 
 := 
\left(
\begin{array}{l l l l}
 a^{2 \, p^{e_1}} 
 & 0 
 & 0
 & b^{2 \, p^{e_1}} \\
 0 
 & 1 
 & 0  
 & 0  \\
 a^{p^{e_1}} \, c^{p^{e_1}} 
 & 0 
 & 1 
 & b^{p^{e_1}} \, d^{p^{e_1}}  \\
c^{2 \, p^{e_1}} 
 & 0 
 & 0
 & d^{2 \, p^{e_1}} 
\end{array}
\right) . 
\end{align*}
We can define a morphism $\sigma^+ : \SL(2, k) \to \SL(4, k)$ as 
\begin{align*}
\renewcommand{\arraystretch}{1.5} 
\sigma^+ 
\left(
\begin{array}{c c}
 A & B \\
 C & D
\end{array}
\right) 
 := 
\left(
\begin{array}{l l l | l}
 A^2 
 & B^2
 & 0 
 & 0  \\
C^2 
 & D^2 
 & 0
 & 0 \\ 
A \, C 
 & B \, D
 & 1 
 & 0 \\
\hline 
0 
 & 0 
 & 0
 & 1
\end{array}
\right) . 
\end{align*}

\begin{lem}
Let $\sigma^* : \SL(2, k) \to \SL(4, k)$ and  
$\sigma^+ : \SL(2, k) \to \SL(4, k)$ be as above. 
Let $P := P_{2, \, 3} \, P_{3, \, 4} \, P_{2, \, 3} \in \GL(4, k)$. 
Then the following assertions {\rm (1)}, {\rm (2)}, {\rm (3)} hold true: 
\begin{enumerate}[label = {\rm (\arabic*)}]
\item $\Inn_P \circ \sigma^* =\sigma^+ \circ  F^{e_1} $. 

\item $\sigma^+$ is a homomorphism. 

\item $\sigma^*$ is a homomorphism. 
\end{enumerate} 
\end{lem}

\begin{proof}
The proof is straightforward. 
\end{proof}

\begin{lem}
The following assertions {\rm (1)}, {\rm (2)}, {\rm (3)}, {\rm (4)} hold true: 
\begin{enumerate}[label = {\rm (\arabic*)}]
\item 
$
 V(4)^{\sigma^+} = 
 k 
 \left(
 \begin{array}{c}
  0 \\
  0 \\
  1 \\
  0
 \end{array}
 \right) 
  \oplus 
 k 
 \left(
 \begin{array}{c}
  0 \\
  0 \\
  0 \\
  1
 \end{array}
 \right)
 $. 

\item 
$
 W(4)^{\sigma^+} 
 = 
 k \left(
 \begin{array}{c c c c}
  0 & 0 & 0 & 1 
 \end{array}
 \right)
 $. 
 
\item $d(\sigma^+) = (2, 1)$. 

\item $d(\sigma^*) = (2, 1)$. 
\end{enumerate} 
\end{lem}

\begin{proof}
Refer to the proof of Lemma 5.15. 

\end{proof}

\subsubsection{${\rm (XXI)}^*$}

Assume $p = 2$. 
Let $e_1$ be an integer such that 
\[
 e_1 \geq 0 . 
\]
We can define a morphism $\sigma^* : \SL(2, k) \to \SL(4, k)$ as 
\begin{align*}
\renewcommand{\arraystretch}{1.5} 
\sigma^*
\left(
\begin{array}{c c}
 a & b \\
 c & d
\end{array}
\right) 
 := 
\left(
\begin{array}{l l l l}
 a^{2 \,p^{e_1}} 
 & 0 
 & a^{p^{e_1}} \, b^{p^{e_1}}
 & b^{2 \, p^{e_1}} \\
 a^{p^{e_1}} \, c^{p^{e_1}} 
 &  1
 & b^{p^{e_1}} \, c^{p^{e_1}}
 & b^{p^{e_1}} \, d^{p^{e_1}} \\
 0
 & 0 
 & 1
 & 0 \\
 c^{2 \, p^{e_1}} 
 & 0 
 & c^{p^{e_1}} \, d^{p^{e_1}} 
 & d^{2 \, p^{e_1}} 
\end{array}
\right) . 
\end{align*}
We can define a morphism $\sigma^+ : \SL(2, k) \to \SL(4, k)$ as 
\begin{align*}
\renewcommand{\arraystretch}{1.5} 
\sigma^+ 
\left(
\begin{array}{c c}
 A & B \\
 C & D
\end{array}
\right) 
 := 
\left(
\begin{array}{l l l l}
 1 
 & A \, C 
 & B \, D
 & B \, C \\
 0 
 & A^2
 & B^2 
 & A \, B \\ 
 0 
 & C^2 
 & D^2 
 & C \, D  \\
 0 
 & 0 
 & 0  
 & 1  
\end{array}
\right) . 
\end{align*}

\begin{lem}
Let 
\[
P_1 := 
\left(
\begin{array}{c c c c}
 1 & 0 & 0 & 1 \\
 0 & 1 & 0 & 0 \\
 0 & 0 & 1 & 0 \\
 1 & 0 & 0 & 0 
\end{array}
\right) 
\in \GL(4, k) 
\]
and let $P_2 := P_{3, \, 4} \, P_{1, \, 2} \in \GL(4, k)$. 
Let $\sigma^*_{\rm (V)^*}$ and $\sigma^+_{\rm (V)^*}$ respectively denote 
the homomorphisms $\sigma^*$ and $\sigma^+$ given in $\rm (V)^*$. 
Then the following assertions {\rm (1)}, {\rm (2)}, {\rm (3)} hold true: 
\begin{enumerate}[label = {\rm (\arabic*)}]
\item $\Inn_{P_1} \circ \sigma^+ = \sigma^+_{\rm (V)^*}$.

\item $\Inn_{P_2  \, P_1} \circ \sigma^* = \sigma^*_{\rm (V)^*}$. 
So, $\sigma^* \sim \sigma^*_{\rm (V)^*}$.

\item $\sigma^*$ is a homomorphism. 
\end{enumerate} 
\end{lem}

\begin{proof}
The proof is straightforward. 

\end{proof}

\subsubsection{${\rm (XXII)}^*$}

Let $e_1$ be an integer such that 
\[ 
e_1 \geq 0 . 
\]
We can define a morphism $\sigma^* : \SL(2, k) \to \SL(4, k)$ as 
\begin{align*}
\renewcommand{\arraystretch}{1.5} 
\sigma^*
\left(
\begin{array}{c c}
 a & b \\
 c & d
\end{array}
\right) 
 := 
\left(
\begin{array}{l l l l}
 a^{p^{e_1}} 
 & 0 
 & b^{p^{e_1}}
 & 0 \\
 0 
 &  a^{p^{e_1}} 
 & 0 
 & b^{p^{e_1}} \\
 c^{p^{e_1}} 
 & 0 
 & d^{p^{e_1}} 
 & 0 \\
0
 & c^{p^{e_1}} 
 & 0 
 & d^{p^{e_1}} 
\end{array}
\right) . 
\end{align*}

\begin{lem}
Let $\sigma^* : \SL(2, k) \to \SL(4, k)$ be as above. 
Let $\sigma_{\rm (XV)^*}^* : \SL(2, k) \to \SL(4, k)$ denote the homomorphism $\sigma^*$ given in 
$\rm (XV)^*$. Assume $e_2 = e_3 = e_1$. 
Let $P := P_{3, \, 4} \in \GL(4, k)$. 
Then the following assertions {\rm (1)} and {\rm (2)} hold true: 
\begin{enumerate}[label = {\rm (\arabic*)}]
\item 
$
\Inn_{P} \circ \sigma^* 
= 
\sigma_{\rm (XV)^*}^*$. 
So, $\sigma^* \sim \sigma_{\rm (XV)^*}^*$. 

\item $\sigma^*$ is a homomorphism. 
\end{enumerate} 
\end{lem}

\begin{proof}
The proof is straightforward. 

\end{proof}

\subsubsection{${\rm (XXIV)}^*$}

Let $e_2$ be an integer such that 
\[ 
e_2 \geq 0 . 
\]
We can define a morphism $\sigma^* : \SL(2, k) \to \SL(4, k)$ as 
\begin{align*}
\renewcommand{\arraystretch}{1.5} 
\sigma^*
\left(
\begin{array}{c c}
 a & b \\
 c & d
\end{array}
\right) 
 := 
\left(
\begin{array}{l l l l}
 a^{p^{e_2}} 
 & 0 
 & 0
 & b^{p^{e_2}} \\
 0 
 & 1 
 & 0
 & 0  \\
 0
 & 0 
 & 1
 & 0 \\
c^{p^{e_2}} 
 & 0 
 & 0
 & d^{p^{e_2}} 
\end{array}
\right) . 
\end{align*}
We can define a morphism $\sigma^+ : \SL(2, k) \to \SL(4, k)$ as 
\begin{align*}
\renewcommand{\arraystretch}{1.5} 
\sigma^+ 
\left(
\begin{array}{c c}
 A & B \\
 C & D
\end{array}
\right) 
 := 
\left(
\begin{array}{l l | l l}
 A 
 & B 
 & 0 
 & 0 \\
 C 
 & D
 & 0 
 & 0 \\ 
\hline 
 0 
 & 0 
 & 1 
 & 0  \\
 0 
 & 0 
 & 0
 & 1  
\end{array}
\right) . 
\end{align*}

%

\begin{lem}
Let $\sigma^* : \SL(2, k) \to \SL(4, k)$ and  
$\sigma^+ : \SL(2, k) \to \SL(4, k)$ be as above. 
Let $P := P_{3, \, 4} \, P_{2, \, 3} \in \GL(4, k)$. 
Then the following assertions {\rm (1)}, {\rm (2)}, {\rm (3)} hold true: 
\begin{enumerate}[label = {\rm (\arabic*)}]
\item $\Inn_P \circ \sigma^* =\sigma^+ \circ  F^{e_2}$. 

\item $\sigma^+$ is a homomorphism. 

\item $\sigma^*$ is a homomorphism. 
\end{enumerate} 
\end{lem}

\begin{proof}
The proof is straightforward. 
\end{proof}

\begin{lem}
The following assertions {\rm (1)}, {\rm (2)}, {\rm (3)}, {\rm (4)} hold true: 
\begin{enumerate}[label = {\rm (\arabic*)}]
\item 
$
 V(4)^{\sigma^+} = 
 k 
 \left(
 \begin{array}{c}
  0 \\
  0 \\
  1 \\
  0
 \end{array}
 \right) 
  \oplus 
 k 
 \left(
 \begin{array}{c}
  0 \\
  0 \\
  0 \\
  1
 \end{array}
 \right) 
 $. 

\item 
$
 W(4)^{\sigma^+} 
 = 
  k \left(
 \begin{array}{c c c c}
  0 & 0 & 1 & 0 
 \end{array}
 \right)
  \oplus 
 k \left(
 \begin{array}{c c c c}
  0 & 0 & 0 & 1 
 \end{array}
 \right)
 $. 
 
\item $d(\sigma^+) = (2, 2)$. 

\item $d(\sigma^*) = (2, 2)$. 
\end{enumerate} 
\end{lem}

\begin{proof}
The proof is straightforward. 
\end{proof}

\subsubsection{${\rm (XXVI)}^*$}

We can define a homomorphism $\sigma^* : \SL(2, k) \to \SL(4, k)$ as 
\begin{align*}
\renewcommand{\arraystretch}{1.5} 
\sigma^*
\left(
\begin{array}{c c}
 a & b \\
 c & d
\end{array}
\right) 
 := 
I_4 . 
\end{align*}

\begin{lem}
The following assertions {\rm (1)}, {\rm (2)}, {\rm (3)} hold true: 
\begin{enumerate}[label = {\rm (\arabic*)}]
\item 
$
 V(4)^{\sigma^*} = 
 V(4)
 $. 

\item 
$
 W(4)^{\sigma^*} 
 = 
 W(4)
 $. 
 
\item $d(\sigma^*) = (4, 4)$. 
\end{enumerate} 
\end{lem}

\begin{proof}
The proof is straightforward. 
\end{proof}

\subsection{An overlapping classification of homomorphisms from $\SL(2, k)$ to $\SL(4, k)$}

\begin{thm}
Let 
\[
\text{
$\nu = {\rm I}$, \;  $\rm II$, \;  $\rm IV$, \;  $\rm V$, \;  $\rm VII$, \;  
$\rm IX$, \;  $\rm XI$, \;  
$\rm XV$, \;  $\rm XIX$, \;  $\rm XXI$, \;  $\rm XXII$, \;  
$\rm XXIV$, \;  $\rm XXVI$. 
}
\]
Let $(\varphi^*, \omega^*)$ be a pair of the form $(\nu)$. 
Let $\psi^* : \rB(2, k) \to \SL(4, k)$ be a  
homomorphism defined by $\psi^* := \psi_{\varphi^*, \, \omega^*} \circ \jmath^{-1}$. 
Then the following assertions {\rm (1)} and {\rm (2)} hold true: 
\begin{enumerate}[label = {\rm (\arabic*)}]
\item Let $\sigma^* : \SL(2, k) \to \SL(4, k)$ be the homomorphism 
of the form $(\nu^*)$. 
Then $\sigma^* \circ \imath_{\rB(2, k)} = \psi^*$.

\item There exists a unique homomorphism $\widehat{\sigma} : \SL(2, k) \to \SL(4, k)$ 
such that $\widehat{\sigma}  \circ \imath_{\rB(2, k)} = \psi^*$. 
\end{enumerate} 
\end{thm}

\begin{proof}

(1) By the construction of $\sigma^*$, we have 
\begin{align*} 
\sigma^*
\left(
\begin{array}{c c}
 a & b \\
 0 & d
\end{array}
\right) 
 = 
\omega^*(a) \, \varphi^*\left(\frac{b}{a}\right)  
 = 
\varphi^*(a \, b) \, \omega^*(a) 
 = 
(\psi_{\varphi^*, \, \omega^*} \circ \jmath^{-1} ) 
\left(
\begin{array}{c c}
 a & b \\
 0 & d
\end{array}
\right) . 
\end{align*}

(2) The existence of $\widehat{\sigma} $ follows from the above assertion (1). 
Let $\widehat{\phi}^+ : \G_a \to \SL(4, k)$, 
$\widehat{\omega} : \G_m \to \SL(4, k)$, 
$\widehat{\phi}^- : \G_a \to \SL(4, k)$ be the homomorphisms defined by 
\[
 {\widehat{\phi}}^+(t) 
 :=
 \widehat{\sigma} 
 \left(
 \begin{array}{c c}
  1 & t \\
  0 & 1
 \end{array}
 \right) , \qquad 
 {\widehat{\omega}}(u) 
 :=
 \widehat{\sigma} 
 \left(
 \begin{array}{c c}
  u & 0 \\
  0 & u^{-1}
 \end{array}
 \right) , \qquad 
 {\widehat{\phi}}^-(s) 
 :=
 \widehat{\sigma} 
 \left(
 \begin{array}{c c}
  1 & 0 \\
  s & 1
 \end{array}
 \right) . 
\]
Since $\widehat{\sigma}$ is an extension of $\psi^*$, we have 
\[
 \varphi^* =  {\widehat{\phi}}^+ , \qquad 
 \omega^* =  {\widehat{\omega}} . 
\]
We know from Lemmas 4.1 -- 4.26 in Subsection 4.1 that 
\[
\phi^- =  {\widehat{\phi}}^- . 
\]
Hence we have $\sigma^* = \widehat{\sigma}$. 
\end{proof}

\begin{thm}
Let $\sigma : \SL(2, k) \to \SL(4, k)$ be a homomorphism. 
Then there exists a homomorphism $\sigma^* : \SL(2, k) \to \SL(4, k)$ 
satisfying the following conditions {\rm (i)} and {\rm (ii)}: 
\begin{enumerate}[label = {\rm (\roman*)}]
\item $\sigma$ and $\sigma^*$ are equivalent, i.e., $\sigma \sim \sigma^*$. 

\item $\sigma^*$ has one of the forms 
$\rm (I)^*$, $\rm (II)^*$, $\rm (IV)^*$, $\rm (V)^*$, $\rm (VII)^*$, 
$\rm (IX)^*$, $\rm (XI)^*$, 
$\rm (XV)^*$, $\rm (XIX)^*$, 
$\rm (XXIV)^*$, $\rm (XXVI)^*$. 
\end{enumerate} 
\end{thm}

\begin{center} 

\renewcommand{\arraystretch}{1.5} 
\begin{tabular}{| p{5em} | p{5em} | p{5em} | p{5em} | p{5em} |}
\hline 
 $p = 2$ 
 & $ p = 3$ 
 & $p \geq 5$ 
 & $p \geq 2$ 
 & $d$ \\
\hline 
 & 
 & $\rm (I)^*$
 &  
 & $(0, 0)$ \\
\hline  
 & $\rm (II)^*$ 
 & 
 & 
 & $(0, 0)$ \\
\hline 
 $\rm (IV)^*$
 & $\rm (IV)^*$ 
 & $\rm (IV)^*$
 & $\rm (IV)^*$ 
 & $(0, 0)$ \\ 
\hline 
 $\rm (V)^*$
 & 
 & 
 & 
 & $(1, 1)$ \\
\hline 
 & $\rm (VII)^*$
 &  
 & 
 & $(0, 0)$ \\
\hline  
 & $\rm (IX)^*$
 & $\rm (IX)^*$ 
 & 
 & $(1, 1)$ \\
\hline 
 $\rm (XI)^*$
 & 
 & 
 & 
 & $(1, 2)$ \\
\hline
 $\rm (XV)^*$
 & $\rm (XV)^*$ 
 & $\rm (XV)^*$
 & $\rm (XV)^*$ 
 & $(0, 0)$ \\
\hline
 $\rm (XIX)^*$
 & 
 & 
 & 
 & $(2, 1)$ \\
\hline
 $\rm (XXIV)^*$
 & $\rm (XXIV)^*$ 
 & $\rm (XXIV)^*$ 
 & $\rm (XXIV)^*$ 
 & $(2, 2)$ \\
\hline
 $\rm (XXVI)^*$
 & $\rm (XXVI)^*$ 
 & $\rm (XXVI)^*$
 & $\rm (XXVI)^*$ 
 & $(4, 4)$ \\
\hline 
$7$ types 
 & $7$ types
 & $6$ types 
 & $4$ types 
 \\
\cline{1-4} 
\end{tabular}

\end{center}

\begin{proof}
There exists a regular matrix $P$ of $\GL(4, k)$ such that 
$\Inn_P \circ \sigma$ is antisymmetric (see Lemma 1.20 (1)). 
Let $\sigma' : = \Inn_P \circ \sigma$ and consider the homomorphisms 
\begin{align*}
 \varphi_{\sigma'} : \G_a \to \SL(4, k),  \qquad 
 \omega_{\sigma'} : \G_m \to \SL(4, k) , \qquad 
 \varphi_{\sigma'}^- : \G_a \to \SL(4, k) . 
\end{align*}
So, $\varphi_{\sigma'} \in \cU_4$ and $\omega_{\sigma'} \in \Omega(4)$. 
We know from Theorem 3.1 (2) that there exists a pair $(\varphi^{*}, \omega^{*})$ 
of $\cU_4 \times \Omega(4)$ such that 
the following conditions {\rm (a)} and {\rm (b)} hold true: 
\begin{enumerate}[label = {\rm (\alph*)}]
\item $(\varphi_{\sigma'}, \omega_{\sigma'}) \sim (\varphi^{*}, \omega^{*})$. 

\item $(\varphi^{*}, \omega^{*})$ has one of the forms 
\text{\rm (I) -- (XXVI)}. 
\end{enumerate} 
Note that the homomorphism $\psi_{\varphi^*, \, \omega^*} \circ \jmath$ is 
extendable. 
In fact, letting $\psi := \sigma \circ \imath_{\rB(2, k)}$, we have the equivalences 
\[
\begin{array}{ c c l c c}
 \sigma 
 & \sim 
 &
 \Inn_P \circ \sigma = \sigma' , 
 & \\
 \psi 
 & \sim 
 &
 \Inn_P \circ \psi = \psi_{\varphi_{\sigma'}, \, \omega_{\sigma'}} \circ \jmath 
 & \sim 
 & 
\psi_{\varphi^*, \, \omega^*} \circ \jmath
\end{array}
\]
and we can apply Lemma 2.3 to the homomorphisms 
$\psi$ and $\psi_{\varphi^*, \, \omega^*} \circ \jmath$. 
We know from Lemmas 4.1 -- 4.26 that 
$\sigma^*$ has one of the forms $(\nu)$, where 
\[
\text{
$\nu = {\rm I}$, \;  $\rm II$, \;  $\rm IV$, \;  $\rm V$, \;  $\rm VII$, \;  
$\rm IX$, \;  $\rm XI$, \;  
$\rm XV$, \;  $\rm XIX$, \;  $\rm XXI$, \;  $\rm XXII$, \;  
$\rm XXIV$, \;  $\rm XXVI$. 
}
\]
Let $\sigma^* : \SL(2, k) \to \SL(4, k)$ be a homomorphism 
such that $\sigma^* \circ \imath_{\rB(2, k)} =  \psi_{\varphi^*, \, \omega^*} \circ \jmath$. 
We know from Theorem 5.25 (1) that 
$\sigma^*$ has one of the forms $(\nu^*)$. 
We can delete the forms $\rm (XXI)^*$ and $\rm (XXII)^*$ from the forms $(\nu)^*$ 
(see Lemmas 5.20 and 5.21). 
So, $\sigma^*$ has one of the forms 
$\rm (I)^*$, $\rm (II)^*$, $\rm (IV)^*$, $\rm (V)^*$, $\rm (VII)^*$, 
$\rm (IX)^*$, $\rm (XI)^*$, 
$\rm (XV)^*$, $\rm (XIX)^*$, 
$\rm (XXIV)^*$, $\rm (XXVI)^*$. 
In addition, we can show that $\sigma$ is equivalent to 
$\sigma^*$ (see Theorem 5.25 (2)). 
\end{proof}

\section{The classification of homomorphisms from $\SL(2, k)$ to 
$\SL(4, k)$}

\subsection{Homomorphisms from $\SL(2, k)$ to $\SL(4, k)$}

\subsubsection{$\rm (I)^\sharp$}

Assume $p \geq 5$. 
For all integer $e_1 \geq 0$, we can define a homomorphism 
$\sigma_{{\rm (I)^\sharp}, \; e_1} : \SL(2, k) \to \SL(4, k)$ as 
\begin{align*}
\renewcommand{\arraystretch}{1.5} 
& 
\sigma_{{\rm (I)^\sharp}, \; e_1} 
\left(
\begin{array}{c c}
 a & b \\
 c & d
\end{array}
\right) \\
& := 
\left(
\begin{array}{l l l l}
 a^{3 \, p^{e_1}}  
 & a^{2 \, p^{e_1}}  \, b^{p^{e_1}}  
 & a^{p^{e_1}} \, b^{2 \, p^{e_1}}  
 & b^{3 \, p^{e_1}}  \\
3 \, a^{2 \, p^{e_1}} \, c^{p^{e_1}}  
 & a^{p^{e_1}} \, (a^{p^{e_1}} \, d^{p^{e_1}} + 2 \, b^{p^{e_1}} \, c^{p^{e_1}}) 
 & b^{p^{e_1}} \, (2 \, a^{p^{e_1}} \, d^{p^{e_1}} + b^{p^{e_1}} \, c^{p^{e_1}}) 
 & 3 \, b^{2 \, p^{e_1}} \, d^{p^{e_1}} \\
3 \, a^{p^{e_1}} \, c^{2 \, p^{e_1}}  
 & c^{p^{e_1}} \,  (2 \, a^{p^{e_1}} \, d^{p^{e_1}} + b^{p^{e_1}} \, c^{p^{e_1}}) 
 & d^{p^{e_1}} \, (a^{p^{e_1}} \, d^{p^{e_1}} + 2 \, b^{p^{e_1}} \, c^{p^{e_1}}) 
 & 3 \, b^{p^{e_1}} \, d^{2 \, p^{e_1}} \\
 c^{3 \, p^{e_1}} 
 & c^{2 \, p^{e_1}} \, d^{p^{e_1}} 
 & c^{p^{e_1}} \, d^{2 \, p^{e_1}} 
 & d^{3 \, p^{e_1}}   
\end{array}
\right) . 
\end{align*}

\begin{lem}
Let $e_1 \geq 0$ and let $\sigma^*$ be the homomorphism given in $\rm (I)^*$. 
Let 
\[
P := \diag(1, \; 1, \; 2, \; 6) \in \GL(4, k). 
\]
Then $\Inn_P \circ \sigma^* = \sigma_{{\rm (I)^\sharp}, \; e_1}$. 
\end{lem}

\begin{proof}
The proof is straightforward. 
\end{proof}

We can define a homomorphism 
$\omega_{{\rm (I)^\sharp}, \; e_1} : \G_m \to \SL(4, k)$ as 
\[
\omega_{{\rm (I)^\sharp}, \; e_1}
:= 
\sigma_{{\rm (I)^\sharp}, \; e_1} \circ \imath_{\rB(2, k)} \circ \imath_2'. 
\] 

\begin{lem}
We have 
\[
\omega_{{\rm (I)^\sharp}, \; e_1} 
\; \sim \; 
\omega_{3 \, p^{e_1}, \; \; p^{e_1}, \; \; p^{ - e_1}, \; \; - 3 \, p^{e_1}} . 
\]
\end{lem}

\begin{proof}
The homomorphism $\omega_{{\rm (I)^\sharp}, \; e_1} : \G_m \to \SL(4, k)$ 
is equivlent to a homomorphism $\omega : \G_m \to \SL(4, k)$ induced from 
$\sigma^*$ given in $\rm (I)^*$, i.e., 
$ \omega^* := \sigma^* \circ \circ \imath_{\rB(2, k)} \circ \imath_2'$. 
\end{proof}

\subsubsection{$\rm (II)^\sharp$}

Assume $p = 3$. 
For all integer $e_1 \geq 0$, we can define a homomorphism 
$\sigma_{{\rm (II)^\sharp}, \; e_1} : \SL(2, k) \to \SL(4, k)$ as 
\begin{align*}
\renewcommand{\arraystretch}{1.5} 
\sigma_{{\rm (II)^\sharp}, \; e_1} 
\left(
\begin{array}{c c}
 a & b \\
 c & d
\end{array}
\right) 
 := 
\left(
\begin{array}{l l | l l}
 a^{3 \, p^{e_1}}  
 & b^{3 \, p^{e_1}} 
 & a^{2 \, p^{e_1}} \, b^{p^{e_1}}  
 & \frac{1}{2} \, a^{p^{e_1}} \, b^{2 \, p^{e_1}} \\
c^{3 \, p^{e_1}} 
 & d^{3 \, p^{e_1}}  
 & c^{2 \, p^{e_1}} \, d^{p^{e_1}}
 & \frac{1}{2} \, c^{p^{e_1}} \, d^{2 \, p^{e_1}}  \\ 
\hline 
0 
 & 0  
 & a^{p^{e_1}} 
 & b^{p^{e_1}} \\
0 
 & 0 
 & c^{p^{e_1}}
 & d^{p^{e_1}} 
\end{array}
\right) . 
\end{align*}

\begin{lem} 
Let $e_1 \geq 0$, 
let $\sigma^*$ be the homomorphism given in $\rm (II)^*$ 
and let $P := P_{3,\, 4} \, P_{2,\, 3} \in \GL(4, k)$. 
Then $\Inn_P \circ \sigma^* = \sigma_{{\rm (II)^\sharp}, \; e_1}$. 
\end{lem}

\begin{proof}
The proof is straightforward. 
\end{proof}

We can define a homomorphism 
$\omega_{{\rm (II)^\sharp}, \; e_1)} : \G_m \to \SL(4, k)$ as 
\[
\omega_{{\rm (II)^\sharp}, \; e_1}
:= 
\sigma_{{\rm (II)^\sharp}, \; e_1} \circ \imath_{\rB(2, k)} \circ \imath_2'. 
\] 

\begin{lem}
We have 
\[
\omega_{{\rm (II)^\sharp}, \; e_1} 
\; \sim \; 
\omega_{p^{e_1 + 1}, \; \; p^{e_1}, \; \; p^{ - e_1}, \; \;  p^{- e_1 - 1}} . 
\]
\end{lem}

\begin{proof}
The proof is straightforward. 
\end{proof}

\subsubsection{$\rm (IV)^\sharp$}

For all integers $e_1$ and $e_2$ satisfying  
\[ 
 e_2 > e_1 \geq 0 , 
\]
we can define a homomorphism 
$\sigma_{{\rm (IV)^\sharp}, \; (e_1, \, e_2)} : \SL(2, k) \to \SL(4, k)$ as 
\begin{align*}
\renewcommand{\arraystretch}{1.5} 
\sigma_{{\rm (IV)^\sharp}, \; (e_1, \, e_2)} 
\left(
\begin{array}{c c}
 a & b \\
 c & d
\end{array}
\right) 
& := 
\renewcommand{\arraystretch}{1.5} 
\left(
\begin{array}{l l l l}
 a^{p^{e_2}} \cdot a^{p^{e_1}} 
 & a^{p^{e_2}} \cdot b^{p^{e_1}} 
 & b^{p^{e_2}} \cdot a^{p^{e_1}} 
 & b^{p^{e_2}} \cdot b^{p^{e_1}} \\
  a^{p^{e_2}} \cdot c^{p^{e_1}} 
 & a^{p^{e_2}} \cdot d^{p^{e_1}} 
 & b^{p^{e_2}} \cdot c^{p^{e_1}} 
 & b^{p^{e_2}} \cdot d^{p^{e_1}} \\
 c^{p^{e_2}} \cdot a^{p^{e_1}} 
 & c^{p^{e_2}} \cdot b^{p^{e_1}} 
 & d^{p^{e_2}} \cdot a^{p^{e_1}} 
 & d^{p^{e_2}} \cdot b^{p^{e_1}} \\
 c^{p^{e_2}} \cdot c^{p^{e_1}} 
 & c^{p^{e_2}} \cdot d^{p^{e_1}} 
 & d^{p^{e_2}} \cdot c^{p^{e_1}} 
 & d^{p^{e_2}} \cdot d^{p^{e_1}}
\end{array}
\right) \\
& = 
\renewcommand{\arraystretch}{1.5} 
\left(
\begin{array}{c c}
a^{p^{e_2}}
\left(
\begin{array}{c c}
 a^{p^{e_1}} & b^{p^{e_1}} \\
 c^{p^{e_1}} & d^{p^{e_1}} 
\end{array}
\right)
& 
b^{p^{e_2}}
\left(
\begin{array}{c c}
 a^{p^{e_1}} & b^{p^{e_1}} \\
 c^{p^{e_1}} & d^{p^{e_1}} 
\end{array}
\right) \\ [1.5em] 
c^{p^{e_2}}
\left(
\begin{array}{c c}
 a^{p^{e_1}} & b^{p^{e_1}} \\
 c^{p^{e_1}} & d^{p^{e_1}} 
\end{array}
\right)
& 
d^{p^{e_2}}
\left(
\begin{array}{c c}
 a^{p^{e_1}} & b^{p^{e_1}} \\
 c^{p^{e_1}} & d^{p^{e_1}} 
\end{array}
\right)
\end{array}
\right) . 
\end{align*}

\begin{lem} 
Let $e_1$, $e_2$ be integers satisfying $e_2 > e_1 \geq 0$ 
and let $\sigma^*$ be the homomorphism given in $\rm (IV)^*$. 
Then $\sigma^* = \sigma_{{\rm (IV)^\sharp}, \; (e_1, \, e_2)}$. 
\end{lem}

\begin{proof}
The proof is straightforward. 
\end{proof}

We can define a homomorphism 
$\omega_{{\rm (IV)^\sharp}, \; (e_1, \, e_2)} : \SL(2, k) \to \SL(4, k)$ as 
\[
 \omega_{{\rm (IV)^\sharp}, \; (e_1, \, e_2)} 
= 
\sigma_{{\rm (IV)^\sharp}, \; (e_1, \, e_2)} \circ \imath_{\rB(2, k)} \circ \imath_2'. 
\]

\begin{lem}
We have 
\[
\omega_{{\rm (IV)^\sharp}, \; (e_1, \, e_2)}  
\; \sim \; 
\omega_{p^{e_1} + p^{e_2}, 
\; \;  p^{e_2} - p^{e_1}, 
\; \;   - (p^{e_2} - p^{e_1}), 
\; \;  - ( p^{e_1} + p^{e_2})} . 
\]
\end{lem}

\begin{proof}
The proof is straightforward. 
\end{proof}

\subsubsection{$\rm (V)^\sharp$}

Assume $p = 2$. 
For all integer $e_1 \geq 0$, we can define a homomorphism 
$\sigma_{{\rm (V)^\sharp}, \; e_1} : \SL(2, k) \to \SL(4, k)$ as 
\begin{align*}
\renewcommand{\arraystretch}{1.5} 
\sigma_{{\rm (V)^\sharp}, \; e_1} 
\left(
\begin{array}{c c}
 a & b \\
 c & d
\end{array}
\right) 
 := 
\left(
\begin{array}{l l l l}
1 
 & 0 
 & 0 
 & 0 \\
a^{p^{e_1}} \, b^{p^{e_1}} 
 & a^{p^{e_1 + 1}}   
 & b^{p^{e_1 + 1}}
 & 0  \\ 
c^{p^{e_1}} \, d^{p^{e_1}} 
 & c^{p^{e_1 + 1}} 
 & d^{p^{e_1 + 1}}
 & 0  \\
b^{p^{e_1}} \, c^{p^{e_1}} 
 & a^{p^{e_1}} \, c^{p^{e_1}} 
 & b^{p^{e_1}} \, d^{p^{e_1}} 
 & 1 
\end{array}
\right) . 
\end{align*}

\begin{lem} 
Let $e_1 \geq 0$, 
let $\sigma^*$ be the homomorphism given in $\rm (V)^*$ 
and let $P := P_{3, \, 4} \, P_{1, \, 2} \in \GL(4, k)$. 
Then $\Inn_P \circ \sigma^* = \sigma_{{\rm (V)^\sharp}, \; e_1}$. 
\end{lem}

\begin{proof}
The proof is straightforward. 
\end{proof}

We can define a homomorphism 
$\omega_{{\rm (V)^\sharp}, \; e_1} : \G_m \to \SL(4, k)$ as 
\[
\omega_{{\rm (V)^\sharp}, \; e_1} 
= 
\sigma_{{\rm (V)^\sharp}, \; e_1} \circ \imath_{\rB(2, k)} \circ \imath_2'. 
\]

\begin{lem}
We have 
\[
\omega_{{\rm (V)^\sharp}, \; e_1} 
\; \sim \; 
\omega_{p^{e_1 + 1}, 
\; \; 0, 
\; \; 0,
\; \;  - p^{e_1 + 1}} . 
\]
\end{lem}

\begin{proof}
The proof is straightforward. 
\end{proof}

\subsubsection{$\rm (VII)^\sharp$}

Assume $p = 3$. 
For all integer $e_1 \geq 0$, we can define a homomorphism 
$\sigma_{{\rm (VII)^\sharp}, \; e_1} : \SL(2, k) \to \SL(4, k)$ as 
\begin{align*}
\renewcommand{\arraystretch}{1.5} 
\sigma_{{\rm (VII)^\sharp}, \; e_1} 
\left(
\begin{array}{c c}
 a & b \\
 c & d
\end{array}
\right) 
 := 
\left( 
\begin{array}{c c | c c}
 a^{p^{e_1}} 
 & b^{p^{e_1}} 
 & \frac{1}{2} \, a^{2 \, p^{e_1}} \, c^{p^{e_1}}
 & \frac{1}{2} \, b^{2 \, p^{e_1}} \, d^{p^{e_1}} \\
 c^{p^{e_1}}
  & d^{p^{e_1}} 
  & a^{p^{e_1}} \, c^{2 \, p^{e_1}} 
  & b^{p^{e_1}} \, d^{2 \, p^{e_1}} \\
\hline 
0 & 0 & a^{3 \, p^{e_1}} & b^{3 \, p^{e_1}} \\
0 & 0 & c^{3 \, p^{e_1}} & d^{3 \, p^{e_1}} 
\end{array}
\right) . 
\end{align*}

\begin{lem} 
Let $e_1 \geq 0$, 
let $\sigma^*$ be the homomorphism given in $\rm (VII)^*$. 
The following assertions {\rm (1)} and {\rm (2)} hold true: 
\begin{enumerate}[label = {\rm (\arabic*)}]
\item Letting $P := P_{1,\, 2} \, P_{2,\, 3} \in \GL(4, k)$, we have 
$\Inn_P \circ \sigma^* = \sigma_{{\rm (VII)^\sharp}, \; e_1}$. 

\item ${^\tau}( \sigma_{\rm (VII)^\sharp, \; e_1} )^\tau = \sigma_{\rm (II)^\sharp, \; e_1}$. 
\end{enumerate} 
\end{lem}

\begin{proof}
The proof is straightforward. 
\end{proof}

We can define a homomorphism 
$\omega_{{\rm (VII)^\sharp}, \; e_1} : \G_m \to \SL(4, k)$ as 
\[
\omega_{{\rm (VII)^\sharp}, \; e_1} 
= 
\sigma_{{\rm (VII)^\sharp}, \; e_1} \circ \imath_{\rB(2, k)} \circ \imath_2'. 
\]

\begin{lem}
We have 
\[
\omega_{{\rm (VII)^\sharp}, \; e_1} 
\; \sim \; 
\omega_{p^{e_1 + 1}, 
\; \; p^{e_1}, 
\; \; - p^{e_1},
\; \; - p^{e_1 + 1}} . 
\]
\end{lem}

\begin{proof}
The proof is straightforward. 
\end{proof}

\subsubsection{$\rm (IX)^\sharp$}

Assume $p \geq 3$. 
For all integer $e_1 \geq 0$, we can define a homomorphism 
$\sigma_{{\rm (IX)^\sharp}, \; e_1} : \SL(2, k) \to \SL(4, k)$ as 
\begin{align*}
\renewcommand{\arraystretch}{1.5} 
\sigma_{{\rm (IX)^\sharp}, \; e_1} 
\left(
\begin{array}{c c}
 a & b \\
 c & d 
\end{array}
\right) 
 := 
\left(
\begin{array}{l l l | l}
 a^{2 \, p^{e_1}} 
 & a^{p^{e_1}} \, b^{p^{e_1}} 
 & b^{2 \, p^{e_1}}  
 & 0  \\
2 \, a^{p^{e_1}} \, c^{p^{e_1}} 
 & a^{p^{e_1}} \, d^{p^{e_1}} + b^{p^{e_1}} \, c^{p^{e_1}}    
 & 2 \, b^{p^{e_1}} \, d^{p^{e_1}}   
 & 0  \\ 
c^{2 \, p^{e_1}} 
 & c^{p^{e_1}} \, d^{p^{e_1}} 
 & d^{2 \, p^{e_1}} 
 & 0 \\ 
\hline 
0 
 & 0 
 & 0 
 & 1 
\end{array}
\right) . 
\end{align*}

\begin{lem} 
Let $e_1 \geq 0$, 
let $\sigma^*$ be the homomorphism given in $\rm (IX)^*$ 
and let $P := P_{3,\, 4} \in \GL(4, k)$. 
Then $\Inn_P \circ \sigma^* = \sigma_{{\rm (IX)^\sharp}, \; e_1}$. 
\end{lem}

\begin{proof}
The proof is straightforward. 
\end{proof}

We can define a homomorphism 
$\omega_{{\rm (IX)^\sharp}, \; e_1} : \G_m \to \SL(4, k)$ as 
\[
\omega_{{\rm (IX)^\sharp}, \; e_1} 
= 
\sigma_{{\rm (IX)^\sharp}, \; e_1} \circ \imath_{\rB(2, k)} \circ \imath_2'. 
\]

\begin{lem}
We have 
\[
\omega_{{\rm (IX)^\sharp}, \; e_1} 
\; \sim \; 
\omega_{2 \, p^{e_1}, 
\; \; 0, 
\; \; 0,
\; \; - 2 \, p^{e_1}} . 
\]
\end{lem}

\begin{proof}
The proof is straightforward. 
\end{proof}

\subsubsection{$\rm (XI)^\sharp$}

Assume $p = 2$. For all integer $e_1 \geq 0$, we can define a homomorphism 
$\sigma_{{\rm (XI)^\sharp}, \; e_1} : \SL(2, k) \to \SL(4, k)$ as 
\begin{align*}
\renewcommand{\arraystretch}{1.5} 
\sigma_{{\rm (XI)^\sharp}, \; e_1} 
\left(
\begin{array}{c c}
 a & b \\
 c & d
\end{array}
\right) 
 := 
\left(
\begin{array}{l l l | l}
 a^{2 \, p^{e_1}}  
 & b^{2 \, p^{e_1}} 
 & a^{p^{e_1}} \, b^{p^{e_1}} 
 & 0  \\
c^{2 \, p^{e_1}}  
 & d^{2 \, p^{e_1}}  
 & c^{p^{e_1}} \, d^{p^{e_1}}  
 & 0 \\ 
0 
 & 0  
 & 1 
 & 0 \\
\hline 
0 
 & 0 
 & 0
 & 1
\end{array}
\right) . 
\end{align*}

\begin{lem} 
Let $e_1 \geq 0$, 
let $\sigma^*$ be the homomorphism given in $\rm (XI)^*$ 
and let $P := P_{3,\, 4} \,  P_{2,\, 3} \in \GL(4, k)$. 
Then $\Inn_P \circ \sigma^* = \sigma_{{\rm (XI)^\sharp}, \; e_1}$. 
\end{lem}

\begin{proof}
The proof is straightforward. 
\end{proof}

We can define a homomorphism 
$\omega_{{\rm (XI)^\sharp}, \; e_1} : \G_m \to \SL(4, k)$ as 
\[
\omega_{{\rm (XI)^\sharp}, \; e_1} 
= 
\sigma_{{\rm (XI)^\sharp}, \; e_1} \circ \imath_{\rB(2, k)} \circ \imath_2'. 
\]

\begin{lem}
We have 
\[
\omega_{{\rm (XI)^\sharp}, \; e_1} 
\; \sim \; 
\omega_{2 \, p^{e_1}, 
\; \; 0, 
\; \; 0,
\; \; - 2 \, p^{e_1}} . 
\]
\end{lem}

\begin{proof}
The proof is straightforward. 
\end{proof}

\subsubsection{$\rm (XV)^\sharp$}

For all integers $e_2$ and $e_3$ satisfying 
\[
 e_2 \geq e_3 \geq 0 , 
\]
we can define a homomorphism 
$\sigma_{{\rm (XV)^\sharp}, \; (e_2, \, e_3)} : \SL(2, k) \to \SL(4, k)$ as 
\begin{align*}
\renewcommand{\arraystretch}{1.5} 
\sigma_{{\rm (XV)^\sharp}, \; (e_2, \, e_3)} 
\left(
\begin{array}{c c}
 a & b \\
 c & d
\end{array}
\right) 
  := 
\left( 
\begin{array}{c c | c c}
a^{p^{e_2}} & b^{p^{e_2}} & 0 & 0\\
c^{p^{e_2}} & d^{p^{e_2}} & 0 & 0 \\
\hline 
 0 & 0 & a^{p^{e_3}} & b^{p^{e_3}} \\
 0 & 0 & c^{p^{e_3}} & d^{p^{e_3}}  
\end{array}
\right) . 
\end{align*}

\begin{lem} 
Let $e_2 \geq e_3 \geq 0$, 
let $\sigma^*$ be the homomorphism given in $\rm (XV)^*$. 
Then the following assertions {\rm (1)} and {\rm (2)} hold true: 
\begin{enumerate}[label = {\rm (\arabic*)}]
\item  Letting $P := P_{3,\, 4} \,  P_{2,\, 3} \in \GL(4, k)$, 
we have $\Inn_P \circ \sigma^* = \sigma_{{\rm (XV)^\sharp}, \; (e_2, \, e_3)} $. 

\item Letting 
\[
\renewcommand{\arraystretch}{1.5}
 Q := \left(
\begin{array}{c | c}
 O_2 & I_2 \\
\hline 
 I_2 & O_2
\end{array}
\right) \in \GL(4, k) , 
\]
we have 
$\Inn_Q \circ {^\tau}( \sigma_{{\rm (XV)^\sharp}, \; (e_2, \, e_3)} )^\tau 
 = \sigma_{{\rm (XV)^\sharp}, \; (e_2, \, e_3)}$. 
\end{enumerate} 
\end{lem}

\begin{proof}
The proof is straightforward. 
\end{proof}

We can define a homomorphism 
$\omega_{{\rm (XV)^\sharp}, \; (e_2, \, e_3)} : \G_m \to \SL(4, k)$ as 
\[
\omega_{{\rm (XV)^\sharp}, \; (e_2, \, e_3)} 
= 
\sigma_{{\rm (XV)^\sharp}, \; (e_2 , \ , e_3)} \circ \imath_{\rB(2, k)} \circ \imath_2'. 
\]

\begin{lem}
We have 
\[
\omega_{{\rm (XV)^\sharp}, \; (e_2, \, e_3)} 
\; \sim \; 
\omega_{p^{e_2}, 
\; \; p^{e_3}, 
\; \; - p^{e_3},
\; \; - p^{e_2}} . 
\]
\end{lem}

\begin{proof}
The proof is straightforward. 
\end{proof}

\subsubsection{$\rm (XIX)^\sharp$}

Assume $p = 2$. 
For all integer $e_1 \geq 0$, we can define a homomorphism 
$\sigma_{{\rm (XIX)^\sharp}, \; e_1} : \SL(2, k) \to \SL(4, k)$ as 
\begin{align*}
\renewcommand{\arraystretch}{1.5} 
\sigma_{{\rm (XIX)^\sharp}, \; e_1}
\left(
\begin{array}{c c}
 a & b \\
 c & d
\end{array}
\right) 
 := 
\left(
\begin{array}{l l l | l}
 a^{2 \, p^{e_1}} 
 & b^{2 \, p^{e_1}} 
 & 0 
 & 0  \\
 c^{2 \, p^{e_1}} 
 & d^{2 \, p^{e_1}} 
 & 0
 & 0 \\ 
a^{p^{e_1}} \, c^{p^{e_1}} 
 & b^{p^{e_1}} \, d^{p^{e_1}}
 & 1 
 & 0 \\
\hline 
0 
 & 0 
 & 0
 & 1
\end{array}
\right) . 
\end{align*}

\begin{lem} 
Let $e_1 \geq 0$, 
let $\sigma^*$ be the homomorphism given in $\rm (XIX)^*$ 
and let $P := P_{2,\, 3} \, P_{3,\, 4} \,  P_{2,\, 3} \in \GL(4, k)$. 
Then $\Inn_P \circ \sigma^* = \sigma_{{\rm (XIX)^\sharp}, \; e_1}$. 
\end{lem}

\begin{proof}
The proof is straightforward. 
\end{proof}

We can define a homomorphism 
$\omega_{{\rm (XIX)^\sharp}, \; e_1} : \G_m \to \SL(4, k)$ as 
\[
\omega_{{\rm (XIX)^\sharp}, \; e_1} 
= 
\sigma_{{\rm (XIX)^\sharp}, \; e_1} \circ \imath_{\rB(2, k)} \circ \imath_2'. 
\]

\begin{lem}
We have 
\[
\omega_{{\rm (XIX)^\sharp}, \; e_1} 
\; \sim \; 
\omega_{2 \, p^{e_1}, 
\; \; 0, 
\; \; 0,
\; \; - 2 \, p^{e_1}} . 
\]
\end{lem}

\begin{proof}
The proof is straightforward. 
\end{proof}

\subsubsection{$\rm (XXIV)^\sharp$}

For all integer $e_2 \geq 0$, we can define a homomorphism 
$\sigma_{{\rm (XXIV)^\sharp}, \; e_2} : \SL(2, k) \to \SL(4, k)$ as 
\begin{align*}
\renewcommand{\arraystretch}{1.5} 
\sigma_{{\rm (XXIV)^\sharp}, \; e_2} 
\left(
\begin{array}{c c}
 a & b \\
 c & d
\end{array}
\right) 
 := 
\left(
\begin{array}{l l | l | l}
 a^{p^{e_2}} 
 & b^{p^{e_2}} 
 & 0 
 & 0 \\
 c^{p^{e_2}} 
 & d^{p^{e_2}}
 & 0 
 & 0 \\ 
\hline 
 0 
 & 0 
 & 1 
 & 0  \\
\hline 
 0 
 & 0 
 & 0
 & 1  
\end{array}
\right) . 
\end{align*}

\begin{lem} 
Let $e_1 \geq 0$, 
let $\sigma^*$ be the homomorphism given in $\rm (XXIV)^*$ 
and let $P := P_{3,\, 4} \,  P_{2,\, 3} \in \GL(4, k)$. 
Then $\Inn_P \circ \sigma^* = \sigma_{{\rm (XXIV)^\sharp}, \; e_1}$. 
\end{lem}

\begin{proof}
The proof is straightforward. 
\end{proof}

We can define a homomorphism 
$\omega_{{\rm (XXIV)^\sharp}, \; e_2} : \G_m \to \SL(4, k)$ as 
\[
\omega_{{\rm (XXIV)^\sharp}, \; e_2} 
= 
\sigma_{{\rm (XXIV)^\sharp}, \; e_2} \circ \imath_{\rB(2, k)} \circ \imath_2'. 
\]

\begin{lem}
We have 
\[
\omega_{{\rm (XXIV)^\sharp}, \; e_2} 
\; \sim \; 
\omega_{p^{e_2}, 
\; \; 0, 
\; \; 0,
\; \; - p^{e_2}} . 
\]
\end{lem}

\begin{proof}
The proof is straightforward. 
\end{proof}

\subsubsection{$\rm (XXVI)^\sharp$}

We can define a homomorphism 
$\sigma_{{\rm (XXVI)^\sharp}} : \SL(2, k) \to \SL(4, k)$ as 
\begin{align*}
\renewcommand{\arraystretch}{1.5} 
\sigma_{{\rm (XXVI)^\sharp}} 
\left(
\begin{array}{c c}
 a & b \\
 c & d
\end{array}
\right) 
 := I_4 .
\end{align*}

\begin{lem}
Let $\sigma^*$ be the homomorphism given in ${\rm (XXVI)^*}$. 
Then $\sigma_{\rm (XXVI)^*} = \sigma_{{\rm (XXVI)^\sharp}} $. 
\end{lem}

\begin{proof}
The proof is straightforward. 
\end{proof}

We can define a homomorphism 
$\omega_{{\rm (XXVI)^\sharp}} : \G_m \to \SL(4, k)$ as 
\[
\omega_{{\rm (XXVI)^\sharp}} 
= 
\sigma_{{\rm (XXVI)^\sharp}} \circ \imath_{\rB(2, k)} \circ \imath_2'. 
\]

\begin{lem}
We have 
\[
\omega_{{\rm (XXVI)^\sharp}} 
=  
\omega_{0, 
\; \; 0, 
\; \; 0,
\; \; 0 } . 
\]
\end{lem}

\begin{proof}
The proof is straightforward. 
\end{proof}

\subsection{Equivalences of homomorphisms from $\SL(2, k)$ to $\SL(4, k)$}

\begin{lem}
The following assertions {\rm (1)}, {\rm (2)}, {\rm (3)} hold true: 
\begin{enumerate}[label = {\rm (\arabic*)}, itemsep = 12pt] 
\item In the case where $p = 2$, we have the following: 
\begin{enumerate}[label = {\rm (\roman*)}, leftmargin = 5em]
\item[$\rm (IV)^\sharp$] For all integers $e_1$ and $e_2$ satisfying $e_2 > e_1 \geq 0$, we have 
\begin{fleqn}[10em] 
\begin{align*}
 d\bigl( \, \sigma_{{\rm (IV)^\sharp}, \; (e_1, \, e_2) } \, \bigr) 
= 
(0, 0) . 
\end{align*}
\end{fleqn}

\item[$\rm (V)^\sharp$] For all integer $e_1 \geq 0$, we have 
\begin{fleqn}[10em] 
\begin{align*}
 d\bigl( \, \sigma_{{\rm (V)^\sharp}, \; e_1 } \, \bigr) 
= 
(1, 1) . 
\end{align*}
\end{fleqn}

\item[$\rm (XI)^\sharp$] For all integer $e_1 \geq 0$, we have 
\begin{fleqn}[10em] 
\begin{align*}
 d\bigl( \, \sigma_{{\rm (XI)^\sharp}, \; e_1 } \, \bigr) 
= 
(1, 2) . 
\end{align*}
\end{fleqn}

\item[$\rm (XV)^\sharp$] For all integers $e_2$ and $e_3$ satisfying $e_2 \geq e_3 \geq 0$, we have 
\begin{fleqn}[10em] 
\begin{align*}
 d\bigl( \, \sigma_{{\rm (XV)^\sharp}, \; (e_2, \, e_3) } \, \bigr) 
= 
(0, 0) . 
\end{align*}
\end{fleqn}

\item[$\rm (XIX)^\sharp$] For all integer $e_1 \geq 0$, we have 
\begin{fleqn}[10em] 
\begin{align*}
 d\bigl( \, \sigma_{{\rm (XIX)^\sharp}, \; e_1 } \, \bigr) 
= 
(2, 1) . 
\end{align*}
\end{fleqn}

\item[$\rm (XXIV)^\sharp$] For all integer $e_2 \geq 0$, we have 
\begin{fleqn}[10em] 
\begin{align*}
 d\bigl( \, \sigma_{{\rm (XXIV)^\sharp}, \; e_2 } \, \bigr) 
= 
(2, 2) . 
\end{align*}
\end{fleqn}

\item[$\rm (XXVI)^\sharp$] We have 
\begin{fleqn}[10em] 
\begin{align*}
 d\bigl( \, \sigma_{{\rm (XXVI)^\sharp} } \, \bigr) 
= 
(4, 4) . 
\end{align*}
\end{fleqn}
\end{enumerate}

\item In the case where $p = 3$, we have the following: 

\begin{enumerate}[label = {\rm (\roman*)}, leftmargin = 5em]
\item[$\rm (II)^\sharp$] For all integer $e_1 \geq 0$, we have 
\begin{fleqn}[10em] 
\begin{align*}
 d\bigl( \, \sigma_{{\rm (II)^\sharp}, \; e_1 } \, \bigr) 
= 
(0, 0) . 
\end{align*}
\end{fleqn}

\item[$\rm (IV)^\sharp$] For all integers $e_1$ and $e_2$ satisfying $e_2 > e_1 \geq 0$, we have 
\begin{fleqn}[10em] 
\begin{align*}
 d\bigl( \, \sigma_{{\rm (IV)^\sharp}, \; (e_1, \, e_2) } \, \bigr) 
= 
(0, 0) . 
\end{align*}
\end{fleqn}

\item[$\rm (VII)^\sharp$] For all integer $e_1 \geq 0$, we have 
\begin{fleqn}[10em] 
\begin{align*}
 d\bigl( \, \sigma_{{\rm (VII)^\sharp}, \; e_1 } \, \bigr) 
= 
(0, 0) . 
\end{align*}
\end{fleqn}

\item[$\rm (IX)^\sharp$] For all integer $e_1 \geq 0$, we have 
\begin{fleqn}[10em] 
\begin{align*}
 d\bigl( \, \sigma_{{\rm (IX)^\sharp}, \; e_1 } \, \bigr) 
= 
(1, 1) . 
\end{align*}
\end{fleqn}

\item[$\rm (XV)^\sharp$] For all integers $e_2$ and $e_3$ satisfying $e_2 \geq e_3 \geq 0$, we have 
\begin{fleqn}[10em] 
\begin{align*}
 d\bigl( \, \sigma_{{\rm (XV)^\sharp}, \; (e_2, \, e_3) } \, \bigr) 
= 
(0, 0) . 
\end{align*}
\end{fleqn}

\item[$\rm (XXIV)^\sharp$] For all integer $e_2 \geq 0$, we have 
\begin{fleqn}[10em] 
\begin{align*}
 d\bigl( \, \sigma_{{\rm (XXIV)^\sharp}, \; e_2 } \, \bigr) 
= 
(2, 2) . 
\end{align*}
\end{fleqn}

\item[$\rm (XXVI)^\sharp$] We have 
\begin{fleqn}[10em] 
\begin{align*}
 d\bigl( \, \sigma_{{\rm (XXVI)^\sharp} } \, \bigr) 
= 
(4, 4) . 
\end{align*}
\end{fleqn}

\end{enumerate}

\item In the case where $p \geq 5$, we have the following: 

\begin{enumerate}[label = {\rm (\roman*)}, leftmargin = 5em]
\item[$\rm (I)^\sharp$] For all integer $e_1 \geq 0$, we have 
\begin{fleqn}[10em] 
\begin{align*}
 d\bigl( \, \sigma_{{\rm (I)^\sharp}, \; e_1 } \, \bigr) 
= 
(0, 0) . 
\end{align*}
\end{fleqn}

\item[$\rm (IV)^\sharp$] For all integers $e_1$ and $e_2$ satisfying $e_2 > e_1 \geq 0$, we have 
\begin{fleqn}[10em] 
\begin{align*}
 d\bigl( \, \sigma_{{\rm (IV)^\sharp}, \; (e_1, \, e_2) } \, \bigr) 
= 
(0, 0) . 
\end{align*}
\end{fleqn}

\item[$\rm (IX)^\sharp$] For all integer $e_1 \geq 0$, we have 
\begin{fleqn}[10em] 
\begin{align*}
 d\bigl( \, \sigma_{{\rm (IX)^\sharp}, \; e_1 } \, \bigr) 
= 
(1, 1) . 
\end{align*}
\end{fleqn}

\item[$\rm (XV)^\sharp$] For all integers $e_2$ and $e_3$ satisfying $e_2 \geq e_3 \geq 0$, we have 
\begin{fleqn}[10em] 
\begin{align*}
 d\bigl( \, \sigma_{{\rm (XV)^\sharp}, \; (e_2, \, e_3) } \, \bigr) 
= 
(0, 0) . 
\end{align*}
\end{fleqn}

\item[$\rm (XXIV)^\sharp$] For all integer $e_2 \geq 0$, we have 
\begin{fleqn}[10em] 
\begin{align*}
 d\bigl( \, \sigma_{{\rm (XXIV)^\sharp}, \; e_2 } \, \bigr) 
= 
(2, 2) . 
\end{align*}
\end{fleqn}

\item[$\rm (XXVI)^\sharp$] We have 
\begin{fleqn}[10em] 
\begin{align*}
 d\bigl( \, \sigma_{{\rm (XXVI)^\sharp} } \, \bigr) 
= 
(4, 4) . 
\end{align*}
\end{fleqn}
\end{enumerate}

\end{enumerate} 
\end{lem}

\begin{proof}
(1) See Lemmas 5.7, \, 5.9, \, 5.15, \, 5.17, \, 5.19, \, 5.23, \, 5.24. 

(2) See Lemmas 5.5, \, 5.7, \, 5.11, \, 5.13, \, 5.17, \, 5.23, \, 5.24. 

(3) See Lemmas 5.3, \, 5.7, \, 5.13, \, 5.17, \, 5.23, \, 5.24. 

\end{proof}

Assume $p = 5$ and for all integer $e_1 \geq 0$,

We can define a homomorphism 
$\omega_{{\rm (IV)^\sharp}, \; (e_1', \, e_2')} : \G_m \to \SL(4, k)$ be the homomorphism 
defined by  
\[
\omega_{{\rm (IV)^\sharp}, \; (e_1', \, e_2')}(u) 
:= 
\sigma_{{\rm (IV)^\sharp}, \; (e_1', \, e_2')}
\left(
\begin{array}{c c}
 u & 0 \\
 0 & u^{-1}
\end{array}
\right) . 
\] 

\begin{lem}
The following assertions {\rm (1)}, {\rm (2)}, {\rm (3)} hold true: 
\begin{enumerate}[label = {\rm (\arabic*)}]
\item In the case where $p = 2$, for all integers $e_1$, $e_2$, $e_2'$, $e_3'$ satisfying 
$e_2 > e_1 \geq 0$ and $e_2' \geq e_3' \geq 0$, 
the homomorphisms  
$\sigma_{{\rm (IV)^\sharp}, \; (e_1, \, e_2)}$ and 
$\sigma_{{\rm (XV)^\sharp}, \; (e_2', \, e_3')}$ are not 
equivalent.

\item In the case where $p = 3$, we have the following: 
\begin{enumerate}[label = {\rm (\roman*)}] 
\item For all integers $e_1$, $e_1'$, $e_2'$ satisfying 
$e_1 \geq 0$ and $e_2' > e_1' \geq 0$, 
the homomorphisms 
$\sigma_{{\rm (II)^\sharp}, \; e_1}$ and $\sigma_{{\rm (IV)^\sharp}, \; (e_1', \, e_2')}$ 
are not equivalent.

\item For all integers $e_1$, $e_1'$ satisfying $e_1 \geq 0$ and $e_1' \geq 0$, 
the homomorphisms 
$\sigma_{{\rm (II)^\sharp}, \; e_1}$ and $\sigma_{{\rm (VII)^\sharp}, \; e_1'}$ are not 
equivalent.

\item For all integer $e_1$, $e_2$, $e_3$ satisfying 
$e_1 \geq 0$ and $e_2 \geq e_3 \geq 0$, 
the homomorphisms 
$\sigma_{{\rm (II)^\sharp}, \; e_1}$ and $\sigma_{{\rm (XV)^\sharp}, \; (e_2, \, e_3)}$ 
are not equivalent.

\item For all integers $e_1$, $e_2$, $e_1'$ satisfying 
$e_2 > e_1 \geq 0$ and $e_1' \geq 0$, 
the homomorphisms  
$\sigma_{{\rm (IV)^\sharp}, \; (e_1, \, e_2)}$ and 
$\sigma_{{\rm (VII)^\sharp}, \; e_1'}$ are not 
equivalent.

\item For all integers $e_1$, $e_2$, $e_2'$, $e_3'$ satisfying 
$e_2 > e_1 \geq 0$ and $e_2' \geq e_3' \geq 0$, 
the homomorphisms  
$\sigma_{{\rm (IV)^\sharp}, \; (e_1, \, e_2)}$ and 
$\sigma_{{\rm (XV)^\sharp}, \; (e_2', \, e_3')}$ are not 
equivalent.

\item For all integers $e_1$, $e_2$, $e_3$ satisfying 
$e_1 \geq 0$ and $e_2 \geq e_3 \geq 0$, 
the homomorphisms  
$\sigma_{{\rm (VII)^\sharp}, \; e_1}$ and 
$\sigma_{{\rm (XV)^\sharp}, \; (e_2, \, e_3)}$ are not 
equivalent. 
\end{enumerate}

\item  In the case where $p \geq 5$, we have the following:
\begin{enumerate}[label = {\rm (\roman*)}]
\item 
For all integers $e_1$, $e_1'$, $e_2'$ satisfying 
$e_1 \geq 0$ and $e_2' > e_1' \geq 0$, 
the homomorphisms 
$\sigma_{{\rm (I)^\sharp}, \; e_1}$ and $\sigma_{{\rm (IV)^\sharp}, \; (e_1', \, e_2')}$ 
are not equivalent.

\item For all integer $e_1$, $e_2$, $e_3$ satisfying 
$e_1 \geq 0$ and $e_2 \geq e_3 \geq 0$, 
the homomorphisms 
$\sigma_{{\rm (I)^\sharp}, \; e_1}$ and $\sigma_{{\rm (XV)^\sharp}, \; (e_2, \, e_3)}$ 
are not equivalent.

\item For all integer $e_1$, $e_2$, $e_2'$, $e_3'$ satisfying 
$e_2 > e_1 \geq 0$ and $e_2' \geq e_3' \geq 0$, 
the homomorphisms $\sigma_{{\rm (IV)^\sharp}, \; (e_1, \, e_2)}$ 
and $\sigma_{{\rm (XV)^\sharp}, \; (e_2', \, e_3')}$ are not equivalent. 
\end{enumerate} 
\end{enumerate}
\end{lem}

\begin{proof}
(1) Suppose to the contrary that 
there exist integers $e_1$, $e_2$, $e_2'$, $e_3'$ satisfying 
\[
e_2 > e_1 \geq 0 , \qquad e_2' \geq e_3' \geq 0 
\]
and the homomorphisms  
$\sigma_{{\rm (IV)^\sharp}, \; (e_1, \, e_2)}$ and 
$\sigma_{{\rm (XV)^\sharp}, \; (e_2', \, e_3')}$ are equivalent. 
So, the homomorphisms 
$\omega_{{\rm (IV)^\sharp}, \; (e_1, \, e_2)}$ and 
$\omega_{{\rm (XV)^\sharp}, \; (e_2', \, e_3')}$ are equivalent. 
By Lemmas 6.6 and 6.16, we have $(p^{e_1} + p^{e_2}, \, p^{e_2} - p^{e_1}) = (p^{e_2'}, \, p^{e_3'})$. 
Thus we have 
\begin{align*}
\left\{
 \begin{array}{r @{\,} l @{\qquad} r} 
 2 \, p^{e_2} & = p^{e_2'} + p^{e_3'} & \textcircled{\scriptsize 1} \\ 
 2 \, p^{e_1} & = p^{e_2'} - p^{e_3'} & \textcircled{\scriptsize 2}
\end{array}
\right.
\end{align*}
By $\textcircled{\scriptsize 1}$, we have $2 \, p^{e_2} = p^{e_3'} \, ( p^{e_2' - e_3'} + 1 )$. 
Note that $e_2' = e_3'$. In fact, suppose to the contrary that $e_2' > e_3'$, 
then $e_3' \geq e_2$ and $2 = p^{e_3' - e_2} \, ( p^{e_2' - e_3'} + 1 )$, 
which implies $p^{e_2' - e_3'} + 1 \geq 3$. 
This is a contradiction. 
By $\textcircled{\scriptsize 2}$, we have $2 \, p^{e_1} = 0$. 
This is a contradiction.

(2) (i) Suppose to the contrary that there exist integers $e_1$, $e_1'$, $e_2'$ such that 
\[
 e_1 \geq 0, \qquad e_2' > e_1' \geq 0 
\]
and the homomorphisms 
$\sigma_{{\rm (II)^\sharp}, \; e_1}$ and $\sigma_{{\rm (IV)^\sharp}, \; (e_1', \, e_2')}$ 
are equivalent. 
Thus $\omega_{{\rm (II)^\sharp}, \; e_1}$ and $\omega_{{\rm (IV)^\sharp}, \; (e_1', \, e_2')}$ 
are equivalent. By Lemmas 6.4 and 6.6, we have 
\[
\left\{
\begin{array}{c @{\,} l }
  p^{e_1 + 1} & = p^{e_1'} + p^{e_2'} , \\
 p^{e_1} & = p^{e_2'} - p^{e_1'} . 
\end{array}
\right.
\]
Summing the above two equalities, we have 
\[
 4 \, p^{e_1} = 2 \, p^{e_2'} , 
\]
which implies $p = 2$. This contradicts $p = 3$.

(ii) Suppose to the contrary that there exist integers $e_1$ and $e_1'$ such that 
\[
 e_1 \geq 0, \qquad e_1' \geq 0
\]
and the homomorphisms 
$\sigma_{{\rm (II)^\sharp}, \; e_1}$ and $\sigma_{{\rm (VII)^\sharp}, \; e_1'}$ are 
equivalent. 
Then there exists a regular matrix 
\[
\renewcommand{\arraystretch}{1.5} 
P = 
\left(
\begin{array}{c | c}
 P_1 & P_2 \\
\hline 
 P_3 & P_4 
\end{array} 
\right) 
\in \GL(4, k)
\qquad \bigl( \, P_1, P_2, P_3, P_4 \in \Mat(2, k) \, \bigr)
\]
such that 
\begin{align*}
& 
\renewcommand{\arraystretch}{1.5} 
\left(
\begin{array}{c c | c c}
 A^{3 \, p^{e_1}} 
 & B^{3 \, p^{e_1}}  
 & A^{2 \, p^{e_1}} \, B^{p^{e_1}} 
 & \frac{1}{2} \, A^{p^{e_1}} \, B^{2 \, p^{e_1}} \\
 C^{3 \, p^{e_1}} 
 & D^{3 \, p^{e_1}} 
 & C^{2 \, p^{e_1}} \, D^{p^{e_1}} 
 & \frac{1}{2} \, C^{p^{e_1}} \, D^{2 \, p^{e_1}} \\
\hline 
 0 & 0 & A^{p^{e_1}} & B^{p^{e_1}} \\
 0 & 0 & C^{p^{e_1}} & D^{p^{e_1}} 
\end{array}
\right)
\left(
\begin{array}{c | c}
 P_1 & P_2 \\
\hline 
 P_3 & P_4 
\end{array} 
\right) \\
& \qquad = 
\renewcommand{\arraystretch}{1.5} 
\left(
\begin{array}{c | c}
 P_1 & P_2 \\
\hline 
 P_3 & P_4 
\end{array} 
\right)
\left( 
\begin{array}{c c | c c}
 A^{p^{e_1'}} 
 & B^{p^{e_1'}} 
 & \frac{1}{2} \, A^{2 \, p^{e_1'}} \, C^{p^{e_1'}}
 & \frac{1}{2} \, B^{2 \, p^{e_1'}} \, D^{p^{e_1'}} \\
 C^{p^{e_1'}}
  & D^{p^{e_1'}} 
  & A^{p^{e_1'}} \, C^{p^{e_1'}} 
  & B^{p^{e_1'}} \, D^{p^{e_1'}} \\
\hline 
0 & 0 & A^{3 \, p^{e_1'}} & B^{3 \, p^{e_1'}} \\
0 & 0 & C^{3 \, p^{e_1'}} & D^{3 \, p^{e_1'}} 
\end{array}
\right) . 
\end{align*}

Letting 
\[
\left(
\begin{array}{c c}
 A & B \\
 C & D
\end{array}
\right) 
 := 
\left(
\begin{array}{c c}
 1 & 0 \\
 1 & 1 
\end{array} 
\right) 
\in \SL(2, k) 
\]
and comparing the $(1, 1)$-th block and the $(2, 1)$-th block of 
tboth sides of the equality, we have 
\begin{align*}
\left\{ 
\begin{array}{l @{\qquad}  l}
\left(
\begin{array}{c c}
 1 & 0 \\
 1 & 1 
\end{array}
\right) 
P_1 
 + 
 \left( 
 \begin{array}{c c}
  0 & 0 \\
  1 & \frac{1}{2} 
 \end{array}
 \right)
P_3 
 = 
 P_1
\left(
\begin{array}{c c}
1 & 0 \\
1 & 1
\end{array} 
\right) 
& \textcircled{\scriptsize 1} \\  [1.5em] 
 \left( 
 \begin{array}{c c}
  1 & 0 \\
  1 & 1
 \end{array}
 \right)
P_3 
 = 
 P_3
\left(
\begin{array}{c c}
1 & 0 \\
1 & 1
\end{array} 
\right) 
& 
\textcircled{\scriptsize 2} 
\end{array} 
\right. 
\end{align*}

Letting 
\[
\left(
\begin{array}{c c}
 A & B \\
 C & D
\end{array}
\right) 
 := 
\left(
\begin{array}{c c}
 1 & 1 \\
 0 & 1 
\end{array} 
\right) \in \SL(2, k) 
\]
and comparing the $(1, 1)$-th block and the $(2, 1)$-th block of 
tboth sides of the equality, we have 
\begin{align*}
\left\{ 
\begin{array}{l @{\qquad} l}
\left(
\begin{array}{c c}
 1 & 1 \\
 0 & 1 
\end{array}
\right) 
P_1 
 + 
 \left( 
 \begin{array}{c c}
  1 & \frac{1}{2} \\
  0 & 0
 \end{array}
 \right)
P_3 
 = 
 P_1
\left(
\begin{array}{c c}
1 & 1 \\
0 & 1
\end{array} 
\right) 
& 
\textcircled{\scriptsize 3} \\ [1.5em] 
 \left( 
 \begin{array}{c c}
  1 & 1 \\
  0 & 1
 \end{array}
 \right)
P_3 
 = 
 P_3
\left(
\begin{array}{c c}
1 & 1 \\
0 & 1
\end{array} 
\right) 
& 
\textcircled{\scriptsize 4} 
\end{array} 
\right. 
\end{align*}

Write 
\[
 P_1 = 
 \left(
 \begin{array}{c c}
  x & y \\
  z & w
 \end{array}
 \right) , \qquad 
 P_3 = 
 \left( 
 \begin{array}{c c}
  s & t \\
  u & v
 \end{array}
 \right) 
 \qquad 
 (\,  x, y, z, w, s, t, u, v \in k   \,) . 
\]
By $\textcircled{\scriptsize 2}$, we have 
\[
\left(
\begin{array}{c c}
 s & t \\
 s + u & t + v
\end{array}
\right) 
 = 
\left( 
\begin{array}{c c}
 s + t & t \\
 u + v & v
\end{array} 
\right) . 
\]
So, $t = 0$ and $s = v$. 
Thus
\[
 P_3 = 
 \left(
 \begin{array}{c c}
   s & 0 \\
   u & s 
 \end{array}
 \right) . 
\]

By $\textcircled{\scriptsize 4}$, we have 
\[
\left(
\begin{array}{c c}
 s + u & s \\
 u & s 
\end{array}
\right) 
 = 
\left( 
\begin{array}{c c}
 s & s \\
 u & u + s 
\end{array} 
\right) . 
\]
So, $u = 0$. 
Thus
\[
 P_3 = 
 \left(
 \begin{array}{c c}
   s & 0 \\
   0 & s 
 \end{array}
 \right) . 
\]

By $\textcircled{\scriptsize 1}$, we have 
\[
\left(
\begin{array}{c c}
 x & y \\
 x + z + s & y + w + \frac{1}{2} s 
\end{array}
\right) 
 = 
\left( 
\begin{array}{c c}
 x + y & y \\
 z + w & w
\end{array} 
\right) . 
\]
So, $y = 0$, $x + s = w$ and $y + \frac{1}{2} s = 0$. 
Therefore, $s = 0$ and $x = w$.  
Thus
\[
 P_1 = 
 \left(
 \begin{array}{c c}
   x & 0 \\
   z & x 
 \end{array}
 \right) , 
 \qquad  
 P_3 = 
 \left(
 \begin{array}{c c}
   0 & 0 \\
   0 & 0 
 \end{array}
 \right) . 
\]

By $\textcircled{\scriptsize 3}$, we have 
\[
\left(
\begin{array}{c c}
 x + z & x \\
 z & x 
\end{array}
\right) 
 = 
\left( 
\begin{array}{c c}
 x & x \\
 z & z + x
\end{array} 
\right) . 
\]
So, $z = 0$. 
Thus 
\[
 P_1 = 
 \left(
 \begin{array}{c c}
   x & 0 \\
   0 & x 
 \end{array}
 \right) . 
\]

Comparing the $(1, 1)$-th block of both sides of the equality, 
we hafve 
\[
\left(
\begin{array}{c c}
 A^{3 \, p^{e_1}} & B^{3 \, p^{e_1}} \\
 C^{3 \, p^{e_1}} & D^{3 \, p^{e_1}} 
\end{array}
\right)
\left(
\begin{array}{c c}
 x & 0 \\
 0 & x
\end{array} 
\right) 
 = 
 \left(
\begin{array}{c c}
 x & 0 \\
 0 & x
\end{array} 
\right) 
\left(
\begin{array}{c c}
 A^{p^{e_1'}} & B^{p^{e_1'}} \\
 C^{p^{e_1'}} & D^{p^{e_1'}} 
\end{array}
\right) . 
\]
So, $A^{3 \, p^{e_1}} \, x = x \, A^{p^{e_1'}}$. 
If $e_1' \ne e_1 + 1$, we have $x = 0$. 
Thus $P_1 = O$. 
So,    
\[
 P =
 \renewcommand{\arraystretch}{1.5}  
 \left(
 \begin{array}{c | c}
  O & P_2 \\
  \hline 
   O & P_4
 \end{array} 
 \right) . 
\]
This contradicts the condition that $P$ is regular. 

Now, we have $e_1' = e_1 + 1$. 
Comparing the $(2, 2)$-th blocks of both sides of the equality, 
we have 
\[
\left(
\begin{array}{c c}
 A^{p^{e_1}} & B^{p^{e_1}} \\
 C^{p^{e_1}} & D^{p^{e_1}} 
\end{array}
\right)
P_4 
 = 
 P_4 
\left(
\begin{array}{c c}
 A^{3 \, p^{e_1'}} & B^{3 \, p^{e_1'}} \\
 C^{3 \, p^{e_1'}} & D^{3 \, p^{e_1'}} 
\end{array}
\right) .  
\]
Thus 
\[
\left(
\begin{array}{c c}
 A^{p^{e_1}} & B^{p^{e_1}} \\
 C^{p^{e_1}} & D^{p^{e_1}} 
\end{array}
\right)
P_4 
 = 
 P_4 
\left(
\begin{array}{c c}
 A^{p^{e_1 + 2}} & B^{p^{e_1 + 2}} \\
 C^{p^{e_1 + 2}} & D^{p^{e_1 + 2}} 
\end{array}
\right) .  
\]
Write 
\[
 P_4 = 
 \left(
 \begin{array}{c c}
  \alpha & \beta \\
  \gamma & \delta 
 \end{array}
 \right) 
 \qquad (\, \alpha, \beta, \gamma, \delta \in k  \,) . 
\]

Letting 
\[
\left(
\begin{array}{c c}
 A & B \\
 C & D
\end{array}
\right) 
:= 
\left( 
\begin{array}{r r}
 0 & - 1 \\
 1 & 0 
\end{array} 
\right) \in \SL(2, k) ,  
\]
we have 
\[
\left(
\begin{array}{r r}
 - \gamma & - \delta \\
 \alpha & \beta
\end{array}
\right) 
 = 
\left( 
\begin{array}{r r}
 \beta & - \alpha \\
 \delta & - \gamma 
\end{array}
\right) , 
\]
we have $\delta = \alpha$ and $\gamma = - \beta$. 
Thus 
\[
 P_4 
  = 
 \left( 
  \begin{array}{r r}
   \alpha & \beta \\
   - \beta & \alpha 
  \end{array}
 \right) . 
\]

Letting 
\[
\left(
\begin{array}{c c}
 A & B \\
 C & D
\end{array}
\right) 
:= 
\left( 
\begin{array}{c c}
 1 & 1 \\
 0 & 1 
\end{array} 
\right) \in \SL(2, k) ,  
\]
we have 
\[
\left(
\begin{array}{r r}
 \alpha - \beta & \beta + \alpha \\
  - \beta & \alpha
\end{array}
\right) 
 = 
\left( 
\begin{array}{r r}
 \alpha & \alpha + \beta \\
 - \beta & - \beta + \alpha
\end{array}
\right) . 
\]
So, $\beta = 0$ and $\alpha = 0$. 
Thus $P_4 = O$. 
So, 
\[
 P =
 \renewcommand{\arraystretch}{1.5}  
 \left(
 \begin{array}{c | c}
  P_1 & P_2 \\
  \hline 
   O & O 
 \end{array} 
 \right) . 
\]
This contradicts the condition that $P$ is regular.

(iii) Suppose to the contrary that there exist integers $e_1$, $e_2$, $e_3$ satisfying 
\[
e_1 \geq 0, \qquad e_2 \geq e_3 \geq 0
\]
and the homomorphisms  
$\sigma_{{\rm (II)^\sharp}, \; e_1}$ and 
$\sigma_{{\rm (XV)^\sharp}, \; (e_2, \, e_3)}$ are equivalent. 
Let $\varphi_{{\rm (II)^\sharp}, \, e_1} : \G_a \to \SL(4, k)$ and 
$\varphi_{{\rm (XV)^\sharp}, \; (e_2, \, e_3)} : \G_a \to \SL(4, k)$ 
be homomorphisms defined by 
\[
\varphi_{{\rm (II)^\sharp}, \, e_1} (t) 
:= 
\sigma_{{\rm (II)^\sharp}, \; e_1} 
\left(
\begin{array}{c c}
 1 & t \\
 0 & 1 
\end{array}
\right) , 
\qquad 
\varphi_{{\rm (XV)^\sharp}, \; (e_2, \, e_3)} (t) 
 := 
\sigma_{{\rm (XV)^\sharp}, \; (e_2, \, e_3)}
\left(
\begin{array}{c c}
 1 & t \\
 0 & 1 
\end{array}
\right)  . 
\]
Let $V := k^{\oplus 4}$ be the four-dimensional row vector space over $k$, 
let $V^{\varphi_{{\rm (II)^\sharp}, \, e_1} }$ denote the $\varphi_{{\rm (II)^\sharp}, \, e_1} $-fixed 
subspace and 
let $V^{\varphi_{{\rm (XV)^\sharp}, \; (e_2, \, e_3)} }$ 
denote the $\varphi^{\varphi_{{\rm (XV)^\sharp}, \; (e_2, \, e_3)}}$-fixed 
subspace. 
Thus $\dim V^{\varphi_{{\rm (II)^\sharp}, \, e_1} }= \dim V^{\varphi_{{\rm (XV)^\sharp}, \; (e_2, \, e_3)} }$. 
But we can show $\dim V^{\varphi_{{\rm (II)^\sharp}, \, e_1} } = 1$ and 
$\dim V^{\varphi_{{\rm (XV)^\sharp}, \; (e_2, \, e_3)} } = 2$. 
This is a contradiction.

 (iv) Suppose to the contrary that there exist integers $e_1$, $e_2$, $e_1'$ satisfying 
\[
e_2 > e_1 \geq 0, \qquad e_1' \geq 0
\]
and the homomorphisms  
$\sigma_{{\rm (IV)^\sharp}, \; (e_1, \, e_2)}$ and 
$\sigma_{{\rm (VII)^\sharp}, \; e_1'}$ are equivalent. 
So, $\omega_{{\rm (IV)^\sharp}, \; (e_1, \, e_2)}$ and 
$\omega_{{\rm (VII)^\sharp}, \; e_1'}$ are equivalent. 
By Lemmas 6.6 and 6.10, we have 
\[
\left\{
\begin{array}{r @{\,} l}
 p^{e_1} + p^{e_2} & = p^{e_1' + 1} , \\
 p^{e_2} - p^{e_1} & = p^{e_1'} . 
\end{array}
\right.
\]
Summing the above two equalities, we have 
$2 \, p^{e_2} = 4 \, p^{e_1'}$, which implies $p = 2$. 
This contradicts $p = 3$.

(v) See the proof of the above assertion (1).

(vi) Suppose to contrary that there exist integers $e_1$, $e_2$, $e_3$ satisfying 
\[
e_1 \geq 0 , \qquad e_2 \geq e_3 \geq 0 
\]
and the homomorphisms  
$\sigma_{{\rm (VII)^\sharp}, \; e_1}$ and 
$\sigma_{{\rm (XV)^\sharp}, \; (e_2, \, e_3)}$ are equivalent. 
Thus $\sigma_{{\rm (II)^\sharp}, \; e_1}$ and 
$\sigma_{{\rm (XV)^\sharp}, \; (e_2, \, e_3)}$ are equivalent (see Lemmas 1.24, 6.9, 6.15). 
This contradicts (iii).

(3) (i) Suppose to the contrary that there exist integers $e_1$, $e_1'$, $e_2'$ satisfying 
$e_1 \geq 0$ and $e_2' > e_1' \geq 0$, 
the homomorphisms 
$\sigma_{{\rm (I)^\sharp}, \; e_1}$ and $\sigma_{{\rm (IV)^\sharp}, \; (e_1', \, e_2')}$ 
are equivalent. 
So, $\omega_{{\rm (I)^\sharp}, \; e_1}$ and $\omega_{{\rm (IV)^\sharp}, \; (e_1', \, e_2')}$ 
are equivalent. 
By Lemmas 6.2 and 6.6, we have 
\[
\left\{
\begin{array}{c @{\,} l}
 3 \, p^{e_1} & = p^{e_1'} + p^{e_2'} , \\
 p^{e_1} & = p^{e_2'} - p^{e_1'} . 
\end{array}
\right.
\]
Summing the above two equalities, we have 
$4 \, p^{e_1} = 2 \, p^{e_2'}$, which implies $p = 2$. 
This contradicts $p \geq 5$.

(ii) Suppose to contrary that there exist integer $e_1$, $e_2$, $e_3$ satisfying 
\[
e_1 \geq 0 , \qquad e_2 \geq e_3 \geq 0 , 
\]
and the homomorphisms 
$\sigma_{{\rm (I)^\sharp}, \; e_1}$ and $\sigma_{{\rm (XV)^\sharp}, \; (e_2, \, e_3)}$ 
are equivalent. 
So, $\omega_{{\rm (I)^\sharp}, \; e_1}$ and $\omega_{{\rm (XV)^\sharp}, \; (e_2, \, e_3)}$ 
are equivalent. 
By Lemmas 6.2 and 6.16, we have $(3 \, p^{e_1}, \, p^{e_1}) = (p^{e_2}, \, p^{e_3})$. 
Therefore $p = 3$. This contradicts $p \geq 5$.

(iii) See the proof of the above assertion (1).

\end{proof}

\begin{lem}
The following assertions hold true:
\begin{enumerate}
\item[$\rm (I)^\sharp$] 
Let $e_1$ and $e_1'$ be integers satisfying $e_1 \geq 0$ and $e_1' \geq 0$. 
If $e_1 \ne e_1'$, then $\sigma_{{\rm (I)^\sharp}, \; e_1}$ 
and $\sigma_{{\rm (I)^\sharp}, \; e_1'}$ are not equivalent.

\item[$\rm (II)^\sharp$] 
Let $e_1$ and $e_1'$ be integers satisfying $e_1 \geq 0$ and $e_1' \geq 0$. 
If $e_1 \ne e_1'$, then $\sigma_{{\rm (II)^\sharp}, \; e_1}$ 
and $\sigma_{{\rm (II)^\sharp}, \; e_1'}$ are not equivalent.

\item[$\rm (IV)^\sharp$] Let $e_1$, $e_2$, $e_1'$, $e_2'$ be integers 
satisfying $e_2 > e_1 \geq 0$ and $e_2' > e_1' \geq 0$. 
If $(e_1, e_2) \ne (e_1', e_2')$, then $\sigma_{{\rm (IV)^\sharp}, \; (e_1, \, e_2)}$ 
and $\sigma_{{\rm (IV)^\sharp}, \; (e_1', \, e_2')}$ are not equivalent.

\item[$\rm (V)^\sharp$] 
Let $e_1$ and $e_1'$ be integers satisfying $e_1 \geq 0$ and $e_1' \geq 0$. 
If $e_1 \ne e_1'$, then $\sigma_{{\rm (V)^\sharp}, \; e_1}$ 
and $\sigma_{{\rm (V)^\sharp}, \; e_1'}$ are not equivalent.

\item[$\rm (VII)^\sharp$] 
Let $e_1$ and $e_1'$ be integers satisfying $e_1 \geq 0$ and $e_1' \geq 0$. 
If $e_1 \ne e_1'$, then $\sigma_{{\rm (VII)^\sharp}, \; e_1}$ 
and $\sigma_{{\rm (VII)^\sharp}, \; e_1'}$ are not equivalent.

\item[$\rm (IX)^\sharp$] 
Let $e_1$ and $e_1'$ be integers satisfying $e_1 \geq 0$ and $e_1' \geq 0$. 
If $e_1 \ne e_1'$, then $\sigma_{{\rm (IX)^\sharp}, \; e_1}$ 
and $\sigma_{{\rm (IX)^\sharp}, \; e_1'}$ are not equivalent.

\item[$\rm (XI)^\sharp$] 
Let $e_1$ and $e_1'$ be integers satisfying $e_1 \geq 0$ and $e_1' \geq 0$. 
If $e_1 \ne e_1'$, then $\sigma_{{\rm (XI)^\sharp}, \; e_1}$ 
and $\sigma_{{\rm (XI)^\sharp}, \; e_1'}$ are not equivalent.

\item[$\rm (XV)^\sharp$] 
Let $e_2$, $e_3$, $e_2'$, $e_3'$ be integers satisfying 
$e_2 \geq e_3 \geq 0$ and $e_2' \geq e_3' \geq 0$. 
If $(e_2, e_3) \ne (e_2', e_3')$, then $\sigma_{{\rm (XV)^\sharp}, \; (e_2, \, e_3)}$ 
and $\sigma_{{\rm (XV)^\sharp}, \; (e_2', \, e_3')}$ are not equivalent.

\item[$\rm (XIX)^\sharp$] 
Let $e_1$ and $e_1'$ be integers satisfying $e_1 \geq 0$ and $e_1' \geq 0$. 
If $e_1 \ne e_1'$, then $\sigma_{{\rm (XIX)^\sharp}, \; e_1}$ 
and $\sigma_{{\rm (XIX)^\sharp}, \; e_1'}$ are not equivalent.

\item[$\rm (XXIV)^\sharp$] 
Let $e_2$ and $e_2'$ be integers satisfying $e_2 \geq 0$ and $e_2' \geq 0$. 
If $e_2 \ne e_2'$, then $\sigma_{{\rm (XXIV)^\sharp}, \; e_2}$ 
and $\sigma_{{\rm (XXIV)^\sharp}, \; e_2'}$ are not equivalent. 
\end{enumerate} 
\end{lem}

\begin{proof}
$\rm (I)^\sharp$  If $\sigma_{{\rm (I)^\sharp}, \; e_1}$ and 
$\sigma_{{\rm (I)^\sharp}, \; e_1'}$ are equivalent, 
then $\omega_{{\rm (I)^\sharp}, \; e_1}$ and 
$\omega_{{\rm (I)^\sharp}, \; e_1'}$ are equivalent. 
By Lemma 6.2, we have $3 \, p^{e_1} = 3 \, p^{e_1'}$ and $p^{e_1} = p^{e_1'}$, 
which implies $e_1 = e_1'$.

$\rm (II)^\sharp$  If $\sigma_{{\rm (II)^\sharp}, \; e_1}$ and 
$\sigma_{{\rm (II)^\sharp}, \; e_1'}$ are equivalent, 
then $\omega_{{\rm (II)^\sharp}, \; e_1}$ and 
$\omega_{{\rm (I)^\sharp}, \; e_1'}$ are equivalent. 
By Lemma 6.4, we have $p^{e_1 + 1} = p^{e_1' + 1}$ and $p^{e_1} = p^{e_1'}$, 
which implies $e_1 = e_1'$.

$\rm (IV)^\sharp$  If $\sigma_{{\rm (IV)^\sharp}, \; (e_1, \, e_2)}$ and 
$\sigma_{{\rm (IV)^\sharp}, \; (e_1', \, e_2')}$ are equivalent, 
then $\omega_{{\rm (IV)^\sharp}, \; (e_1, \, e_2)}$ and 
$\omega_{{\rm (IV)^\sharp}, \; (e_1', \, e_2')}$ are equivalent. 
By Lemma 6.6, we have $p^{e_1} + p^{e_2} = p^{e_1'} + p^{e_2'}$ and 
$p^{e_2} - p^{e_1} = p^{e_2'} - p^{e_1'}$, 
which implies $e_2 = e_2'$ and $e_1 = e_1'$.

$\rm (V)^\sharp$  If $\sigma_{{\rm (V)^\sharp}, \; e_1}$ and 
$\sigma_{{\rm (V)^\sharp}, \; e_1'}$ are equivalent, 
then $\omega_{{\rm (V)^\sharp}, \; e_1}$ and 
$\omega_{{\rm (V)^\sharp}, \; e_1'}$ are equivalent. 
By Lemma 6.8, we have $p^{e_1 + 1} = p^{e_1' + 1}$, which implies $e_1 = e_1'$.

$\rm (VII)^\sharp$  If $\sigma_{{\rm (VII)^\sharp}, \; e_1}$ and 
$\sigma_{{\rm (VII)^\sharp}, \; e_1'}$ are equivalent, 
then $\omega_{{\rm (VII)^\sharp}, \; e_1}$ and 
$\omega_{{\rm (VII)^\sharp}, \; e_1'}$ are equivalent. 
By Lemma 6.10, we have $p^{e_1 + 1} = p^{e_1' + 1}$ and $p^{e_1} = p^{e_1'}$, 
which implies $e_1 = e_1'$.

$\rm (IX)^\sharp$  If $\sigma_{{\rm (IX)^\sharp}, \; e_1}$ and 
$\sigma_{{\rm (IX)^\sharp}, \; e_1'}$ are equivalent, 
then $\omega_{{\rm (IX)^\sharp}, \; e_1}$ and 
$\omega_{{\rm (IX)^\sharp}, \; e_1'}$ are equivalent. 
By Lemma 6.12, $2 \, p^{e_1} = 2 \, p^{e_1'}$, which implies $e_1 = e_1'$.

$\rm (XI)^\sharp$  If $\sigma_{{\rm (XI)^\sharp}, \; e_1}$ and 
$\sigma_{{\rm (XI)^\sharp}, \; e_1'}$ are equivalent, 
then $\omega_{{\rm (XI)^\sharp}, \; e_1}$ and 
$\omega_{{\rm (XI)^\sharp}, \; e_1'}$ are equivalent. 
By Lemma 6.14, we have $2 \, p^{e_1} = 2 \, p^{e_1'}$, which implies $e_1 = e_1'$.

$\rm (XV)^\sharp$  If $\sigma_{{\rm (XV)^\sharp}, \; (e_2, \, e_3)}$ and 
$\sigma_{{\rm (XV)^\sharp}, \; (e_2', \, e_3')}$ are equivalent, 
then $\omega_{{\rm (XV)^\sharp}, \; (e_2, \, e_3)}$ and 
$\omega_{{\rm (XV)^\sharp}, \; (e_2', \, e_3')}$ are equivalent. 
By Lemma 6.16, we have $p^{e_2} = p^{e_2'}$ and $p^{e_3} = p^{e_3'}$, 
which implies $e_2 = e_2'$ and $e_3 = e_3'$.

$\rm (XIX)^\sharp$  If $\sigma_{{\rm (XIX)^\sharp}, \; e_1}$ and 
$\sigma_{{\rm (XIX)^\sharp}, \; e_1'}$ are equivalent, 
then $\omega_{{\rm (XIX)^\sharp}, \; e_1}$ and 
$\omega_{{\rm (XIX)^\sharp}, \; e_1'}$ are equivalent. 
By Lemma 6.18, we have  $2 \, p^{e_1} = 2 \, p^{e_1'}$, which implies $e_1 = e_1'$.

$\rm (XXIV)^\sharp$  If $\sigma_{{\rm (XXIV)^\sharp}, \; e_2}$ and 
$\sigma_{{\rm (XXIV)^\sharp}, \; e_2'}$ are equivalent, 
then $\omega_{{\rm (XXIV)^\sharp}, \; e_2}$ and 
$\omega_{{\rm (XXIV)^\sharp}, \; e_2'}$ are equivalent. 
By Lemma 6.20, we have $p^{e_2} = p^{e_2'}$, which implies $e_2 = e_2'$. 
\end{proof}

\subsection{The classification of homomorphisms from 
$\SL(2, k)$ to $\SL(4, k)$}

For 
\[
\nu = {\rm I}, \quad {\rm II}, \quad {\rm IV}, \quad {\rm V}, \quad {\rm VII}, \quad 
{\rm IX}, \quad {\rm XI}, \quad {\rm XV}, \quad {\rm XIX}, \quad {\rm XXIV}, 
\quad {\rm XXVI}, 
\]
we denote by $S_{(\nu)^\sharp}$ the set of all equivalence classes 
of homomorphisms from $\SL(2, k)$ to $\SL(4, k)$ of the form $(\nu)^\sharp$, 
i.e., 
\[
 S_{(\nu)^\sharp} 
  := 
  \{ \, [ \, \sigma \, ] \in \Hom(\SL(2, k), \SL(4, k)) / \sim  \; \mid 
  \text{$\sigma$ has the form $(\nu)^\sharp$} \, \} . 
\]

Let 
\begin{align*}
\left\{
\begin{array}{r @{\,} l }
 \Lambda_1 & :=  \{ \, (e_1, e_2) \in \Z^2 \mid e_2 > e_1 \geq 0 \,\}, \\
 \Lambda_2 & := \{ \, e \in \Z \mid e \geq 0 \,\}, \\ 
 \Lambda_3 & := \{ \ (e_2, e_3) \in \Z^2 \mid e_2 \geq e_3 \geq 0 \,\} , \\
 \Lambda_4  & := \{ \, I_4 \, \} . 
\end{array}
\right. 
\end{align*}

\begin{thm}
The following assertions {\rm (1)}, {\rm (2)}, {\rm (3)} hold true: 
\begin{enumerate}[label = {\rm (\arabic*)}]
\item In the case where $p = 2$, we have the following 
natural one-to-one corredpondences: 
\begin{fleqn}[6em]   
\begin{align*}
& \Hom(\SL(2, k), \SL(4, k)) / \sim \\
 & \qquad \cong \;   
  S_{\rm (IV)^\sharp} \, \sqcup \, S_{\rm (V)^\sharp} \, \sqcup \, S_{\rm (XI)^\sharp} \, \sqcup \, 
  S_{\rm (XV)^\sharp} \, \sqcup \, S_{\rm (XIX)^\sharp} \, \sqcup \, S_{\rm (XXIV)^\sharp} \, \sqcup \, 
  S_{\rm (XXVI)^\sharp} \\
   & \qquad \cong  
 \Lambda_1 \, \sqcup \, 
\bigl( \,  \Lambda_2 \, \sqcup \, \Lambda_2 \, \sqcup \, \Lambda_2 \, \sqcup \, \Lambda_2 \, \bigr)
\, \sqcup \, 
\Lambda_3 \, \sqcup \, \Lambda_4 . 
\end{align*} 
\end{fleqn}

\item In the case where $p = 3$, we have the following natural 
one-to-one correspondences:  
\begin{fleqn}[6em] 
\begin{align*}
& \Hom(\SL(2, k), \SL(4, k)) / \sim \\
 & \qquad \cong \;   
  S_{\rm (II)^\sharp} \, \sqcup \, S_{\rm (IV)^\sharp} \, \sqcup \, S_{\rm (VII)^\sharp} \, \sqcup \, 
  S_{\rm (IX)^\sharp} \, \sqcup \, S_{\rm (XV)^\sharp} \, \sqcup \, S_{\rm (XXIV)^\sharp} \, \sqcup \, 
  S_{\rm (XXVI)^\sharp} \\
  & \qquad \cong  
 \Lambda_1 \, \sqcup \, 
\bigl( 
\,  \Lambda_2 \, \sqcup \, \Lambda_2 \, \sqcup \, \Lambda_2 \, \sqcup \, \Lambda_2 \, 
\bigr)
\, \sqcup \, 
\Lambda_3 \, \sqcup \, \Lambda_4 . 
\end{align*} 
\end{fleqn}

\item In the case where $p \geq 5$, we have the following natural 
one-to-one correspondences:  
\begin{fleqn}[6em] 
\begin{align*}
& \Hom(\SL(2, k), \SL(4, k)) / \sim \\
 & \qquad \cong \;   
  S_{\rm (I)^\sharp} \, \sqcup \, S_{\rm (IV)^\sharp} \, \sqcup \, S_{\rm (IX)^\sharp} \, \sqcup \, 
  S_{\rm (XV)^\sharp} \, \sqcup \, S_{\rm (XXIV)^\sharp} \, \sqcup \, 
  S_{\rm (XXVI)^\sharp} \\
  & \qquad \cong  
 \Lambda_1 \, \sqcup \, 
\bigl( 
\,  \Lambda_2 \, \sqcup \, \Lambda_2 \, \sqcup \, \Lambda_2  \, 
\bigr)
\, \sqcup \, 
\Lambda_3 \, \sqcup \, \Lambda_4 . 
\end{align*} 
\end{fleqn} 

\end{enumerate} 
\end{thm}

\begin{proof}
See Theorem 5.26 and Lemmas 6.23, 6.24, 6.25.

\end{proof}

\section{Indecomposable decompositions of homomorphisms from $\SL(2, k)$ 
to $\SL(n, k)$, where $1 \leq n \leq 4$}

Given a homomorphism $\sigma : \SL(2, k) \to \SL(n, k)$, 
we can regard $V(n)$ as an $\SL(2, k)$-module, 
where $V(n)$ is the $n$-dimensional column vector space over $k$. 
We say that $\sigma$ is {\it indecomposable} if 
$V(n)$ is an indecomposable $\SL(2, k)$-module.

\subsection{$1 \leq n \leq 3$} 

Let $1 \leq n \leq 3$. 
In the following, we define homomorphisms $\SL(2, k) \to \SL(n, k)$: 
\begin{enumerate}[label = {\rm (\arabic*)}]

\item Assume $n = 1$. 
Let $\sigma^{(1)} : \SL(2, k) \to \SL(1, k)$ be the homomorphism 
defined by 
\begin{fleqn}[10em]
\begin{align*}
\sigma^{(1)}(A) = I_1 .
\end{align*}
\end{fleqn}

\item Assume $n = 2$. 
\begin{enumerate}[label = {\rm (2.\arabic*)}]
\item 
For all integer $e \geq 0$, we can define a homomorphism 
$  \sigma^{(2.1)}_e  : \SL(2, k) \to \SL(2, k)$ as 
\begin{fleqn}[10em]
\begin{align*}
\renewcommand{\arraystretch}{1.5} 
 \sigma^{(2.1)}_e 
 \left(
 \begin{array}{c c}
  a & b \\
  c & d
 \end{array}
 \right)
 =
 \left(
 \begin{array}{c c}
  a^{p^e} & b^{p^e} \\
  c^{p^e} & d^{p^e} 
 \end{array}
 \right) . 
 \end{align*}
 \end{fleqn}

\item We can define a homomorphism 
$ \sigma^{(2.2)}  : \SL(2, k) \to \SL(2, k)$ as 
\begin{fleqn}[10em]
\begin{align*}
\renewcommand{\arraystretch}{1.5} 
 \sigma^{(2.2)} 
 \left(
 \begin{array}{c c}
  a & b \\
  c & d
 \end{array}
 \right)
 =
I_2 . 
\end{align*}
\end{fleqn} 
\end{enumerate}

\item Assume $n = 3$.  
\begin{enumerate}[label = {\rm (3.\arabic*)}]
\item  In the case where $p = 2$, we define homomorphisms 
$\SL(2, k) \to \SL(3, k)$, as follows:  
\begin{enumerate}[label = {\rm (\alph*)}]
\item For all integer $e \geq 0$, we can define a homomorphism 
$\sigma^{(3.1.{\rm a})}_e  : \SL(2, k) \to \SL(3, k)$ as 
\begin{fleqn}[10em]
\begin{align*}
\renewcommand{\arraystretch}{1.5} 
 \sigma^{(3.1.{\rm a})}_e 
 \left(
 \begin{array}{c c}
  a & b \\
  c & d
 \end{array}
 \right)
 =
\left(
\begin{array}{c c | c}
 a^{p^{e + 1}} &  b^{p^{e + 1}} & 0 \\
 c^{p^{e + 1}} &  d^{p^{e + 1}} & 0 \\ 
 \hline 
 a^{p^e} \, c^{p^e}  & b^{p^e} \, d^{p^e} & 1 
\end{array} 
\right) . 
\end{align*}
\end{fleqn}

\item  For all integer $e \geq 0$, we can define a homomorphism 
$\sigma^{(3.1.{\rm b})}_e  : \SL(2, k) \to \SL(3, k)$ as 
\begin{fleqn}[10em]
\begin{align*}
\renewcommand{\arraystretch}{1.5} 
 \sigma^{(3.1.{\rm b})}_e 
\left(
\begin{array}{c c}
 a & b \\
 c & d 
\end{array}
\right) 
 = 
\left(
\begin{array}{c c | c}
 a^{p^{e + 1}} &  b^{p^{e + 1}} & a^{p^e} \, b^{p^e} \\
 c^{p^{e + 1}} &  d^{p^{e + 1}} & c^{p^e} \, d^{p^e} \\ 
 \hline 
 0 & 0 & 1 
\end{array} 
\right) . 
\end{align*}
\end{fleqn}

\item  For all integer $e \geq 0$, we can define a homomorphism 
$\sigma^{(3.1.{\rm c})}_e  : \SL(2, k) \to \SL(3, k)$ as 
\begin{fleqn}[10em]
\begin{align*}
\renewcommand{\arraystretch}{1.5} 
 \sigma^{(3.1.{\rm c})}_e 
\left(
\begin{array}{c c}
 a & b \\
 c & d 
\end{array}
\right) 
 = 
\left(
\begin{array}{c c | c}
 a^{p^{e}} &  b^{p^{e}} & 0 \\
 c^{p^{e}} &  d^{p^{e}} & 0 \\ 
 \hline 
 0 & 0 & 1 
\end{array} 
\right) . 
\end{align*}
\end{fleqn}

\item We can define a homomorphism 
$\sigma^{(3.1.{\rm d})} : \SL(2, k) \to \SL(3, k)$ as 
\begin{fleqn}[10em]
\begin{align*}
\renewcommand{\arraystretch}{1.5} 
 \sigma^{(3.1.{\rm d})} 
\left(
\begin{array}{c c}
 a & b \\
 c & d 
\end{array}
\right) 
 = 
I_3 . 
\end{align*}
\end{fleqn}

\end{enumerate}

\item In the case where $p \geq 3$, we define homomorphisms 
$\SL(2, k) \to \SL(3, k)$, as follows:  
\begin{enumerate}[label = {\rm (\alph*)}]
\item For all integer $e \geq 0$, we can define a homomorphism 
$\sigma^{(3.2.{\rm a})}_e  : \SL(2, k) \to \SL(3, k)$ as 
\begin{fleqn}[10em]
\begin{align*}
\renewcommand{\arraystretch}{1.5} 
 \sigma^{(3.2.{\rm a})}_e 
 \left(
 \begin{array}{c c}
  a & b \\
  c & d
 \end{array}
 \right)
 =
\left(
\begin{array}{c c c}
 a^{2 \, p^{e}} & a^{p^{e}} \, b^{p^{e}}  & b^{2 \, {p^{e}}}  \\
 2 \, a^{p^{e}} \, c^{p^{e}} & a^{p^{e}} \, d^{p^{e}} + b^{p^{e}} \, c^{p^{e}} &  2 \, b^{p^{e}} \, d^{p^{e}} \\
 c^{2 \, p^{e}} & c^{p^{e}} \, d^{p^{e}}  & d^{2 \, p^{e}}  
\end{array}
\right) . 
\end{align*}
\end{fleqn}

\item  For all integer $e \geq 0$, we can define a homomorphism 
$\sigma^{(3.2.{\rm b})}_e  : \SL(2, k) \to \SL(3, k)$ as 
\begin{fleqn}[10em]
\begin{align*}
\renewcommand{\arraystretch}{1.5} 
 \sigma^{(3.2.{\rm b})}_e 
\left(
\begin{array}{c c}
 a & b \\
 c & d 
\end{array}
\right) 
 = 
\left(
\begin{array}{c c | c}
 a^{p^{e}} &  b^{p^{e}} & 0 \\
 c^{p^{e}} &  d^{p^{e}} & 0 \\ 
 \hline 
 0 & 0 & 1 
\end{array} 
\right) . 
\end{align*}
\end{fleqn}

\item We can define a homomorphism 
$\sigma^{(3.2.{\rm c})} : \SL(2, k) \to \SL(3, k)$ as 
\begin{fleqn}[10em]
\begin{align*}
\renewcommand{\arraystretch}{1.5} 
 \sigma^{(3.2.{\rm c})} 
\left(
\begin{array}{c c}
 a & b \\
 c & d 
\end{array}
\right) 
 = 
I_3 . 
\end{align*}
\end{fleqn}

\end{enumerate}

\end{enumerate}

\end{enumerate}

\begin{lem}
Let $1 \leq n \leq 3$ and let $\sigma : \SL(2, k) \to \SL(n, k)$ be a homomorphism. 
Then the following assertions {\rm (1)}, {\rm (2)}, {\rm (3)} hold true: 
\begin{enumerate}[label = {\rm (\arabic*)}]
\item Assume $n = 1$. 
Then $\sigma$ coincides with the trivial homomorphism $\sigma^{(1)}$.

\item Assume $n = 2$.  
the homomorphism $\sigma : \SL(2, k) \to \SL(3, k)$ 
is equivalent to one of the 
homomorphisms $\sigma^{(2.1)}_e$ $(e \geq 0)$ and $\sigma^{(2.2)}$.

\item Assume $n = 3$. 
\begin{enumerate}[label = {\rm (3.\arabic*)}]
\item In the case where $p = 2$, 
the homomorphism $\sigma : \SL(2, k) \to \SL(3, k)$ 
is equivalent to one of the 
homomorphisms $\sigma^{(3.1.{\rm a})}_e$ $(e \geq 0)$, 
$\sigma^{(3.1.{\rm b})}_e$ $(e \geq 0)$,  
$\sigma^{(3.1.{\rm c})}_e$ $(e \geq 0)$, 
$\sigma^{(3.1.{\rm d})}$. 

\item In the case where $p \geq 3$, 
the homomorphism $\sigma : \SL(2, k) \to \SL(3, k)$ 
is equivalent to one of the 
homomorphisms $\sigma^{(3.2.{\rm a})}_e$ $(e \geq 0)$, 
$\sigma^{(3.2.{\rm b})}_e$ $(e \geq 0)$,  
$\sigma^{(3.2.{\rm c})}$.

\end{enumerate} 
\end{enumerate} 
\end{lem}

\begin{proof}
(1) The proof is straightforward.

(2) The proof is an exercise to the reader.

(3) We can prove this assertion (cf. \cite[Section 4]{Tanimoto 2022}). 
\end{proof}

\begin{lem}
The following assertions {\rm (1)}, {\rm (2)}, {\rm (3)} hold true: 
\begin{enumerate}[label = {\rm (\arabic*)}, itemsep = 6pt]
\item Assume $n = 1$. 
Then $d \bigl( \, \sigma^{(1)} \, \bigr) = (1, 1) $.

\item Assume $n = 2$. 
For all integer $e \geq 0$, we have the following: 
\begin{enumerate}[label = {\rm (2.\arabic*)}]
\item $d \bigl( \, \sigma^{(2.1)}_e \, \bigr) = (0, 0)$. 

\item $d \bigl( \, \sigma^{(2.2)} \, \bigr) = (2, 2)$. 
\end{enumerate}

\item Assume $n = 3$. 
\begin{enumerate}[label = {\rm (3.\arabic*)}]
\item In the case where $p = 2$, for all integer $e \geq 0$, we have the following: 
\begin{enumerate}[label = {\rm (\alph*)}]

\item $d \bigl( \, \sigma^{(3.1.{\rm a})}_e \, \bigr) = (1, 0)$.

\item $d \bigl( \, \sigma^{(3.1.{\rm b})}_e \, \bigr) = (0, 1)$.

\item $d \bigl( \, \sigma^{(3.1.{\rm c})}_e \, \bigr) = (1, 1)$.

\item $d \bigl( \, \sigma^{(3.1.{\rm d})} \, \bigr) = (3, 3)$.
\end{enumerate}

\item In the case where $p \geq 3$, for all integer $e \geq 0$, we have the following: 
\begin{enumerate}[label = {\rm (\alph*)}]
\item $d \bigl( \, \sigma^{(3.2.{\rm a})}_e \, \bigr) = (0, 0)$. 

\item $d \bigl( \, \sigma^{(3.2.{\rm b})}_e \, \bigr) = (1, 1)$.

\item $d \bigl( \, \sigma^{(3.2.{\rm c})} \, \bigr) = (3, 3)$. 
\end{enumerate} 
\end{enumerate} 
\end{enumerate} 
\end{lem}

\begin{proof} 
The proof is straightforward. 

\end{proof}

\begin{thm}
The following assertions {\rm (1)}, {\rm (2)}, {\rm (3)} hold true: 
\begin{enumerate}[label = {\rm (\arabic*)}, itemsep = 6pt]
\item Assume $n = 1$. 
Then $\sigma^{(1)}$ is indecomposable.

\item Assume $n = 2$. 
For all integer $e \geq 0$, we have the following: 
\begin{enumerate}[label = {\rm (2.\arabic*)}] 
\item $\sigma^{(2.1)}_e$ is indecomposable.

\item $\sigma^{(2.2)}$ has the following indecompsable decomposition: 
\begin{fleqn}[10em] 
\begin{align*}
 \sigma^{(2.2)} = \sigma^{(1)} \oplus \sigma^{(1)} . 
\end{align*}
\end{fleqn} 
\end{enumerate}

\item Assume $n = 3$.  
\begin{enumerate}[label = {\rm (3.\arabic*)}]
\item In the case where $p = 2$, for all integer $e \geq 0$, 
we have the following: 
\begin{enumerate}[label = {\rm (\alph*)}] 

\item $\sigma^{\rm (3.1.a)}_e$ is indecompsable. 

\item $\sigma^{ {\rm (3.1.b)}  }_e$ is indecompsable.

\item $\sigma^{ {\rm (3.1.c)}  }_e$ has the following indecomposable decomposition: 
\begin{fleqn}[10em] 
\begin{align*}
 \sigma^{ {\rm (3.1.c)}  }_e = \sigma^{(2.1)}_{e} \oplus \sigma^{(1)} . 
\end{align*}
\end{fleqn}

\item $\sigma^{{\rm (3.1.d)} } $ has the following indecomposable decomposition: 
\begin{fleqn}[10em] 
\begin{align*}
\sigma^{{\rm (3.1.d)} }  = \sigma^{(1)} \oplus \sigma^{(1)} \oplus \sigma^{(1)} . 
\end{align*}
\end{fleqn}

\end{enumerate}

\item In the case where $p \geq 3$, for all integer $e \geq 0$, we have the following: 
\begin{enumerate}[label = {\rm (\alph*)}]  
\item $\sigma^{{\rm (3.2.a)} }_e $ is indecompsable.

\item $\sigma^{{\rm (3.2.b)} }_e $ has the following indecomposable decomposition: 
\begin{fleqn}[10em] 
\begin{align*}
\sigma^{{\rm (3.2.b)} }_e  = \sigma^{(2.1)}_e \oplus \sigma^{(1)} . 
\end{align*}
\end{fleqn}

\item $\sigma^{{\rm (3.2.c)} } $ has the following indecomposable decomposition: 
\begin{fleqn}[10em] 
\begin{align*}
 \sigma^{{\rm (3.2.c)} }  = \sigma^{(1)} \oplus \sigma^{(1)} \oplus \sigma^{(1)} .  
\end{align*}
\end{fleqn}

\end{enumerate}

\end{enumerate} 
\end{enumerate} 
\end{thm}

\begin{proof}
The proof is straightforwad (use Lemma 7.2 for assertion (2.1), 
assertions (3.1) (a), (b), assertion (3.2) (a)).

\end{proof}

\subsection{$n = 4$}

\begin{thm} 
The follwowing assertions {\rm (1)}, {\rm (2)}, {\rm (3)} hold true: 
\begin{enumerate}[label = {\rm (\arabic*)}, itemsep = 12pt]
\item In the case where $p = 2$, we have the following: 
\begin{enumerate}[label = {\rm (\arabic*)}]
\item[$\rm (IV)^\sharp$] For all integers $e_1$ and $e_2$ satisfying 
$e_2 > e_1 \geq 0$, the homomorphism $\sigma_{{\rm (IV)^\sharp}, \; (e_1, \, e_2)}$ 
is indecomposable.

\item[$\rm (V)^\sharp$] For all integer $e_1 \geq 0$, 
the homomorphism $\sigma_{{\rm (V)^\sharp}, \; e_1}$ is indecomposable.

\item[$\rm (XI)^\sharp$] For all integer $e_1 \geq 0$, 
the homomorphism $\sigma_{{\rm (XI)^\sharp,} \; e_1}$ has the following 
indecomposable decomposition: 
\begin{fleqn}[10em] 
\begin{align*}
\sigma_{{\rm (XI)^\sharp}, \; e_1} = \sigma^{\rm (3.1.b)}_{e_1} \oplus \sigma^{(1)} . 
\end{align*}
\end{fleqn}

\item[$\rm (XV)^\sharp$] For all integers $e_2$ and $e_3$ satisfying 
$e_2 \geq e_3 \geq 0$, the homomorphism $\sigma_{{\rm (XV)^\sharp}, \; (e_2, \, e_3)}$ has the following indecomposable decomposition: 
\begin{fleqn}[10em] 
\begin{align*}
\sigma_{{\rm (XV)^\sharp}, \; (e_2, \, e_3)}
  = \sigma^{(2.1)}_{e_2} \oplus \sigma^{(2.1)}_{e_3}  . 
\end{align*}
\end{fleqn}

\item[$\rm (XIX)^\sharp$] For all integer $e_1 \geq 0$, 
the homomorphism $\sigma_{{\rm (XIX)^\sharp}, \; e_1}$ has the following indecomposable decomposition: 
\begin{fleqn}[10em] 
\begin{align*}
 \sigma_{{\rm (XIX)^\sharp}, \; e_1} 
 = \sigma^{\rm (3.1.a)}_{e_1} \oplus \sigma^{(1)} . 
\end{align*}
\end{fleqn}

\item[$\rm (XXIV)^\sharp$] For all integer $e_2 \geq 0$, 
the homomorphism $\sigma_{{\rm (XXIV)^\sharp}, \; e_2}$ has the following indecomposable decomposition: 
\begin{fleqn}[10em] 
\begin{align*}
  \sigma_{{\rm (XXIV)^\sharp}, \; e_2} 
  = \sigma^{\rm (2.1)}_{e_2} \oplus \sigma^{(1)} \oplus \sigma^{(1)}. 
\end{align*}
\end{fleqn}

\item[$\rm (XXVI)^\sharp$] The homomorphism $\sigma_{\rm (XXVI)^\sharp}$ has the following indecomposable decomposition: 
\begin{fleqn}[10em] 
\begin{align*}
  \sigma_{\rm (XXVI)^\sharp} 
  = \sigma^{(1)} \oplus \sigma^{(1)} \oplus \sigma^{(1)} \oplus \sigma^{(1)} . 
\end{align*}
\end{fleqn} 

\end{enumerate}

\item In the case where $p = 3$, we have the following: 
\begin{enumerate}[label = {\rm (\arabic\sharp)}]
\item[$\rm (II)^\sharp$] For all integer $e_1 \geq 0$, 
the homomorphism $\sigma_{{\rm (II)^\sharp}, \; e_1}$ is indecomposable.

\item[$\rm (IV)^\sharp$] For all integers $e_1$ and $e_2$ satisfying 
$e_2 > e_1 \geq 0$, the homomorphism $\sigma_{{\rm (IV)^\sharp}, \; (e_1, \, e_2)}$ 
is indecomposable.

\item[$\rm (VII)^\sharp$] For all integer $e_1 \geq 0$, 
the homomorphism $\sigma_{{\rm (VII)^\sharp}, \; e_1}$ is indecomposable.

\item[$\rm (IX)^\sharp$] For all integer $e_1 \geq 0$, 
the homomorphism $\sigma_{{\rm (IX)^\sharp}, \; e_1}$ 
has the following indecomposable decomposition: 
\begin{fleqn}[10em] 
\begin{align*}
 \sigma_{{\rm (IX)^\sharp}, \; e_1}
  = \sigma^{(\rm 3.2.a)}_{e_1} \oplus \sigma^{(1)} . 
\end{align*}
\end{fleqn}

\item[$\rm (XV)^\sharp$] For all integers $e_2$ and $e_3$ satisfying 
$e_2 \geq e_3 \geq 0$, the homomorphism $\sigma_{{\rm (XV)^\sharp}, \; (e_2, \, e_3)}$ has the following indecomposable decomposition: 
\begin{fleqn}[10em] 
\begin{align*}
\sigma_{{\rm (XV)^\sharp}, \; (e_2, \, e_3)}
  = \sigma^{(2.1)}_{e_2} \oplus \sigma^{(2.1)}_{e_3}  . 
\end{align*}
\end{fleqn}

\item[$\rm (XXIV)^\sharp$] For all integer $e_2 \geq 0$, 
the homomorphism $\sigma_{{\rm (XXIV)^\sharp}, \; e_2}$ has the following indecomposable decomposition: 
\begin{fleqn}[10em] 
\begin{align*}
  \sigma_{{\rm (XXIV)^\sharp}, \; e_2} 
  = \sigma^{\rm (2.1)}_{e_2} \oplus \sigma^{(1)} \oplus \sigma^{(1)}. 
\end{align*}
\end{fleqn}

\item[$\rm (XXVI)^\sharp$] The homomorphism $\sigma_{\rm (XXVI)^\sharp}$ has the following indecomposable decomposition: 
\begin{fleqn}[10em] 
\begin{align*}
  \sigma_{\rm (XXVI)^\sharp} 
  = \sigma^{(1)} \oplus \sigma^{(1)} \oplus \sigma^{(1)} \oplus \sigma^{(1)} . 
\end{align*}
\end{fleqn}

\end{enumerate}

\item In the case where $p \geq 5$, we have the following: 
\begin{enumerate}[label = {\rm (\arabic\sharp)}]
\item[$\rm (I)^\sharp$] For all integer $e_1 \geq 0$, 
the homomorphism $\sigma_{{\rm (I)^\sharp}, \; e_1}$ is indecomposable.

\item[$\rm (IV)^\sharp$] For all integers $e_1$ and $e_2$ satisfying 
$e_2 > e_1 \geq 0$, the homomorphism $\sigma_{{\rm (IV)^\sharp}, \; (e_1, \, e_2)}$ 
is indecomposable.

\item[$\rm (IX)^\sharp$] For all integer $e_1 \geq 0$, 
the homomorphism $\sigma_{{\rm (IX)^\sharp}, \; e_1}$ 
has the following indecomposable decomposition: 
\begin{fleqn}[10em] 
\begin{align*}
 \sigma_{{\rm (IX)^\sharp}, \; e_1}
  = \sigma^{(\rm 3.2.a)}_{e_1} \oplus \sigma^{(1)} . 
\end{align*}
\end{fleqn}

\item[$\rm (XV)^\sharp$] For all integers $e_2$ and $e_3$ satisfying 
$e_2 \geq e_3 \geq 0$, the homomorphism $\sigma_{{\rm (XV)^\sharp}, \; (e_2, \, e_3)}$ has the following indecomposable decomposition: 
\begin{fleqn}[10em] 
\begin{align*}
\sigma_{{\rm (XV)^\sharp}, \; (e_2, \, e_3)}
  = \sigma^{(2.1)}_{e_2} \oplus \sigma^{(2.1)}_{e_3}  . 
\end{align*}
\end{fleqn}

\item[$\rm (XXIV)^\sharp$] For all integer $e_2 \geq 0$, 
the homomorphism $\sigma_{{\rm (XXIV)^\sharp}, \; e_2}$ has the following indecomposable decomposition: 
\begin{fleqn}[10em] 
\begin{align*}
  \sigma_{{\rm (XXIV)^\sharp}, \; e_2} 
  = \sigma^{\rm (2.1)}_{e_2} \oplus \sigma^{(1)} \oplus \sigma^{(1)}. 
\end{align*}
\end{fleqn}

\item[$\rm (XXVI)^\sharp$] The homomorphism $\sigma_{\rm (XXVI)^\sharp}$ has the following indecomposable decomposition: 
\begin{fleqn}[10em] 
\begin{align*}
  \sigma_{\rm (XXVI)^\sharp}
   = \sigma^{(1)} \oplus \sigma^{(1)} \oplus \sigma^{(1)} \oplus \sigma^{(1)}. 
\end{align*}
\end{fleqn}

\end{enumerate}

\end{enumerate} 
\end{thm}

\begin{proof}
\quad 

\begin{enumerate}[label = {\rm (\arabic*)}]
\item  Assertions $\rm (XI)^\sharp$, $\rm (XV)^\sharp$, $\rm (XIX)^\sharp$, 
$\rm (XXIV)^\sharp$, $\rm (XXVI)^\sharp$ are clear. 
\begin{enumerate}
\item[$\rm (IV)^\sharp$]   
Suppose to the contrary that $\sigma_{{\rm (IV)^\sharp}, \; (e_1, \, e_2)}$ is decomposable 
for some $e_2 > e_1 \geq 0$. 
Since $d(\sigma_{{\rm (IV)^\sharp}, \; (e_1, \, e_2)}) = (0, 0)$, 
we must have 
\[
 \sigma_{{\rm (IV)^\sharp}, \; (e_1, \, e_2)}
 \sim \sigma_{e_2'}^{(2.1)} \oplus \sigma_{e_3'}^{(2.1)} 
\]
for some $e_2' \geq e_3' \geq 0$. 
We have a contradiction (see Lemma 6.24 (1)).

\item[$\rm (V)^\sharp$]  
Since $d( \sigma_{{\rm (V)^\sharp}, \; e_1} ) = (1, 1)$, we can show that 
$\sigma_{{\rm (V)^\sharp}, \; e_1}$ is indecompsable (see Lemma 7.2). 
\end{enumerate}

\item Assertions $\rm (IX)^\sharp$, $\rm (XV)^\sharp$, $\rm (XXIV)^\sharp$, 
$\rm (XXVI)^\sharp$ are clear. 
\begin{enumerate}
\item[$\rm (II)^\sharp$]  See Lemmas 6.23 (2) and 6.24 (2) (iii).  

\item[$\rm (IV)^\sharp$] See Lemmas 6.23 (2) and 6.24 (2) (v). 

\item[$\rm (VII)^\sharp$] See Lemmas 6.23 (2) and 6.24 (2) (vi). 
\end{enumerate}

\item Assertions $\rm (IX)^\sharp$, $\rm (XV)^\sharp$, $\rm (XXIV)^\sharp$, 
$\rm (XXVI)^\sharp$ are clear. 
\begin{enumerate}
\item[$\rm (I)^\sharp$]  See Lemmas 6.23 (3) and 6.24 (3) (ii). 

\item[$\rm (IV)^\sharp$]  See Lemmas 6.23 (3) and 6.24 (3) (iii). 
\end{enumerate} 

\end{enumerate}

\end{proof}

\end{document}